\documentclass[a4paper,final,texthight24cm,11pt]{article}
\usepackage{amssymb}
\usepackage{graphicx}
\usepackage{amsmath}

\def\func#1{\mathop{\rm #1}}%
\def\QTR#1#2{{\csname#1\endcsname #2}}%(gp) Is this the best?

\begin{document}

\title{Existence and A Priori Estimates for Euclidean Gibbs States}
\author{S. Albeverio, Yu. Kondratiev, T. Pasurek, and M. R\"{o}ckner}
\maketitle

\begin{abstract}
\noindent We prove a priori estimates and, as sequel, existence of Euclidean
Gibbs states for quantum lattice systems. For this purpose we develop a new
analytical approach, the main tools of which are: first, a characterization
of the Gibbs states in terms of their Radon--Nikodym derivatives under shift
transformations as well as in terms of their logarithmic derivatives through
integration by parts formulae, and second, the choice of appropriate
Lyapunov functionals describing stabilization effects in the system. The
latter technique becomes applicable since on the basis of the integration by
parts formulae the Gibbs states are characterized as solutions of an
infinite system of partial differential equations. Our existence result
generalize essentially all previous ones. In particular, superquadratic
growth of the interaction potentials is allowed and $N$-particle
interactions for $N\in \mathbb{N}\cup \{\infty \}$ are included. We also
develop abstract frames both for the necessary single spin space analysis
and for the lattice analysis apart from their applications to our concrete
models. Both types of general results obtained \ in these two frames should
be also of their own interest in infinite dimensional analysis.\medskip
\medskip

\noindent MATHEMATICAL\ SUBJECT\ CLASSIFICATION (2000):

\noindent Primary: 82B10;\quad Secondary: 46G12, 60H30\medskip

KEY WORDS\ AND\ PHRASES:\quad quantum lattice systems; Euclidean Gibbs
states; smooth measures on infinite dimensional spaces and their logarithmic
derivatives; integration by parts formulae; Lyapunov functionals\newpage
\end{abstract}

\tableofcontents

\section{Introduction}

\noindent \textbf{I. Background, methods and purposes} \ This paper is
concerned with models of quantum anharmonic lattice systems. In statistical
physics they are commonly viewed as models for \emph{quantum crystals} and
are closely related to (Euclidean) quantum lattice field theory with
continuous time (cf. [DLP78, GJ81]). A mathematical description of
equilibrium properties of quantum systems might be given in terms of their%
\emph{\ Gibbs states}. However, a corresponding general concept of quantum
Gibbs states (including their rigorous definition) as that based on the
standard algebraic approach [BrRo81] still remains open. Thus we shall take
the \emph{Euclidean approach.} This approach is conceptually analogous to
the well-known Euclidean strategy in quantum field theory (see e.g. [Sim74,
GJ81]. This analogy was pointed out and first developed in [AH-K75] (see
also the recent developments in [AKKR01,02]). It transforms the problem of
constructing quantum Gibbs states as functionals on the algebra of local
observables into the problem of studying certain \emph{Euclidean Gibbs
measures} $\mu $ on the loop lattice $\Omega _{\beta }:=[C(S_{\beta })]^{%
\mathbb{Z}^{d}}.$ Here $\beta >0$ has the meaning of the inverse temperature
and $S_{\beta }\cong \lbrack 0,\beta ]$ is a circle with Lebesgue measure $%
d\tau .$ For a more detailed discussion of the relations between quantum and
Euclidean Gibbs states we refer to the Appendix (Sect.\thinspace 8) below.

As compared with classical lattice systems, the situation with Euclidean
Gibbs measures of quantum anharmonic systems is much more rich and
complicated, since for last systems the single spin spaces (e.g., $%
C(S_{\beta })$ or $L_{\beta }^{R}(S_{\beta }),$ $R\geq 1$) themselves are 
\emph{infinite dimensional} and their topological features should be taken
into account properly.

On the other hand, Gibbs measures, usually defined by local specifications
through the Dobrushin--Lanford--Ruelle equations (see, e.g., [Pr76, Ge88])
can also be described via their \emph{Radon--Nikodym derivatives }under
shift transformations of the configuration spaces (cf. Theorem 4.6 in
[AKR97b] and its extension in Proposition 4.2 below). If the interaction
potentials are smooth enough, this flow characterization is in turn
equivalent to the characterization of Gibbs states as differentiable
measures satisfying integration by parts formulas with prescribed \emph{%
logarithmic derivatives }(cf. Proposition 4.9). As is typically for
non-compact spin spaces, we have to restrict ourselves to the set $\mathcal{G%
}_{t}$ of tempered Gibbs measures $\mu ,$ which we specify by some natural
support condition (see Definition 3.9). The aim of this paper is to show
that such alternative descriptions in a \emph{direct analytical way} give
both the existence of Euclidean Gibbs measures $\mu \in \mathcal{G}_{t}$
(cf. Main Theorem I in Subsect.\thinspace 2.3) and a priori estimates on
their correlations in terms of parameters of the interaction (cf. Main
Theorem II in Subsect.\thinspace 2.3).

The essential ingredient of the proofs is that the characterization of Gibbs
measures via integration by parts gives the possibility to deal with them as
with the solutions of an infinite system of first order PDE's 
\begin{equation}
\partial _{h_{i}}\mu =b_{i}\cdot \mu ,\quad i\in \mathbb{N}.  \tag{1.1}
\end{equation}%
Here $h_{i}$ are some admissible directions forming an orthonormal basis in
the tangent Hilbert space $\mathcal{H}:=l^{2}(\mathbb{Z}^{d})\otimes
L^{2}(S_{\beta })$ and $b_{i}$ are the corresponding partial logarithmic
derivatives of (all) measures $\mu \in \mathcal{G}_{t}$ along $h_{i}.$ This
enables us to employ the Lyapunov function method (in a similar way as in
finite dimensional PDE's) in order to get a priori moment estimates on $\mu
\in \mathcal{G}^{t}$. However, when dealing with the integration by parts
formulas (1.1), the main difficulty one encounters is that one does not know
in advance whether $b_{i}\in L^{1}(\mu )$ for any fixed direction $h_{i}$
and measure $\mu \in \mathcal{G}_{t}$. This problem is successfully overcome
by a special choice of test functions $f$ to which we can correctly apply
both sides of the distributional identity (1.1). The local Gibbs
specifications also satisfy the same integration by parts formulas, from
which we deduce moment estimates\emph{\ uniformly} in volume. The latter is
crucial for our proof of the existence result for Euclidean Gibbs measures $%
\mu $, i.e., $\mathcal{G}_{t}\neq \emptyset $.

It should be mentioned that this approach has first been realized for
classical lattice systems in [AKRT99,00]. However, the abstract scheme, in
the form suggested in those papers, is not directly applicable (see also
[BR00, BRW01]) to the present case. The reason is that for Euclidean Gibbs
states we have to perform not only \textquotedblright \emph{lattice analysis}%
\textquotedblright , but need also a separate and quite different
\textquotedblright \emph{single spin space analysis}\textquotedblright .

Some results on the existence of Euclidean Gibbs measures, concerning
specific classes of anharmonic interactions and using various techniques,
have been already known before (see Subsection 2.6 for the references and a
detailed discussion of the previous \textquotedblright state of the
art\textquotedblright\ of the problem). But we emphasize that our technique
is completely different and rather elementary, provided one has the
integration by parts description of Gibbs measures (or uses it as a
definition). We also would like to point out an advantage of our approach
over the so-called stochastic dynamics method, since now we can avoid the
extremely difficult ergodicity problems for the related infinite dimensional
stochastic evolutional equations (cf. the corresponding discussion in
Subsect.\thinspace 5.3). We would also like to stress that our approach for
the first time gives existence of Euclidean Gibbs states for many-particle
interactions with infinite radius and superquadratic growth, which has not
been possible to obtain by other known methods. In addition, we obtain
useful and precise information on the support properties of all tempered
Euclidean Gibbs states.

For the reason of notation simplicity only, in this paper we decided to
restrict ourselves to the systems of \emph{one-dimensional} (i.e., \emph{%
polarized}) oscillators. However, all the results obtained naturally extend
to the multi-dimensional case. A possible way of such generalization is
described in [AKPR04].\medskip

\noindent \textbf{II. Structure and contents of the paper}\quad Along with
the Introduction (\emph{Section 1}) and the Appendix (\emph{Section 8}), the
paper consists of two parts, namely, \textbf{Part I}: {\large %
\textquotedblright Gibbs setting\textquotedblright }\textbf{\ }(\emph{%
Sections 2--5})\textbf{\ }and \textbf{Part II}: {\large \textquotedblright
Abstract setting\textquotedblright }\ (\emph{Sections 6, 7}).\medskip

\textbf{Part I} is devoted to general aspects of the theory of Euclidean
Gibbs measures and their characterization through integration by parts
formulas.

\emph{Section 2:}\quad In order to present main ideas, we start with some
particular classes of quantum lattice systems, \textquotedblright anharmonic
crystals\textquotedblright . Then we formulate our main Theorems I--III on
the existence and a priori estimates for tempered Euclidean Gibbs measures
and compare them with the previous results obtained by other methods.

\emph{Section 3:}\quad We introduce the general models of quantum lattice
systems and discuss the assumptions on the interactions (which are not
necessarily translation invariant and possibly of infinite range). Next, we
describe in detail the corresponding Gibbsian formalism and define the set $%
\mathcal{G}_{t}$ of all tempered Euclidean Gibbs measures $\mu $ on the loop
lattice $\Omega .$

\emph{Section 4:}\quad We begin with the alternative description of $\mu \in 
\mathcal{G}_{t}$ in terms of their Radon--Nikodym derivatives under shifts
on the configuration space (see Proposition 4.1). Then we prove an
equivalent description of $\mu \in \mathcal{G}_{t}$ as differentiable
measures satisfying (in the distributional sense by pairing with proper test
functions $f$) the integration by parts formula (1.1) with the given
logarithmic derivatives $b_{i}$ (see Proposition 4.9).

\emph{Section 5:}\quad Here we demonstrate the applications of the
(IbP)-formula to the study of $\mu \in \mathcal{G}_{t}.$ Also we discuss the
validity of the same flow and integration by parts characterizations for
local Gibbs specifications. And finally, we give elementary proofs of our
Main Theorems I--III, but under some Hypotheses (\textbf{H}) and (\textbf{H}$%
_{loc}$), which are proved to always hold for our systems subsequently in
Section 7.\medskip

In\textbf{\ Part II} we transfer the problem in a general framework of
symmetrizing measures on Banach lattices.

\emph{Section 6:}\quad Here we perform the \textquotedblright \emph{single
spin space analysis}\textquotedblright\ on a general Banach space $X$, which
later on will be taken as the loop space $C(S_{\beta })$. More precisely, we
study symmetrizing measures $\mu \in \mathcal{M}^{b}(X)$ on a single Banach
space $X$ (i.e., those measures satisfying integration by parts formulas
with some prescribed partial logarithmic derivatives $b_{i})$. We suppose
that the logarithmic derivatives have a linear component $A$ (being a
positive selfadjoint operator in some tangent Hilbert space $H\supseteq X$)
and the nonlinear one $F$ (possessing certain coercivity properties w.r.t. $%
H $). Developing a Lyapunov functional method in this situation, we derive a
priori estimates on the moments of $\mu \in \mathcal{M}^{b}(X).$

\emph{Section 7:}\quad We further enrich the abstract setting of the
previous section by adding an extra \textquotedblright \emph{lattice
structure}\textquotedblright . So, in order to include the case of Euclidean
Gibbs measures, we consider symmetrizing measures $\mu \in \mathcal{M}^{b}(%
\mathcal{X})$ on Banach lattices\ $\mathcal{X}:=X^{\mathbb{Z}^{d}}.$ The
required a priori estimates on $\mu \in \mathcal{M}^{b}(\mathcal{X})$ are
formulated as Theorems 7.1 and 7.3. At the end of this Section we come back
to the Euclidean Gibbs measures and (on the basis of Theorems 7.1 and 7.3
just proved) verify the validity of the Hypotheses (\textbf{H}) and (\textbf{%
H}$_{\mathtt{loc}}$) for them.\medskip

In the Appendix (\emph{Section 8}) we enclose a discussion on the Euclidean
approach to quantum Gibbs states.\medskip

Finally we mention that the results of the paper have been announced in
[AKPR01a,01b,04] and presented in various talks since December 2000 during
seminars or conferences, e.g., in Berlin, Kiev, Moscow, Oberwolfach, and
Pisa.

\bigskip

\begin{center}
{\Huge Part I: Gibbs setting}
\end{center}

\section{Particular models of quantum anharmonic crystals}

In order to fix the main ideas and make the reader more familiar with the
topic, in this section we give a self-contained presentation of some
particular models of\emph{\ quantum lattice systems} (\emph{QLS}, for
shorthand). So, in Subsects.\thinspace 2.1--2.3 we concentrate on the
simplest case of the\emph{\ translation invariant} system with the\emph{\
harmonic pair interaction\ }between \emph{nearest neighbors} only (cf. \emph{%
Model I}). Following the Euclidean approach, in Subsect.\thinspace 2.2 we
define the corresponding quantum Gibbs states as classical Gibbs measures,
but with \emph{infinite dimensional} single spin (i.e., \emph{loop}) spaces.
In Subsect.\thinspace 2.3 we formulate our main results on the existence and
a priori estimates for tempered Gibbs measures on loop lattices. Further
modifications that one needs in order to treat the case of the pair
interaction potentials of \emph{superquadratic growth} (\emph{Model II})
resp. of possibly \emph{infinite range} (\emph{Model III }) are briefly
described in Subsect.\thinspace 2.4 resp.\thinspace 2.5. Finally, in
Subsect.\thinspace 2.6 we discuss fundamental problems and basic methods in
the study\ of Euclidean Gibbs states, as well as compare our results with
those previously obtained by other authors.

\subsection{Model I: harmonic pair interaction}

We start with the following simplest and sufficiently popular in the
literature model of a quantum crystal. Let $\mathbb{Z}^{d}(\subset \mathbb{R}%
^{d})$ be the $d$-dimensional integer lattice with the Euclidean distance $%
\left\vert k-j\right\vert ,\ k,j\in \mathbb{Z}^{d}$. We consider an infinite
system of interacting quantum particles performing one-dimensional (i.e.,
polarized) oscillations with displacements $q_{k}\in \mathbb{R}$ around
their equilibrium positions at points $k\in \mathbb{Z}^{d}.$ Each particle
individually is described by the quantum mechanical Hamiltonian 
\begin{equation}
\mathbb{H}_{k}:=-\frac{1}{2\mathfrak{m}}\frac{d^{2}}{dq_{k}^{2}}+\frac{a^{2}%
}{2}q_{k}^{2}+V(q_{k})  \tag{2.1}
\end{equation}%
acting in the (physical) Hilbert state space $\mathcal{H}_{k}:=L^{2}(\mathbb{%
R},dq_{k}).$ Here $\mathfrak{m}$ $(=\mathfrak{m}_{ph}/\hbar ^{2})>0$ is the
(reduced) mass of the particles and $a^{2}>0$ is their rigidity w.r.t. the
harmonic oscillations. Concerning the anharmonic self-interaction potential,
we suppose that $V\in C^{2}(\mathbb{R\rightarrow R}),$ i.e.,\emph{\ }twice
continuously differentiable, and, moreover, that it satisfies the following
growth condition:\smallskip

\noindent \textbf{Assumption }$\mathbf{(V}_{\mathbf{0}}\mathbf{)}$\textbf{:}
\ \emph{There exist some constants }$P>2$ \emph{and }$K_{V},C_{V}>0$ \emph{%
such that for all }$q\in \mathbb{R}$%
\begin{equation*}
K_{V}^{-1}|q|^{P-l}-C_{V}\leq (\text{sgn}q)^{l}\cdot V^{(l)}(q)\leq
K_{V}|q|^{P-l}+C_{V},\text{ \ }l=0,1,2.
\end{equation*}

Next, we add the harmonic nearest-neighbor interaction $W(q_{k},q_{k^{\prime
}}):=J(q_{k}-q_{k^{\prime }})^{2}$ with the intensity $J>0,$ the sum being
taken over all (unordered) pairs $\langle k,k^{\prime }\rangle $ in $\mathbb{%
Z}^{d}$ such that $|k-k^{\prime }|=1.$ The whole system is then described by
a heuristic Hamiltonian of the form 
\begin{equation}
\mathbb{H}:=-\frac{1}{2\mathfrak{m}}\sum_{k\in \mathbb{Z}^{d}}\frac{d^{2}}{%
dq_{k}^{2}}+\frac{a^{2}}{2}\sum_{k\in \mathbb{Z}^{d}}q_{k}^{2}+\sum_{k\in 
\mathbb{Z}^{d}}V(q_{k})+J\sum\limits_{\langle k,k^{\prime }\rangle \subset 
\mathbb{Z}^{d}}(q_{k}-q_{k^{\prime }})^{2}.  \tag{2.2}
\end{equation}%
Actually, the infinite-volume Hamiltonian (2.2) cannot be defined directly
as a mathematical object and is represented by the local (i.e., indexed by
finite volumes $\Lambda \subset \mathbb{Z}^{d})$ Hamiltonians%
\begin{equation*}
\mathbb{H}_{\Lambda }:=\sum_{k\in \Lambda }H_{k}+J\sum\limits_{\langle
k,k^{\prime }\rangle \subset \Lambda }(q_{k}-q_{k^{\prime }})^{2}
\end{equation*}%
(as self-adjoint and lower bounded Schr\"{o}dinger operators) acting in the
Hilbert spaces $\mathcal{H}_{\Lambda }:=\otimes _{k\in \Lambda }\mathcal{H}%
_{k}$. \smallskip

Lattice systems of the above type (as well as their generalizations studied
in Sect.\thinspace 3 below) are commonly viewed in quantum statistical
physics as mathematical models of a crystalline substance (for more physical
background see e.g. [DLP79, FM99, MVZ00, AKKR02]). The study of such systems
is especially motivated by the reason, that they provide a mathematically
rigorous and physically quite realistic description for the important
phenomenon of phase transitions (i.e., non-uniqueness of Gibbs states). So,
if the potential $V$ has a double-well shape, in the large mass limit $%
\mathfrak{m\rightarrow \infty }$ the QLS (2.2) may undergo (ferroelectic)
structural phase transitions connected with the appearance of macroscopic
displacements of particles for low temperatures $\beta ^{-1}<\beta
_{cr}^{-1}(\mathfrak{m})$ (for a mathematical justification of this effect
see e.g. [BaK91, HM00]). \smallskip

\noindent \textbf{Remark~2.1.} \textbf{(i)} \ In fact, from the potential $V$
one can always extract the quadratic term $a^{2}q^{2}/2$ with a small $%
a^{2}>0$, so that $(\mathbf{V}_{\mathbf{0}})$ is still true.\textbf{\ }%
Typical potentials satisfying Assumption $\mathbf{(V}_{\mathbf{0}}\mathbf{)}$
are polynomials of even degree and with a positive leading coefficient, i.e.,%
\begin{equation}
P(q):=V(q):=\sum_{1\leq l\leq 2n}b_{l}q^{l}\text{ \ with }b_{2n}>0\text{ \
and \ }n\geq 2.  \tag{2.3}
\end{equation}%
In this case one speaks about so-called \emph{ferromagnetic} $P(\phi )$\emph{%
--models}, which also can naturally be looked upon as lattice
discretizations of quantum $P(\phi )$--fields (cf. [Sim74, GJ81]). Let us
also mention a special choice in (2.3), when%
\begin{equation}
P(q):=\sum_{0\leq l\leq n}b_{2l}q^{2l}\text{ \ with }b_{2l}\geq 0\text{ \
for all \ }2\leq l\leq n.  \tag{2.4}
\end{equation}%
Since $b_{2}\in \mathbb{R}$ can be a large negative number, such polynomials
may have arbitrary deep double wells. \smallskip

\noindent \textbf{Notation 2.2.}\ \ Throughout the paper we shall use the
following notation. For a set $\Lambda \subset \mathbb{Z}^{d}$, we denote by 
$|\Lambda |$ its \emph{cardinality}, by $diam$ $\Lambda :=\sup_{k,k^{\prime
}\in \Lambda }|k-k^{\prime }|$ its \emph{diameter}, by $\Lambda ^{c}:=%
\mathbb{Z}^{d}\backslash \Lambda $ its \emph{complement}, and by $\partial
\Lambda :=\left\{ k^{\prime }\in \Lambda ^{c}\text{\thinspace }\left\vert 
\text{\thinspace }\exists k\in \Lambda ,\text{ }\left\vert k-k^{\prime
}\right\vert =1\right. \right\} $ its \emph{boundary}. In particular, $%
\partial k:=\left\{ k^{\prime }\in \mathbb{Z}^{d}\text{\thinspace }%
\left\vert \text{\thinspace }\left\vert k-k^{\prime }\right\vert =1\right.
\right\} $ is the set of all neighbors of $k$ consisting of $2d$ points. We
write $\Lambda \Subset \mathbb{Z}^{d}$ whenever $1\leq |\Lambda |<\infty $.
As usual, $\Lambda \nearrow \mathbb{Z}^{d}$ means the limit as $N\rightarrow
\infty $ along any increasing sequence of volumes $\Lambda ^{(N)}\subset
\Lambda ^{(N+1)}\Subset \mathbb{Z}^{d}$ such that $\bigcup_{N\in \mathbb{N}%
}\Lambda ^{(N)}=\mathbb{Z}^{d}.$

\subsection{Definition of Euclidean Gibbs measures}

As was already mentioned in the Introduction, we take the\emph{\ }Euclidean
approach based on a path space representation for quantum Gibbs states.
Since such approach involves intricate relations between quantum statistical
mechanics and stochastic processes, for convenience of the non-expert reader
we enclose a more extended discussion of the related topics in the Appendix.
Now we briefly describe the corresponding Euclidean Gibbsian formalism just
for the concrete class of quantum lattice systems (2.2); for all necessary
details presented already in the context of the general model (3.1) we refer
the reader to Section 3 below.\smallskip

Let $S_{\beta }\cong \lbrack 0,\beta ]$ be a circle of length $\beta ,$
where we fix positive $0<\beta :=(kT)^{-1}$ having the meaning of inverse
temperature. For each $k\in \mathbb{Z}^{d},$ denote by 
\begin{gather*}
C_{\beta }^{n+\alpha }:=C^{n+\alpha }(S_{\beta }\rightarrow \mathbb{R}),%
\text{ \ \ }n\in \mathbb{N}\cup \{0\},\text{ }\alpha \in (0,1), \\
L_{\beta }^{r}:=L^{r}(S_{\beta }\rightarrow \mathbb{R},d\tau ),\text{ \ \ }%
r\geq 1,
\end{gather*}%
the standard Banach spaces of all (H\"{o}lder) continuous resp. integrable
(w.r.t. Lebesgue measure $d\tau $) functions, i.e., loops, $\omega
_{k}:S_{\beta }\rightarrow \mathbb{R}.$ In particular, $C_{\beta }$ with the 
$\sup $-norm $|\cdot |_{C_{\beta }}$ will be treated as the \emph{single
spin space}, whereas $L_{\beta }^{2}$ with the inner product $(\cdot ,\cdot
)_{L_{\beta }^{2}}:=|\cdot |_{L_{\beta }^{2}}^{2}$ as the\emph{\ }Hilbert
space\ tangent to\emph{\ }$C_{\beta }$.

As the \emph{configuration space} for the infinite volume system we define
the\emph{\ }product loop space%
\begin{equation}
\Omega :=[C_{\beta }]^{\mathbb{Z}^{d}}=\left\{ \omega =(\omega _{k})_{k\in 
\mathbb{Z}^{d}}\left\vert \omega :S_{\beta }\rightarrow \mathbb{R}^{\mathbb{Z%
}^{d}},\text{ }\omega _{k}\in C_{\beta }\right. \right\}  \tag{2.5}
\end{equation}%
equipped with the corresponding\ Borel\emph{\ }$\sigma $\emph{-}algebra $%
\mathcal{B}(\Omega ).$ Let $\mathcal{M}(\Omega )$ denote the set of all
probability measures on $(\Omega ,\mathcal{B}(\Omega )).$ Next, we introduce
the subset of (\emph{\textquotedblright exponentially
increasing\textquotedblright }) \emph{tempered\ configurations}%
\begin{equation}
\Omega _{(e)t}:=\left\{ \omega \in \Omega \text{ }\left\vert \text{ }\forall
\delta \in (0,1):\text{ }||\omega ||_{-\delta }:=\left[ \sum\nolimits_{k\in 
\mathbb{Z}^{d}}e^{-\delta |k|}|\omega _{k}|_{L_{\beta }^{2}}^{2}\right] ^{%
\frac{1}{2}}<\infty \right. \right\}  \tag{2.6}
\end{equation}%
and respectively the subset of\ tempered\emph{\ }measures supported by $%
\Omega _{(e)t}\in \mathcal{B}(\Omega )$, i.e., 
\begin{equation}
\mathcal{M}_{(e)t}:=\left\{ \mu \in \mathcal{M}(\Omega )\,\left\vert \,\mu
\left( \Omega _{(e)t}\right) =1\right. \right\} .  \tag{2.7}
\end{equation}%
\medskip

Heuristically, the Euclidean Gibbs measures $\mu $ we are interested in have
the following representation 
\begin{equation}
d\mu (\omega ):=Z^{-1}\exp \left\{ -\mathcal{I}(\omega )\right\}
\prod\limits_{k\in \mathbb{Z}^{d}}d\gamma _{\beta }(\omega _{k}),  \tag{2.8}
\end{equation}%
where $Z$ is the normalization factor, $\gamma _{\beta }$ is a centered
Gaussian measure on $(C_{\beta },$ $\mathcal{B}(C_{\beta }))$ with
correlation operator $A_{\beta }^{-1},$ and $A_{\beta }:=-\mathfrak{m}\Delta
_{\beta }+a^{2}\mathbf{1}$ is the shifted Laplace--Beltrami operator on the
circle $S_{\beta }$. Respectively the map%
\begin{equation}
\Omega \ni \omega \longmapsto \mathcal{I}(\omega ):=\int_{S_{\beta }}\left[
\sum_{k\in \mathbb{Z}^{d}}V(\omega _{k})+\sum_{<k,k^{\prime }>\subset 
\mathbb{Z}^{d}}W(\omega _{k},\omega _{k^{\prime }})\right] d\tau  \tag{2.9}
\end{equation}%
might be viewed as a potential energy functional describing an interacting
system of loops $\omega _{k}\in C_{\beta }$ indexed by $k\in \mathbb{Z}^{d}.$
Of course it is impossible to use this presentation for $\mu $ literally,
since the series in (2.9) do not converge in any sense. We recall that,
relying on the \emph{Dobrushin--Lanford--Ruelle (DLR) formalism }(cf. [Do70,
Pr76, Ge88]), a rigorous meaning can be given to the measures $\mu $ as
random fields on $\mathbb{Z}^{d}$ with a prescribed family of local
specifications $\left\{ \pi _{\Lambda }\right\} _{\Lambda \Subset \mathbb{Z}%
^{d}}$ in the following way:

For $\Lambda \Subset \mathbb{Z}^{d},$ $\pi _{\Lambda }$ is defined as a
probability kernel on $(\Omega ,\mathcal{B}(\Omega ))$: for all $\Delta \in 
\mathcal{B}(\Omega )$ and $\xi \in \Omega $%
\begin{equation}
\pi _{\Lambda }(\Delta |\xi ):=Z_{\Lambda }^{-1}(\xi )\int_{\Omega _{\Lambda
}}\exp \left\{ -\mathcal{I}_{\Lambda }(\omega |\xi )\right\} \mathbf{1}%
_{\Delta }(\omega _{\Lambda },\xi _{\Lambda ^{c}})
\text{$\prod\limits_{k\in\Lambda}$}
d\gamma _{\beta }(\omega _{k})  \tag{2.10}
\end{equation}%
(where $\mathbf{1}_{\Delta }$ denotes the indicator on $\Delta $). Here $%
Z_{\Lambda }(\xi )$ is the normalization factor and 
\begin{equation}
\mathcal{I}_{\Lambda }(\omega |\xi ):=\int_{S_{\beta }}\left[ \sum_{k\in
\Lambda }V(\omega _{k})+\sum\limits_{<k,k^{\prime }>\subset \Lambda
}W(\omega _{k},\omega _{k^{\prime }})+\sum_{k\in \Lambda ,\text{ }k^{\prime
}\in \Lambda ^{c}}W(\omega _{k},\xi _{k^{\prime }})\right] d\tau  \tag{2.11}
\end{equation}%
is the interaction in the volume $\Lambda $ under the boundary condition $%
\xi _{\Lambda ^{c}}:=(\xi _{k^{\prime }})_{k^{\prime }\in \Lambda ^{c}}.$
Obviously, $\inf_{\omega \in \Omega }\mathcal{I}_{\Lambda }(\omega |\xi
)>-\infty $ and the RHS in (2.11) makes sense for the potentials $V,W$ we
deal here with. An important point is the consistency property for $\{\pi
_{\Lambda }\}_{\Lambda \Subset \mathbb{Z}^{d}}$:\emph{\ }for all $\Lambda
\subset \Lambda ^{\prime }\Subset \mathbb{Z}^{d},$\emph{\ }$\xi \in \Omega $
and $\Delta \in \mathcal{B}(\Omega )$%
\begin{equation*}
(\pi _{\Lambda ^{\prime }}\pi _{\Lambda })(\Delta |\xi ):=\int_{\Omega }\pi
_{\Lambda ^{\prime }}(d\omega |\xi )\pi _{\Lambda }(\Delta |\omega )=\pi
_{\Lambda ^{\prime }}(\Delta |\xi ).
\end{equation*}

\noindent \textbf{Definition 2.3. \ }\emph{A probability measure} $\mu $\ 
\emph{on} $(\Omega ,\mathcal{B}(\Omega ))$\ \emph{is called\ Euclidean Gibbs
measure for the specification} $\{\pi _{\Lambda }\}_{\Lambda \Subset \mathbb{%
Z}^{d}}$ \emph{(corresponding to the quantum lattice system (2.2) at inverse
temperature} $\beta >0$\emph{) if it satisfies the DLR} \emph{equilibrium
equations: for all }$\Lambda \Subset \mathbb{Z}^{d}$\ \emph{and} $\Delta \in 
\mathcal{B}(\Omega )$%
\begin{equation}
\mu \pi _{\Lambda }(\Delta ):=\int_{\Omega }\mu (d\omega )\pi _{\Lambda
}(\Delta |\omega )=\mu (\Delta ).  \tag{2.12}
\end{equation}

Fixing $\beta >0$, let $\mathcal{G}$ denote the set of all such measures $%
\mu .$ We shall be concerned with the subset $\mathcal{G}_{(e)t}$ of \emph{%
tempered }Gibbs measures supported by $\Omega _{(e)t},$ i.e.,%
\begin{equation}
\mathcal{G}_{(e)t}:=\mathcal{G}\cap \mathcal{M}_{(e)t}=\left\{ \mu \in 
\mathcal{G}\,\left\vert \,\mu \left( \Omega _{(e)t}\right) =1\right.
\right\} .  \tag{2.13}
\end{equation}%
\smallskip \noindent \textbf{Remark~2.4. (i)} \ Later it will be instructive
to compare our results on quantum systems with the analogous classical ones.
The large-mass limit $\mathfrak{m\rightarrow \infty }$ (or $\hbar
\rightarrow 0$) of model (2.2) gives us an infinite system of interacting
classical particles moving in the external field $V.$ Such system is
described by the potential energy functional%
\begin{equation}
H_{cl}(q)=\frac{a^{2}}{2}\sum_{k\in \mathbb{Z}^{d}}q_{k}^{2}+\sum_{k\in 
\mathbb{Z}^{d}}V(q_{k})+\sum\limits_{\langle k,k^{\prime }\rangle \subset 
\mathbb{Z}^{d}}W(q_{k},q_{k^{\prime }})  \tag{2.14}
\end{equation}%
on the configuration space $\Omega _{cl}:=\mathbb{R}^{\mathbb{Z}^{d}}\ni
\{q_{k}\}_{k\in \mathbb{Z}^{d}}:=q$ (see [AKKR02]). Again, the formal
Hamiltonian (2.14) does not make sense itself and is represented by the
local Hamiltonians%
\begin{equation}
\mathbb{H}_{cl,\Lambda }(q|y):=\frac{a^{2}}{2}\sum_{k\in \Lambda
}q_{k}^{2}+\sum_{k\in \Lambda }V(q_{k})+\sum\limits_{<k,k^{\prime }>\subset
\Lambda }W(q_{k},q_{k^{\prime }})+\sum_{k\in \Lambda ,\text{ }k^{\prime }\in
\Lambda ^{c}}W(q_{k},y_{k^{\prime }})  \tag{2.15}
\end{equation}%
in the volumes $\Lambda \Subset \mathbb{Z}^{d}$ given the boundary
conditions $y\in \Omega _{cl}$. The corresponding Gibbs states $\mu \in 
\mathcal{G}_{cl}$ at inverse temperature $\beta >0$ are defined as
probability measures on $\Omega _{cl}$ satisfying the DLR equations $\mu \pi
_{\Lambda }=\mu ,$ $\Lambda \Subset \mathbb{Z}^{d},$ with the family of
local specifications%
\begin{equation}
\pi _{\Lambda }(\Delta |y):=Z_{\Lambda }^{-1}(y)\int_{\mathbb{R}^{\Lambda
}}\exp \left\{ -\beta H_{cl,\Lambda }(q|y)\right\} \mathbf{1}_{\Delta
}(q_{\Lambda },y_{\Lambda ^{c}})
\text{$\prod\limits_{k\in \Lambda }$}
dq_{k},\text{ \ 
}\Delta \in \mathcal{B}(\Omega _{cl}),\text{ \ }y\in \Omega _{cl}. 
\tag{2.16}
\end{equation}%
Starting from the pioneering papers [LP76, Ro77, COPP78, BH-K82], such
unbounded spin systems have been under intensive investigation in classical
statistical mechanics (for recent developments see, e.g., [AKRT00, BH00,
Yo01]).\smallskip

\textbf{(ii)} \ Our definition of temperedness (as well as its modification
to the classical systems (2.14) with $|q_{k}|$ substituting $|\omega
_{k}|_{L_{\beta }^{2}}$) is less restrictive (and simpler) than those
usually used in the literature (for comparison, see e.g. [COPP78, BH-K82]).
So, obviously, $\Omega _{(e)t}\supseteq \Omega _{(s)t}$ resp. $\mathcal{M}%
_{(e)t}\supseteq \mathcal{M}_{(s)t},$ where the subsets of all (\emph{%
\textquotedblright slowly increasing\textquotedblright })\emph{\ tempered}
configurations resp.\emph{\ }measures are defined by 
\begin{gather}
\Omega _{(s)t}:=\left\{ \omega \in \Omega \text{ }\left\vert \text{ }\exists
p=p(\omega )>0:\text{ }||\omega ||_{-p}:=\left[ \sum\nolimits_{k\in \mathbb{Z%
}^{d}}(1+|k|)^{-2p}|\omega _{k}|_{L_{\beta }^{2}}^{2}\right] ^{\frac{1}{2}%
}<\infty \right. \right\} ,  \notag \\
\mathcal{M}_{(s)t}:=\left\{ \mu \in \mathcal{M}(\Omega )\ \left\vert \text{ }%
\exists p=p(\mu )>0:\text{ }||\omega ||_{-p}<\infty \text{ \ for }\mu -\text{%
a.e. }\omega \in \Omega \text{ }\right. \right\} .  \tag{2.17}
\end{gather}%
Moreover, $\mathcal{M}_{(s)t}$ contains all measures $\mu \in \mathcal{M}%
(\Omega )$ satisfying the following condition in terms of their moment
sequence: 
\begin{equation*}
\exists p=p(\mu )>0:\text{ }\sum\nolimits_{k\in \mathbb{Z}%
^{d}}(1+|k|)^{-2p}E_{\mu }|\omega _{k}|_{L_{\beta }^{2}}^{2}<\infty ,
\end{equation*}%
in particular, those having the so-called Ruelle support (see Remark 5.11
(ii) below). Here and further on we write%
\begin{equation*}
E_{\mu }f:=\int fd\mu
\end{equation*}%
for any $\mu $-integrable function $f.\smallskip $

\textbf{(iii)} \ If it does not lead to the reader's confusion, as a rule in
the concrete model setting we shall just use\emph{\ the standard notation} $%
\Omega _{t},$ $\mathcal{M}_{t}$ and $\mathcal{G}_{t}$ by omitting all the
additional sub- and superscripts (as, e.g., here $(e)$ or $(s)$).

\subsection{Formulation of the main results}

Now we are ready to present our results on the existence and a priori
estimates for Euclidean Gibbs measures corresponding to the QLS (2.2). We
assume that all the conditions on the interaction potentials imposed in
Subsect.\thinspace 2.1 are fulfiled without mentioning them again in the
formulations of our statements.\smallskip

\noindent \textbf{Main} \textbf{Theorem I} (Existence of Tempered Gibbs
States)\textbf{.} \ \emph{For all values of the mass} $\mathfrak{m}>0$ \emph{%
and the inverse temperature} $\beta >0:$%
\begin{equation*}
\mathcal{G}_{t}\neq \emptyset .
\end{equation*}%
\smallskip

\noindent \textbf{Main Theorem II} (A Priori Estimates on Tempered Gibbs
States)\textbf{. }\ \emph{Every} $\mu \in \mathcal{G}_{t}$ \emph{is
supported by the set of H\"{o}lder} \emph{loops\ }$\bigcap_{0\leq \alpha <%
\frac{1}{2}}[C_{\beta }^{\alpha }]^{\mathbb{Z}^{d}}.$ \emph{Moreover, for
all }$Q\geq 1$\emph{\ and\ }$\alpha \in \lbrack 0,\frac{1}{2})$%
\begin{equation}
\sup_{\mu \in \mathcal{G}_{t}}\sup_{k\in \mathbb{Z}^{d}}\int_{\Omega
}|\omega _{k}|_{C_{\beta }^{\alpha }}^{Q}d\mu (\omega )<\infty .  \tag{2.18}
\end{equation}%
\smallskip

\noindent \textbf{Corollary from Theorem II. } \emph{The set }$\mathcal{G}%
_{t}$ \emph{is compact w.r.t. the topology of weak convergence of measures
on any of the spaces }$[C_{\beta }^{\alpha }]^{\mathbb{Z}^{d}}$, $0\leq
\alpha <1/2,$\emph{\ equipped by the system of seminorms} $|\omega
_{k}|_{C_{\beta }^{\alpha }\text{ }},$ $k\in \mathbb{Z}^{d}.$\smallskip

Actually, the existence result for $\mu \in \mathcal{G}_{t}$ immediately
follows by Prokhorov's tightness criterion from the moment estimate (2.19)
for the family $\pi _{\Lambda }(d\omega |\xi =0),$ which holds uniformly in
volume $\Lambda \Subset \mathbb{Z}^{d}$. Moreover, the a priori estimates
for the probability kernels $\pi _{\Lambda }(d\omega |\xi )$ of the local
specification, subject to the fixed boundary condition $\xi \in \Omega _{t},$
stated in Theorem III below, are also of independent interest and have
various applications.\smallskip

\noindent \textbf{Main Theorem III} (Moment Estimates Uniformly in Volume)%
\textbf{. }\ \emph{Let} \emph{us} \emph{fix any boundary condition} $\xi \in
\Omega .$ \emph{Then} \emph{for all }$\delta >0,$ $\alpha \in \lbrack 0,%
\frac{1}{2})$ \emph{and} $Q\geq 1$ \emph{hold}$:$%
\begin{gather}
\sup_{k\in \mathbb{Z}^{d}}|\xi _{k}|_{C_{\beta }^{\alpha }}<\infty
\Longrightarrow \sup_{\Lambda \Subset \mathbb{Z}^{d}}\sup_{k\in \mathbb{Z}%
^{d}}\int_{\Omega }|\omega _{k}|_{C_{\beta }^{\alpha }}^{Q}\pi _{\Lambda
}(d\omega |\xi )=:C_{Q,\xi }<\infty ,  \tag{2.19} \\
\xi \in \Omega _{t}\Longrightarrow \sup_{\Lambda \Subset \mathbb{Z}%
^{d}}{\displaystyle\sum\limits_{k\in \Lambda }}e^{-\delta |k|}\int_{\Omega }|\omega
_{k}|_{C_{\beta }^{\alpha }}^{Q}\pi _{\Lambda }(d\omega |\xi )=:C_{Q,\xi
}^{\prime }<\infty .  \tag{2.20}
\end{gather}%
\smallskip

It should be mentioned that all the above statements extend to quantum
lattice systems with general many-particle interactions (cf.
Subsects.\thinspace 2.5, 2.6 and 3.1--3 below). The proof of Theorems I--III
will be given in Subsect.\thinspace 5.2. A detailed analysis of the
\textquotedblright \emph{state of} \emph{the art}\textquotedblright\ of the
problems dealt with in the literature and a comparison of our results with
the previous ones obtained by other authors will be performed in
Subsect.\thinspace 2.6.

\subsection{Model II: pair interaction of superquadratic growth}

Here we briefly discuss how to modify the previous setting in order to
include many-particle interaction potentials of \emph{superquadratic growth}%
. Namely, let us consider the following generalization of the QLS (2.2)
described by a heuristic Hamiltonian of the form%
\begin{equation}
\mathbb{H}:=-\frac{1}{2\mathfrak{m}}\sum_{k\in \mathbb{Z}^{d}}\frac{d^{2}}{%
dq_{k}^{2}}+\frac{a^{2}}{2}\sum_{k\in \mathbb{Z}^{d}}q_{k}^{2}+\sum_{k\in 
\mathbb{Z}^{d}}V(q_{k})+\sum\limits_{\langle k,k^{\prime }\rangle \subset 
\mathbb{Z}^{d}}W(q_{k},q_{k^{\prime }}).  \tag{2.21}
\end{equation}%
The one-particle potential $V\in C^{2}(\mathbb{R\rightarrow R)}$ satisfies
the same Assumption $(\mathbf{V}_{\mathbf{0}})$ as in Subsect.\thinspace
2.1, i.e., has asymptotic behaviour at infinity as a polynomial of order $%
P>2 $. Concerning the pair potential, we suppose that $W\in C^{2}(\mathbb{R}%
^{2}\rightarrow \mathbb{R)}$ has respectively at most polynomial growth of
the order $R<P$:\smallskip

\noindent \textbf{Assumption }$\mathbf{(W}_{\mathbf{0}}\mathbf{)}$\textbf{:}
\ \emph{There exist constants }$R\in \lbrack 2,P)$ \emph{and }$K_{W},C_{W}>0$
\emph{such that for all }$q,q^{\prime }\in \mathbb{R}$%
\begin{equation*}
|\partial _{q}^{(l)}W(q,q^{\prime })|\leq K_{W}\left( |q|^{R-l}+|q^{\prime
}|^{R-l}\right) +C_{W},\quad l=0,1,2.
\end{equation*}%
\smallskip

\noindent \textbf{Remark\thinspace 2.5.} \ A trivial example for pair
potentials which fit $\mathbf{(W}_{\mathbf{0}}\mathbf{)}$ are the
polynomials $W(q,q^{\prime }):=\sum_{l=0}^{2r}(q-q^{\prime })^{l}$ of even
degree $2r<P.$ In other words, our assumptions mean the so-called \emph{%
lattice} \emph{stabilization},\emph{\ }when the pair interaction is
dominated by the single-particle one. The case $P=R$ is also allowed, but it
needs a more accurate analysis which will be performed in Sect.\thinspace 3
below (in this respect see the general Assumptions $(\mathbf{V})$, $(\mathbf{%
W})$ on the interaction there).\smallskip

As compared with the initial QLS model (2.2), the only principal difference
in dealing with its generalization (2.21) is that we should proper change
the notion of temperedness. Now we define the subset of \emph{tempered}\ 
\emph{configurations} by%
\begin{equation}
\Omega _{(e)t}^{R}:=\left\{ \omega \in \Omega \text{ }\left\vert \text{ }%
\forall \delta \in (0,1):\text{ }||\omega ||_{-\delta ,R}:=\left[
\sum\nolimits_{k\in \mathbb{Z}^{d}}e^{-\delta |k|}|\omega _{k}|_{L_{\beta
}^{R}}^{2}\right] ^{\frac{1}{2}}<\infty \right. \right\} .  \tag{2.24}
\end{equation}%
It makes sense to consider the largest subset of such type by taking the
smallest possible value of the parameter $R\geq 2$ describing the polynomial
growth of the pair potential $W(q,q^{\prime })$ in Assumption $\mathbf{(W}_{%
\mathbf{0}}\mathbf{)}$. Note that for $R=2$ we just repeat the previous
definition (2.6), i.e., $\Omega _{(e)t}=\Omega _{(e)t}^{R=2}.$ Then all our
main statements presented in Subsect.\thinspace 2.3 \emph{remain true}, as
soon as in their formulation we substitute the single spin space $L_{\beta
}^{2}$ by $L_{\beta }^{R}$ and respectively specify the subset $\mathcal{G}%
_{t}:=\mathcal{G}_{(e)t}^{R}$ of\emph{\ tempered\ Gibbs measures} as those
supported by $\Omega _{t}:=\Omega _{(e)t}^{R}.$ Let us stress that (even in
the case of translation invariant interactions we now deal with) we cannot
guarantee that (outside the uniqueness regime) any tempered Gibbs measure
will be invariant w.r.t. lattice translations. So, the above set $\mathcal{G}%
_{(e)t}^{R}$ is in a certain sense the largest one so that for any of its
points $\mu $ we are technically able to get moment estimates like (2.18) 
\emph{uniformly} w.r.t. the lattice parameter $k\in \mathbb{Z}^{d}$.

\subsection{Model III: pair interaction of infinite range}

A further generalization of the QLS models (2.2) and (2.21) concerns the
case of\emph{\ not necessarily translation-invariant} pair interaction of
possibly \emph{infinite range}. Namely, let us consider the model described
by a heuristic Hamiltonian of the form%
\begin{equation}
\mathbb{H}:=-\frac{1}{2\mathfrak{m}}\sum_{k\in \mathbb{Z}^{d}}\frac{d^{2}}{%
dq_{k}^{2}}+\frac{a^{2}}{2}\sum_{k\in \mathbb{Z}^{d}}q_{k}^{2}+\sum_{k\in 
\mathbb{Z}^{d}}V_{k}(q_{k})+\sum\limits_{\{k,k^{\prime }\}\subset \mathbb{Z}%
^{d}}W_{\{k,k^{\prime }\}}(q_{k},q_{k^{\prime }}).  \tag{2.25}
\end{equation}%
The one-particle potentials $V_{k}\in C^{2}(\mathbb{R\rightarrow R)}$
satisfy the same Assumption $(\mathbf{V}_{\mathbf{0}})$ as before, but with
the constants $P>2$ and $K_{V},C_{V}>0$ which are \emph{uniform} for all $%
k\in \mathbb{Z}^{d}.$ The two-particle interactions (taken over all
unordered pairs $\{k,k^{\prime }\}\subset \mathbb{Z}^{d}$, $k\neq k^{\prime
} $) are given by symmetric functions $W_{\{k,k^{\prime }\}}\in C^{2}(%
\mathbb{R}^{2}\rightarrow \mathbb{R)}$ satisfying the following growth
condition:

\noindent \textbf{Assumption }$\mathbf{(W}_{\mathbf{0}}^{\ast }\mathbf{)}$%
\textbf{:} \ \emph{There exist some constants }$2\leq R<P$ \emph{and }$%
J_{k,j}\geq 0$ \emph{such that for all }$\{k,j\}\subset \mathbb{Z}^{d}$ 
\emph{and} $q,q^{\prime }\in \mathbb{R}$%
\begin{equation*}
|\partial _{q}^{(l)}W_{\{k,k^{\prime }\}}(q,q^{\prime })|\leq J_{k,k^{\prime
}}\left( 1+|q|^{R-l}+|q^{\prime }|^{R-l}\right) ,\quad l=0,1,2.
\end{equation*}

For the matrix $\mathbf{J}:=\{J_{k,k^{\prime }}\}$ we allow different rates
of decay (for instance, polynomial or exponential) when the distance $|k-j|$
between the points of the lattice gets large:

\noindent \textbf{Assumption }$\mathbf{(J}_{\mathbf{0}}\mathbf{)}$\textbf{:} 
$\mathbf{(i)}$ \ \emph{For all }$p\geq 0$ \emph{\ holds}%
\begin{equation*}
||\mathbf{J}||_{p}:=\sup_{k\in \mathbb{Z}^{d}}\left\{ \sum_{j\in \mathbb{Z}%
^{d}\backslash \{k\}}J_{k,j}(1+|k-j|)^{p}\right\} <\infty ,
\end{equation*}

\emph{or, even stronger,\smallskip }

$\mathbf{(ii)}$ \emph{There exist some} $\delta >0$\emph{\ such that} 
\begin{equation*}
||\mathbf{J}||_{\delta }:=\sup_{k\in \mathbb{Z}^{d}}\left\{ \sum_{j\in 
\mathbb{Z}^{d}\backslash \{k\}}J_{k,j}e^{\delta |k-j|}\right\} <\infty .
\end{equation*}%
\smallskip

Again, we first need to choose the proper notion of the temperedness, which
entirely depends on the decay rate of the pair interaction. A new issue
caused by the \emph{infinite range} of the interaction is that one has to
check (cf. Lemma 3.7 below) that the probability kernels $\pi _{\Lambda
}(d\omega |\xi )$ are well defined for all boundary conditions $\xi \in
\Omega _{t}$. So, in view of Assumption $\mathbf{(J}_{\mathbf{0}}\mathbf{)(i)%
}$\textbf{, }we define the subset $\Omega _{(s)t}^{R}\subset \Omega
_{(e)t}^{R}$ of (\emph{\textquotedblright slowly
increasing\textquotedblright }) \emph{tempered}\ configurations by%
\begin{equation}
\Omega _{(s)t}^{R}:=\left\{ \omega \in \Omega \text{ }\left\vert \text{ }%
\exists p=p(\omega )>0:\text{ }||\omega ||_{-p,R}:=\left[ \sum\nolimits_{k%
\in \mathbb{Z}^{d}}(1+|k|)^{-2p}|\omega _{k}|_{L_{\beta }^{R}}^{2}\right] ^{%
\frac{1}{2}}<\infty \right. \right\} .  \tag{2.24}
\end{equation}%
Respectively, we introduce the subset of tempered Gibbs measures 
\begin{equation*}
\mathcal{G}_{(s)t}^{R}:=\left\{ \mu \in \mathcal{G}\text{ }\left\vert \text{ 
}\exists p=p(\mu )>d:\text{ }||\omega ||_{-p,R}<\infty \text{ }\forall
\omega \in \Omega \text{ }(\func{mod}\mu )\right. \right\} .
\end{equation*}%
Then our main Theorems I and II about existence and a priori estimates for
the tempered Euclidean Gibbs measures \emph{remain true}, provided in their
formulation one substitutes the single spin spaces $L_{\beta }^{2}$ by $%
L_{\beta }^{R}$ and, respectively, $\Omega _{t}$ by $\Omega _{(s)t}^{R}$ and 
$\mathcal{G}_{t}$ by $\mathcal{G}_{(s)t}^{R}$. In the formulation of
Theorems III describing the properties of the probability kernels $\pi
_{\Lambda }(d\omega |\xi )$ one also needs obvious changes, which we shall
discuss later in Subsect.\thinspace 5.2.2. On the other hand, if we want to
deal with the larger subset $\mathcal{G}_{t}^{R}\supset \mathcal{G}%
_{(s)t}^{R}$ and completely keep the previous setup of the QLS\ Model II, we
should correspondingly impose the stronger Assumption $\mathbf{(J}_{\mathbf{0%
}}\mathbf{)(ii)}$ on the decay of matrix $\mathbf{J}$.

\subsection{Comments on Main Theorems \textbf{I}-\textbf{III}}

\subsubsection{Basic problems and known results}

>From a general viewpoint of the\emph{\ }theory of \emph{Markov} \emph{random
fields}, the following are fundamental problems in the study of Euclidean
Gibbs measures on loop lattices:\smallskip \newline
\textbf{I. Existence problem.} \ The initial step in any study of Gibbs
measures is to check whether the set $\mathcal{G}_{t}$ is nonempty. However,
as is typical for systems with noncompact (in our case, even
infinite-dimensional) spin spaces, the existence of $\mu \in \mathcal{G}_{t}$
stated by Theorem I is not evident at all. An important observation in this
respect is that, under natural assumptions on the interaction, any
accumulating point of the family $\pi _{\Lambda },$ $\Lambda \Subset \mathbb{%
Z}^{d},$ is certainly Gibbs. Depending on the specific class of quantum
lattice models one deals with, the required convergence $\pi _{\Lambda
^{(N)}}\rightarrow \mu ,$ $\Lambda ^{(N)}\nearrow \mathbb{Z}^{d},$ and thus
the existence of\emph{\ }$\mu \in \mathcal{G}_{t},$ can be proved by the
following main methods listed below:

\begin{description}
\item[(i)] \textbf{General Dobrushin's criterion for existence of Gibbs
distributions} [Do68,70]: The validity of the sufficient conditions of the 
\emph{Dobrushin existence theorem} for some classical unbounded spin systems
(including $P(\phi )$-lattice models, cf. Remark 2.1) has been verified e.g.
in [COPP78, BH-K82, Sin82] (however, under assumptions on the interaction
potentials more restrictive than $\mathbf{(V}_{0}\mathbf{)}$ and $\mathbf{(W}%
_{0}\mathbf{)}$). Contrary to the classical case, the same problem for
quantum lattice systems was not covered at all by any previous work and will
be treated in one of our subsequent papers.

\item[(ii)] \textbf{Ruelle's} \textbf{technique of superstability estimates }%
(cf. the original papers [Rue69, LP76] and resp. [PaY94,95] for its
extension to the quantum case\textbf{\ }): This technique otherwise requires
that the interaction is \emph{translation invariant} and (according to the
so-called regularity condition) the many-particle potentials have \emph{at
most quadratic growth} (i.e., ($\mathbf{W}_{\mathbf{0}}$) holds with $R=2$).
As was shown in [PaY94], for the subclass of boundary conditions $\xi \in
\Omega _{(ss)t}\subset \Omega _{(s)t}$ (for instance, such that $\sup_{k\in 
\mathbb{Z}^{d}}|\xi _{k}|_{L_{\beta }^{2}}<\infty $) the family of
probability kernels $\pi _{\Lambda }(d\omega |\xi ),$ $\Lambda \Subset 
\mathbb{Z}^{d},$ $\Lambda \nearrow \mathbb{Z}^{d},$ is tight (in the sense
of local weak convergence on $\Omega $) and has at least one accumulation
point $\mu $ from the subset of superstable Gibbs measures $\mathcal{G}%
_{(ss)t}\subset \mathcal{G}_{(s)t}$ (for the corresponding definitions see
Remark 5.14 (ii) below). This technique also shows that any $\mu \in 
\mathcal{G}_{(ss)t}$ is a priori of sub-Gaussian growth$.$

\item[(iii)] \textbf{Cluster expansions} is one of the most powerful methods
for the study of Gibbs fields, but it works only in a \emph{perturbative
regime}, i.e., when an effective parameter of\emph{\ }the\emph{\ }%
interaction is small. In particular, various versions of this technique
imply both existence and also uniqueness (but in some weaker than the DLR
sense) of the associated infinite volume Gibbs distributions\emph{\ }(see
e.g. the early works [AH-K75, GRS75, GJ81] and the recent developments in
[PaY95, AKMR98, FM99, MRZ00, MVZ00]).

\item[(iv)] \textbf{Method of correlations inequalities} involves more
detailed information about the structure of the interaction (for instance,
whether they are ferromagnetic or convex, see Remark 2.1). Starting from a
number of correlations inequalities (such as \emph{FKG, GKS, Lebowitz,
Brascamp-Lieb} etc.) commonly known for classical lattice systems, by a
lattice approximation technique (similar to that one used in Euclidean field
theory) one can extend them to the quantum case (cf. [AKK98, OS99,
AKKR01,02]).

\item[(v)] \textbf{Method of\ reflection positivity} (as a part of \textbf{%
(iv)}) applies to \emph{translation invariant }systems of type (2.23) with 
\emph{near-neighbours pair interactions }(i.e., when $V_{k}:=V,$ $%
W_{\{k,k^{\prime }\}}:=W$, and $W=0$ if $|k-k^{\prime }|>1$). For the
general description of the method and its applications to classical lattice
systems we refer the reader e.g. to [Shl86]. The proper modification of this
technique for the QLS (2.23) gives the existence of so-called \emph{periodic}
Gibbs states (for their definition see [AKKR02] as well as Remark 3.10 (iv)
below). Moreover, the reflection positivity method can also be used to study
phase transitions in such models with the double-well anharmonicity $V.$
This has been implemented under certain conditions (e.g., in the dimension $%
d\geq 3$ and for large enough $\beta ,\mathfrak{m}>>1$) in [DLP79, PK87,
BaK91, HM00].

\item[\textbf{(vi)}] \textbf{Method of stochastic dynamics} (also referred
to in quantum physics as \textquotedblright \textbf{stochastic quantization}%
\textquotedblright ; see e.g. [Fu91, DPZ96, AKRT01] and the related
bibliography therein): In this method the Gibbs measures are directly
treated as invariant (more precise, reversible) distributions for the
so-called Glauber or Langevin stochastic dynamics. However, some additional
technical assumptions are required on the interaction (among them \emph{at
most} \emph{quadratic growth} of the pair potentials $W_{\{k,k^{\prime
}\}}(q,q^{\prime })$) in order to ensure the solvability of the
corresponding stochastic evolution equations in infinite dimensions (not to
mention the extremely difficult ergodicity problem for them). This method
has been first applied in [AKRT01] to prove existence of Euclidean Gibbs
states for the particular QLS model (2.2) (see also the discussion in
Subsect.\thinspace 5.3 below).\medskip
\end{description}

\noindent \textbf{II. A priori estimates for measures in\ }$\mathcal{G}_{t}$%
\textbf{.} \ In turn, Theorem II above contributes to the fundamental
problem of getting \emph{uniform}\textbf{\ }\emph{estimates\ }on correlation
functionals of Gibbs measures in terms of parameters of the interaction.
This problem was initially posed for classical lattice systems in [COPP78,
BH-K82] and is closely related with the compactness of the set of tempered
Gibbs states (cf. Corollary~after Theorem II in Subsect.\thinspace 2.3); we
refer also to [AKRT00] for a detailed discussion of the classical lattice
case. There are very few results in the literature about a priori
integrability properties of tempered Gibbs measures on loop or path spaces
(see, for instance, [Iw85, Fu91, OS99] for the case of Euclidean $P(\phi
)_{1}$-fields and resp. [AKRT01] for the case of quantum anharmonic
crystals). All of them are based on the method of stochastic dynamics
already mentioned above. It is worth noting that the other methods (cf.
Items \textbf{I. (i)--(v)}) give also some estimates on limit points for $%
\pi _{\Lambda },$ $\Lambda \Subset \mathbb{Z}^{d},$ but not uniformly and
not for all $\mu \in \mathcal{G}_{t}.$ Besides, the finiteness of the
moments of the measures $\mu \in \mathcal{G}_{t}$ is also useful for the
study of Gibbs measures by means of the associated\emph{\ }Dirichlet
operators $\mathbb{H}_{\mu }$ in the spaces $L^{p}(\mu ),\ p\geq 1,$ (this
is known as the Holley--Stroock approach [HS76, AKR97a,b]). In particular,
by [AKR97a,b] $\mu $ is an \emph{extreme point} (or pure phase) in $\mathcal{%
G}_{t},$ if and only if the corresponding Markov semigroup $\exp (-t\mathbb{H%
}_{\mu }),$ $t\geq 0$, is \emph{ergodic} in $L^{2}(\mu )$ (which extends the
famous result in [HS76] related to the Ising model); see also a related
discussion in Subsect.\thinspace 5.3 below.\medskip

In this paper we shall be mainly focused on the first two problems described
above, namely existence and a priori estimates for $\mu \in \mathcal{G}_{t}.$
But, in order to give the reader a greater insight into the subject, we also
mention the next \emph{important directions}:\medskip

\noindent \textbf{III. Uniqueness problem.} \ The validity of the sufficient
conditions of the \emph{Dobrushin uniqueness criterion} (see [Do70, F\"{o}%
82]) for the QLS's (2.19) with \emph{pair interactions of at most quadratic
growth} has been first verified in [AKRT97a,b]. In doing so, the
coefficients of Dobrushin's matrix were estimated by means of log-Sobolev
inequalities proved on the single loop spaces $L_{\beta }^{2}$ and the
uniqueness of $\mu \in \mathcal{G}_{t}$ was established for small values of
the inverse temperature $\beta \in (0,\beta _{0}),$ but under conditions
independent of the particle mass $\mathfrak{m}>0$. For a \emph{special class}
of\emph{\ }ferromagnetic models with the polynomial self-interaction (2.4),
these results have been essentially improved in the recent series of papers
[AKKR01--03]. The latter papers propose a new technique which combines the
classical ideas of [LP76, BH-K82] based on the use of FKG and GKS
correlation inequalities with the spectral analysis of one-site oscillators
(2.1) specific for the quantum case. The strongest result of such type
obtained in [AKKR03] establishes the uniqueness of $\mu \in \mathcal{G}_{t}$
in the small-mass domain $\mathfrak{m}\in (0,\mathfrak{m}_{0})$ \emph{%
uniformly at all values of} $\beta >0.$ This provides a mathematical
justification for the well-known physical phenomenon that structural phase
transition for a given mass $\mathfrak{m}>0$ can be suppressed not only by
thermal fluctuations (i.e., high temperatures $\beta ^{-1}>\beta _{cr}^{-1}$%
), but for the light particles (with $\mathfrak{m}<\mathfrak{m}_{cr}$) also
by the quantum fluctuations (i.e., tunneling in a double-well potential)
simultaneously at all temperatures $\beta >0$.\smallskip

\noindent \textbf{IV. Phase transitions. \ }There are basically \emph{two
general methods} for proving existence of phase transitions (i.e.,
non-uniqueness of $\mu \in \mathcal{G}_{t}$) for \emph{low} \emph{%
temperatures} $\beta ^{-1}$, namely, the \emph{reflection positivity} (for $%
d\geq 3$) and the \emph{energy-entropy (Peierls-type) argument} (for $d\geq 2
$). However, in practice their successful applications to quantum lattice
systems have been limited so far to ferromagnetic $P(\phi )$-models (cf.
[DLP79, PK87, BaK91, HM00] resp. [GJS75, AKRe98]). The first method (already
mentioned in Item I.(v)) enables one to prove the positivity of \emph{a} 
\emph{long-range order parameter} $\lim_{|\Lambda |\rightarrow \infty
}E_{\mu _{per,\Lambda }}\left[ \sum_{k\in \Lambda }\omega _{k}(\tau )\right]
^{2}/|\Lambda |^{2}$ for large enough $\mathfrak{m}>\mathfrak{m}_{0}$ and $%
\beta >\beta _{0}(\mathfrak{m}_{0})$ via the so-called infrared (Gaussian)
bounds on two-point correlation functions $E_{\mu _{per,\Lambda }}\omega
_{k}(\tau )\omega _{k^{\prime }}(\tau )$ w.r.t. the local Gibbs measures $%
\mu _{per,\Lambda }$ with periodic boundary conditions. The second method
has originally been discovered (as the so-called Peierls argument) for the
Ising model and further developed to apply to rather general classical spin
systems (now known as the Pirogov--Sinai contour method), cf. [Sin82]. Its
quantum modification was firstly implemented in [GJS75] to the study of
phase transition in the $(\varphi ^{4})_{2}-$model of Euclidean field theory
and then in [AKRe98] to its lattice approximation (2.2). Following the idea
of that papers, one defines a \textquotedblright \emph{collective spin
variable}\textquotedblright\ $\sigma _{k}:=$sign$\int_{S_{\beta }}\omega
_{k}(\tau )d\tau $ taking values $\pm 1$ and a long-range parameter as the
correlation function $<\sigma _{k}\sigma _{k^{\prime }}>:=\lim_{|\Lambda
|\rightarrow \infty }E_{\mu _{\Lambda }^{per}}\sigma _{k}\sigma _{k^{\prime
}}.$ The existence of long-range behaviour, and hence phase transition,
follows from the estimate $<\sigma _{k}\sigma _{k^{\prime }}>\geq 1/2$ valid
for large enough values of $\mathfrak{m}$ and $\beta $.\medskip 

\noindent \textbf{V. Euclidean ground states. \ }Of special interest for
quantum systems is the case of \emph{zero absolute temperature}, i.e., $%
\beta =\infty $, which is technically more complicated and less studied in
the literature. In particular, it involves an important problem of the
operator realization of the formal Hamiltonian (2.2) in quantum mechanics
(cf. [AH-K75, BeK94]). The corresponding Gibbs measures $\mu \in \mathcal{G}%
_{gr}$ on the \textquotedblright path lattice\textquotedblright\ $[C(\mathbb{%
R})]^{\mathbb{Z}^{d}},$ so-called \emph{Euclidean ground states}, also allow
the DLR-description, but through a family of local specifications $\pi
_{I\times \Lambda }$ indexed by \textquotedblright
time-space\textquotedblright\ windows $I\times \Lambda $ with $I\Subset 
\mathbb{R},$ $\Lambda \Subset \mathbb{Z}^{d},$ cf. [MRZ00]. A principal
difference with the previous case $0<\beta <\infty $ is that now there is
not available any such (independent from boundary conditions $\xi $)
reference measure, so that all $\pi _{I\times \Lambda }(d\omega |\xi )$ are
defined as its Gibbs modifications. So far, there are very few rigorous
results about\textbf{\ }Gibbs measures\textbf{\ }on the path space $[C(%
\mathbb{R})]^{\mathbb{Z}^{d}}$, which all are mainly related to the
existence problem. A recent progress in this direction was achieved in the
series of papers [FM99, MRZ00, MVZ00], where the limit measures $%
\lim_{\Lambda \nearrow \mathbb{Z}^{d}}\lim_{I\nearrow \mathbb{R}}\pi
_{I\times \Lambda }\in \mathcal{G}_{gr}$ for the $P(\phi )-$lattice models
(2.2) have been constructed through cluster expansions w.r.t. the small mass
parameter $\mathfrak{m}<<1$. At the same time, for fixed $\Lambda \Subset 
\mathbb{Z}^{d},$ the corresponding unique Gibbs measures $\mu _{gr,\Lambda
}:=\lim_{I\nearrow \mathbb{R}}\pi _{I\times \Lambda }$ on the path space $[C(%
\mathbb{R})]^{\Lambda }$ are well-known as the $P(\phi )_{1}-$processes and
can be looked upon as a special case of Euclidean field theory in
space-dimension zero (cf. [Iw85, OS99]). Besides it should be noted that the
Gibbs measures on the path space $[C(\mathbb{R})]^{\mathbb{Z}^{d}}$ also
appear in a natural way as weak solutions for SDE's in $\mathbb{Z}^{d}$
[Deu87, MRZ00].

\subsubsection{Comparison with the previous results}

We emphasize that our Main Theorems I and II improve essentially all the
known results on the existence and a priori estimates for tempered Euclidean
Gibbs measures presented in Subsect. 2.6.1 above. As already stressed in the
Introduction, in order to prove our statements formulated in
Subsect.\thinspace 2.3, we shall propose \emph{a new technique}, which
completely differs from those listed under Items I.(i)--(v) and relies on
the alternative description of $\mu \in \mathcal{G}$ via integration by
parts. Then our main theorems will follow immediately from corresponding
results on symmetrizing measures in the abstract framework of
Sects.\thinspace 6, 7. Moreover, our technique obviously extends (cf.
Sect.\thinspace 3) to \emph{general many-particle interactions} (not
necessarily translation invariant and possibly having superquadratic growth,
unbounded order and infinite range), which were not covered at all by any
previous work. On the other hand, our approach is conceptually more
straightforward and technically easier in comparison to the stochastic
dynamics method mentioned under I.(vi) above.\smallskip

This alternative approach has been first realized in [AKRT99,00], however in
the much simpler situation of classical lattice systems like (2.14) with
finite dimensional spins. But the concrete technique suggested in those
papers does not apply to loop spaces, so that a proper (highly non-trivial)
modification for the quantum case is necessary. The reason is that in the
quantum case we have to do not only a \textquotedblright \emph{lattice} 
\emph{analysis}\textquotedblright\ (depending on the properties of the
interaction potentials $V$, $W$), but also an additional \textquotedblright 
\emph{single} \emph{spin space analysis}\textquotedblright\ (taking into
account the spectral properties of the elliptic operator $A_{\beta })$%
.\smallskip

It should also be mentioned, that in the recent preprint [Ha01] some
(deterministic) version of integration by parts for local specifications has
been used to prove existence of Gibbs measures relative to Brownian motion
on the path space $C(\mathbb{R}\rightarrow \mathbb{R}^{d})$. The study of
such Gibbsian (in general non Markovian) processes has been initiated in
[OS99]. As a special case they include the so-called $P(\phi )_{1}$%
-processes as Gibbs distributions corresponding to a single quantum particle
at zero temperature, i.e., $\beta =\infty $ (see e.g. [Iw85]). Finally, let
us notice that our method can be also modified to apply to the case of zero
absolute temperature, i.e., $\beta =\infty $, and corresponding Euclidean
ground states on the \textquotedblright path lattice\textquotedblright\ $[C(%
\mathbb{R})]^{\mathbb{Z}^{d}}$(cf. Item \textbf{IV} in Subsect. 2.6.1). This
case is under present investigation.

\section{A general model of quantum anharmonic crystals}

In this section we describe in detail the Euclidean Gibbsian setup for our
main model (3.1), which obviously includes all the particular models
introduced in Section\thinspace 2 and hence will be further referred to as
the \emph{general QLS model}.

\subsection{Assumptions on the interaction}

In the subsequent we shall consider the following system of quantum
anharmonic oscillators on $\mathbb{Z}^{d}$ with general (\emph{not
necessarily translation-invariant}) many-particle interaction (\emph{%
possibly having unbounded order and infinite range}), which is described by
the heuristic infinite-dimensional Hamiltonian 
\begin{equation}
\mathbb{H}:=-\frac{1}{2\mathfrak{m}}\sum_{k\in \mathbb{Z}^{d}}\frac{d^{2}}{%
dq_{k}^{2}}+\frac{a^{2}}{2}\sum_{k\in \mathbb{Z}^{d}}q_{k}^{2}\mathbb{+}%
\sum\limits_{k\in \mathbb{Z}^{d}}V_{k}(q_{k})+\sum_{M=2}^{N}\sum\limits_{%
\{k_{1},...,k_{M}\}\subset \mathbb{Z}^{d}}W_{\{k_{1},...,k_{M}%
\}}(q_{k_{1}},...,q_{k_{M}}).  \tag{3.1}
\end{equation}%
\bigskip

\noindent \textbf{Notation 3.1.}\ \ Below we shall distinguish between the
following notation: $(k_{1},...,k_{M})\in (\mathbb{Z}^{d})^{M}$ will stand
for \emph{ordered} \emph{sets (=} \emph{sequences) of length }$M,$ and $%
\{k_{1},...,k_{M}\}\subset \mathbb{Z}^{d}$ for \emph{unordered sets
consisting of }$M$ \emph{distinct} \emph{points}.\smallskip

\noindent \textbf{Definition 3.2. \ }\emph{We specify Assumptions }$(\mathbf{%
W}_{\mathbf{i-iii}}),$ $(\mathbf{J})$ and $(\mathbf{V}_{\mathbf{i-v}})$ 
\emph{on the system (3.1) as follows:}

\begin{enumerate}
\item[\textbf{(}$\mathbf{W}$\textbf{)}] The $M$-particle interaction
potentials (taken over all sets $\{k_{1},...,k_{M}\}\subset \mathbb{Z}^{d}$
with finite $M\in \{2,...,N\}$ and possibly infinite $N\in \mathbb{N\cup
\{+\infty \})}$ are given by twice continuously differentiable \emph{%
symmetric} functions $W_{\{k_{1},...,k_{M}\}}\in C^{2}(\mathbb{R}%
^{M}\rightarrow \mathbb{R})$ satisfying a \emph{polynomial growth}
condition. The latter means that there exist $R\geq 2$ and symmetric
matrices 
\begin{equation*}
\mathbf{J}_{M}=\{J_{k_{1},...,k_{M}}\}_{(k_{1},...,k_{M})\in \mathbb{Z}%
^{d^{M}}},\quad J_{k_{1},...,k_{M}}\geq 0,
\end{equation*}%
such that\ for all $M\leq N,$ $\{k_{1},...,k_{M}\}\subset \mathbb{Z}^{d}$
and $q_{1},...,q_{M}\in \mathbb{R}$%
\begin{gather}
|W_{\{k_{1},...,k_{M}\}}(q_{1},...,q_{M})|\leq
J_{k_{1},...,k_{M}}(1+\sum\limits_{l=1}^{M}|q_{l}|)^{R},  \tag{\QTR{bf}{i}}
\\
|\partial _{1}W_{\{k_{1},...,k_{M}\}}(q_{1},...,q_{M})|\leq
J_{k_{1},...,k_{M}}(1+\sum\limits_{l=1}^{M}|q_{l}|)^{R-1}, 
\tag{\QTR{bf}{ii}} \\
|\partial _{1}\partial _{l}W_{\{k_{1},...,k_{M}\}}(q_{1},...,q_{M})|\leq
J_{k_{1},...,k_{M}}(1+\sum\limits_{l=1}^{M}|q_{l}|)^{R-2}, 
\tag{\QTR{bf}{iii}}
\end{gather}%
where $\partial _{l}$ denotes derivative w.r.t. the variable $q_{l},$ $1\leq
l\leq M.\ $Without loss of generality, we put 
\begin{equation*}
J_{k_{1},...,k_{M}}=0\text{\quad if }k_{l_{1}}=k_{l_{2}}\text{ for some }%
1\leq l_{1}<l_{2}\leq M.
\end{equation*}
\end{enumerate}

\begin{description}
\item[\textbf{(}$\mathbf{J}$\textbf{)}] The matrices $\mathbf{J}%
_{M}=\{J_{k_{1},...,k_{M}}\}_{k_{1},...,k_{M}\in \mathbb{Z}^{d}},$ $%
M=2,...,N,$ are \emph{fastly decreasing}, that is for any $p\geq 0$%
\begin{equation*}
||\mathbf{J}||_{p}:=\sum_{M=2}^{N}M^{R}||\mathbf{J}_{M}||_{p}<\infty ,
\end{equation*}%
where we define the following seminorms on $\mathbb{R}^{(\mathbb{Z}%
^{d})^{M}} $%
\begin{gather*}
||\mathbf{J}_{M}||_{p}:=\sup_{k_{1}\in \mathbb{Z}^{d}}\left\{
\sum\limits_{\{k_{2},...,k_{M}\}\subset \mathbb{Z}^{d}}J_{k_{1},...,k_{M}}%
\left( 1+\sum\limits_{l=1}^{M}|k_{1}-k_{l}|\right) ^{p}\right\} , \\
\text{e.g.,\quad }||\mathbf{J}_{M}||_{0}:=\sup_{k_{1}\in \mathbb{Z}%
^{d}}\left\{ \sum\limits_{\{k_{2},...,k_{M}\}\subset \mathbb{Z}%
^{d}}J_{k_{1},...,k_{M}}\right\} .
\end{gather*}

\item[\textbf{(}$\mathbf{V}$\textbf{)}] The anharmonic self-interactions are
given by twice continuously differentiable functions $V_{k}\in C^{2}(\mathbb{%
R\rightarrow R})$ which satisfy the following \emph{coercivity }estimates%
\emph{\ }with \emph{some fixed} $K_{1},K_{2}>0,$ \emph{small enough} $%
K_{3}>0 $ (cf. Lemma 3.6 and Theorem$\,$7.6)$,$ \emph{arbitrarily small} $%
K_{4}>0,$ and corresponding $L_{1}(K_{1}),...,L_{4}(K_{4})>0$: 
\begin{equation*}
V_{k}^{\prime }(q)\cdot q\geq \max \left\{ 
\begin{array}{r}
K_{1}^{-1}(|V_{k}^{\prime }(q)|+|V_{k}^{\prime \prime }(q)|-L_{1}), \\ 
K_{2}^{-1}(|V_{k}^{\prime \prime }(q)\cdot q|-L_{2}), \\ 
K_{3}^{-1}(|q|^{R}-L_{3}), \\ 
K_{4}^{-1}(q^{2}-L_{4})%
\end{array}%
\right\} 
\begin{array}{r}
\qquad (\text{\textbf{i}}) \\ 
(\text{\textbf{ii}}) \\ 
(\text{\textbf{iii}}) \\ 
(\text{\textbf{iv}})%
\end{array}%
\end{equation*}%
and the \emph{growth} condition with some fixed $K_{0},L_{0}>0$%
\begin{equation*}
|V_{k}^{\prime \prime }(q)|\leq K_{0}(|V_{k}^{\prime
}(q)|+|q|^{R-1})+L_{0},\quad \quad \quad \quad \qquad \qquad (\text{\textbf{v%
}})
\end{equation*}%
\emph{uniformly} for all $k\in \mathbb{Z}^{d}$ and $q\in \mathbb{R}%
.\smallskip $
\end{description}

\noindent \textbf{Remark 3.3. (i) \ }Formally, no bound on the growth at
infinity of the one-particle potentials is\textbf{\ }directly\textbf{\ }%
imposed. But, it is easy to show that ($\mathbf{V}_{\text{\textbf{iii}}}$)
implies that $V_{k}$ growth strongly enough: for any $0<\sigma
<(K_{3}R)^{-1} $ there exists $C_{k}:=C_{k}(\sigma )\in \mathbb{R}$ such that%
\begin{equation}
V_{k}(q)\geq \sigma |q|^{R}+C_{k}\text{, \ \ }\forall q\in \mathbb{R}. 
\tag{3.2}
\end{equation}%
On the other hand, Assumption ($\mathbf{V}_{\text{\textbf{v}}}$) garantees
(by Gronwall's inequality) that with the necessity there exists $%
C_{k}^{\prime }\in \mathbb{(}0,\mathbb{\infty )}$ such that%
\begin{equation}
|V_{k}^{\prime }(q)|\leq C_{k}^{\prime }(1+|q|^{R})\exp K_{0}|q|\text{, \ }%
\forall q\in \mathbb{R}.  \tag{3.3}
\end{equation}%
Typical examples of $V_{k}$ satisfying Assumptions ($\mathbf{V}_{\text{%
\textbf{i--v}}}$) are linear exponential functionals or polynomials of even
degree, i.e., 
\begin{equation*}
e^{\lambda q}+e^{-\lambda q},\quad b_{2n}q^{2n}+...+b_{1}q+b_{0},\quad q\in 
\mathbb{R},
\end{equation*}%
with $\lambda \neq 0,$ $b_{2n}>0$ and $2n\geq R,$ $n\in \mathbb{N}$, as well
as their products and sums.\smallskip

\textbf{(ii)} \ Coercivity assumptions on potentials like ($\mathbf{V}_{%
\text{\textbf{iii}}}$) are standardly used in mathematical physics,
especially when one studies stability properties of dynamical systems (for
more concrete applications to the infinite dimensional SDE's see, e.g.,
[DPZ96]). If $R>2,$ then Assumption ($\mathbf{V}_{\text{\textbf{iv}}}$)
itself becomes superfluous as a trivial sequel of ($\mathbf{V}_{\text{%
\textbf{iii}}}$). Moreover, Assumptions ($\mathbf{V}_{\mathbf{i-}\text{%
\textbf{v}}}$) hold just with arbitrary small $K_{0},...,K_{4}>0$ for every
of the particular QLS models described in Sect.\thinspace 2.\smallskip

(\textbf{iii) \ }We say that the interaction in (3.1) is\emph{\ local,} if
it has \emph{bounded order} $M\in \mathbb{N}$ and \emph{finite range} $\rho
\in (0,\infty ).$ This means that all $W_{\{k_{1},...,k_{M}\}}=0$ whenever $%
M>N$ or $diam\{k_{1},...,k_{M}\}>\rho .$ On the other hand, the interaction
is said to be\emph{\ translation invariant }whenever $V_{k_{0}}=V$ and $%
W_{\{k_{1},...,k_{M}\}}=W_{\{k_{1}+k_{0},...,k_{M}+k_{0}\}}$ for all $%
\{k_{0},k_{1},...,k_{M}\}\Subset \mathbb{Z}^{d}$ with $2\leq M\leq N$%
.\smallskip

\textbf{(iv) \ }The simplest and exactly solvable case of the QLS (3.1) is
the so-called called \emph{harmonic systems} with $N=2$ and $%
W_{\{k,k^{\prime }\}}(q_{1},q_{2}):=a_{k,k^{\prime }}q_{1}q_{2}.$ Here $%
a_{k,k^{\prime }}$ are the elements of the \textquotedblright \emph{%
dynamical matrix}\textquotedblright\ $\mathbf{D}=(a_{k,k^{\prime
}})_{k,k^{\prime }\in \mathbb{Z}^{d}},$ which is usually supposed to be
symmetric, bounded and strictly positive in the Hilbert space $l^{2}(\mathbb{%
Z}^{d}).$ Obviously, this interaction satisfies the above Assumptions ($%
\mathbf{W}$) and ($\mathbf{J}$) as soon as $||\mathbf{D}||_{p}<\infty ,$\ $%
\forall p\geq 0.$ For a detailed study of the harmonic QLS see, e.g.,
[Ga02].\smallskip

\textbf{(v) \ }It would be worth to give here a nontrivial example of
potentials $W_{\{k_{1},...,k_{M}\}},$ $M\in \mathbb{N},$ satisfying
Assumptions ($\mathbf{W}$) and ($\mathbf{J}$) with $N=\infty .$ Let 
\begin{equation*}
W_{\{k_{1},...,k_{M}\}}:=C_{k_{1},...,k_{M}}\left(
1+\sum_{\{k_{l_{1}},k_{l_{2}}\}\subset
\{k_{1},...,k_{M}\}}|q_{l_{1}}-q_{l_{1}}|^{2}\right) ^{R/2}
\end{equation*}%
where 
\begin{gather*}
C_{k_{1},...,k_{M}}:=C_{M}^{-1}\cdot \exp \left\{ -2\sigma
\sum_{\{k_{l_{1}},k_{l_{2}}\}\subset
\{k_{1},...,k_{M}\}}|k_{l_{1}}-k_{l_{2}}|\right\} , \\
C_{M}:=M^{R+\sigma +2}R\sum_{n\in \mathbb{N}}(2n+1)^{dM}\exp (-\sigma
n)<\infty ,
\end{gather*}%
and $\sigma >0$ is arbitrary. Then even a very rough estimate for $%
J_{k_{1},...,k_{M}}:=M^{2}RC_{k_{1},...,k_{M}}$ is enough to show that%
\begin{gather*}
||\mathbf{J}_{M}||_{p}\leq C_{M}^{-1}\cdot \sum_{n\in \mathbb{N}%
}(2n+1)^{dM}(n+1)^{p}\exp (-2\sigma n)\leq C_{p}M^{-(R+\delta +1)} \\
\text{with \ }C_{p}:=\sum_{n\in \mathbb{N}}(n+1)^{p}\exp (-\sigma n)<\infty ,
\end{gather*}%
and thus 
\begin{equation*}
||\mathbf{J}||_{p}:=\sum_{M=2}^{N}M^{R}||\mathbf{J}_{M}||_{p}\leq
C_{p}\sum_{M=2}^{N}M^{-(\delta +1)}<\infty ,\quad \forall p\geq 0.
\end{equation*}

\subsection{Loop spaces and the support of Euclidean Gibbs measures}

Let us fix some finite value of the inverse temperature $\beta >0.$ In this
subsection we define the corresponding \emph{temperature loop lattices}
which describe configurations of the infinite volume system with Hamiltonian
(3.1). It should be noted that, in contrary to the classical lattice systems
(2.14), in the present situation even the single spin spaces are themselves
infinite dimensional (e.g., nonreflexive and nonsmooth Banach spaces) and
their topological features should be taken properly into account (cf.
Subsect.\thinspace 3.1). Although some notation and definitions has been
already used in Sect.\thinspace 2 before, for the convenience of the reader
we recall them in the context of the general model (3.1).

\subsubsection{One-particle loop spaces}

Let $S_{\beta }$ be a circle of length $\beta ,$ which will be considered as
a compact Riemannian manifold with Lebesgue measure $d\tau $ as a volume
element and distance%
\begin{equation*}
\rho (\tau ,\tau ^{\prime }):=\min \{|\tau -\tau ^{\prime }|,\text{ }\beta
-|\tau -\tau ^{\prime }|\},\quad \tau ,\tau ^{\prime }\in S_{\beta }.
\end{equation*}%
We define some standard spaces of functions (i.e., loops) $\upsilon
:S_{\beta }\rightarrow \mathbb{R}.$ Namely, let $C_{\beta }^{m+\alpha
}:=C^{m+\alpha }(S_{\beta }),$ $m\in \mathbb{N}\cup \{0\},$ $0\leq \alpha
<1, $ denote the Banach space of all continuous loops on $S_{\beta }$ whose $%
m$-th derivative is $\alpha $-H\"{o}lder continuous, which is endowed with
the finite norm 
\begin{equation*}
|\upsilon |_{C_{\beta }^{m+\alpha }}:=\sup\limits_{\tau \in S_{\beta
}^{{}}}\sum\limits_{n=0}^{m}|\upsilon ^{(n)}(\tau )|+\sup\limits_{\tau ,\tau
^{\prime }\in S_{\beta }^{{}},\tau \neq \tau ^{\prime }}\frac{|\upsilon
^{(m)}(\tau )-\upsilon ^{(m)}(\tau ^{\prime })|}{|\tau -\tau ^{\prime
}|^{\alpha }}.
\end{equation*}%
Let $L_{\beta }^{r}:=L^{R}(S_{\beta })$ resp. $W_{\beta
}^{r,q}:=W^{r,q}(S_{\beta }),$ $r\geq 1,$ $q\in \mathbb{R},$ be the Lebesgue
resp. Sobolev spaces (with $L_{\beta }^{r}:=W_{\beta }^{r,q=0})$ relative to
the measure $d\tau .$ These spaces can be viewed as completion of $C_{\beta
}^{\infty }$ for the norms 
\begin{equation*}
|\upsilon |_{L_{\beta }^{r}}:=\left[ \int_{S_{\beta }^{{}}}|\upsilon (\tau
)|^{r}d\tau \right] ^{1/r}\text{\quad resp.\quad }|\upsilon |_{W_{\beta
}^{r,q}}:=|(-d^{2}/d\tau ^{2}+\mathbf{1)}^{q/2}\upsilon |_{L_{\beta }^{R}}.
\end{equation*}%
Below we will use the following well-known embeddings: 
\begin{gather}
W_{\beta }^{2,q}\underrightarrow{\subset }_{compact}\ L_{\beta }^{r},\quad q>%
\frac{1}{2}-\frac{1}{r},  \notag \\
W_{\beta }^{2,1}\underrightarrow{\subset }_{compact}\ C_{\beta }^{\alpha
^{\prime }}\underrightarrow{\subset }_{compact}\ C_{\beta }^{\alpha },\quad
0\leq \alpha <\alpha ^{\prime }<\frac{1}{2}.  \tag{3.4}
\end{gather}%
However, the reader is warned that for all positive noninteger numbers $%
\alpha <\alpha ^{\prime }$ the H\"{o}lder spaces $C_{\beta }^{\alpha
^{\prime }}$ in (3.4) are not separable and the embeddings $C_{\beta
}^{\alpha ^{\prime }}\underrightarrow{\subset }C_{\beta }^{\alpha }$ are not
dense.

At every site $k\in \mathbb{Z}^{d}$ of the lattice as the \emph{single spin
space} resp. \emph{tangent space} we define the Banach space $C_{\beta
}:=C_{\beta }^{\alpha =0}$ of all real-valued \emph{continuous loops} on $%
S_{\beta }$ with the sup-norm $|\cdot |_{C_{\beta }}$ resp.\ the Hilbert
space $H:=L_{\beta }^{2}$ of all \emph{square} \emph{integrable loops} with
the $L^{2}$-norm $|\cdot |_{H}=(\cdot ,\cdot )_{H}^{1/2}$. For the
corresponding Borel $\sigma $-algebras we have $C_{\beta }\in \mathcal{B}%
(L_{\beta }^{r})$ and $\mathcal{B}(C_{\beta })=\mathcal{B}(L_{\beta
}^{r})\cap C_{\beta }$. This follows from Kuratowski's theorem (cf. [Pa67,
p.\thinspace 21, Theorem 3.9]), that says 
\begin{subequations}
\begin{gather}
X_{2}\in \mathcal{B}(X_{1})\text{ \emph{and} }\mathcal{B}(X_{2})=\mathcal{B}%
(X_{1})\cap X_{2}\text{ \emph{for all Polish spaces} }X_{1}\text{\emph{\ and}
}X_{2}  \notag \\
\text{\emph{for which there is a measurable embedding} }X_{2}\ 
\underrightarrow{\subset }\ X_{1}.  \tag{3.5}
\end{gather}%
On the other hand, $C_{\beta }^{\alpha }\in \mathcal{B}(C_{\beta })$ for all 
$\alpha \geq 0,$ which can be easily proved by showing the measurability of $%
|\cdot |_{C_{\beta }^{\alpha }}$ (cf. [RS75, Subsect.\thinspace X.11]

\subsubsection{Spaces of sequences over $\mathbb{Z}^{d}$}

In order to describe the behaviour of the system (3.1) when $|k|\rightarrow
\infty ,$ we introduce some spaces of real-valued functions (i.e.,
sequences) on a lattice.

As usual, $\mathbb{R}^{\mathbb{Z}^{d}}$ (resp. $\mathbb{R}_{0}^{\mathbb{Z}%
^{d}}$) stands for the set of all (resp. its subset of finite) sequences
over $\mathbb{Z}^{d}.$ By $l^{q}(\gamma ),$ $q\geq 1,$ we denote a Banach
space of all sequences which are summable w.r.t. the given weight $\gamma
=\{\gamma _{k}\}_{k\in \mathbb{Z}^{d}}$, $\gamma _{k}>0$: 
\end{subequations}
\begin{equation}
l^{q}(\gamma ):=l^{q}(\mathbb{Z}^{d};\gamma ):=\left\{ x=(x_{k})_{k\in 
\mathbb{Z}^{d}}\in \mathbb{R}^{\mathbb{Z}^{d}}\left\vert |x|_{l^{q}(\gamma
)}:=\left[ \sum_{k\in \mathbb{Z}^{d}}\gamma _{k}x_{k}^{q}\right]
^{1/q}<\infty \right. \right\} .  \tag{3.6}
\end{equation}%
We shall mainly use the following two systems of weights $\gamma
_{p}=\{\gamma _{p,k}\}_{k\in \mathbb{Z}^{d}}$ resp. $\gamma _{\delta
}=\{\gamma _{\delta ,k}\}_{k\in \mathbb{Z}^{d}}$ indexed by $p\in \mathbb{R}$
resp. $\delta \in \mathbb{R}$:%
\begin{equation}
\gamma _{p,k}:=(1+|k|)^{p},\text{ \ }\gamma _{\delta ,k}:=\exp (\delta |k|),%
\text{ \ \ }k\in \mathbb{Z}^{d}.  \tag{3.7}
\end{equation}%
\smallskip

In particular, $l^{2}:=l^{2}(\mathbb{Z}^{d};\gamma \equiv 1)$ will be the
standard space of square summable sequences over $\mathbb{Z}^{d}$ with the
inner product $(\cdot ,\cdot )_{l^{2}}:=(\cdot ,\cdot )_{0}$ and the natural
orthonormal basis $e_{k}:=(\delta _{k-j})_{j\in \mathbb{Z}^{d}}$, $k\in 
\mathbb{Z}^{d}.$ For convenience, in every Hilbert space $l^{2}(\gamma )$ we
fix the orthonormal basis $e_{\gamma ,k}:=\gamma _{k}^{-1/2}e_{k},$ $k\in 
\mathbb{Z}^{d},$ so that%
\begin{equation*}
|x|_{l^{2}(\gamma )}^{2}=\sum_{k\in \mathbb{Z}^{d}}\gamma
_{k}^{-1}x_{k}^{2}=\sum_{k\in \mathbb{Z}^{d}}(x,e_{\gamma ,k})_{0}^{2},\quad
x\in l^{2}(\gamma ).
\end{equation*}%
Choosing $l^{2}(\mathbb{Z}^{d})$ as the tangent Hilbert space, we next
define its rigging%
\begin{equation*}
\mathcal{E}(\mathbb{Z}^{d})\subset \mathcal{S}(\mathbb{Z}^{d})\subset l^{2}(%
\mathbb{Z}^{d})\subset \mathcal{S}^{\prime }(\mathbb{Z}^{d})\subset \mathcal{%
E}^{\prime }(\mathbb{Z}^{d})
\end{equation*}%
by the following pairs of mutually dual nuclear spaces$:$%
\begin{gather}
\mathcal{S}(\mathbb{Z}^{d}):=\text{pr}\lim\limits_{p=1,2,...}l^{2}(\gamma
_{p})\text{, \ \ }\mathcal{S}^{\prime }(\mathbb{Z}^{d}):=\text{ind}%
\lim\limits_{p=1,2,...}l^{2}(\gamma _{-p}),  \notag \\
\mathcal{E}(\mathbb{Z}^{d}):=\text{ind}\lim\limits_{\delta >0}l^{2}(\gamma
_{\delta })\text{, \ \ }\mathcal{E}^{\prime }(\mathbb{Z}^{d}):=\text{pr}%
\lim\limits_{\delta >0}l^{2}(\gamma _{-\delta })\text{.}  \tag{3.8}
\end{gather}%
In particular, $\mathcal{S}(\mathbb{Z}^{d})$ and $\mathcal{S}^{\prime }(%
\mathbb{Z}^{d})$ are well known as the Schwartz spaces of fastly decreasing
resp. slowly increasing sequences over $\mathbb{Z}^{d}$. Obviously, one has
the Hilbert--Schmidt embeddings for\ $p^{\prime }>p+d$ resp. $\delta
^{\prime }>\delta $:%
\begin{gather}
O_{p}^{p^{\prime }}:l^{2}(\gamma _{p^{\prime }})\rightarrow l^{2}(\gamma
_{p}),\text{ \ \ }||O_{p}^{p^{\prime }}||_{HS}=\left[ \sum\nolimits_{k\in 
\mathbb{Z}^{d}}(1+|k|)^{-(p^{\prime }-p)}\right] ^{1/2}<\infty ,  \notag \\
O_{\delta }^{\delta ^{\prime }}:l^{2}(\gamma _{\delta ^{\prime
}})\rightarrow l^{2}(\gamma _{\delta }),\text{ \ \ \ }||O_{\delta }^{\delta
^{\prime }}||_{HS}=\left[ \sum\nolimits_{k\in \mathbb{Z}^{d}}\exp (-(\delta
^{\prime }-\delta )|k|)\right] ^{1/2}<\infty .  \tag{3.9}
\end{gather}%
Moreover, 
\begin{equation*}
\mathcal{S}^{\prime }(\mathbb{Z}^{d}):=\bigcup_{p\geq 1}l^{2}(\gamma
_{-p})=\bigcup_{p\geq 1}l^{1}(\gamma _{-p}),\text{ \ }\mathcal{E}^{\prime }(%
\mathbb{Z}^{d}):=\bigcap_{\delta >0}l^{2}(\gamma _{-\delta
})=\bigcap_{\delta >0}l^{1}(\gamma _{-\delta }),
\end{equation*}%
with the equivalence of the corresponding systems of weighted $l^{2}$ and $%
l^{1}$-norms:%
\begin{equation}
|x|_{l^{2}(\gamma _{-2p})}\leq |x|_{l^{1}(\gamma _{-p})}\leq
||O_{0}^{p}||_{HS}|x|_{l^{2}(\gamma _{-p})},\text{ \ }|x|_{l^{2}(\gamma
_{-2\delta })}\leq |x|_{l^{1}(\gamma _{-\delta })}\leq ||O_{0}^{\delta
}||_{HS}|x|_{l^{2}(\gamma _{-\delta })}.  \tag{3.10}
\end{equation}%
Let us note that $\Omega _{cl}:=\mathbb{R}^{\mathbb{Z}^{d}}$ is the
configuration space for the classical lattice systems like (2.14), whereas
its subspaces $\mathcal{S}^{\prime }(\mathbb{Z}^{d})$ and\ $\mathcal{E}%
^{\prime }(\mathbb{Z}^{d})$ are commonly used to describe the support
properties of tempered Gibbs measures $\mu \in \mathcal{G}_{cl}$ (see e.g.
[AKRT00]).

\subsubsection{Loop lattices}

Repeating the definition (2.5) above, as the \emph{configuration space} for
the infinite volume system (3.1) we introduce the \emph{"temperature loop\
lattice"}%
\begin{equation}
\Omega :=[C_{\beta }]^{\mathbb{Z}^{d}}=\{\omega =(\omega _{k})_{k\in \mathbb{%
Z}^{d}}\ |\ \omega :S_{\beta }\rightarrow \mathbb{R}^{\mathbb{Z}^{d}},\text{ 
}\omega _{k}\in C_{\beta }\}.  \tag{3.11}
\end{equation}%
We consider $\Omega $ as a Polish (i.e., complete separable metrizable)
space with the \emph{product topology}, i.e., the weakest topology such that
all finite volume projections%
\begin{equation*}
\Omega \ni \omega \mapsto \mathbb{P}_{\Lambda }\omega :=\omega _{\Lambda
}:=(\omega _{k})_{k\in \Lambda }\in \lbrack C_{\beta }]^{\Lambda }=:\Omega
_{\Lambda },\quad \Lambda \Subset \mathbb{Z}^{d},
\end{equation*}%
are continuous. This topology is generated by the system of seminorms $%
|\omega _{k}|_{C_{\beta }},$ $k\in \mathbb{Z}^{d},$ or, what is the same, by
any one of the (mutually equivalent) metrics%
\begin{equation*}
\rho _{-p}(\omega ,\omega ^{\prime }):=\left[ \sum_{k\in \mathbb{Z}%
^{d}}(1+|k|)^{-p}\left( \frac{|\omega _{k}-\omega _{k}^{\prime }|_{C_{\beta
}}^{2}}{1+|\omega _{k}-\omega _{k}^{\prime }|_{C_{\beta }}^{2}}\right) %
\right] ^{1/2},\text{ \ }p>d.
\end{equation*}%
Next, we provide $\Omega $ with the corresponding\emph{\ Borel\ }$\sigma $%
\emph{-}algebra $\mathcal{B}(\Omega ),$ i.e., the smallest $\sigma $-algebra
containing all open sets. It is a well-known fact for product spaces (cf.
e.g., [CTV85, Proposition 1.4]) that $\mathcal{B}(\Omega )$ also coincides
with the $\sigma $-algebra generated by cylinder sets 
\begin{equation*}
\left\{ \omega \in \Omega \text{ }\left\vert \text{ }\omega _{\Lambda }\in
\Delta _{\Lambda }\right. \right\} ,\text{ \ }\Delta _{\Lambda }\in \mathcal{%
B}(\Omega _{\Lambda }),\text{ }\Lambda \Subset \mathbb{Z}^{d}\text{.}
\end{equation*}%
By $\mathcal{M}(\Omega )$ we, as usual, shall mean the family of all
probability measures on $(\Omega ,\mathcal{B}(\Omega )).$ Depending on the
questions under discussion, $\mathcal{M}(\Omega )$ itself will be endowed
with the \emph{topologies} either\emph{\ of weak }or\emph{\ }of\emph{\
locally weak convergence}, which turn it in a Polish space. It is worth
recalling that a sequence $\{\mu _{n}\}_{n\in \mathbb{N}}\subset \mathcal{M}%
(\Omega )$ \emph{weakly} (resp. \emph{locally weakly}) \emph{converges} to $%
\mu \in \mathcal{M}(\Omega )$ iff $E_{\mu _{n}}f\rightarrow E_{\mu }f,$ as $%
n\rightarrow \infty ,$ for all bounded continuos (and resp. local, i.e.,
such that $f=f(\mathbb{P}_{\Lambda _{0}})$ with some $\Lambda _{0}\Subset 
\mathbb{Z}^{d}$) functions $f:\Omega \rightarrow \mathbb{R}.\smallskip $

Because of the possibly infinite radius of interaction and unboundedness of
interaction potentials, we will also need some \emph{subspaces of tempered
configurations }(cf. the related discussion in Subsect.\thinspace 2.6
above).\smallskip

In general, if $X$ is one of the single loop spaces from Subsect.\thinspace
3.2.1, then $l^{q}(\gamma ;X):=l^{q}(\mathbb{Z}^{d}\rightarrow X;\gamma )$
will stand for the corresponding space of weighted sequences over $\mathbb{Z}%
^{d}$ with values in $X.$ In particular, for $\gamma :=\gamma _{p}$ given by
(3.7), we shall use the following scales of Banach spaces%
\begin{equation}
\mathcal{L}_{-p}^{r}:=l^{2}(\mathbb{Z}^{d}\rightarrow L_{\beta }^{r};\gamma
_{-p})=\left\{ \omega \in \lbrack L_{\beta }^{r}]^{\mathbb{Z}^{d}}\left\vert 
\text{ }||\omega ||_{-p,r}^{2}:=\sum_{k\in \mathbb{Z}^{d}}(1+|k|)^{-p}|%
\omega _{k}|_{L_{\beta }^{r}}^{2}<\infty \right. \right\} ,  \tag{3.12}
\end{equation}%
\begin{equation}
\mathcal{C}_{-p}^{\alpha }:=l^{2}(\mathbb{Z}^{d}\rightarrow C_{\beta
}^{\alpha };\gamma _{-p})=\left\{ \omega \in \Omega \left\vert \text{ }%
||\omega ||_{-p,\alpha }^{2}:=\sum_{k\in \mathbb{Z}^{d}}(1+|k|)^{-p}|\omega
_{k}|_{C_{\beta }^{\alpha }}^{2}<\infty \right. \right\} ,  \tag{3.13}
\end{equation}%
indexed by $p\in \mathbb{R},$ $r\geq 1$ and $\alpha \in \lbrack 0,1/2)$. As
a special case of (3.12) for $H:=L_{\beta }^{2}$ and $p=0$, the Hilbert
space 
\begin{equation}
\mathcal{H}:=l^{2}(\mathbb{Z}^{d}\rightarrow L_{\beta }^{2})=\left\{ \omega
\in (L_{\beta }^{2})^{\mathbb{Z}^{d}}\left\vert <\omega ,\omega >_{\mathcal{H%
}}=||\omega ||_{\mathcal{H}}^{2}:=\sum_{k\in \mathbb{Z}^{d}}|\omega
_{k}|_{H}^{2}<\infty \right. \right\}  \tag{3.14}
\end{equation}%
will be treated as the \emph{tangent space }to $\Omega .$

Next we define the \emph{support spaces} of Euclidean Gibbs measures 
\begin{equation}
\Omega _{-p}^{r}:=\Omega \cap \mathcal{L}_{-p}^{r},\quad p>d,\text{ }r\geq 1,
\tag{3.15}
\end{equation}%
as locally convex Polish spaces with the topology generated by the system of
seminorms $||\omega ||_{-p,r}$ and $|\omega _{k}|_{C_{\beta }},$ $k\in $ $%
\mathbb{Z}^{d},$ or equivalently, by the metric 
\begin{equation*}
\rho _{-p,r}(\omega ,\omega ^{\prime }):=\left[ \sum_{k\in \mathbb{Z}%
^{d}}(1+|k|)^{-p}\left( |\omega _{k}-\omega _{k}^{\prime }|_{L_{\beta
}^{r}}^{2}+\frac{|\omega _{k}-\omega _{k}^{\prime }|_{C_{\beta }}^{2}}{%
1+|\omega _{k}-\omega _{k}^{\prime }|_{C_{\beta }}^{2}}\right) \right]
^{1/2}.
\end{equation*}%
Respectively, for every $k\in \mathbb{Z}^{d},$%
\begin{equation}
\Omega _{-p;k}^{r}:=\left\{ \omega \in \mathcal{L}_{-p}^{r}\left\vert
\,\omega _{k}\in C_{\beta }\right. \right\}  \tag{3.16}
\end{equation}%
will be a Banach space with the norm%
\begin{equation*}
||\omega ||_{-p,r;k}:=||\omega ||_{-p,r}+|\omega _{k}|_{C_{\beta }}.
\end{equation*}%
Again, all the above spaces are equipped with their Borel $\sigma $-algebras
(coinciding, due to (3.5), with the $\sigma $-algebras induced on them by $%
\mathcal{B}(\Omega )).\smallskip $

In Subsect.\thinspace 3.3.3 we shall start with the most extended definition
of tempered Gibbs measures $\mu \in \mathcal{G}_{(s)t}^{R}$ for the QLS
(3.1) as those supported on the subset of \emph{\textquotedblright slowly
increasing\textquotedblright }\ configurations%
\begin{equation}
\Omega _{(s)t}^{R}:=\bigcup\limits_{p\geq 1}\Omega _{-p}^{R}=\left\{ \omega
\in \Omega \text{ }\left\vert \text{ }\exists p=p(\omega )>0:\text{ }%
||\omega ||_{-p,R}<\infty \right. \right\} .  \tag{3.17}
\end{equation}%
Here the parameter $R\geq 2$ describes a possible order of polynomial growth
of the interaction potentials $W_{\{k_{1},...,k_{M}\}}$ according to
Assumption $(\mathbf{W})$. But actually, as will be shown in the proof of
our Main Theorem II (cf. Corollary 5.6 (ii) below), the initial assumption 
\begin{equation*}
\mu (\Omega _{-p_{0}}^{R})=1\text{ for some }p_{0}>d
\end{equation*}%
implies further regularity and temperedness properties for any $\mu \in 
\mathcal{G}_{(s)t}^{R},$ namely, that 
\begin{equation*}
\mu (\mathcal{C}_{-p}^{\alpha })=1\text{ for all }p>d\text{ and }0\leq
\alpha <1/2.
\end{equation*}%
$\smallskip $

We recall that the above definition of temperedness (3.17) has been already
used in the context of the particular QLS model (2.23) with the pair
interaction of infinite range (cf. Subsect.\thinspace 2.5). To get a more
precise setting for the models (2.2) and (2.21) with the nearest neighbour
interaction (cf. Subsects.\thinspace 2.2--2.4), we should modify
(3.15)--(3.17) by taking instead of $\mathcal{L}_{-p}^{R}$ another weighted
spaces 
\begin{equation}
\mathcal{L}_{-\delta }^{R}:=l^{2}(\mathbb{Z}^{d}\rightarrow L_{\beta
}^{R};\gamma _{-\delta })=\left\{ \omega \in \lbrack L_{\beta }^{R}]^{%
\mathbb{Z}^{d}}\left\vert \text{ }||\omega ||_{-\delta ,r}^{2}:=\sum_{k\in 
\mathbb{Z}^{d}}\exp (-\delta |k|)|\omega _{k}|_{L_{\beta }^{R}}^{2}<\infty
\right. \right\} .  \tag{3.18}
\end{equation}%
Then respectively%
\begin{equation}
\Omega _{-\delta }^{R}:=\Omega \cap \mathcal{L}_{-\delta }^{R},\quad \delta
>0,\text{ }r\geq 1,  \tag{3.19}
\end{equation}%
and (as was already defined by (2.17) and (2.22))%
\begin{equation}
\Omega _{(e)t}^{R}:=\bigcap\limits_{\delta >0}\Omega _{-\delta }^{R}=\left\{
\omega \in \Omega \text{ }\left\vert \text{ }||\omega ||_{-\delta ,R}<\infty
,\text{ \ }\forall \delta >0\right. \right\} .  \tag{3.20}
\end{equation}%
Actually, $\Omega _{(e)t}^{R}\supset \Omega _{(s)t}^{R}$ is the largest
support set for $\mu \in \mathcal{G}_{t}$ considered so far (cf. Remark
2.4.(ii)). \smallskip

If there is no confusion, we shall shorten the above notation just by
writing $\Omega _{t},$ $\mathcal{M}_{t}$ or $\mathcal{G}_{t}$.

\subsection{DLR approach to Euclidean Gibbs measures}

According to [AH-K75], the Euclidean measure $\mu $ corresponding to a Gibbs
state $G_{\beta }$ of the quantum lattice system with Hamiltonian (3.1) has
the \emph{heuristic representation}%
\begin{equation}
d\mu (\omega ):=\frac{1}{Z}\exp \left\{ -\mathcal{I}(\omega )\right\}
\prod\limits_{k\in \mathbb{Z}^{d}}d\gamma _{\beta }(\omega _{k}).  \tag{3.21}
\end{equation}%
Here $Z$ is the normalization factor, $\gamma _{\beta }$ is the canonical
distribution on $C_{\beta }$ of the oscillator bridge process of length $%
\beta ,$ and%
\begin{equation}
\mathcal{I}(\omega ):=\int_{S_{\beta }}\left[ \sum\limits_{M=2}^{N}\sum%
\limits_{\{k_{1},...,k_{M}\}\subset \mathbb{Z}^{d}}W_{\{k_{1},...,k_{M}\}}(%
\omega _{k_{1}}(\tau ),...,\omega _{k_{m}}(\tau ))+\sum\limits_{k\in \mathbb{%
Z}^{d}}V_{k}(\omega _{k}(\tau )\right] d\tau  \notag
\end{equation}%
is the formal Euclidean action functional of the system in the infinite
volume $\mathbb{Z}^{d}.$ In full analogy with classical statistical
mechanics, a rigorous meaning can be given to the measure $\mu $ by the 
\emph{DLR formalism} as a Gibbs distribution on the lattice $\mathbb{Z}^{d}$%
, but with the infinite-dimensional single spin (= loop) space $C_{\beta }$.
So, in Subsects.\thinspace 3.3.1--2 we shall firstly study the corresponding
family of local specifications $\pi _{\Lambda },\ \Lambda \Subset \mathbb{Z}%
^{d}$. Thereafter, in Subsect.\thinspace 3.3.3 we define the set $\mathcal{G}
$ of all Gibbs measures $\mu $ on the loop lattice $\Omega $ (as the
solutions to the DLR equations $\mu \,\pi _{\Lambda }=\mu ,\ \forall \Lambda
\Subset \mathbb{Z}^{d})$ and its subset $\mathcal{G}_{t}:=\mathcal{G}%
_{(s)t}^{R}$ of \emph{tempered} (as those with physical relevance) Gibbs
measures supported by $\Omega _{t}:=\Omega _{(s)t}^{R}$.

\subsubsection{One-particle Euclidean measures}

Let us assume that a potential $V\in C_{b,loc}^{1}(\mathbb{R})$ satisfies
the coercivity condition ($\mathbf{V}_{\text{\textbf{iii}}}$) in Definition
3.2. Our aim here is to identify the following one-particle distribution 
\begin{equation}
d\sigma _{\beta }(\upsilon ):=\frac{1}{Z_{k}}\exp \left\{ -\int_{S_{\beta
}^{{}}}V(\upsilon (\tau ))d\tau \right\} d\gamma _{\beta }(\upsilon ) 
\tag{3.22}
\end{equation}%
\smallskip as a probability measure on the loop space $C_{\beta }$ and
collect some its properties to be used later on.\smallskip

We begin with the detailed construction of the basic Gaussian measure $%
\gamma _{\beta }$, which corresponds to a single harmonic oscillator of the
mass $\mathfrak{m}>0$ and rigidity $a>0$. So, in the Hilbert space $%
H:=L_{\beta }^{2}$ let us consider a positive self-adjoint operator (which
describes the quantum character of the system) $A_{\beta }:=-\mathfrak{m}%
\Delta _{\beta }+a^{2}\mathbf{1}$ with the domain $\mathcal{D}$($A_{\beta
}):=W_{\beta }^{2,2}$. Here $\Delta _{\beta }$ is the usual
Laplace--Beltrami operator on the circle $S_{\beta }$ (= maximal extension
in $H$ of the differential expression $d^{2}\varphi /d\tau ^{2}$, $\varphi
\in $ $C_{\beta }^{\infty }$). Let us recall that the operator $A_{\beta }$
has discrete spectrum%
\begin{equation}
\lambda _{n}=\left( \frac{2\pi }{\beta }n\right) ^{2}\mathfrak{m}+a^{2},%
\text{\quad }n\in \mathbb{Z}\text{,}  \tag{3.23}
\end{equation}%
and a complete orthonormal system of trigonometric functions%
\begin{equation}
\varphi _{n}(\tau )=\left\{ 
\begin{array}{cc}
\left( \frac{1}{\beta }\right) ^{\frac{1}{2}}, & n=0\text{ }, \\ 
\left( \frac{2}{\beta }\right) ^{\frac{1}{2}}\cos \frac{2\pi }{\beta }n\tau ,
& n=1,2,...\text{ }, \\ 
-\left( \frac{2}{\beta }\right) ^{\frac{1}{2}}\sin \frac{2\pi }{\beta }n\tau
, & n=-1,-2,...\text{ }.%
\end{array}%
\right.  \tag{3.24}
\end{equation}%
Obviously, 
\begin{equation}
\sup_{n\in \mathbb{Z}}|\varphi _{n}|_{L_{\beta }^{r}}:=\kappa _{r}^{{}}:=2^{%
\frac{1}{2}}\beta ^{\frac{1}{r}-\frac{1}{2}},\text{ \ \ }1\leq r\leq \infty ,%
\text{ \ \ (i.e., }\kappa _{\infty }^{{}}:=(2\beta ^{-1})^{\frac{1}{2}}. 
\tag{3.25}
\end{equation}%
For any $f\in L_{\beta }^{1}$ we define its Fourier and Ces\`{a}ro partial
sums in the standard way (cf. [Ed82]) by 
\begin{equation}
\mathbb{S}_{K}(f):=\sum\limits_{|n|\leq K}(f,\varphi _{n})_{\beta
}^{{}}\varphi _{n}\,,\quad \mathbb{M}_{K}(f):=\frac{1}{K+1}%
\sum\limits_{L=0}^{K}\mathbb{S}_{L}(f),\quad K\in \mathbb{N}\cup \{0\}. 
\tag{3.26}
\end{equation}%
Since, by the Riesz theorem, $\{\varphi _{n}\}_{n\in \mathbb{Z}}$ is a
Schauder basis in every $L_{\beta }^{r},\ 1<r<\infty ,$ there exists a
finite constant $\varsigma _{r}>0$ such that $\forall f\in L_{\beta }^{r}$%
\begin{equation}
\sup_{K\in \mathbb{N}}\left\vert \mathbb{S}_{K}(f)\right\vert _{L_{\beta
}^{r}}\leq \varsigma _{r}|f|_{L_{\beta }^{r}}\text{\quad and\quad }%
\lim_{K\rightarrow \infty }\left\vert f-\mathbb{S}_{K}(f)\right\vert
_{L_{\beta }^{r}}=0.  \tag{3.27}
\end{equation}%
On the other hand, by the Fej\'{e}r theorem, $\forall f\in C_{\beta }$%
\begin{equation}
\sup_{K\in \mathbb{N}}\left\vert \mathbb{M}_{K}(f)\right\vert _{C_{\beta
}}\leq |f|_{C_{\beta }}\text{\quad and\quad }\lim_{K\rightarrow \infty
}\left\vert f-\mathbb{M}_{K}(f)\right\vert _{C_{\beta }}=0.  \tag{3.28}
\end{equation}%
>From (3.24) it follows, in particular, that the set of all trigonometric
polynomials 
\begin{equation*}
T_{\beta }:=lin\{\varphi _{n}\}_{n\in \mathbb{Z}}\subset C_{\beta }^{\infty }
\end{equation*}%
is a domain of essential self-adjointness for $A_{\beta }^{{}}.$ As is well
known, the corresponding semigroup $P_{t}=\exp (-tA_{\beta }^{{}}),$ $t\geq
0,$ is Markovian in all spaces $C_{\beta },$ $L_{\beta }^{r},$ $r\geq 1,$
and one has $(t,\tau )\longmapsto (P_{t}f)(\tau )\in C^{\infty }((0,\infty
)\times S_{\beta }^{{}})$ for any $f\in L_{\beta }^{r}.\smallskip $

Furthermore, from (3.24) the resolvent $A_{\beta }^{-1}$ is of trace class
in $L_{\beta }^{2}$ (moreover, $Tr_{H}A_{\beta }^{\alpha -1}<\infty ,$ $%
\forall \alpha <1/2$). The corresponding\emph{\ Green function} (i.e.,
integral kernel of $A_{\beta }^{-1})$ 
\begin{equation*}
\mathfrak{G}(\tau ,\tau ^{\prime }):=\mathfrak{G}_{\tau }^{{}}(\tau ^{\prime
}):=(A_{\beta }^{-1}\delta _{\tau }^{{}})(\tau ^{\prime })
\end{equation*}%
is given by 
\begin{multline}
\mathfrak{G}(\tau ,\tau ^{\prime })=\sum_{n\in \mathbb{Z}}\lambda
_{n}^{-1}\varphi _{n}(\tau )\varphi _{n}(\tau ^{\prime })=\frac{1}{\beta
a^{2}}+\frac{2}{\beta }\sum_{n\in \mathbb{N}}^{\infty }\frac{\cos 2\pi
n\beta ^{-1}(\tau -\tau ^{\prime })}{\left( 2\pi n\beta ^{-1}\right) ^{2}%
\mathfrak{m}+a^{2}}  \notag \\
=\frac{\mathfrak{g}}{2}\left( e^{-\frac{a}{\sqrt{\mathfrak{m}}}(\beta -|\tau
-\tau ^{\prime }|_{\mathbb{R}}^{{}})}+e^{-\frac{a}{\sqrt{\mathfrak{m}}}|\tau
-\tau ^{\prime }|_{\mathbb{R}}^{{}}}\right) ,\quad \tau ,\tau ^{\prime }\in
S_{\beta }^{{}}.  \tag{3.29}
\end{multline}%
Here, for the sake of convenience, we introduce a parameter (which will be
relevant, e.g., in Subsect. 5.2.1)%
\begin{equation}
\mathfrak{g}:=\left[ 2a\sqrt{\mathfrak{m}}\left( 1-e^{-\frac{a}{\sqrt{%
\mathfrak{m}}}\beta }\right) \right] ^{-1}.  \tag{3.30}
\end{equation}%
Then obviously, 
\begin{equation*}
Tr_{H}A_{\beta }^{-1}=\int_{S_{\beta }}\mathfrak{G}(\tau ,\tau )d\tau =\beta 
\mathfrak{G}_{0}^{{}}(0)\leq \frac{1}{a^{2}}+\frac{\beta ^{2}}{\mathfrak{m}}.
\end{equation*}%
Also we note the following regularity properties of $\mathfrak{G}_{\tau }\in 
\mathcal{D}(A_{\beta }^{\frac{1}{2}}):=W_{\beta }^{2,1}$ resulting from
(3.29): 
\begin{equation}
|\mathfrak{G}_{\tau }|_{C_{\beta }}=\mathfrak{G}_{0}(0)\leq \mathfrak{g}%
,\quad |\mathfrak{G}_{\tau }-\mathfrak{G}_{\tau ^{\prime }}|_{C_{\beta
}}\leq \mathfrak{g}\frac{a}{\sqrt{\mathfrak{m}}}\rho (\tau ,\tau ^{\prime
}),\quad \forall \tau ,\tau ^{\prime }\in S_{\beta }.  \tag{3.31}
\end{equation}

Let now $\gamma _{\beta }$ be the Gaussian measure on $(H,$ $\mathcal{B(}H))$
with zero mean value and correlation operator $A_{\beta }^{-1}$, which is
given uniquely by its Fourier transform%
\begin{equation*}
\int_{H}\exp i(\varphi ,\upsilon )_{H}d\gamma _{\beta }(\upsilon )=\exp
\left\{ -\frac{1}{2}(A_{\beta }^{-1}\varphi ,\varphi )_{H}\right\} ,\text{%
\quad }\varphi \in H.
\end{equation*}%
Actually, the set $C_{\beta }^{\alpha },$ $0\leq \alpha <1/2,$ of H\"{o}lder
continuous loops has full measure, i.e., $\gamma _{\beta }^{{}}(C_{\beta
}^{\alpha })=1$, and the measure $\gamma _{\beta }^{{}}$ on the space $%
(C_{\beta }^{{}},\mathcal{B}(C_{\beta }^{{}}))$ can be viewed as the
canonical realization of the well-known oscillator bridge process of length $%
\beta $ (see, e.g., [Sim79, p.43]). Hence, the abstract Fernique's theorem
(see, e.g. [DeuS89, Theorem 1.3.24]) applied to $\gamma _{\beta }$
immediately gives us that for any $\alpha \in \lbrack 0,1/2)$ there exists $%
\lambda _{0}:=\lambda _{0}(\alpha )>0$ such that 
\begin{equation}
\int_{C_{\beta }}\exp \left\{ \lambda |\upsilon |_{C_{\beta }^{\alpha
}}^{2}\right\} d\gamma _{\beta }(\upsilon )<\infty ,\ \ \forall \lambda \in
\lbrack 0,\lambda _{0}].  \tag{3.32}
\end{equation}%
\smallskip

Now, we can define by the Feynman--Kac formula (3.22) the probability
measure $\sigma _{\beta }$ on $(C_{\beta }^{{}},\mathcal{B}(C_{\beta
}^{{}})) $. By Remark 3.3.(i), the coercivity Assumption ($\mathbf{V}_{\text{%
\textbf{iii}}}$) implies the polynomial growth of $V$ at the infinity and,
hence, its semiboundedness from below, i.e., $\inf_{\mathbb{R}}V\geq
C>-\infty .$ Therefore 
\begin{equation*}
0<Z:=\int_{C_{\beta }}\exp \left\{ -\int_{S_{\beta }}V(\upsilon (\tau
))d\tau \right\} d\gamma _{\beta }^{{}}(\upsilon )<\exp (-\beta C)<\infty
\end{equation*}%
and the definition (3.22) makes sense. Moreover, the measure $\sigma _{\beta
}$ has all moments of the form%
\begin{equation}
\int_{C_{\beta }^{{}}}\exp \left\{ \lambda |\upsilon |_{L_{\beta
}^{R}}^{R}\right\} \left( 1+|\upsilon |_{C_{\beta }^{\alpha }}^{Q}\right)
d\sigma _{\beta }(\upsilon )<\infty \text{, \ }0\leq \lambda <(K_{3}R)^{-1},%
\text{ }0\leq \alpha <\frac{1}{2},\text{ }Q\geq 1,  \tag{3.33}
\end{equation}%
due to (3.2) and the corresponding property (3.32) of the Gaussian measure $%
\gamma _{\beta }$.\smallskip

An important in the subsequent observation is that $\sigma _{\beta }$ is
quasi-invariant w.r.t. shifts 
\begin{equation*}
\upsilon \rightarrow \upsilon +\theta \varphi _{n},\text{\quad }\theta \in 
\mathbb{R},\text{ }n\in \mathbb{Z},
\end{equation*}%
with the \emph{Radon--Nikodym derivatives}%
\begin{equation}
\frac{d\sigma _{\beta }(\upsilon +\theta \varphi _{n})}{d\sigma _{\beta
}(\upsilon )}=:a_{\theta \varphi _{n}}(\upsilon )=a_{\theta \varphi
_{n}}^{A_{\beta }}(\upsilon )a_{\theta \varphi _{n}}^{V}(\upsilon ),\text{ \ 
}  \tag{3.34}
\end{equation}%
where 
\begin{gather}
a_{\theta \varphi _{n}}^{A_{\beta }}(\upsilon ):=\exp \left\{ -\theta
(A_{\beta }\varphi _{n},\upsilon )_{H}-\frac{\theta ^{2}}{2}(A_{\beta
}\varphi _{n},\varphi _{n})_{H}\right\} ,  \notag \\
a_{\theta \varphi _{n}}^{V}(\upsilon ):=\exp \left\{ -\int_{S_{\beta
}^{{}}}[V(\upsilon (\tau )+\theta \varphi _{n}(\tau ))-V(\upsilon (\tau
)]d\tau \right\} .  \tag{3.35}
\end{gather}

If, moreover, $V\in C_{b,loc}^{2}(\mathbb{R}),$ then for all $\theta \in 
\mathbb{R}$ exists%
\begin{equation}
\frac{d}{d\theta }a_{\theta \varphi _{n}}(\upsilon )=a_{\theta \varphi
_{n}}(\upsilon )\left[ (A_{\beta }\varphi _{n},\upsilon +\theta \varphi
_{n})_{H}+(V^{\prime }(\upsilon +\theta \varphi _{n}),\varphi _{n})_{H}%
\right]  \tag{3.36}
\end{equation}%
and the functions%
\begin{gather}
(\upsilon ,\theta )\longmapsto a_{\theta \varphi _{n}}^{{}}(\upsilon ),\ 
\frac{d}{d\theta }a_{\theta \varphi _{n}}(\upsilon )  \notag \\
\text{\emph{are uniformly Lipschitz-continuous }}\emph{on\ balls\ in\ }%
C_{\beta }\times \mathbb{R}.  \tag{3.37}
\end{gather}%
As a special case of (3.36) when $\theta =0$, we define the so-called \emph{%
partial logarithmic derivatives} $b_{\varphi _{n}}:C_{\beta }\rightarrow 
\mathbb{R}$ of the measure $\sigma _{\beta }$ along the basis directions $%
\varphi _{n},$ i.e.,%
\begin{equation}
b_{\varphi _{n}}(\upsilon ):=\left( a_{\theta \varphi _{n}}(\upsilon
)\right) _{\theta =0}^{\prime }=(A_{\beta }\varphi _{n},\upsilon
)_{H}+(V^{\prime }(\upsilon ),\varphi _{n})_{H}.  \tag{3.38}
\end{equation}%
Note that by construction always $a_{\theta \varphi _{n}}\in L^{1}(\sigma
_{\beta }).$ In contrast, the global integrability properties of $b_{\varphi
_{n}}$ are \emph{not known} in advance. In order to ensure them, we need
additional assumptions on the asymptotic behaviour of potential, for
instance, that $V^{\prime }$ satisfies the polynomial growth condition 
\begin{equation*}
|V^{\prime }(q)|\leq C^{\prime }(1+|q|^{R}),\text{ \ }q\in \mathbb{R}.
\end{equation*}%
Then from the moment estimate (3.33) one can easily conclude that for small
enough $\theta _{0}\in (0,\infty )$%
\begin{equation}
\sup_{0\leq \theta \leq \theta _{0}}\left\vert (d/d\theta )a_{\theta \varphi
_{n}}\right\vert \in L^{Q}(\sigma _{\beta }),\text{ and hence also }%
b_{\varphi _{n}}\in L^{Q}(\sigma _{\beta }).  \tag{3.39}
\end{equation}%
If (3.39) holds at least for $Q=1,$ then (3.36) in turn implies (by
Lebesgue's theorem applied to (4.23) below) the following \emph{integration
by parts formula}%
\begin{equation}
\int\nolimits_{C_{\beta }}(d/d\theta )f(\upsilon +\theta \varphi
_{n})d\sigma _{\beta }(\upsilon )=-\int\nolimits_{C_{\beta }}f(\upsilon
)b_{\varphi _{n}}(\upsilon )d\sigma _{\beta }(\upsilon )  \tag{3.40}
\end{equation}%
valid on all smooth cylinder functions $f(\upsilon ):=f_{L}((\upsilon
,\varphi _{1})_{H},...,(\upsilon ,\varphi _{L})_{H})$ with $f_{L}\in
C_{b}^{1}(\mathbb{R}^{L})$ and $L\in \mathbb{N}$. In this case one says that
the measure $\sigma _{\beta }$ is \emph{differentiable}\ (e.g. in the
well-known sense of Fomin and of Skorohod, cf. [Bo97]) along vectors $%
\varphi _{n}$ with the logarithmic derivatives $b_{\varphi _{n}}\in
L^{1}(\sigma _{\beta }).$

\subsubsection{Local specification}

The specification $\pi =\{\pi _{\Lambda }\}_{\Lambda \Subset \mathbb{Z}^{d}}$
(corresponding to the Hamiltonian (3.1)) is defined as a family of measure
kernels 
\begin{equation*}
\mathcal{B}(\Omega )\times \Omega \ni (\Delta ,\xi )\rightarrow \pi
_{\Lambda }(\Delta |\xi )\in \lbrack 0,1]
\end{equation*}%
in the following way:%
\begin{equation}
\pi _{\Lambda }(\Delta |\xi ):=\left\{ 
\begin{array}{ll}
Z_{\Lambda }^{-1}(\xi ){\displaystyle \int_{\Omega _{\Lambda }}}\exp \left\{ -\mathcal{I}%
_{\Lambda }^{W}(\omega |\xi )\right\} \mathbf{1}_{\Delta }(\omega _{\Lambda
},\xi _{\Lambda ^{c}})\prod\limits_{k\in \Lambda }d\sigma _{\beta ,k}(\omega
_{k}), & \xi \in \Omega _{t} \\ 
0, & \xi \in \Omega \backslash \Omega _{t}.%
\end{array}%
\right.  \tag{3.41}
\end{equation}%
(where $\mathbf{1}_{\Delta }$ denotes the indicator on $\Delta $). Here $%
\sigma _{\beta ,k}$ are the one-particle Euclidean measures on $C_{\beta }$
given by (3.22) with $V:=V_{k},$ 
\begin{equation}
Z_{\Lambda }(\xi ):=\int_{\Omega _{\Lambda }^{{}}}\exp \left\{ -\mathcal{I}%
_{\Lambda }^{W}(\omega |\xi )\right\} \prod\limits_{k\in \Lambda }d\sigma
_{\beta }(\omega _{k})  \tag{3.42}
\end{equation}%
is the normalization factor (the so-called partition function), and 
\begin{equation}
\mathcal{I}_{\Lambda }^{W}(\omega |\xi
):=\sum\limits_{M=2}^{N}\sum\limits_{L=1}^{M}\sum\limits_{\substack{ %
\{k_{1},...,k_{L}\}\subset \Lambda  \\ \{k_{L+1},...,k_{M}\}\subset \Lambda
^{c}}}\int_{S_{\beta }}W_{\{k_{1},...,k_{M}\}}(\omega _{k_{1}},...,\omega
_{k_{L}},\xi _{k_{L+1}},...,\xi _{k_{M}})d\tau  \tag{3.43}
\end{equation}%
is the many-particle interaction in the volume $\Lambda $ under the boundary
condition $\xi \in \Omega _{t}:=\Omega _{(s)t}^{R}$ (by convention $%
\{k_{L+1},...,k_{M}\}=\varnothing $ if $L\geq M$.)\smallskip

More precisely, $\pi $ is a $\left\{ \mathcal{B}(\Omega _{\Lambda
^{c}})\right\} _{\Lambda \Subset \mathbb{Z}^{d}}\mathcal{-}$\emph{%
specification} in the sense of [Pr76, Ge88], that means the following
properties hold for all $\Lambda \subset \Lambda ^{\prime }\Subset \mathbb{Z}%
^{d}$:

\begin{description}
\item[($\mathbf{S}_{\text{\textbf{i}}}$)] $\pi _{\Lambda }^{{}}(\Omega |\xi
) $ \emph{is either} $0$ \emph{or} $1$\emph{\ for all} $\xi \in \Omega ;$

\item[($\mathbf{S}_{\text{\textbf{ii}}}$)] $\pi _{\Lambda }^{{}}(\Delta
|\cdot )$\emph{\ is} $\mathcal{B}(\Omega _{\Lambda ^{c}}^{{}})-$\emph{%
measurable for all} $\Delta \in \mathcal{B}(\Omega );$

\item[($\mathbf{S}_{\text{\textbf{iii}}}$)] $\int_{\Omega }f(\omega
)g(\omega )\pi _{\Lambda }(d\omega |\xi )=f(\xi )\int_{\Omega }g(\omega )\pi
_{\Lambda }(d\omega |\xi )$ \emph{for all bounded and} $\mathcal{B}(\Omega
_{\Lambda ^{c}}^{{}})$ (\emph{resp}. $\mathcal{B}(\Omega )$) -- \emph{%
measurable functions} $f$ (\emph{resp}. $g$) on $\Omega ;$

\item[($\mathbf{S}_{\text{\textbf{iv}}}$)] \emph{Consistency:} $\pi
_{\Lambda ^{\prime }}^{{}}=\pi _{\Lambda ^{\prime }}^{{}}\pi _{\Lambda
}^{{}} $ \emph{where the kernel} $\pi _{\Lambda ^{\prime }}^{{}}\pi
_{\Lambda }^{{}} $ \emph{is defined by} 
\begin{equation*}
(\pi _{\Lambda ^{\prime }}^{{}}\pi _{\Lambda }^{{}})(\Delta |\xi
):=\int_{\Omega }\pi _{\Lambda ^{\prime }}(d\omega |\xi )\pi _{\Lambda
}(\Delta |\omega ),\text{\quad }\forall \text{ }\xi \in \Omega ,\text{ }%
\forall \Delta \in \mathcal{B}(\Omega ).
\end{equation*}
\end{description}

Actually, by (3.41) every $\pi _{\Lambda }(d\omega |\xi )$\ is concentrated
on configurations of the form $\omega =(\omega _{\Lambda },\xi _{\Lambda
^{c}})\in \Omega _{t}$ whenever $\xi \in \Omega _{t}$. For each $\Lambda
\Subset \mathbb{Z}^{d}$ and $\xi \in \Omega _{t},$ it is also reasonable to
consider the so-called \emph{Gibbs distribution} \emph{in the volume} $%
\Lambda $ \emph{with the boundary condition} $\xi _{\Lambda ^{c}}$ as the
probability measure 
\begin{equation}
\mu _{\Lambda }(d\omega _{\Lambda }^{{}}|\xi _{\Lambda ^{c}}):=Z_{\Lambda
}^{-1}(\xi )\exp \left\{ -\mathcal{I}_{\Lambda }^{W}(\omega |\xi )\right\}
\prod\limits_{k\in \Lambda }d\sigma _{\beta ,k}(\omega _{k})  \tag{3.44}
\end{equation}%
on the loop space $\Omega _{\Lambda }^{{}}=[C_{\beta }^{{}}]^{\Lambda }$.
Equivalently, $\mu _{\Lambda }^{{}}(d\omega _{\Lambda }^{{}}|\xi _{\Lambda
^{c}})$ is the projection of $\pi _{\Lambda }^{{}}(d\omega |\xi )$ onto $%
\Omega _{\Lambda }^{{}},$ i.e., $\mu _{\Lambda }^{{}}(\cdot |\xi _{\Lambda
^{c}})=\pi _{\Lambda }^{{}}(\cdot |\xi )\circ \mathbb{P}_{\Lambda
}^{-1}.\smallskip $

Because of the possibly infinite range of interaction, some verification
(see Lemmas 3.4 and 3.6 below) is needed of whether definitions (3.41) and
(3.44) make sense, i.e., axioms ($\mathbf{S}_{\mathbf{i-iii}}$) hold. The
validity of the consistency axiom ($\mathbf{S}_{\mathbf{iv}}$) for our model
follows then from the additive structure of the functional $\mathcal{I}%
_{\Lambda }^{W}$ (cf. [Ge88]).

Before proceeding further, here we collect a few technical estimates on the
matrix $\mathbf{J}$ to be repeatedly used below:\smallskip \newline
\textbf{Lemma 3.4. (i)} \ \emph{Let us define a symmetric matrix }$\mathbf{%
\tilde{J}}=\{\tilde{J}_{k,j}\}_{k,j\in \mathbb{Z}^{d}}$\emph{\ with the
nonnegative entries} 
\begin{equation}
\tilde{J}_{k,j}:=\sum_{M=2}^{N}M^{R}\sum\limits_{\substack{ %
\{k_{1},...,k_{M}\}\subset \mathbb{Z}^{d}  \\ k_{1}:=k,\text{ }k_{2}:=j}}%
J_{k_{1},...,k_{M}}\quad (\text{\emph{i.e., }}\tilde{J}_{k,j}:=0\text{ \ 
\emph{for } }k=j).  \tag{3.45}
\end{equation}%
\emph{Then} $\forall p\geq 0$%
\begin{equation}
||\mathbf{\tilde{J}}||_{p}:=\sup_{k\in \mathbb{Z}^{d}}\left\{
\sum\limits_{j\in \mathbb{Z}^{d}}\tilde{J}_{k,j}\left( 1+|k-j|\right)
^{p}\right\} \leq \sum_{M=2}^{N}M^{R}||\mathbf{J}_{M}||_{p}=||\mathbf{J}%
||_{p}<\infty .  \tag{3.46}
\end{equation}

\textbf{(ii)} \ \emph{Let us introduce one more system of seminorms} \emph{on%
} $\mathbb{R}^{(\mathbb{Z}^{d})^{M}}$\emph{\ by} 
\begin{equation}
|||\mathbf{J}_{M}|||_{p}:=\sup_{k_{1}\in \mathbb{Z}^{d}}\left\{
\sum\limits_{\{k_{2},...,k_{M}\}\subset \mathbb{Z}^{d}}J_{k_{1},...,k_{M}}%
\sum\limits_{l=1}^{M}(1+|k_{1}-k_{l}|)^{p}\right\} ,\emph{\quad }p\geq 0%
\emph{.}  \tag{3.47}
\end{equation}%
\emph{Then} $||\mathbf{J}_{M}||_{p}$\emph{\ and }$|||\mathbf{J}_{M}|||_{p}$ 
\emph{are equivalent, i.e.,} 
\begin{equation}
|||\mathbf{J}_{M}|||_{0}=M||\mathbf{J}_{M}||_{0}\quad \text{\emph{and}\quad }%
\frac{1}{M}|||\mathbf{J}_{M}|||_{p}\leq ||\mathbf{J}_{M}||_{p}\leq M^{p-1}|||%
\mathbf{J}_{M}|||_{p},\emph{\quad }p\geq 1.  \tag{3.48}
\end{equation}

\textbf{(iii)} \ \emph{Given any} $\Lambda \Subset \mathbb{Z}^{d},$ \emph{%
let us define seminorms} 
\begin{equation}
|||\mathbf{J}_{M}|||_{\Lambda ,p}:=\sum\limits_{L=1}^{M}\sum\limits 
_{\substack{ \{k_{1},...,k_{L}\}\subset \Lambda  \\ \{k_{L+1},...,k_{M}\}%
\subset \Lambda ^{c}}}J_{k_{1},...,k_{M}}\sum%
\limits_{l=1}^{M}(1+|k_{l}|)^{p},\quad p\geq 0.  \tag{3.49}
\end{equation}%
\emph{Then} 
\begin{equation}
|||\mathbf{J}_{M}|||_{\Lambda ,p}\leq C(M,\Lambda ,p)\cdot ||\mathbf{J}||_{p}%
\text{ \ \emph{with \ }}C(M,\Lambda ,p):=M|\Lambda |\cdot \sup_{k\in \Lambda
}\{(1+|k|)^{p}\}.  \tag{3.50}
\end{equation}

\textbf{Proof:} \textbf{(i)} \ (3.46) follows by a direct calculation from
assumption ($\mathbf{J}$).

\textbf{(ii)} \ (3.48) is readily apparent from (3.47) and the inequality 
\begin{equation}
\frac{1}{M}\sum\limits_{l=1}^{M}|q_{l}|^{p}\leq \left(
\sum\limits_{l=1}^{M}|q_{l}|\right) ^{p}\leq
M^{p-1}\sum\limits_{l=1}^{M}|q_{l}|^{p}  \tag{3.51}
\end{equation}%
which is valid for any $p\geq 1$ and $q_{1},...,q_{M}\in \mathbb{R}.$

\textbf{(iii)} \ (3.50) follows from ($\mathbf{J}$) and (3.48) by the
estimate%
\begin{multline}
|||\mathbf{J}_{M}|||_{\Lambda ,p}\leq |\Lambda |\sup_{k_{1}\in \Lambda
}\{(1+|k_{1}|)^{p}\}\sup_{k_{1}\in \Lambda }\left\{
\sum\limits_{\{k_{2},...,k_{M}\}\subset \mathbb{Z}^{d}}J_{k_{1},...,k_{M}}%
\sum\limits_{l=1}^{M}(1+|k_{1}-k_{l}|)^{p}\right\}  \notag \\
\leq |\Lambda |\sup_{k\in \Lambda }\{(1+|k|)^{p}|||\mathbf{J}_{M}|||_{p}. 
\tag{3.52}
\end{multline}%
$\blacksquare $\smallskip \newline
\textbf{Lemma 3.5.} \ \emph{Along with }($\mathbf{W}_{\text{\textbf{i}}}$)%
\emph{\ and }($\mathbf{J}$)\emph{, let also Assumption }($\mathbf{V}_{\text{%
\textbf{iii}}}$)\emph{\ hold} \emph{with} \emph{some fixed, but} \emph{small
enough} $0<K_{3}<K_{3}^{(0)}$ \emph{(e.g., satisfying (3.55) below)}$.$\emph{%
\ Then} 
\begin{equation}
0<Z_{\Lambda }(\xi )<\infty \quad \text{\emph{for all}}\ \xi \in \Omega _{t}%
\text{ \emph{and} }\Lambda \Subset \mathbb{Z}^{d},  \notag
\end{equation}%
\emph{i.e., the right-hand-sides in (3.39) and (3.42) make sense}.\smallskip

\textbf{Proof:} \ Let $\xi \in \Omega _{-p}^{R}$ for some $p>d.$ Due to ($%
\mathbf{W}_{\text{\textbf{i}}}$), ($\mathbf{J}$) and (3.49)--(3.51) we have
that 
\begin{multline}
\sum\limits_{M=2}^{N}\sum\limits_{L=1}^{M}\sum\limits_{\substack{ %
\{k_{1},...,k_{L}\}\subset \Lambda  \\ \{k_{L+1},...,k_{M}\}\subset \Lambda
^{c}}}\int_{S_{\beta }^{{}}}|W_{\{k_{1},...,k_{M}\}}(\omega
_{k_{1}},...,\omega _{k_{M}})|d\tau  \notag \\
\leq \sum\limits_{M=2}^{N}(M+1)^{R-1}\sum\limits_{L=1}^{M}\sum\limits 
_{\substack{ \{k_{1},...,k_{L}\}\subset \Lambda  \\ \{k_{L+1},...,k_{M}\}%
\subset \Lambda ^{c}}}J_{k_{1},...,k_{M}}\left( \beta
+\sum\limits_{l=1}^{M}|\omega _{k_{l}}|_{L_{\beta }^{R}}^{R}\right)  \notag
\\
\leq \left[ \sum\limits_{M=2}^{N}(M+1)^{R}||\mathbf{J}_{M}||_{pR}\right]
|\Lambda |\cdot \sup_{k\in \Lambda }\{(1+|k|)^{pR}\}\left( \beta +||\omega
||_{-p,R}^{R}\right)  \notag \\
\leq 2^{R}||\mathbf{J}||_{pR}|\Lambda |\cdot \sup_{k\in \Lambda
}\{(1+|k|)^{pR}\}\left( \beta +||\omega ||_{-p,R}^{R}\right) .  \tag{3.53}
\end{multline}%
Thus, returning to the notation (3.39)--(3.41), 
\begin{multline}
|\mathcal{I}_{\beta ,\Lambda }^{W}(\omega |\xi )|\leq
\sum\limits_{M=2}^{N}(M+1)^{R-1}||\mathbf{J}_{M}||_{0}\sum_{k\in \Lambda
}|\omega _{k}|_{L_{\beta }^{R}}^{R}  \notag \\
+2^{R}||\mathbf{J}||_{pR}|\Lambda |\cdot \sup_{k\in \Lambda }\left\{
(1+|k|)^{pR}\right\} \left( \beta +||\xi _{\Lambda ^{c}}||_{-p,R}^{R}\right)
.  \tag{3.54}
\end{multline}%
Combining (3.54) and (3.33), we can get even more than stated in Lemma 3.5.
Namely, if%
\begin{equation}
0<K_{3}<K_{3}^{(0)}:=(2^{R-2}R\cdot ||\mathbf{J}||_{0})^{-1},  \tag{3.55}
\end{equation}%
then the local partition function $Z_{\Lambda }(\xi )$ is uniformly positive
and bounded on all finite radius balls%
\begin{equation*}
B_{-p}^{R}(\rho ):=\{\xi \in \mathcal{L}_{-p}^{R}|\ ||\xi ||_{-p,R}\leq \rho
\},\text{ \ }0<\rho <\infty ,
\end{equation*}%
i.e., it holds%
\begin{equation}
0<\inf_{\xi \in B_{-p}^{R}(\rho )}Z_{\Lambda }(\xi )\leq \sup_{\xi \in
B_{-p}^{R}(\rho )}Z_{\Lambda }(\xi )<\infty .  \tag{3.56}
\end{equation}%
$\blacksquare $\smallskip

\noindent \textbf{Lemma 3.6.} \ \emph{Under the assumptions of Lemma 3.5,
the specification }$\pi =\{\pi _{\Lambda }\}_{\Lambda \Subset \mathbb{Z}%
^{d}} $ \emph{possesses the following regularity properties: for any} $%
\Lambda \Subset \mathbb{Z}^{d},$\emph{\ }$p>d$%
\begin{equation}
(\omega ,\xi )\mapsto \mathcal{I}_{\Lambda }(\omega |\xi )\ \text{\emph{is\
uniformly\ Lipschitz-continuous\ on\ balls\ in }}\mathcal{L}_{-p}^{R}\times 
\mathcal{L}_{-p}^{R},  \tag{3.57}
\end{equation}%
\emph{and hence for any bounded measurable} $f:\Omega \rightarrow \mathbb{R}$%
\begin{equation}
\xi \mapsto (\pi _{\Lambda }f)(\xi ):={\displaystyle\int_{\Omega }}f(\omega )\pi _{\Lambda
}(d\omega |\xi )\text{\ \emph{is uniformly Lipschitz-continuous on balls in }%
}\mathcal{L}_{-p}^{R}.  \tag{3.58}
\end{equation}

\textbf{Proof:} \ (3.57) is a consequence of the estimate%
\begin{multline*}
\sum\limits_{M=2}^{N}\sum\limits_{L=1}^{M}\sum\limits_{\substack{ %
\{k_{1},...,k_{L}\}\subset \Lambda  \\ \{k_{L+1},...,k_{M}\}\subset \Lambda
^{c}}}\int_{S_{\beta }^{{}}}|W_{\{k_{1},...,k_{M}\}}(\omega
_{k_{1}},...,\omega _{k_{M}})-W_{\{k_{1},...,k_{M}\}}(\omega
_{k_{1}}^{\prime },...,\omega _{k_{M}}^{\prime })|d\tau \\
\leq \sum\limits_{M=2}^{N}\sum\limits_{L=1}^{M}\sum\limits_{\substack{ %
\{k_{1},...,k_{L}\}\subset \Lambda  \\ \{k_{L+1},...,k_{M}\}\subset \Lambda
^{c}}}J_{k_{1},...,k_{M}}\int\limits_{S_{\beta
}^{{}}}\sum\limits_{l=1}^{M}|\omega _{k_{l}}-\omega _{k_{l}}^{\prime
}|\left( 1+\sum\limits_{l=1}^{M}|\omega
_{k_{l}}|+\sum\limits_{l=1}^{M}|\omega _{k_{l}}^{\prime }|\right) ^{R-1}d\tau
\\
\leq 3^{R}||\mathbf{J}||_{pR}^{R}|\Lambda |\cdot \sup_{k\in \Lambda
}\{(1+|k|)^{pR}\}||\omega -\omega ^{\prime }||_{-p,R}^{{}}\left( \beta
^{(1-R^{-1})}+||\omega ||_{-p,R}^{R-1}+||\omega ^{\prime
}||_{-p,R}^{R-1}\right) ,
\end{multline*}%
which one can easily get from ($\mathbf{W}_{\text{\textbf{ii}}}$), ($\mathbf{%
J}$), H\"{o}lder's inequality and claim (iii) of Lemma 3.4. Thereafter,
(3.58) immediately follows from (3.42), (3.56) and (3.57).

\noindent $\blacksquare \smallskip $

\noindent \textbf{Remark 3.7.} \ In fact, the proof of Lemma 3.6 shows that,
at least for small enough $\lambda ,K_{3}>0$ such that%
\begin{equation}
K_{3}<\left[ R(\lambda +2^{R-2}||\mathbf{J}||_{0})\right] ^{-1}=:K_{3}^{(%
\lambda )}<K_{3}^{(0)},  \tag{3.59}
\end{equation}%
the measures $\pi _{\Lambda }(d\omega |\xi )$ certainly have all moments of
the form 
\begin{equation}
\int_{\Omega }\exp \left\{ \lambda |\omega _{k}|_{L_{\beta
}^{R}}^{R}\right\} \left( 1+|\omega _{k}|_{C_{\beta }^{\alpha }}^{Q}\right)
d\pi _{\Lambda }(d\omega |\xi )<\infty \text{, \ \ }0\leq \alpha <\frac{1}{2}%
,\text{ }Q\geq 1.  \tag{3.60}
\end{equation}%
So, (3.60) holds for all $\lambda >0$ for the particular QLS models
described in Sect.\thinspace 2 (cf. Remark 3.3 (ii)).

\subsubsection{Definition of tempered Euclidean Gibbs measures}

In this subsection assume that ($\mathbf{W}_{\text{\textbf{i}}}$), ($\mathbf{%
J}$) and ($\mathbf{V}_{\text{\textbf{iii}}}$) in Definition 2.2
hold.\smallskip

\noindent \textbf{Definition 3.8.} \ \emph{A probability measure }$\mu $%
\emph{\ on }$(\Omega ,\mathcal{B}(\Omega ))$\emph{\ is called\ Euclidean
Gibbs state for the specification }$\pi =\{\pi _{\Lambda }^{{}}\}_{\Lambda
\Subset \mathbb{Z}^{d}}$ \emph{(corresponding to the lattice system (3.1) at
inverse temperature }$\beta >0$\emph{) if it satisfies the DLR equilibrium
equations:}%
\begin{equation}
\mu \pi _{\Lambda }=\mu \text{,}\emph{\ \ }\forall \Lambda \Subset \mathbb{Z}%
^{d},  \tag{3.61}
\end{equation}%
\emph{where the measures }$\mu \pi _{\Lambda }$ \emph{are defined by} 
\begin{equation}
\mu \pi _{\Lambda }(\Delta ):=\int_{\Omega }\mu (d\omega )\pi _{\Lambda
}(\Delta |\omega )\text{,\quad }\forall \Delta \in \mathcal{B}(\Omega ). 
\tag{3.62}
\end{equation}%
\smallskip

In other words, $\mu $\ is a measure on $(\Omega ,\mathcal{B}(\Omega ))$
such that its regular conditional distributions given $\omega _{k},$ $k\in
\Lambda ^{c},$ coincide $\mu -$a.e. with the prescribed local specifications 
$\pi _{\Lambda }$ for all $\Lambda \Subset \mathbb{Z}^{d}.$ Let $\mathcal{G}$
denote the set of all Gibbs measures for the system (3.1) at a fixed inverse
temperature $\beta >0$. Applying (3.62) to $B:=(\Omega _{t})^{c}$, it
follows from (3.41) that $\mu (\Omega _{t})=1,$ i.e., any $\mu \in \mathcal{G%
}$ is supported by $\Omega _{t}:=\Omega _{(s)t}^{R}:=\cup _{p>0}\Omega
_{-p}^{R}$. From technical reasons (which become more clear after proof of
our main Hypotheses $\mathbf{(H)}$ in Sect.\thinspace 7), we shall mainly
restrict our consideration to \textit{\emph{tempered }}Gibbs measures $\mu
\in \mathcal{G}_{t}:=\mathcal{G}_{(s)t}^{R}$ supported on the whole spaces $%
\Omega _{-p}^{R}$, $p>d$.\textit{\smallskip }\newline
\noindent \textbf{Definition 3.9.}\emph{\ \ A probability measure }$\mu $%
\emph{\ on }$(\Omega ,\mathcal{B}(\Omega ))$ \emph{is called\
(exponentially) tempered, i.e., }$\mu \in \mathcal{M}_{(e)t}^{R}$, \emph{if}%
\begin{equation}
\forall \delta >0:\text{ \ \ }\mu \left( \Omega _{-\delta }^{R}\right) =1, 
\tag{3.63}
\end{equation}%
\emph{and, moreover, (slowly) tempered, i.e., }$\mu \in \mathcal{M}%
_{(s)t}^{R}\subseteq \mathcal{M}_{(e)t}^{R}$, \emph{if} 
\begin{equation}
\exists p=p(\mu )>d:\text{ \ \ }\mu \left( \Omega _{-p}^{R}\right) =1. 
\tag{3.64}
\end{equation}%
\emph{Respectively we define the subsets }$\mathcal{G}_{(s)t}^{R}\subseteq 
\mathcal{G}_{(e)t}^{R}$ \emph{of}\textit{\ \emph{tempered Gibbs measures by }%
}%
\begin{equation}
\mathcal{G}_{(e)t}^{R}:=\mathcal{G}\cap \mathcal{M}_{(e)t}^{R},\text{ \ \ }%
\mathcal{G}_{(s)t}^{R}:=\mathcal{G}\cap \mathcal{M}_{(s)t}^{R}.  \tag{3.65}
\end{equation}%
\smallskip

\noindent \textbf{Remark} \textbf{3.10. (i) \ }In fact (cf. e.g. [Ge88,
Theorem 1.33]), under conditions imposed any tempered $\mu \in \mathcal{G}$
is \emph{completely determined} just by the family of its \emph{one-site
conditional distributions} $\pi _{\Lambda _{k}}$ with $\Lambda _{k}:=\{k\},$ 
$k\in \mathbb{Z}^{d}$. More precisely, a probability measure $\mu $\ on\emph{%
\ }$(\Omega _{(s)t}^{R},\mathcal{B}(\Omega _{(s)t}^{R}))$ is Gibbs for\emph{%
\ }$\pi =\{\pi _{\Lambda }\}_{\Lambda \Subset \mathbb{Z}^{d}}$ iff $\mu \pi
_{\Lambda _{k}}=\mu $ for all $k\in \mathbb{Z}^{d}.\smallskip $

\textbf{(ii)} \ Under our assumptions every cluster point $\mu _{\ast }$ of
the family $\{\pi _{\Lambda }(d\omega |\xi )\}_{\Lambda \Subset \mathbb{Z}%
^{d},\text{ }\xi \in \Omega _{-p}^{R}}$ w.r.t. to the topology of \emph{weak}
\emph{convergence} of probability measures on the given Polish space $\Omega
_{-p}^{R},$ $p>d$, is surely Gibbs, i.e., $\mu _{\ast }\in \mathcal{G}_{t}.$
Indeed, let $\pi _{\Lambda ^{(K)}}(d\omega |\xi ^{(K)})\overset{w}{%
\rightarrow }\mu _{\ast }(d\omega )$ as $\Lambda ^{(K)}\nearrow \mathbb{Z}%
^{d},$ $K\rightarrow \infty ,$ and let us choose arbitrary $\Lambda
_{0}\Subset \mathbb{Z}^{d}$ and $f\in C_{b}(\Omega _{-p}^{R})$. Since by
Lemma 3.6 also $\pi _{\Lambda _{0}}f\in C_{b}(\Omega _{-p}^{R}),$ one can
directly pass to the limit $K\rightarrow \infty $ in the consistency
condition ($\mathbf{S}_{\mathbf{iv}}$) written for $\Lambda ^{(K)}\supseteq
\Lambda _{0}$%
\begin{equation}
\int_{\Omega }\int_{\Omega }f(\omega ^{\prime })\pi _{\Lambda _{0}}(d\omega
^{\prime }|\omega )\pi _{\Lambda ^{(K)}}(d\omega |\xi ^{(K)})=\int_{\Omega
}f(\omega )\pi _{\Lambda ^{(K)}}(d\omega |\xi ^{(K)}),  \tag{3.66}
\end{equation}%
and thus get the DLR equation 
\begin{equation}
\int_{\Omega }\int_{\Omega }f(\omega ^{\prime })\pi _{\Lambda _{0}}(d\omega
^{\prime }|\omega )\mu _{\ast }(d\omega )=\int_{\Omega }f(\omega )\mu _{\ast
}(d\omega ).  \tag{3.67}
\end{equation}%
\smallskip

\textbf{(iii)} \ Analogously one can show that the subset $\{\mu \in 
\mathcal{G}_{t}|$ $\mu (\Omega _{-p}^{R})=1\}$ is \emph{closed }in the
topology of weak convergence of probability measures on $\Omega _{-p}^{R}$
with $p>d.$ Since the embedding $\mathcal{C}_{-p}^{\alpha }\underset{%
\longrightarrow }{\subset }\Omega _{-p}^{R}$ is continuous, $\{\mu \in 
\mathcal{G}_{t}|$ $\mu (\mathcal{C}_{-p}^{\alpha })=1\}$ is thus closed in
the topology of weak convergence on $\mathcal{C}_{-p}^{\alpha }$ with $0\leq
\alpha <1/2$. For the particular QLS$\ $models I, II the same result is
obviously valid just in the product spaces $[C_{\beta }^{\alpha }]^{\mathbb{Z%
}^{d}}$ with $p>d.$ \smallskip

\textbf{(iv)} \ For the \emph{translation invariant} systems there naturally
arises a question how to construct such $\mu \in \mathcal{G}_{t}$ which are
also invariant w.r.t. to the discrete group of lattice translations $\mathbb{%
Z}_{0}^{d}$, i.e., 
\begin{equation*}
\mu (\Delta )=\mu (T_{j}\Delta ),\text{ \ }\forall j\in \mathbb{Z}_{0}^{d},%
\text{ }\forall \Delta \in \mathcal{B}(\Omega ),
\end{equation*}%
where $T_{j}$ is the shift transformation, $(T_{j}\omega )_{k}:=\omega
_{k+j},$ $k\in \mathbb{Z}^{d}$. With this aim one standardly uses the
so-called local Gibbs distributions with \emph{periodic boundary conditions }%
(cf. e.g. [AKKR02]). Consider any box $\Lambda $ of the form 
\begin{equation}
\Lambda :=\times _{\nu =1}^{d}[a_{\nu },a_{\nu }+l_{\nu }]\cap \mathbb{Z}^{d}%
\text{ with }a=(a_{\nu })_{\nu =1}^{d}\in \mathbb{Z}^{d}\text{, \ }l=(l_{\nu
})_{\nu =1}^{d}\in \mathbb{N}_{0}^{d},  \tag{3.68}
\end{equation}%
and let $\mathfrak{T}(\Lambda )\cong \mathbb{Z}^{d}/l\mathbb{Z}^{d}$ be the
torus obtained by identifying its opposite walls. Respectively, we define
the $\Lambda $-periodic continuation $\tilde{\omega}_{\Lambda }\in \Omega $
of the projection $\omega _{\Lambda }\in \Omega _{\Lambda }$ by $\tilde{%
\omega}_{\Lambda }:=(\omega _{j(k)})_{k\in \mathbb{Z}^{d}},$ where $j(k)\in 
\mathfrak{T}(\Lambda )$ is the unique element such that $j(k)\cong k$ $(%
\func{mod}$ $l).$ Next, we can define the $\Lambda $-periodic modification
of the interaction%
\begin{equation}
\mathcal{I}_{per,\Lambda }^{W}(\omega
):=\sum\limits_{M=2}^{N}\sum\limits_{\{k_{1},...,k_{M}\}\subset \Lambda
}\int_{S_{\beta }}W_{\{k_{1},...,k_{M}\}}(\omega _{j(k_{1})}(\tau
),...,\omega _{j(k_{M})}(\tau ))d\tau   \tag{3.69}
\end{equation}%
and the associated with it periodic local Gibbs distribution 
\begin{equation}
\mu _{\mathfrak{T}(\Lambda )}(d\omega _{\emph{T}(\Lambda )}):=Z_{\mathfrak{T}%
(\Lambda )}^{-1}\exp \left\{ -\mathcal{I}_{per,\Lambda }^{W}(\omega
)\right\} \prod\limits_{j\in \mathfrak{T}(\Lambda )}d\sigma _{\beta }(\omega
_{j})  \tag{3.70}
\end{equation}%
as a probability measure on $\Omega _{\mathfrak{T}(\Lambda )}:=\mathbb{R}^{%
\mathfrak{T}(\Lambda )}.$ Suppose now that a sequence $\mu _{\mathfrak{T}%
(\Lambda ^{(K)})}$ converges locally weakly on $\Omega $ as $\Lambda
^{(K)}\nearrow \mathbb{Z}^{d},$ $K\rightarrow \infty ,$ to some $\mu _{\ast
}\in \mathcal{M}(\Omega _{t}),$ that means 
\begin{equation}
\lim_{K\rightarrow \infty }\int_{\mathfrak{T}(\Lambda _{n})}f(\omega )\mu _{%
\mathfrak{T}(\Lambda ^{(K)})}(d\omega _{T(\Lambda ^{(K)})})=\int_{\Omega
_{t}}f(\omega )\mu _{\ast }(d\omega )  \tag{3.71}
\end{equation}%
for any local $f=f(\mathbb{P}_{\Lambda _{0}})\in C_{b}(\Omega )$. Assuming
again that $W_{\{k_{1},...,k_{M}\}}\in C(\mathbb{R}^{M})$ are local, it is
easy to show that $\mu _{\ast }\in \mathcal{G}_{t}.$ So, by construction one
has the consistency relation%
\begin{equation}
\int_{\mathfrak{T}(\Lambda )}\int_{\Lambda _{0}}f(\omega ^{\prime })\mu
_{\Lambda _{0}}(d\omega _{\Lambda _{0}}^{\prime }|\omega )\mu _{\mathfrak{T}%
(\Lambda ^{(K)})}(d\omega _{T(\Lambda )})=\int_{\mathfrak{T}(\Lambda
^{(K)})}f(\omega )\mu _{\mathfrak{T}(\Lambda )}(d\omega _{\mathfrak{T}%
(\Lambda )})  \tag{3.72}
\end{equation}%
in any boxes $\Lambda _{0},\Lambda \Subset \mathbb{Z}^{d}$ such that $%
\Lambda ^{(\rho +1)}:=\left\{ k\in \mathbb{Z}^{d}|\text{ dist\thinspace }%
(k,\Lambda _{0})\leq \rho +1\right\} \subset \Lambda .$ Since $\pi _{\Lambda
_{0}}f=:g=g(\mathbb{P}_{\Lambda ^{(K)}})\in C_{b}(\Omega ),$ one can
directly pass in (3.71) to the limit as $\Lambda :=\Lambda ^{(K)}\nearrow 
\mathbb{Z}^{d}$ and thus get the DLR equation (3.67). Note that the required 
$\mathbb{Z}_{0}^{d}-$translation invariance of $\mu _{\ast }$ follows from
the invariance of any $\mu _{\mathfrak{T}(\Lambda )}$ w.r.t. the
translations of the torus $\mathfrak{T}(\Lambda ).\smallskip $

\section{Flow and integration by parts characterization of Euclidean Gibbs
measure}

Analogously to what has been already done for classical lattice systems in
[AKRT99], in this paper we would like to develop an\emph{\ alternative
approach to the construction and study of Euclidean Gibbs measures} for
quantum lattice systems. The main ingredients of this approach are: first,
the\emph{\ flow}\textit{\ }\emph{characterization of }$\mu \in \mathcal{G}$ 
\emph{in terms of Radon--Nikodym derivatives} w.r.t. shift transformations
of the configuration space\emph{\ }$\Omega $ and, second, the
characterization (resulting from the previous one) in terms of their
logarithmic derivatives via corresponding\emph{\ integration by parts
formulas}. If the interaction potentials are differentiable (as they are in
our case), both characterizations are equivalent. Such alternative (to the
usual through the DLR equations) characterizations of Gibbs measures have
long been known for a number of specific models in statistical mechanics and
field theory (see, e.g., [Roy75,77, HS76, Fr\"{o}77, Fri 82, Deu87, Ki95]).
However, in rather full generality the flow description for both classical
and quantum Gibbs measures has first been proved in [AKR97a,b]. In
particular, [AKR97b, Sect.\thinspace 4.2] contains the flow characterization
for Euclidean Gibbs measures of quantum lattice systems (3.1) with harmonic
pair interactions. In this section we extend the latter result to systems
with many-particle interactions and give the complete characterization of $%
\mu \in \mathcal{G}_{t}$ in terms of their Radon--Nikodym and logarithmic
derivatives. Also we observe that the local Gibbs specifications $\pi
_{\Lambda },$ $\Lambda \Subset \mathbb{Z}^{d},$ also satisfy the same flow
and integration by parts descriptions, which later will be crucial for our
proof of the existence of $\mu \in \mathcal{G}_{t}$.

As an application of this alternative description of Euclidean Gibbs states,
a direct analytic proof of our main Theorems I and II using the
corresponding (IbP)-formulas will be presented in Subsect.\thinspace 4.6
under Hypotheses \textbf{(H) }and \textbf{(H}$_{\mathtt{loc}}$\textbf{)}
(which in turn will later be proved to be satisfied, cf. Subsect.\thinspace
7.3).

We assume throughout Section 4 that conditions \textbf{(W)}, \textbf{(J)}, 
\textbf{(V) }in Definition 3.2 always hold without any additional reference
to them.

\subsection{Flow description of Euclidean Gibbs measures}

We start with the flow description of $\mu \in \mathcal{G}_{t}$ in terms of
their \emph{\textquotedblright shift\textquotedblright --Radon--Nikodym
derivatives }$a_{\theta h_{i}},$ $\theta \in \mathbb{R},$ along some set of
admissible directions $h_{i}$, $i\in \mathcal{I},$ whose linear span is
dense in $\Omega _{t}:=\Omega _{(s)t}^{R}.$

As already said in Subsect.\thinspace 3.2.3, we consider 
\begin{equation*}
\mathcal{H}:=l^{2}\otimes L_{\beta }^{2}\cong l^{2}(\mathbb{Z}%
^{d}\rightarrow L_{\beta }^{2})
\end{equation*}%
with the inner product $<\omega ,\omega >_{\mathcal{H}}:=||\omega ||_{%
\mathcal{H}}^{2}$ as a \emph{tangent Hilbert space} to $\Omega .$ For the
remainder of the paper, we fix the concrete \emph{orthonormal} \emph{basis}
in $\mathcal{H}$ (for shorthand, we denote it by $bas(\mathcal{H}%
):=\{h_{i}\}_{i\in \mathbb{Z}^{d+1}}$) consisting of the vectors%
\begin{equation}
h_{i}:=e_{k}\otimes \varphi _{n}=(\delta _{k-k^{\prime }}\varphi
_{n})_{k^{\prime }\in \mathbb{Z}^{d}},\text{ \ \ }i=(k,n)\in \mathbb{Z}%
^{d+1},\ k\in \mathbb{Z}^{d},\ n\in \mathbb{Z}.  \tag{4.1}
\end{equation}%
Here $e_{k}:=(\delta _{k-k^{\prime }})_{k^{\prime }\in \mathbb{Z}^{d}},$ and 
$\varphi _{n}\in C_{\beta }^{\infty }$ are the eigenvectors of the operator $%
A$ in $H:=L_{\beta }^{2}$ given by (3.24).

For $i=(k,n)\in \mathbb{Z}^{d+1}$ and $\theta \in \mathbb{R},$ we define the
corresponding \emph{relative many-particle interaction} by 
\begin{multline}
\mathcal{I}_{rel}^{W}(\omega |\theta h_{i}):=  \notag \\
\sum\limits_{M=2}^{N}\sum\limits_{\substack{ k_{1}:=k  \\ %
\{k_{1},...,k_{M}\}\subset \mathbb{Z}^{d}}}\int\limits_{S_{\beta }^{{}}}%
\left[ W_{\{k_{1}..,k_{M}\}}(\omega _{k_{1}}+\theta \varphi _{n},\omega
_{k_{2}},...,\omega _{k_{M}})-W_{\{k_{1},...,k_{M}\}}(\omega
_{k_{1}},...,\omega _{k_{m}})\right] d\tau .  \tag{4.2}
\end{multline}%
(In the notation of (3.41) it heuristically equals to $\mathcal{I}_{\mathbb{Z%
}^{d}}^{W}(\omega +\theta h_{i})-\mathcal{I}_{\mathbb{Z}^{d}}^{W}(\omega )$%
). Since $(\mathbf{W}_{\text{\textbf{i}}})$ and $(\mathbf{J})$ hold,
arguments similar to those applied in the proof of Lemmas 3.5 and 3.6 yield
the following properties of the function $(\theta ,\omega )\mapsto \mathcal{I%
}_{rel}^{W}(\omega |\theta h_{i})$ for every $i\in \mathbb{Z}^{d}$:%
\begin{gather}
\text{\emph{The sum in (4.2) with the integrands replaced by their absolute
values}}  \notag \\
\text{\emph{converges for all }}\theta \in \mathbb{R}\text{, }\omega \in 
\mathcal{L}_{-p}^{R},\text{ \ \emph{and}}  \notag \\
(\theta ,\omega )\mapsto \mathcal{I}_{rel}^{W}(\omega |\theta h_{i})\text{%
\emph{\ is uniformly Lipschitz-continuous on balls in }}\mathbb{R\times }%
\mathcal{L}_{-p}^{R}\text{, }p>d.  \tag{4.3}
\end{gather}%
Moreover, exploiting $(\mathbf{W}_{\text{\textbf{ii}}})$ and using a
standard theorem of analysis about the differentiability of uniformly
convergent series (cf. e.g. [Car67, Theorem\thinspace 3.6.2]), one can show
that there exists%
\begin{equation}
\frac{\partial }{\partial \theta }\mathcal{I}_{rel}^{W}(\omega |\theta
h_{i})=\sum\limits_{M=2}^{N}\sum\limits_{\substack{ k_{1}:=k  \\ %
\{k_{1},...,k_{M}\}\subset \mathbb{Z}^{d}}}\int\limits_{S_{\beta }^{{}}}%
\left[ \partial _{1}W_{\{k_{1}..,k_{M}\}}(\omega _{k_{1}}+\theta \varphi
_{n},\omega _{k_{2}},...,\omega _{k_{M}})\cdot \varphi _{n}\right] (\tau
)\,d\tau  \tag{4.4}
\end{equation}%
with the same regularity properties of the function $(\theta ,\omega
)\longmapsto \frac{\partial }{\partial \theta }\mathcal{I}_{rel}^{W}(\omega
|\theta h_{i})$ as described in (4.3).

Next, let us define the following densities $a_{\theta h_{i}}:\Omega
_{t}\rightarrow \mathbb{R},$ where (cf. (3.35) and (4.2))%
\begin{gather}
a_{\theta h_{i}}:=a_{\theta h_{i}}^{A_{\beta }}\cdot a_{\theta
h_{i}}^{V}\cdot a_{\theta h_{i}}^{W}\text{ \ \ and}  \notag \\
a_{\theta h_{i}}^{A_{\beta }}(\omega ):=a_{\theta \varphi _{n}}^{A_{\beta
}}(\omega _{k}),\text{ \ }a_{\theta h_{i}}^{V}(\omega ):=a_{\theta \varphi
_{n}}^{V_{k}}(\omega _{k}),\text{ \ }a_{\theta h_{i}}^{W}(\omega ):=\exp \{-%
\mathcal{I}_{rel}^{W}(\omega |\theta h_{i})\}.  \tag{4.5}
\end{gather}%
According to the discussion above%
\begin{gather}
(\theta ,\omega )\longmapsto a_{\theta h_{i}}^{W}(\omega ),\ \frac{\partial 
}{\partial \theta }a_{\theta h_{i}}^{W}(\omega )\text{\emph{\ are uniformly
Lipschitz-continuous}}  \notag \\
\text{\emph{on balls in }}\mathbb{R\times }\mathcal{L}_{-p}^{R}\text{, }p>d.
\tag{4.6}
\end{gather}%
Thus, combining (3.37) and (4.6), we get the important the for later use
assertion which we formulate as follows:\smallskip\ 

\noindent \textbf{Lemma 4.1 \ }\emph{For every} $i=(k,n)\in \mathbb{Z}^{d+1}$
\emph{the functions}%
\begin{equation*}
(\theta ,\omega )\longmapsto a_{\theta h_{i}}^{{}}(\omega ),\ \frac{\partial 
}{\partial \theta }a_{\theta h_{i}}(\omega )\emph{\ }
\end{equation*}%
\emph{are uniformly Lipschitz-continuous on balls in }$\mathbb{R\times }%
\Omega _{-p;k}^{R}\emph{,}$ $p>d.$\smallskip

The main result of this subsection is Proposition 4.2\textbf{\ }which
extends the corresponding result of Theorem 4.6 in [AKR97b], where (for
simplicity only) the case of harmonic pair interactions (i.e., the
particular QLS model (1.2)) was treated. The method of the proof proposed in
that paper (see also the earlier papers [Iw85, Fu91, Ro75, Fri82] on related
topics) also applies in our situation with some technical alterations.
However,\emph{\ }in order to keep the exposition self-contained, here we
repeat the key steps of [AKR97b] and present the complete proof.\smallskip 
\newline
\textbf{Proposition 4.2 }(Flow Description of Tempered Gibbs States) \ \emph{%
Let} $\mathcal{M}_{t}^{a}$ \emph{denote the set of all} \emph{probability
measures} $\mu $ \emph{on} $(\Omega ,\mathcal{B}(\Omega ))$ \emph{which
satisfy the temperedness condition (3.64)} \emph{and are} \emph{%
quasi-invariant w.r.t. the shifts} $\omega \longmapsto \omega +\theta h_{i}$%
, $\theta \in \mathbb{R}$, $i=(k,n)\in \mathbb{Z}^{d+1},$ \emph{with the
Radon--Nikodym derivatives} 
\begin{equation}
\frac{d\mu (\omega +\theta h_{i})}{d\mu (\omega )}:=a_{\theta h_{i}}(\omega
).  \tag{4.7}
\end{equation}%
\emph{Then} $\mathcal{G}_{t}=\mathcal{M}_{t}^{a}.\smallskip $

\textbf{Proof: \ }Define $\mathcal{T}:=[T_{\beta }]^{\mathbb{Z}%
_{0}^{d}}:=lin\{h_{i}\}_{i\in \mathcal{I}}\subset \Omega $ as a subspace of
all finite sequences $\eta =(\eta _{k})_{k\in \mathbb{Z}^{d}}$ $(\eta _{k}=0$
when $|k|>K(\eta )>0)$ with components $\eta _{k}\in T_{\beta }:=lin\left\{
\varphi _{n}\right\} _{n\in \mathbb{Z}}.$ By the standard semi-group
argument, $\mathcal{M}_{t}^{a}$ is exactly the set of all $\mu \in \mathcal{M%
}_{t}$ which are $\mathcal{T}-$quasi-invariant with the cocycle 
\begin{equation}
\frac{d\mu (\omega +\eta )}{d\mu (\omega )}:=a_{\eta }(\omega )=a_{\eta
}^{A_{\beta }}(\omega )\cdot a_{\eta }^{V}(\omega )\cdot a_{\eta
}^{W}(\omega )>0,\quad \eta \in \mathcal{T},\text{ }\omega \in \Omega _{t}. 
\tag{4.8}
\end{equation}%
Here 
\begin{equation}
\begin{array}{l}
a_{\eta }^{A_{\beta }}(\omega ):=\exp \left\{ -\sum_{k\in \mathbb{Z}^{d}:\
|k|\leq K(\eta )}\left[ (A_{\beta }\eta _{k},\omega _{k})_{H}+\frac{1}{2}%
(A_{\beta }^{{}}\eta _{k},\eta _{k})_{H}\right] \right\} , \\ 
a_{\eta }^{V}(\omega ):=\exp \left\{ -\sum_{k\in \mathbb{Z}^{d}:\ |k|\leq
K(\eta )}\int\limits_{S_{\beta }^{{}}}[V_{k}(\omega _{k}+\eta
_{k})-V_{k}(\omega _{k})]\,d\tau \right\} , \\ 
a_{\eta }^{W}(\omega ):=\exp \left\{ -\mathcal{I}_{rel}^{W}(\omega |\eta
)\right\} ,%
\end{array}
\tag{4.9}
\end{equation}%
and 
\begin{multline}
\mathcal{I}_{rel}^{W}(\omega |\eta ):=\sum_{M=2}^{N}\sum_{\substack{ %
\{k_{1},...,k_{M}\}\subset \mathbb{Z}^{d}  \\ \exists 1\leq L\leq M:\text{ }%
|k_{L}|\leq K(\eta )}}\int\limits_{S_{\beta }^{{}}}\left[ W_{\{k_{1}..,k_{M}%
\}}(\omega _{k_{1}}+\eta _{k_{1}},...,\omega _{k_{M}}+\eta _{k_{M}})\right. 
\notag \\
\left. -W_{\{k_{1}..,k_{M}\}}(\omega _{k_{1}},...,\omega _{k_{M}})\right]
\,d\tau .  \tag{4.10}
\end{multline}%
As follows from the preceding to Lemma 4.1 discussion, the above definitions
are correct, and for all $k\in \mathbb{Z}^{d},$ $\eta \in \mathcal{T}$ the%
\emph{\ }functions 
\begin{equation}
\mathbb{R}\times \Omega _{t}\ni (\theta ,\omega )\longmapsto a_{\theta \eta
}(\omega )\text{ \ \emph{are continuous and locally bounded.}}  \tag{4.11}
\end{equation}%
\smallskip

\textbf{(i)} $\mathcal{G}_{t}\subseteq \mathcal{M}_{t}^{a}:$ \ This
inclusion is obvious. Suppose that $\mu \in \mathcal{G}_{t},$ i.e., it
satisfies the DLR equations (3.59). Then from the quasi-invariance property
of the probability kernels $\pi _{\Lambda }(d\omega |\xi )$ (see also
Subsect.\thinspace 5.1 below), one has that for any $k\in \Lambda \Subset 
\mathbb{Z}^{d},$ $\eta :=\theta h_{i}\in \mathcal{T}$ and $\Delta \in 
\mathcal{B}(\Omega )$%
\begin{multline*}
\int_{\Omega _{t}}\mathbf{1}_{\Delta }(\omega +\eta )\mu (d\omega
)=\int_{\Omega _{t}}\int_{\Omega _{t}}\mathbf{1}_{\Delta }(\omega +\eta )\pi
_{\Lambda }(d\omega |\xi )\mu (d\xi )= \\
\int_{\Omega _{t}}\int_{\Omega _{t}}\mathbf{1}_{\Delta }(\omega )a_{\eta
}(\omega )\pi _{\Lambda }(d\omega |\xi )\mu (d\xi )=\int_{\Omega _{t}}%
\mathbf{1}_{\Delta }(\omega )a_{\eta }(\omega )\mu (d\omega ).
\end{multline*}%
This relation exactly means that $\mu \in \mathcal{M}_{a}^{t}.\smallskip $

\textbf{(ii)} $\mathcal{M}_{t}^{a}\subseteq \mathcal{G}_{t}:$ \ Keeping the
notation of Subsect.\thinspace 3.3, for every $k\in \mathbb{Z}^{d}$ define $%
\Lambda _{k}:=\{k\},$ $\Lambda _{k}^{c}:=\mathbb{Z}^{d}\backslash \{k\}$ and 
$\omega _{\Lambda _{k}^{c}}:=\mathbb{P}_{\Lambda _{k}^{c}}\omega \in \Omega
_{\Lambda _{k}^{c}}^{{}}.$ Now we start with arbitrary $\mu \in \mathcal{M}%
_{t}^{a},$ and let us disintegrate this measure w.r.t. its projection $\mu
_{\Lambda _{k}^{c}}:=\mu \mathbb{P}_{\Lambda _{k}^{c}}^{-1}$ onto $\Omega
_{\Lambda _{k}^{c}}\,$(see, e.g., [Par67, p.\thinspace 147, Theorem 8.1]):%
\begin{equation}
\mu (d\omega _{k},d\omega _{\Lambda _{k}^{c}})=\nu _{\omega _{\Lambda
_{k}^{c}}}(d\omega _{k})\mu _{\Lambda _{k}^{c}}(d\omega _{\Lambda _{k}^{c}}).
\tag{4.12}
\end{equation}%
Here $\nu _{\omega _{\Lambda _{k}^{c}}}(d\omega _{k})$ are some probability
measures (= regular conditional distributions given $\omega _{\Lambda
_{k}^{c}})$ on $(C_{\beta }^{{}},\mathcal{B}(C_{\beta }^{{}}))$. By the
temperedness condition on $\mu ,$ it is clear that $\mu _{\Lambda
_{k}^{c}}(\Omega _{\Lambda _{k}^{c},t})=1$ where $\Omega _{\Lambda
_{k}^{c},t}:=\mathbb{P}_{\Lambda _{k}^{c}}\Omega _{t}.$ On the other hand,
from (4.8) and (4.12) one can straightforwardly verify (using crucially the
continuity property (4.11) of $a_{\eta _{k}}(\omega )$; cf. [Roy75,
Proposition 3]) the quasi-invariance of the measures $\nu _{\omega _{\Lambda
_{k}^{c}}}^{{}}$ in the following sense: There exists a Borel subset $\Delta
_{k}\subseteq \Omega _{\Lambda _{k}^{c},t}$ of full measure $\mu _{\Lambda
_{k}^{c}},$ i.e., $\mu _{\Lambda _{k}^{c}}(\Delta _{k})=1,$ such that for
any $\omega _{\Lambda _{k}^{c}}\in \Delta _{k}$ and $\eta _{k}\in T_{\beta }$%
\begin{equation}
\frac{d\nu _{\omega _{\Lambda _{k}^{c}}}(\omega _{k}+\eta _{k})}{d\nu
_{\omega _{\Lambda _{k}^{c}}}(d\omega _{k})}=a_{\eta _{k}}(\omega
_{k},\omega _{\Lambda _{k}^{c}}),\quad \omega _{k}\in C_{\beta }(\func{mod}%
\nu _{\omega _{\Lambda _{k}^{c}}}^{{}}).  \tag{4.13}
\end{equation}

>From now on fix any $\xi _{\Lambda _{k}^{c}}\in \Delta _{k},$ and let us
show that (4.13) implies 
\begin{equation}
\nu _{\xi _{\Lambda _{k}^{c}}}(d\omega _{k})=\mu _{\Lambda _{k}}(d\omega
_{k}|\xi _{\Lambda _{k}^{c}})  \tag{4.14}
\end{equation}%
where $\mu _{\Lambda _{k}}^{{}}(d\omega _{k}^{{}}|\xi _{\Lambda
_{k}^{c}}^{{}})$ is the Gibbs measure in the volume $\Lambda _{k}$ with the
boundary condition $\xi _{\Lambda _{k}^{c}}:=\omega _{\Lambda
_{k}^{c}}^{{}}. $ We recall that according to definition (3.42)%
\begin{multline}
\mu _{\Lambda _{k}}(d\omega _{k}^{{}}|\xi _{\Lambda
_{k}^{c}}^{{}}):=Z_{\Lambda _{k}^{{}}}^{-1}(\xi )\cdot \exp \left\{
-\int_{S_{\beta }}V_{k}(\omega _{k})d\tau \right\}  \notag \\
\times \exp \left\{ \sum_{\substack{ \{k_{1},...,k_{M}\}\subset \mathbb{Z}%
^{d}  \\ k_{1}:=k;\ 2\leq M\leq N}}\int_{S_{\beta
}^{{}}}W_{\{k_{1}..,k_{M}\}}(\omega _{k_{1}},\xi _{k_{2}},...,\xi
_{k_{M}})d\tau \right\} \gamma _{\beta }(d\omega _{k}^{{}}),  \tag{4.15}
\end{multline}%
where $\gamma _{\beta }$ is the Gaussian measure with correlation operator $%
A_{\beta }^{-1}$. But $\nu (d\omega _{k}):=\mu _{\Lambda _{k}}^{{}}(d\omega
_{k}^{{}}|\xi _{\Lambda _{k}^{c}}^{{}})$ is the unique probability measure
on $(C_{\beta },\mathcal{B}(C_{\beta }))$ which satisfies the flow
description (4.13), i.e., 
\begin{equation}
\frac{d\nu (\omega _{k}^{{}}+\eta _{k})}{d\nu (\omega _{k}^{{}})}=a_{\eta
_{k}}(\omega _{k}),\quad \omega _{k}\in C_{\beta }\text{ }(\func{mod}\nu ), 
\tag{4.16}
\end{equation}%
for all $\eta _{k}\in T_{\beta }$ with a cocycle having the explicit form%
\begin{gather}
a_{\eta _{k}}(\omega _{k}):=\exp \left\{ -(A_{\beta }^{{}}\eta _{k},\omega
_{k})_{H}^{{}}-\frac{1}{2}(A_{\beta }^{{}}\eta _{k},\eta
_{k})_{H}^{{}}-\int_{S_{\beta }^{{}}}\left[ V_{k}^{{}}(\omega _{k}+\eta
_{k})-V_{k}^{{}}(\omega _{k})\right] d\tau \right\}  \notag \\
\times \exp \left\{ -\sum_{\substack{ \{k_{1},...,k_{M}\}\subset \mathbb{Z}%
^{d}  \\ k_{1}:=k;\ 2\leq M\leq N}}\int_{S_{\beta }^{{}}}\left[
W_{\{k_{1}..,k_{M}\}}(\omega _{k_{1}}+\eta _{k_{1}},\xi _{k_{2}},...,\xi
_{k_{M}})-W_{\{k_{1}..,k_{M}\}}(\omega _{k_{1}},\xi _{k_{2}},...,\xi
_{k_{M}})\right] d\tau \right\} .  \notag \\
\tag{4.17}
\end{gather}%
To check this, let us introduce the new measure 
\begin{equation}
\sigma (d\omega _{k}):=\exp \left\{ \int_{S_{\beta }^{{}}}\left[
V_{k}(\omega _{k})+\sum_{\substack{ \{k_{1},...,k_{M}\}\subset \mathbb{Z}%
^{d}  \\ k_{1}:=k;\ 2\leq M\leq N}}W_{\{k_{1}..,k_{M}\}}(\omega _{k_{1}},\xi
_{k_{2}},...,\xi _{k_{M}})\right] d\tau \right\} \nu (d\omega _{k}) 
\tag{4.18}
\end{equation}%
which (due to our assumptions on the potentials $V$ and $W$) is at least $%
\sigma $-finite on $(C_{\beta },\mathcal{B}(C_{\beta })).$ By (4.16)--(4.18)
we have that%
\begin{equation}
\frac{d\sigma (\omega _{k}+\eta _{k})}{d\sigma (\omega _{k})}=\exp \left\{
-(A_{\beta }^{{}}\eta _{k},\omega _{k})_{H}^{{}}-\frac{1}{2}(A_{\beta
}^{{}}\eta _{k},\eta _{k})_{H}^{{}}\right\}  \tag{4.19}
\end{equation}%
or, equivalently, 
\begin{equation}
\frac{d\sigma (\omega _{k}+\eta _{k})}{d\sigma (\omega _{k})}=\frac{d\gamma
_{\beta }(\omega _{k}+\eta _{k})}{d\gamma _{\beta }(\omega _{k})}. 
\tag{4.20}
\end{equation}%
We claim that (4.20) implies 
\begin{equation}
\sigma (C_{\beta }^{{}})<\infty \text{\quad and\quad }\frac{\sigma (d\omega
_{k}^{{}})}{\sigma (C_{\beta }^{{}})}=\gamma _{\beta }^{{}}(d\omega
_{k}^{{}}).  \tag{4.21}
\end{equation}%
To this end we exactly repeat the corresponding arguments from the proof of
Theorem 4.6 in [AKR97]. Let $A_{\beta }\varphi _{n}=\lambda _{n}\varphi _{n}$
with any $n\in \mathbb{Z}\backslash \{0\},$ and let us consider the image
measure $\sigma _{n}$ of $\sigma $ under the mapping $\omega _{k}\longmapsto
(\omega _{k},\varphi _{n})_{H}.$ Then, by (4.20) and the product structure
of the right hand side of (4.19), $\sigma _{n}$ is quasi-invariant with
\textquotedblright shift\textquotedblright --Radon--Nikodym derivatives
equal to that of the Gaussian measure $\sqrt{\frac{\lambda _{n}}{2\pi }}\exp
(-\frac{1}{2}\lambda _{n}q^{2})dq:=\gamma _{n}(q)dq$ on $\mathbb{R}.$
Herefrom it is well-known and could be easily verified that $\sigma
_{n}(dq)\thicksim dq.$ Hence, if $\rho :=d\sigma _{n}/dq,$ then for all $%
\theta \in \mathbb{R}$%
\begin{equation*}
\frac{\rho (q+\theta )}{\gamma _{n}(q+\theta )}=\frac{\rho (q)}{\gamma
_{n}(q)},\quad q\in \mathbb{R}\text{ }(\func{mod}\text{ }dq).
\end{equation*}%
The latter means that $\rho /\gamma _{n}\equiv const$ $(\func{mod}$ $dq).$
Thus $\sigma _{n}(dq)=const\cdot \gamma _{n}(q)dq$ and, in particular, $%
\sigma (C_{\beta })=\sigma _{n}(\mathbb{R)<}\infty .$ Since $\sigma $ is
finite and satisfies (6.22), it readily follows from [Roy75, Proposition 4]
that $\sigma (d\omega _{k}^{{}})=\sigma (C_{\beta }^{{}})\cdot \gamma
_{\beta }^{{}}(d\omega _{k}^{{}}).$

Consequently, combining (4.15) and (4.18), we deduce that $\nu (d\omega
_{k})=const\cdot \nu _{\xi _{\Lambda _{k}^{c}}}(d\omega _{k})$ and, since
both are probability measures, they coincide. Thus (4.14) is shown.
Herefrom, noting that each measure from $\mathcal{G}_{t}$\ is fully
determined by $\left\{ \pi _{\Lambda _{k}}\right\} _{k\in \mathbb{Z}^{d}}$
(cf. Remark 3.10 (i)), we get the desired inclusion $\mu \in \mathcal{G}%
_{t}. $

\noindent $\blacksquare \medskip $

\noindent \textbf{Remark 4.3 \ }Actually, the flow characterization in
Proposition 4.2 is true under \emph{minimal assumptions} on the potentials,
which guarantee, (besides the well-definedness of the local specification $%
\pi _{\Lambda }$) merely the continuity and local boundedness of the
functions $\mathbb{R\times }\Omega _{t}\ni (\theta ,\omega )\longmapsto
a_{\theta h_{i}}(\omega )\in \mathbb{R}$ for all $i\in \mathcal{I}.$ For the
quantum systems with \emph{local} interaction like the particular models
(2.2) and (2.23) discussed in Sect.\thinspace 2, the Radon--Nikodym
derivatives $a_{\theta h_{i}}(\omega )$ given by (4.8) are well-defined for
all $\omega \in \Omega $ and hence the flow description (4.7) is valid\emph{%
\ for all} Gibbs measures $\mu \in \mathcal{G}$.\smallskip

\subsection{Smooth functions on $\Omega $}

In applications, however, it is more convenient to use not the flow
characterization itself, but its\emph{\ infinitesimal form} which we proceed
to describe in the next subsections.\smallskip

So, by Proposition 4.2, for any $\mu \in \mathcal{G}_{t}$ and all bounded
measurable functions $f:\Omega \rightarrow \mathbb{R}$%
\begin{equation}
\int_{\Omega }f(\omega )a_{\theta h_{i}}^{{}}(\omega )d\mu (\omega
)=\int_{\Omega }f(\omega -\theta h_{i})d\mu (\omega ),  \tag{4.22}
\end{equation}%
and thus 
\begin{equation}
\lim_{\theta \rightarrow \pm 0}\int_{\Omega }f(\omega )\frac{a_{\theta
h_{i}}^{{}}(\omega )-1}{\theta }d\mu (\omega )=-\lim_{\theta \rightarrow \mp
0}\int_{\Omega }\frac{f(\omega +\theta h_{i})-f(\omega )}{\theta }d\mu
(\omega )  \tag{4.23}
\end{equation}%
provided the above limits exist. Henceforth, we should first introduce some
spaces of differentiable functions on the loop lattice $\Omega $ needed to
specify the meaning of (4.23) and (IbP)-formulas resulting from it.
Actually, Subsects.\thinspace 4.2.1--4 can be viewed as a collection of
standard definitions and facts from convex analysis, which we modify for our
concrete situation.

\subsubsection{Partially differentiable functions}

Let $X$, $Y$ be locally convex spaces and let $\Phi :X\rightarrow Y$. We
recall that the \emph{partial derivatives on the right} resp. \emph{left in
the direction} $h\in X$ of the function $\Phi $ at a point $x\in X$ are
defined by%
\begin{equation}
\partial _{h}^{+}\Phi (x):=\lim_{\theta \rightarrow +0}\frac{\Phi (x+\theta
h)-\Phi (x)}{\theta },\quad \partial _{h}^{-}\Phi (x):=\lim_{\theta
\rightarrow -0}\frac{\Phi (x+\theta h)-\Phi (x)}{\theta }.  \tag{4.24}
\end{equation}%
If the right and left limits in (4.24) coincide, one says that there exists
the corresponding \emph{partial derivative }$\partial _{h}\Phi (x)$ \emph{in
the direction} $h$. For $m\in \mathbb{N}$ and given vectors $%
h_{1},...,h_{m}\in X,$ we denote by $C^{m}(X\rightarrow Y;h_{1},...,h_{m})$
the set of all functions $\Phi \in C(X\rightarrow Y)$ having continuous
partial derivatives $\partial _{h_{m}}...\partial _{h_{1}}\Phi :X\rightarrow
Y.$ Respectively, by $C^{m}(X\rightarrow Y)$, $m\in \mathbb{N}\cup \{\infty 
\mathbb{\}},$ we shall denote the set of all $\Phi \in C(X\rightarrow Y)$
having continuous partial derivatives $\partial _{h_{l}}...\partial
_{h_{1}}\Phi :X\rightarrow Y$ of any order $0\leq l\leq m$ along arbitrary
vectors $h_{1},...,h_{l}\in X$. Then, as usual, 
\begin{equation*}
C_{b,loc}^{m}(X\rightarrow Y)\supseteq C_{b}^{m}(X\rightarrow Y)\supseteq
C_{0}^{m}(X\rightarrow Y)
\end{equation*}%
will mean the subspaces of those $\Phi \in C^{m}(X\rightarrow Y)$ which
satisfy the extra assumptions of\emph{\ }local resp. global boundedness of
all their partial derivatives, or, in addition, of boundedness of their
support (supp$\,\Phi :=\{x\in X$\thinspace $|$\thinspace $\Phi (x)\neq 0\})$%
. For shortness, we shall omit $Y:=\mathbb{R}$ in the corresponding notation
and write, e.g., $C^{m}(X):=C^{m}(X\rightarrow \mathbb{R}).$ It should be
stressed that such type of directional differentiability (actually, even in
a weaker sense, along some total set of $h\in X$) will be quite enough for
our applications, and thus we do not discuss here the more involved notions
of G\^{a}teaux or Fr\'{e}chet derivatives. Starting from Subsect.\thinspace
3.2.3 we will mainly use the functional spaces $C^{m}(\Omega
_{-p}^{R})\supset C^{m}(\Omega _{-p;k}^{R}),$ $k\in \mathbb{Z}^{d},$ (and
their subspaces as described above) with the underlying space $X$ being in
this case $(\Omega _{-p}^{R},\rho _{-p,R})$ resp. $(\Omega _{-p,k}^{R},$ $%
||\cdot ||_{-p,R;k}).$

\subsubsection{The norm-function on $L_{\protect\beta }^{R}$ and $C_{\protect%
\beta }^{{}}$}

It is well-known that the norm-function $x\rightarrow \Phi (x):=|x|_{X}$ on
an arbitrary Banach space $X$ may be not differentiable. But as follows by a
simple convexity argument applied to (4.24), there always exist the partial%
\emph{\ }derivatives $\partial _{h}^{+}|x|_{X}^{{}}$ resp. $\partial
_{h}^{-}|x|_{X}^{{}}$ on the right resp. on the left at every point $x\in X$
and along every vector $h\in X.$ Since $\partial _{h}^{+}|x|_{X}^{{}}$ and $%
\partial _{h}^{-}|x|_{X}^{{}}$ do not necessarily coincide, one introduces
the \emph{subdifferential} $\partial |\cdot |_{X}:X\rightarrow X^{\ast }$ as
the (possibly multivalued) mapping 
\begin{equation}
\partial |x|_{X}:=\left\{ x^{\ast }\in X^{\ast }\left\vert |x^{\ast
}|_{X^{\ast }}=1,\text{ }(x,x^{\ast })=|x|_{X}\right. \right\} ,\quad
\forall x\in X,  \tag{4.25}
\end{equation}%
where $(\cdot ,\cdot )$ denotes the canonical pairing between $X$ and its
dual $X^{\ast }.$ Then for all $x,h\in X$%
\begin{equation}
-|h|_{X}\leq \partial _{h}^{-}|x|_{X}=\min_{x^{\ast }\in \partial
|x|_{X}}(x,x^{\ast })\leq \max_{x^{\ast }\in \partial |x|_{X}}(x,x^{\ast
})=\partial _{h}^{+}|x|_{X}\leq |h|_{X}.  \tag{4.26}
\end{equation}%
Besides, as is easy to see from the definition (4.24), for fixed $h\in X,$
the corresponding derivatives $\partial _{h}^{+}|x|_{X}$ resp. $\partial
_{h}^{-}|x|_{X}$ are semicontinuous (above resp. below) functions of $x\in
X. $

The following two examples will be of special importance for us. Consider
the norm-function in the spaces $X:=L_{\beta }^{R},$ $1<R<\infty .$ It is (Fr%
\'{e}chet) differentiable if $x\neq 0,$ that means $\partial
_{h}^{+}|x|_{L_{\beta }^{R}}^{{}}=\partial _{h}^{-}|x|_{L_{\beta }^{R}}^{{}}$
and $\partial |x|_{L_{\beta }^{R}}$ consists of the unique $x^{\ast
}:=|x|_{L_{\beta }^{R}}^{1-R}|x|^{R-2}x\in L_{\beta }^{R^{\prime }},$ $%
R^{\prime }=R(R-1)^{-1}.$ But this is not the case for $X:=C_{\beta },$ for
which it is well known that the norm-function $|\cdot |_{C_{\beta }}$ is not
(G\^{a}teaux) differentiable everywhere on $C_{\beta }\backslash \{0\}.$
Indeed,%
\begin{gather}
\exists \partial _{h}|x|_{C_{\beta }}=\partial _{h}^{\pm }|x|_{C_{\beta }}%
\text{ \ for }x,h\in C_{\beta }\text{ \ if and only if}  \notag \\
h(\tau )=\pm h(\tau ^{\prime })\text{ \ for all \ }\tau ,\tau ^{\prime }\in
S_{\beta }\text{ \ such that }x(\tau )=\pm x(\tau ^{\prime })=|x|_{C_{\beta
}}.  \tag{4.27}
\end{gather}%
Thus in the relevant calculations we will use either the right or left
derivatives $\partial _{h}^{+}|x|_{C_{\beta }^{{}}}^{{}}$ resp. $\partial
_{h}^{-}|x|_{C_{\beta }^{{}}}^{{}}$ and estimate them through the inequality
(4.26). For much information on related topics in convex analysis we refer,
e.g., to [Dei85].

\subsubsection{Cylinder functions on\textbf{\ }$\Omega $}

Since configurations of the system are described by sequences of continuous
loops over $\mathbb{Z}^{d},$ there are few natural notions of cylinder
functions w.r.t. the single spin space and the lattice structure:\smallskip

\textbf{(i)} \ (\emph{Cylindricity w.r.t. the lattice basis} $%
\{e_{k}\}_{k\in \mathbb{Z}^{d}}$) \ By $\mathcal{F}C^{m}(\Omega ;\mathbb{Z}%
^{d}),$ $m\in \mathbb{N}\cup \{0,+\infty \}$, we shall denote the set of all
local functions $f:\Omega \rightarrow \mathbb{R}$ which can be represented
as $f(\omega )=f_{\Lambda }(\mathbb{P}_{\Lambda }\omega )$ with some $%
f_{\Lambda }\in C^{m}(\Omega _{\Lambda })$ and $\Lambda \Subset \mathbb{Z}%
^{d}.\smallskip $

\textbf{(ii)} \ (\emph{Cylindricity w.r.t. the }$\emph{bas}$\emph{is} $%
\{h_{i}\}_{i\in \mathbb{Z}^{d+1}}$ \emph{in} $\mathcal{H}$) \ By $\mathcal{F}%
C^{m}(\Omega ;\mathcal{H})$ ($\subset \mathcal{F}C^{m}(\Omega ;\mathbb{Z}%
^{d})$) we shall respectively denote the set of all functions $f:\Omega
\rightarrow \mathbb{R}$ of the form%
\begin{equation}
f(\omega )=f_{L}(<\omega ,h_{i_{1}}>_{\mathcal{H}},...,<\omega ,h_{i_{L}}>_{%
\mathcal{H}}),  \tag{4.28}
\end{equation}%
with some $L\in \mathbb{N},$ $f_{L}\in C^{l}(\mathbb{R}^{L})$ and$\
i_{1},...,i_{L}\in \mathbb{Z}^{d+1}$. Then for the corresponding partial
derivatives in directions $h\in \Omega $ we have 
\begin{equation}
\partial _{h}f(\omega )=\sum_{l=1}^{L}\partial _{l}f_{L}(<\omega
,h_{i_{1}}>_{\mathcal{H}},...,<\omega ,h_{i_{L}}>_{\mathcal{H}%
})<h,h_{i_{l}}>_{\mathcal{H}}.  \tag{4.29}
\end{equation}%
Replacing $C^{m}(\mathbb{R}^{L})$ in (4.28) by $C_{b}^{m}(\mathbb{R}^{L})$
(resp. $C_{b,loc}^{m}(\mathbb{R}^{L})$), we get the subsets $\mathcal{F}%
C_{b}^{m}(\Omega ;\mathcal{H})$ (resp. $\mathcal{F}C_{b,loc}^{m}(\Omega ;%
\mathcal{H}))$.\smallskip

\textbf{(iii)} \ (\emph{Space-time cylindricity}) \ Note that for any fixed $%
k\in \mathbb{Z}^{d}$ and $\tau \in S_{\beta }$, 
\begin{equation*}
<\omega ,\delta _{(k,\tau )}>_{\mathcal{H}}:=\omega _{k}(\tau ),\quad \omega
=(\omega _{k})_{k\in \mathbb{Z}^{d}}\in \Omega ,
\end{equation*}%
is a well-defined bounded linear functional $\delta _{(k,\tau )}\in $ $%
\Omega ^{\ast }.$ Hence, we can introduce the set $\mathcal{F}C^{m}(\Omega ;%
\mathbb{Z}^{d}\times S_{\beta })$ (and also its subsets $\mathcal{F}%
C_{b}^{m}(\Omega ;\mathbb{Z}^{d}\times S_{\beta })$, $\mathcal{F}%
C_{b,loc}^{m}(\Omega ;\mathbb{Z}^{d}\times S_{\beta })$) consisting of all
cylinder (w.r.t. the family $\{\delta _{(k,\tau )}\}_{(k,\tau )\in \mathbb{Z}%
^{d}\times S_{\beta }^{{}}}\subset \Omega ^{\ast }$) functions $f:\Omega
\rightarrow \mathbb{R}$ of the form%
\begin{equation}
f(\omega )=f_{L}\left( \omega _{k_{1}}(\tau _{1}),...,\omega _{k_{L}}(\tau
_{L})\right) =f_{L}\left( <\omega ,\delta _{(k_{1},\tau _{1})}>_{\mathcal{H}%
},...,<\omega ,\delta _{(k_{L},\tau _{L})}>_{\mathcal{H}}\right) , 
\tag{4.30}
\end{equation}%
where $L\in \mathbb{N},\ f_{L}\in C^{m}(\mathbb{R}^{L})$ and $%
(k_{1}^{{}},\tau _{1}^{{}}),...,(k_{L}^{{}},\tau _{L}^{{}})\in \mathbb{Z}%
^{d}\times S_{\beta }^{{}}.$ As can be checked by direct calculation, for
any $h\in \Omega $ there exists 
\begin{equation}
\partial _{h}f(\omega )=\sum_{l=1}^{L}\partial _{l}f_{L}\left( \omega
_{k_{1}}(\tau _{1}),...,\omega _{k_{L}}(\tau _{L})\right) \cdot
h_{k_{l}}(\tau _{l}).  \tag{4.31}
\end{equation}%
Obviously, $\mathcal{F}C^{m}(\Omega ;\mathcal{H}),$ $\mathcal{F}C^{m}(\Omega
;\mathbb{Z}^{d}\times S_{\beta })\subset C^{m}(\Omega _{-p}^{R})$ for all $%
p\geq 0,$ $R\geq 1.$

\subsubsection{Approximations by smooth functions}

Since we are dealing with functions on infinite-dimensional vector spaces
(by the way, in our case $\mathcal{F}C_{b}(\Omega ;\mathcal{H})\cap
C_{0}(\Omega _{-p;k}^{-R})=\emptyset $ and the norms $||\cdot ||_{-p,R;k}$
defining the topology in $\Omega _{-p}^{R}$ are not differentiable), a
possible approximation of $f:\Omega \rightarrow \mathbb{R}$ by smooth resp.
boundedly supported functions should be treated carefully. For the sake of
completeness, here we collect several technical assertions (more or less of
common knowledge and even true in a more general setting on locally convex
spaces), which will be used below when extending the (IbP)-formula (4.39) to
suitable classes of differentiable functions on $\Omega .\smallskip $\newline
\textbf{Lemma 4.3 \ }\emph{Let} $\mu $ \emph{be a Borel measure on} $\Omega $
\emph{such that} $\mu (\Omega _{-p}^{R})=1$ \emph{for some} $p>d,$ $R\geq 2.$
\emph{Then the following assertions hold:}

\textbf{(i) (}$C^{1}$-Approximation by Smooth Cylinder Functions) \ \emph{%
Given any }$h_{i}\in bas(\mathcal{H})$ \emph{and} $f\in C_{b}^{1}(\Omega
_{-p}^{R};h_{i})$,\emph{\ in each of the spaces }$\mathcal{F}C^{\infty
}(\Omega ;\mathcal{H})$\emph{\ and }$\mathcal{F}C^{\infty }(\Omega ;\mathbb{Z%
}^{d}\times S_{\beta })$ \emph{there exist corresponding sequences }$%
\{f^{(K)}\}_{K\in \mathbb{N}}\mathcal{\ }$\emph{such that} 
\begin{gather}
\inf f\leq f^{(K)}\leq \sup f,\quad ||\partial _{h_{i}}f^{(K)}||_{L^{\infty
}}^{{}}\leq 2||\partial _{h_{i}}f||_{L^{\infty }}^{{}},\quad \text{\emph{and}%
}  \notag \\
\lim_{K\rightarrow \infty }\left( ||f^{(K)}-f||_{L^{q}(\mu
)}^{{}}+||\partial _{h_{i}}f^{(K)}-\partial _{h_{i}}f||_{L^{q}(\mu
)}^{{}}\right) =0\text{\quad }\emph{for\ all\ }1\leq q<\infty .  \tag{4.32}
\end{gather}

\textbf{(ii)} (Approximation by Boundedly Supported Functions) \ \emph{Given
any }$h\in \Omega _{-p}^{R}$\emph{\ and} $f\in C_{0}(\Omega _{-p;k}^{R})$,%
\emph{\ there exists a sequence} $\{f^{(K)}\}_{K\in \mathbb{N}}\subset
C_{0}^{1}(\Omega _{-p;k}^{R};h)$ \emph{such that} 
\begin{gather}
\inf f\leq f^{(K)}\leq \sup f,\text{\quad }f^{(K)}\underset{K\rightarrow
\infty }{\longrightarrow }f\text{\quad \emph{pointwise on} }\Omega _{-p}^{R}
\notag \\
\text{\emph{and thus\quad }}\lim_{K\rightarrow \infty
}||f^{(K)}-f||_{L^{q}(\mu )}^{{}}=0\text{\quad }\emph{for\ all\ }1\leq
q<\infty .  \tag{4.33}
\end{gather}%
\emph{If, additionally, }$f$ \emph{has a bounded right derivative }$\partial
_{h}^{+}f:\Omega _{-p}^{R}\rightarrow \mathbb{R},$\emph{\ then also} 
\begin{gather}
|\partial _{h}f^{(K)}|\leq \sup |\partial _{h}^{+}f|,\text{\quad }\partial
_{h}^{{}}f^{(K)}\underset{K\rightarrow \infty }{\longrightarrow }\partial
_{h}^{+}f\text{\quad \emph{pointwise on} }\Omega _{-p}^{R}  \notag \\
\text{\emph{and thus\quad }}\lim_{K\rightarrow \infty }||\partial
_{h}f^{(K)}-\partial _{h}^{+}f||_{L^{q}(\mu )}^{{}}=0\text{\quad \emph{for\
all }}1\leq q<\infty .  \tag{4.34}
\end{gather}

\textbf{(iii) (}$L^{q}$\emph{-}Approximation) \ \emph{Given any} $h\in
\Omega _{-p}^{R},$ \emph{each of the sets }$\mathcal{F}C_{b}^{\infty
}(\Omega ;\mathbb{Z}^{d}\times S_{\beta }),$ $\mathcal{F}C_{b}^{\infty
}(\Omega ;\mathcal{H})$ \emph{and} $C_{0}^{1}(\Omega _{-p;k}^{R};h)$ \emph{%
is dense in all spaces} $L^{q}(\mu ),$ $1\leq q<\infty $. \emph{Moreover,} 
\emph{for }$0\leq f\in L^{\infty }(\mu )$ \emph{the corresponding sequences }%
$\{f^{(K)}\}_{K\in \mathbb{N}}$\emph{, }$\lim_{K\rightarrow \infty
}||f^{(K)}-f||_{L^{q}(\mu )}^{{}}=0,$ \emph{can be chosen so that} 
\begin{equation}
0\leq \inf f\leq f^{(K)}\leq \sup f\quad \text{\emph{and}\quad }%
\lim_{K\rightarrow \infty }f^{(K)}=f\text{\quad }(\mu -\text{a.e.}). 
\tag{4.35}
\end{equation}

\textbf{Proof:} \textbf{(i) \ }First, we set 
\begin{equation*}
g^{(K)}:=f(\mathbb{P}_{\Lambda _{K}}(\mathbb{M}_{K}))\in \mathcal{F}%
C_{b}(\Omega ;\mathbb{Z}^{d+1}))\cap C_{b}^{1}(\Omega _{-p}^{R};h_{i}),
\end{equation*}%
where 
\begin{equation*}
(\mathbb{M}_{K}(\omega ))_{k}:=\mathbb{M}_{K}(\omega _{k}),\quad (\mathbb{P}%
_{\Lambda _{K}}\omega )_{k}:=\left\{ 
\begin{tabular}{rr}
$\omega _{k},$ & $k\in \Lambda _{K}$ \\ 
$0,$ & $\text{otherwise}$%
\end{tabular}%
\ \right. ,
\end{equation*}%
and $\Lambda _{K},$ $K\in \mathbb{N},$ are bounded domains exhausting $%
\mathbb{Z}^{d}.$ Based on the properties (3.26)--(3.28) of the Ces\`{a}ro
partial sums $\mathbb{M}_{K}(\omega _{k})$ in the space of continuous loops $%
C_{\beta }$, it is easy to check that $\{g^{(K)}\}_{K\in \mathbb{N}}$
satisfies 
\begin{gather}
\inf f\leq g^{(K)}\leq \sup f,\quad ||\partial _{h_{i}}g^{(K)}||_{L^{\infty
}}\leq ||\partial _{h_{i}}f||_{L^{\infty }}\,,\text{\quad and}  \notag \\
g^{(K)}\rightarrow f,\text{ }\partial _{h_{i}}g^{(K)}\rightarrow \partial
_{h_{i}}f\ \text{pointwise on }\Omega _{-p}^{R}\text{ when }K\rightarrow
\infty .  \tag{4.36}
\end{gather}%
Next, noting that any cylinder function $g\in \mathcal{F}C(\Omega ;\mathcal{H%
})$ is of the form (4.29) with some $g_{L}\in C(\mathbb{R}^{L})$ and using a
standard mollifier argument on $\mathbb{R}^{L},$ for each $g^{(K)}$ defined
above one can construct an approximating sequence (in the sense of (4.36)) $%
\{f^{(N)}\}_{N\in \mathbb{N}}:=\{g^{(K,N)}\}_{N\in \mathbb{N}}\subset 
\mathcal{F}C_{b}^{\infty }(\Omega ;\mathcal{H})$. And finally, we take into
account that for any $\omega \in \Omega $ the Riemann integral $<\omega
,h_{i}>_{\mathcal{H}}:=\int_{S_{\beta }}\omega _{k}(\tau )\varphi _{n}(\tau
)d\tau $ can be approximated by its partial sums in the following way:%
\begin{gather*}
<\omega ,h_{i}>_{\mathcal{H}}=\lim_{M\rightarrow \infty }\frac{\beta }{M}%
\sum_{m=0}^{M-1}\varphi _{n}(\frac{m}{M}\beta )\cdot <\omega ,\delta _{(k,%
\frac{m}{M}\beta )}>_{\mathcal{H}}^{{}}, \\
\partial _{h_{i}}<\omega ,h_{i}>_{\mathcal{H}}=\lim_{M\rightarrow \infty }%
\frac{\beta }{M}\sum_{m=0}^{M-1}\varphi _{n}^{2}(\frac{m}{M}\beta )=|\varphi
_{n}|_{H}^{2}=1, \\
\frac{\beta }{M}\sum_{m=0}^{M-1}\varphi _{n}^{2}(\frac{m}{M}\beta )\leq
\beta \cdot |\varphi _{n}|_{L_{\beta }^{\infty }}^{2}\leq 2,\quad \forall
M,n\in \mathbb{N}.
\end{gather*}%
By Lebesgue's dominated convergence theorem, this latter also enables us to
approximate each $g^{(K,N)}\in \mathcal{F}C_{b}^{\infty }(\Omega ;\mathcal{H}%
)$ by $\{g^{(K,N,M)}\}_{M\in \mathbb{N}}\subset \mathcal{F}C_{b}^{\infty
}(\Omega ;\mathbb{Z}^{d}\times S_{\beta })$ in the sense of (4.32)$%
.\smallskip $

\textbf{(ii)} \ Here one can put, for instance, 
\begin{gather*}
f^{(K)}(\omega ):=K\int_{0}^{1/K}f(\omega +\theta h)d\theta \\
\text{with \ }\partial _{h}f^{(K)}(\omega )=K[f(\omega +\frac{1}{K}%
h)-f(\omega )].
\end{gather*}%
Then it is straightforward to show that the functions $\{f^{(K)}\}_{K\in 
\mathbb{N}}\subset C_{0}^{1}(\Omega _{-p;k}^{R};h)$ satisfy the required
assumptions (4.33) and (4.34).\smallskip

\textbf{(iii)} \ In this respect we recall the regularity property of Borel
measures $\mu \in \mathcal{M}(X)$ on a Polish space $X$, according to which $%
\mu (\Delta )=\sup \{\mu (B)|B\subseteq \Delta ,$ $B$ is closed\}, $\forall
\Delta \in \mathcal{B}(X)$ (cf. [RS72, Subsect.\thinspace IV.4]). In our
situation this standardly implies that the set $C_{b}(\Omega _{-p;k}^{R})$
(and thus, by an obvious cut-off argument, also its subset $C_{0}(\Omega
_{-p;k}^{R})$) is dense in all $L^{q}(\mu ),$ $1\leq q<\infty ,$ as well as
in $L^{\infty }(\mu )$ (the latter is only in the sense of (4.35)). Starting
from this point, in order to construct the desired approximation of $f\in
C_{0}(\Omega _{-p;k}^{R})$ by functions from $\mathcal{F}C_{b}^{\infty
}(\Omega ;\mathcal{H}),$ $\mathcal{F}C_{b}^{\infty }(\Omega ;\mathbb{Z}%
^{d}\times S_{\beta })$ or resp. $C_{0}^{1}(\Omega _{-p;k}^{R};h),$ one can,
e.g., repeat the previous arguments from the proof of Assertions (i) resp.
(ii).

\noindent $\blacksquare $

\subsection{Partial logarithmic derivatives}

Given a probability measure $d\mu (x)=\exp \{-\rho (x)\}dx$ on $\mathbb{R}%
^{n}$ with (smooth enough) density $\rho ,$ the vector field $%
b=(b_{l})_{l=1}^{n}:=-\nabla \rho :\mathbb{R}^{n}\rightarrow \mathbb{R}^{n}$
is usually called a logarithmic gradient of $\mu $. Equivalently, one can
define its components $b_{l}$, the so-called partial logarithmic derivatives
along the basis vectors $e_{l},$ through 
\begin{equation}
b_{l}(x):=\left( \text{ln}\frac{d\mu (x+\theta e_{l})}{d\mu (x)}\right)
_{\theta =0}^{\prime },\quad x=(x_{l})_{l=1}^{n}=\sum_{l=1}^{n}x_{l}e_{l}\in 
\mathbb{R}^{n}.  \tag{4.37}
\end{equation}%
Whereas for Gibbs measures $\mu \in \mathcal{G}_{t}$ the presentation (3.21)
with the Euclidean \textquotedblright density\textquotedblright\ $\mathcal{I}%
(\omega )$ has a heuristic sense only, below we will use their rigorous
description via Radon--Nikodym derivatives (4.7) in order to modify
definition (4.37) to the infinite-dimensional case.\smallskip

Namely, having regard to Lemma 4.1, we define the\emph{\ partial logarithmic
derivative} of (\emph{all}) measures $\mu \in \mathcal{G}_{t}$ along the\
fixed direction $h_{i}=e_{k}\otimes \varphi _{n},$ $i=(k,n)\in \mathbb{Z}%
^{d+1},$ as a mapping $b_{i}:\Omega _{t}\rightarrow \mathbb{R},$%
\begin{equation}
b_{i}(\omega ):=\partial _{h_{i}}a_{h_{i}}(\omega )=(a_{\theta h_{i}}(\omega
))_{\theta =0}^{\prime }=-(A_{\beta }^{{}}\varphi _{n},\omega
_{k})_{H}-(F_{k}^{V,W}(\omega ),\varphi _{n})_{H}.  \tag{4.38}
\end{equation}%
Here $F_{k}^{V,W}:$ $\Omega _{t}\rightarrow L_{\beta }^{R^{\prime }}$ (with $%
R^{-1}+(R^{\prime })^{-1}=1$) is a nonlinear Nemytskii-type operator acting
by 
\begin{multline}
F_{k}^{V,W}(\omega ):=F_{k}^{V_{k}}(\omega _{k})+F_{k}^{W}(\omega )  \notag
\\
:=V_{k}^{\prime }(\omega _{k})+\sum\limits_{M=2}^{N}\sum\limits_{\substack{ %
\{k_{1},...,k_{M}\}\subset \mathbb{Z}^{d}  \\ k_{1}:=k}}\partial
_{1}W_{\{k_{1},...,k_{M}\}}(q_{k_{1}},...,q_{k_{M}})|_{\substack{ %
q_{l}=\omega _{l}  \\ 1\leq l\leq M}}\text{ }.  \tag{4.39}
\end{multline}%
\smallskip \newline
\textbf{Lemma 4.4} \ \emph{For all }$i=(k,n),$ $i^{\prime }=(k^{\prime
},n^{\prime })\in \mathbb{Z}^{d+1}$ \emph{and} $p>d$%
\begin{equation}
F_{k}^{V}\in C_{b,loc}^{1}(C_{\beta }\rightarrow C_{\beta };\varphi
_{n^{\prime }}),\text{ \ }F_{k}^{W}\in C_{b,loc}^{1}(\mathcal{L}%
_{-p}^{R}\rightarrow L_{\beta }^{R^{\prime }};h_{i^{\prime }}),  \tag{4.40}
\end{equation}%
\emph{and thus the partial logarithmic derivatives (4.38) give rise to
smooth mappings} 
\begin{equation}
b_{i}\in C_{b,loc}^{1}(\Omega _{-p;k}^{R};h_{i^{\prime }}),\quad \forall p>d.
\tag{4.41}
\end{equation}

\textbf{Proof:} \ For $V\in C_{b,loc}^{2}(\mathbb{R})$ the properties of the
operators $F_{k}^{V}$ in the space $C_{\beta }$ are well known: namely, they
are (Fr\'{e}chet) differentiable and 
\begin{equation*}
\partial _{\varphi }F_{k}^{V}(\upsilon ):=V_{k}^{\prime \prime }(\upsilon
)\cdot \varphi ,\quad \forall \varphi \in C_{\beta }^{{}}.
\end{equation*}%
Fixed $k,k^{\prime }\in \mathbb{Z}^{d},$ along with $F_{k}^{W}$ let us
introduce one more Nemytskii-type operator $\partial _{k^{\prime
}}F_{k}^{V,W}:\Omega _{t}\rightarrow L_{\beta }^{R^{\prime }}$ by%
\begin{multline}
\partial _{k^{\prime }}F_{k}^{V,W}(\omega ):=\partial _{k^{\prime
}}F_{k}^{V}(\omega _{k})+\partial _{k^{\prime }}F_{k}^{W}(\omega )  \notag \\
:=\left\{ 
\begin{array}{cc}
V_{k}^{\prime \prime }(\omega _{k})+\sum\limits_{M=2}^{N}\sum\limits
_{\substack{ \{k_{1},...,k_{M}\}\subset \mathbb{Z}^{d} \\ k_{1}:=k}}\partial
_{1}^{2}W_{\{k_{1},...,k_{M}\}}(\omega _{k_{1}},...,\omega _{k_{M}})|
_{\substack{ q_{l}=\omega _{l} \\ 1\leq l\leq M}}, & k=k^{\prime }, \\ 
\sum\limits_{M=2}^{N}\sum\limits_{\substack{ \{k_{1},...,k_{M}\}\subset 
\mathbb{Z}^{d} \\ k_{1}:=k,\text{ }k_{2}:=k^{\prime }}}\partial
_{1,2}^{2}W_{\{k_{1},...,k_{M}\}}(\omega _{k_{1}},...,\omega _{k_{M}})|
_{\substack{ q_{l}=\omega _{l} \\ 1\leq l\leq M}}, & k\neq k^{\prime }.%
\end{array}%
\right.   \tag{4.42}
\end{multline}%
Then by ($\mathbf{W}_{\text{\textbf{ii,iii}}}$), ($\mathbf{J}$) and Lemma
3.4 applied to the RHS\ in (4.39) and (4.42) we get that 
\begin{multline}
\max \left\{ |F_{k}^{W}(\omega )|_{L_{\beta }^{R^{\prime }}},|\partial
_{k^{\prime }}F_{k}^{W}(\omega )|_{L_{\beta }^{R^{\prime }}}\right\} \leq
2^{R-1}\left[ \sum\limits_{j\neq k}\tilde{J}_{k,j}|\omega _{j}|_{L_{\beta
}^{R}}^{R-1}+||J||_{0}\left( |\omega _{j}|_{L_{\beta }^{R}}^{R-1}+|\mathbf{1}%
|_{L_{\beta }^{R^{\prime }}}\right) \right]   \notag \\
\leq 2^{R-1}\left[ ||\mathbf{J}||_{p(R-1)}(1+|k|)_{{}}^{p(R-1)}||\omega
||_{-p,R}^{R-1}+||\mathbf{J}||_{0}|\mathbf{1}|_{L_{\beta }^{R^{\prime }}}%
\right]   \tag{4.43}
\end{multline}%
with the matrix $\{\tilde{J}_{k,j}\}_{k,j\in \mathbb{Z}^{d}}$ defined by
(3.45). Hence, obviously, 
\begin{multline}
\max \left\{ |(F_{k}^{W}(\omega ),\varphi _{n})_{H}^{{}}|,\text{ }|(\partial
_{h_{i}}F_{k}^{W}(\omega ),\varphi _{n})_{H}^{{}}|\right\}   \notag \\
\leq 2^{R-1}\kappa _{R}\left[ \sum\nolimits_{j\in \mathbb{Z}^{d}}\tilde{J}%
_{k,j}|\omega _{j}|_{L_{\beta }^{R}}^{R}+||\mathbf{J}||_{0}\left( |\omega
_{k}|_{L_{\beta }^{R}}^{R}+\beta +1\right) \right] .  \tag{4.44}
\end{multline}%
Since the series in the RHS of (4.43) converges uniformly on balls in $%
\mathcal{L}_{-p}^{R},$ this standardly yields (cf., e.g., [Ca67,
Theorem\thinspace 3.6.2]) that%
\begin{equation}
F_{k}^{W}\in C_{b,loc}^{1}(\mathcal{L}_{-p}^{R}\rightarrow L_{\beta
}^{R^{\prime }};h_{i^{\prime }})\quad \text{with\quad }\partial
_{h_{i^{\prime }}}F_{k}^{W}(\omega )=\partial _{k^{\prime }}F_{k}^{W}(\omega
)\cdot \varphi _{n^{\prime }}.  \tag{4.45}
\end{equation}%
Thus, from the previous discussion it is evident that $b_{i^{\prime }}\in
C_{b,loc}^{1}(\Omega _{-p;k}^{R};h_{i^{\prime }})$ with%
\begin{equation}
\partial _{h_{i^{\prime }}}b_{i}=-\delta _{k-k^{\prime }}\left( \lambda
_{n}+V_{k}^{\prime \prime }(\omega _{k}),\varphi _{n}\varphi _{n^{\prime
}}\right) _{H}-\left( \partial _{k^{\prime }}F_{k}^{W}(\omega ),\varphi
_{n}\varphi _{n^{\prime }}\right) _{H}^{{}}.  \tag{4.46}
\end{equation}%
In particular, taking into account (3.25) and (4.43), one has merely the
following bounds on the growth of the logarithmic derivatives:%
\begin{equation*}
|b_{i}(\omega )|\leq \kappa _{\infty }\left\{ \lambda _{n}|\omega
_{k}|_{L_{\beta }^{1}}+|F_{k}^{V,W}(\omega )|_{L_{\beta }^{1}}\right\} \leq
C_{i}\left\{ 1+|\omega _{k}|_{L_{\beta }^{1}}+|V_{k}^{\prime }(\omega
_{k})|_{L_{\beta }^{1}}+||\omega ||_{-p,R}^{R-1}\right\} ,
\end{equation*}%
\begin{equation}
|\partial _{h_{i^{\prime }}}b_{i}(\omega )|\leq \lambda _{n}+\kappa _{\infty
}^{2}|\partial _{k^{\prime }}F_{k}^{V,W}(\omega )|_{L_{\beta }^{1}}\leq
C_{i}^{\prime }\left\{ 1+|V_{k}^{\prime \prime }(\omega _{k})|_{L_{\beta
}^{1}}+||\omega ||_{-p,R}^{R-1}\right\}   \tag{4.47}
\end{equation}%
with some absolute constants $C_{i},C_{i}^{\prime }\in (0,\infty )$
depending, besides $i\in \mathbb{Z}^{d+1}$, also on chosen $p>d.$

\noindent $\blacksquare \smallskip $

Combining Lemmas 4.1 and 4.4, by a straightforward calculation one get\ the
following useful relations between the Radon-Nikodym and logarithmic
derivatives\textbf{.}\smallskip

\noindent \textbf{Corollary 4.5 \ }\emph{Under the assumptions on the
interaction potentials imposed in Definition 3.1}%
\begin{equation}
\partial _{h_{i}}a_{\theta h_{i}}(\omega )=a_{\theta h_{i}}(\omega )\left[
b_{i}(\omega +\theta h_{i})-b_{i}(\omega )\right] ,\text{ \ }\frac{\partial 
}{\partial \theta }a_{\theta h_{i}}(\omega )=a_{\theta h_{i}}(\omega
)b_{i}(\omega +\theta h_{i}),  \tag{4.48}
\end{equation}%
\emph{and thus} $a_{\theta h_{i}}$ \emph{can be recovered from} $b_{i}$ 
\emph{through} 
\begin{equation}
a_{\theta h_{i}}(\omega )=\exp \int\nolimits_{0}^{\theta }b_{i}(\omega
+\vartheta h_{i})d\vartheta .  \tag{4.49}
\end{equation}%
\smallskip \noindent \textbf{Remark 4.6} \textbf{(i)} \ In relation to the
further applications in Subsect.\thinspace 7 below, it is important to note
that the pointwise \emph{coercivity} \emph{and growth} assumptions on the
one-particle potentials $V_{k}$ imply the corresponding properties for the
operators $F_{k}^{V}$ w.r.t. the tangent space $H:=L_{\beta }^{2}.$ Namely,
from ($\mathbf{V}_{\text{\textbf{i--v}}}$) we have uniformly for all $%
i=(k,n)\in \mathbb{Z}^{d+1}$ and $\omega \in \Omega _{-p}^{R}$:%
\begin{equation}
(F_{k}^{V}(\omega _{k}),\omega _{k})_{H}\geq \max \left\{ 
\begin{array}{r}
K_{1}^{-1}\left[ \kappa _{\infty }^{-1}|(F_{k}^{V}(\omega _{k}),\varphi
_{n})_{H}^{{}}|+\kappa _{\infty }^{-2}|(\partial _{n}F_{k}^{V}(\omega
_{k}),\varphi _{n})_{H}^{{}}|-L_{1}\beta \right] \\ 
K_{2}^{-1}\left[ \kappa _{\infty }^{-1}|(\partial _{n}F_{k}^{V}(\omega
_{k}),\omega _{k})_{H}^{{}}|-L_{2}\beta \right] \\ 
K_{1}^{-1}\left[ |F_{k}^{V}(\omega _{k}|_{L_{\beta }^{1}}-L_{1}\beta \right]
\\ 
K_{3}^{-1}[|\omega _{k}|_{L_{\beta }^{R}}^{R}-L_{3}\beta ] \\ 
K_{4}^{-1}\left[ |\omega _{k}|_{H}^{2}-L_{4}\beta \right]%
\end{array}%
\right\}  \tag{4.50}
\end{equation}%
and 
\begin{equation}
|(\partial _{\varphi _{n}}F_{k}^{V}(\omega _{k})|_{L_{\beta }^{R^{\prime
}}}\leq \kappa _{\infty }\left[ K_{0}\left( |F_{k}^{V}(\omega
_{k})|_{L_{\beta }^{R^{\prime }}}+|\omega _{k}|_{L_{\beta
}^{R}}^{R-1}\right) +L_{0}|\mathbf{1}|_{L_{\beta }^{R^{\prime }}}^{{}}\right]
.  \tag{4.51}
\end{equation}%
On the other hand, together with Young's inequality (4.43) implies the
following upper bound for the mappings $F_{k}^{W}$:%
\begin{multline}
\max \left\{ |F_{k}^{W}(\omega )|_{L_{\beta }^{R^{\prime }}},\text{ }%
|(F_{k}^{W}(\omega ),\omega _{k})_{H}^{{}}|,\text{ }\kappa _{\infty
}^{-1}\cdot |\partial _{h_{i}}F_{k}^{W}(\omega )|_{L_{\beta }^{R^{\prime }}},%
\text{ }\kappa _{\infty }^{-1}\cdot |(\partial _{h_{i}}F_{k}^{W}(\omega
),\omega _{k})_{H}^{{}}|\right\}  \notag \\
\leq 2^{R}\left[ \sum\nolimits_{j\in \mathbb{Z}^{d}}\tilde{J}_{k,j}|\omega
_{j}|_{L_{\beta }^{R}}^{R}+||\mathbf{J}||_{0}\left( |\omega _{k}|_{L_{\beta
}^{R}}^{R}+\beta +1\right) \right] .  \tag{4.52}
\end{multline}%
These estimates will be essential, e.g., for the proof of Theorems 7.5 and
7.6 below.\smallskip \newline

\textbf{(ii)} \ Of course, one can also define the \textquotedblright \emph{%
generalized}\textquotedblright\ logarithmic gradient of the measures $\mu
\in \mathcal{G}_{t}$ as a measurable vector field $b:=(b_{k})_{k\in \mathbb{Z%
}^{d}}$ with components 
\begin{equation*}
\Omega _{-p}^{R}\ni \omega \rightarrow b_{k}(\omega ):=\sum\limits_{n\in 
\mathbb{Z}}b_{(k,n)}(\omega )\cdot \varphi _{n}=-A_{\beta }\omega
_{k}-F_{k}^{V,W}(\omega )\in W_{\beta }^{2,-2}.
\end{equation*}%
But even though the important for the sequel coercivity estimate (4.47) for
its nonlinear components $F_{k}^{V,W}$ holds, it can not help us to control
the coercivity properties of the logarithmic gradient $b$ \emph{as a whole }%
w.r.t. \ the given tangent space $\mathcal{H}:=l^{2}(\mathbb{Z}^{d})\otimes
L_{\beta }^{2}$ (in contrast with the technique suggested for classical
Gibbs measures $\mu _{cl}$ on $\mathcal{S}^{\prime }(\mathbb{Z}^{d})$ in
[AKRT99,00]). The reason is that, taken for every $k\in \mathbb{Z}^{d},$ the
bilinear form $(A_{\beta }^{{}}\omega _{k},\omega _{k})_{H}$ is defined only
on a $\mu -$measure zero set of $\omega _{k}\in W_{\beta }^{2,1}$. Thus in
the next Sects.\thinspace 5 and 6 we \emph{separately }have to do a
\textquotedblright \emph{lattice analysis}\textquotedblright , analogous to
the one for classical Gibbs states and depending on the coercivity
properties in $\mathcal{H}$ of the Nemytskii operators $F_{k}^{V,W}$, as
well as an additional \textquotedblright \emph{single spin space analysis}%
\textquotedblright , taking into account the spectral properties in $%
L_{\beta }^{2}$ of the elliptic operator $\emph{A}_{\beta }$ in the linear
part of the logarithmic derivatives (4.38).

\subsection{Integration by parts (IbP) description of Euclidean Gibbs
measures}

Next we shall show that, if the interaction potentials are smooth (as they
are in our case), the flow characterization (3.9) of $\mu \in \mathcal{G}%
_{t} $ is equivalent to their definition as differentiable measures
satisfying (in the distributional sense by pairing both sides with proper
test functions $f)$ the integration by parts (for short,\emph{\ IbP})
formulas 
\begin{equation}
\partial _{h_{i}}\mu (d\omega )=b_{i}(\omega )\mu (d\omega ),\text{ \ }i\in 
\mathbb{Z}^{d+1},  \tag{4.53}
\end{equation}%
with the given by (4.38) logarithmic derivatives $b_{i}$ along basis vectors 
$h_{i}\in bas(\mathcal{H})$.

So, we return to formula (4.23). The principal problem (remaining open until
Theorem 7.6 below) is here that \emph{we do not know a priori whether} $%
b_{i}\in L^{1}(\mu ).$ This information is needed before a standard theory
of differentiable measures (see, e.g., [Bo97]) can be applied. Thus, based
on the previous discussion in Subsect.\thinspace 4.2, we first have to find
a large enough amount of functions $f$ (which, in fact, will be $\mu -$a.e.
differentiable according to Corollary 4.11 (ii)), such that certainly $%
fb_{i}\in L^{\infty }(\mu )$ and passage to the limit in the both sides of
(4.23) is correct.\smallskip

\noindent \textbf{Definition 4.7} \ \emph{Fixing a basis vectors} $h_{i}$, $%
i=(k,n)\in \mathbb{Z}^{d+1},$ \emph{we define }(conventionally in this paper)%
\emph{\ a set of test functions}%
\begin{equation}
C_{b}^{1,\pm }(\Omega _{-p}^{R};h_{i}):=\left\{ f:\Omega
_{-p}^{R}\rightarrow \mathbb{R}\left\vert 
\begin{array}{c}
f(\omega )=g(\omega )\chi _{2}(|\omega _{k}|_{C_{\beta }},||\omega ||_{-p,R})
\\ 
\text{with }\forall g\in C^{1}(\Omega _{-p}^{R};h_{i}),\text{ }\forall \chi
_{2}\in C^{1}(\mathbb{R}^{2}), \\ 
\sup_{\Omega _{-p}^{R}}\left( |f|+|\partial _{h_{i}}^{+}f|+|\partial
_{h_{i}}^{-}f|\right) <\infty%
\end{array}%
\right. \right\} .  \tag{4.54}
\end{equation}%
\emph{Its subsets} $C_{dec}^{1,\pm }(\Omega _{-p}^{R};h_{i})\supset
C_{0}^{1,\pm }(\Omega _{-p}^{R};h_{i})$\emph{\ are defined as below}$:$ 
\emph{Namely,} \emph{let} $C_{dec}^{1,\pm }(\Omega _{-p}^{R};h_{i})$ \emph{%
consist of those }$f\in C_{b}^{1,\pm }(\Omega _{-p}^{R};h_{i})$ \emph{which
satisfy the} \emph{extra decay condition}%
\begin{equation}
\sup_{\omega \in \Omega _{-p}^{R}}\left\vert f(\omega )\left( 1+|\omega
_{k}|_{L_{\beta }^{R^{\prime }}}+|F_{k}(\omega )|_{L_{\beta }^{R^{\prime
}}}\right) \right\vert <\infty .  \tag{4.55}
\end{equation}%
\emph{Respectively,} $C_{0}^{1,\pm }(\Omega _{-p}^{R};h_{i})$ \emph{consists
of those (boundedly supported) functions} $f\in C_{b}^{1,\pm }(\Omega
_{-p}^{R};h_{i})$ \emph{which can be written in the form}%
\begin{equation}
f(\omega )=g(\omega )\cdot \chi _{1}(||\omega ||_{-p,R;k})\text{ \emph{with} 
}\forall g\in C_{b}^{1}(\Omega _{-p}^{R};h_{i}),\text{ }\forall \chi _{1}\in
C_{0}^{1}(\mathbb{R}).  \tag{4.56}
\end{equation}%
\smallskip

Note that, by the chain rule and relation (4.26), for any $f\in C_{0}^{1,\pm
}(\Omega _{-p}^{R};h_{i})$ the derivatives on the right $\partial
_{h_{i}}^{+}f$ and on the left $\partial _{h_{i}}^{-}f$ exist as globally
bounded functions on $\Omega _{-p}^{R}$ (but \emph{not necessarily
continuous, }i.e.,\emph{\ }$\partial _{h_{i}}^{+}f\notin C(\Omega _{-p}^{R})$%
)$.$ If $\mu (\Omega _{-p}^{R})=1,$ then obviously every set $C_{0}^{1,\pm
}(\Omega _{-p}^{R};h_{i})\ $is dense in all $L^{q}(\mu ),$ $q\in \lbrack
1,\infty ),$ (and also in $L^{\infty }(\mu )$ w.r.t. the $\mu -$a.e.
convergence (4.35)) due to the corresponding density property of the sets $%
C_{0}^{1}(\Omega _{-p;k}^{R};h_{i})\subset C_{b}^{1}(\Omega _{-p}^{R})$ (cf.
Lemma 4.3 (iii)).\smallskip 

Having introduced such special classes of test functions, now we are able to
prove some (in fact, preliminary) version of the (IbP)-characterization\ for
Euclidean Gibbs measures given by Proposition 4.8. It is inspired by the
well known fact (cf., e.g., [Sk84, Be85, DaS88]) that every probability
measure $\mu $ on a vector space $X$, which is differentiable along some
direction $h\in X$ with corresponding logarithmic derivative $b_{h}\in
L^{1}(\mu )$, is for sure also quasi-invariant w.r.t. all shifts $x\mapsto
x+\theta h$. However, the new difficulty and the principal difference
compared with the above mentioned papers is that \emph{no assumptions on the
global integrability} of the logarithmic derivatives $b_{h_{i}}$ are imposed
here.\textbf{\ }Moreover, the test functions $f,$ for this reason chosen to
have bounded support, are not better than \emph{partially differentiable on
the left and right}. Instead, we shall crucially use the proper
approximation procedure (cf. Lemma 4.3) and the observation from Lemma 4.1
that $a_{\theta h_{i}}$, $b_{h_{i}}$ are continuous locally bounded
functions on $\Omega _{-p}^{R}.\smallskip $ \newline
\textbf{Proposition 4.8 }((IbP)-Characterization of Tempered Gibbs Measures)
\ \emph{Denote by} $\mathcal{M}_{t}^{b}$ \emph{the set of all probability
measures }$\mu $ \emph{on} $(\Omega ,\mathcal{B}(\Omega ))$ \emph{which
satisfy the temperedness} \emph{condition (3.64)} \emph{with some} \emph{%
fixed} $p=p(\mu )>d$\emph{\ and for any }$i\in \mathbb{Z}^{d+1}$ \emph{the
(IbP)-formula}%
\begin{equation}
\int_{\Omega }\partial _{h_{i}}^{\pm }f(\omega )\,d\mu (\omega
)=-\int_{\Omega }f(\omega )b_{i}(\omega )\,d\mu (\omega )  \tag{4.57}
\end{equation}%
\emph{for all functions} $f\in C_{0}^{1,\pm }(\Omega _{-p}^{R};h_{i}).$ 
\emph{Then} 
\begin{equation}
\mathcal{M}_{t}^{b}=\mathcal{M}_{t}^{a}=\mathcal{G}_{t}.  \tag{4.58}
\end{equation}

\textbf{Proof}:\textbf{\ (i) } $\mathcal{M}_{t}^{a}\subseteq \mathcal{M}%
_{t}^{b}$:\ \ Consider arbitrary $\mu \in \mathcal{M}^{a}$ supported on some 
$\Omega _{-p}^{R}$ with $p>d.$ Then (4.22) holds, in particular, on
functions $f\in C_{0}^{1,\pm }(\Omega _{-p}^{R};h_{i}),$ $i=(k,n)\in \mathbb{%
Z}^{d+1},$ with 
\begin{equation*}
\text{supp}f\subset B_{-p;k}^{R}(\rho )=\left\{ \omega \in \Omega
_{-p}^{R}\left\vert \ ||\omega ||_{-p,R;k}\leq \rho \right. \right\} ,\
0<\rho =\rho _{f}<\infty .
\end{equation*}%
According to definition (4.56), for all such $f$%
\begin{equation*}
\sup_{\omega \in \Omega _{-p}^{R}}\sup_{|\theta |\leq 1}\left\vert \frac{%
f(\omega +\theta h_{i})-f(\omega )}{\theta }\right\vert \leq \sup_{\omega
\in \Omega _{-p}^{R}}|\partial _{i}^{\pm }f|<\infty .
\end{equation*}%
On the other hand, from Lemma 4.1 
\begin{equation*}
\sup_{\omega \in B_{-p,k}^{R}(\rho )}\sup_{|\theta |\leq 1}\left\vert \frac{%
a_{\theta h_{i}}^{{}}(\omega )-1}{\theta }\right\vert <\infty .
\end{equation*}%
Thus, by Lebesgue's dominated convergence theorem one can pass to the limit $%
\theta \rightarrow \pm 0$ in both sides of (4.22) and get the (IbP)-formula
(4.57). This means the desired inclusion $\mu \in \mathcal{M}_{t}^{b}.$%
\smallskip\ 

\textbf{(ii)} \ $\mathcal{M}_{t}^{b}\subseteq \mathcal{M}_{t}^{a}$:\ \ We
claim that each $\mu \in \mathcal{M}_{t}^{b}$ is quasi-invariant w.r.t. the
shifts $\omega \mapsto \omega +\theta h_{i},$ $\theta \in \mathbb{R}$, with
the Radon--Nikodym derivatives 
\begin{equation}
\frac{d\mu (\omega +\theta h_{i})}{d\mu (\omega )}:=\exp
\int\nolimits_{0}^{\theta }b_{i}(\omega +\vartheta h_{i})d\vartheta
:=a_{\theta h_{i}}(\omega ).  \tag{4.59}
\end{equation}%
(the latter identity in (4.59) is due to Corollary 4.7). So, let $\mu
(\Omega _{-p}^{R})=1$ and let (4.57) hold for $i=(k,n)\in \mathbb{Z}^{d+1}$.
Given any $f\in C_{0}^{1,\pm }(\Omega _{-p}^{R};h_{i}),$ let us define a
family of functions indexed by $\theta \in \mathbb{R}:$%
\begin{equation}
f(\theta ,\cdot )\in C_{0}^{1,\pm }(\Omega _{-p}^{R};h_{i}),\quad f(\theta
,\omega ):=f(\omega +\theta h_{i})a_{\theta h_{i}}(\omega ).  \tag{4.60}
\end{equation}%
Denoting for convenience%
\begin{equation*}
I_{f}(\theta ):=\int_{\Omega }f(\theta ,\omega )d\mu (\omega ),
\end{equation*}%
one can check by a direct calculation that $\forall \theta \in \mathbb{R}$ 
\begin{equation}
\exists \frac{d^{\pm }}{d\theta }I_{f}(\theta )=\int_{\Omega }\left[
\partial _{h_{i}}^{\pm }f(\omega +\theta h_{i})+f(\omega +\theta
h_{i})b_{i}(\omega +\theta h_{i})\right] a_{\theta h_{i}}(\omega )d\mu
(\omega ).  \tag{4.61}
\end{equation}%
Substituting the exact expression (4.48) for $\partial _{h_{i}}a_{\theta
h_{i}}(\omega )$ in (4.61) and then applying the (IbP)-formula (4.57) to the
function $f(\theta ,\cdot ),$ we find that $\frac{d^{\pm }}{d\theta }%
I_{f}(\theta )=0.$ In view of the continuity of $\theta \rightarrow
I_{f}(\theta ),$ the latter yields $I_{f}\equiv const,$ i.e.,%
\begin{equation}
\int_{\Omega }f(\omega +\theta h_{i})a_{\theta h_{i}}(\theta ,\omega )d\mu
(\omega )=\int_{\Omega }f(\omega )d\mu (\omega ),\quad \theta \in \mathbb{R}.
\tag{4.62}
\end{equation}%
Herefrom by Fatou's lemma and Lemma 4.3\thinspace (iii) we conclude that%
\begin{equation}
E_{\mu }(a_{\theta h_{i}}^{{}}(\theta ,\cdot ))=\sup \left\{ E_{\mu }f(\cdot
-\theta h_{i})\text{\thinspace }\left\vert \text{\thinspace }f\in
C_{0}^{1,\pm }(\Omega _{-p}^{R};h_{i}),\text{ }0\leq f\leq 1\right. \right\}
=1,  \tag{4.63}
\end{equation}%
as well as that (4.62) holds for all $f\in L^{1}(\mu )$. Since (4.62) and
(4.7) are equivalent, the converse inclusion $\mu \in \mathcal{M}_{a}^{t}$
is also proved.

\noindent $\blacksquare \smallskip $\newline
\textbf{Remark 4.9. \ }The proof of Proposition 4.8 shows the following
relations:\smallskip 

\textbf{(i) \ }The classes of measures $\mu \in \mathcal{M}_{t}$ which
satisfy the flow resp. integration by parts descriptions coincide, i.e. $%
\mathcal{M}_{t}^{a}=\mathcal{M}_{t}^{b},$ as soon as $\forall i\in \mathbb{Z}%
^{d+1}$ there exists $\partial _{\theta }a_{\theta h_{i}}\in C_{b,loc}(%
\mathbb{R}\times \Omega _{t})$ (and thus also $\partial _{\theta }a_{\theta
h_{i}}|_{\theta =0}=:b_{i}\in C_{b,loc}(\Omega _{t})$). To this end, along
with the assumption $V_{k}\in C_{b,loc}^{1}(\mathbb{R}),$ $k\in \mathbb{Z}%
^{d},$ one could impose $(\mathbf{J})$ and $(\mathbf{W}_{\text{\textbf{i,ii}}%
})$ only.\smallskip

\textbf{(ii) \ }Let $\mu $ be a probability measure\textbf{\ }on $\Omega
_{-p}^{R}$ such that the (IbP)-formula (4.57) holds for all $f\in
C_{0}^{1,\pm }(\Omega _{-p}^{R};h_{i})$ with fixed $i=(k,n)\in \mathbb{Z}%
^{d+1}.$ Then it also holds for all $f\in C_{0}^{1,\pm }(\Omega _{-p^{\prime
}}^{R^{\prime }};h_{i})$ (even though $C_{0}^{1,\pm }(\Omega _{-p^{\prime
}}^{R^{\prime }};h_{i})\upharpoonright \Omega _{-p}^{R}\nsubseteq
C_{0}^{1,\pm }(\Omega _{-p}^{R};h_{i})$) provided $\partial _{\theta
}a_{\theta h_{i}}^{{}}\in C_{b,loc}(\mathbb{R}\times \Omega _{-p^{\prime
};k}^{R^{\prime }})$ and\smallskip $\ p^{\prime }\geq p>d/2,$ $R^{\prime
}\geq R.\smallskip $\newline
\textbf{Corollary 4.10 }(Differentiability $\mu -$a.e. of $|\omega
_{k}|_{C_{\beta }^{{}}}$)\textbf{.} \ \emph{Let} $\mu (\Omega _{-p}^{R})=1$ 
\emph{with some }$p>d,$\emph{\ and let the (IbP)-formula (4.56) in the fixed
direction} $h_{i}\in bas(\mathcal{H}),$ $i=(k,n)\in \mathbb{Z}^{d+1},$ \emph{%
hold for all} \emph{functions from }$C_{0}^{1,\pm }(\Omega _{-p}^{R};h_{i})$%
\emph{. Then there exists a Borel set }$\Delta \subset \Omega $ \emph{such
that}%
\begin{equation}
\mu (\Delta )=1\text{ \ and \ }\forall \omega \in \Delta :\text{ }\partial
_{\varphi _{n}}^{+}|\omega _{k}|_{C_{\beta }^{{}}}=\partial _{\varphi
_{n}}^{-}|\omega _{k}|_{C_{\beta }^{{}}}.  \tag{4.64}
\end{equation}%
\emph{Hence}%
\begin{equation}
\forall f\in C_{b}^{1,\pm }(\Omega _{-p}^{R};h_{i})\text{ }\forall \omega
\in \Delta :\text{ }\exists \text{\thinspace }\partial _{h_{i}}^{{}}f(\omega
)=\partial _{h_{i}}^{\pm }f(\omega ),  \tag{4.65}
\end{equation}%
\emph{and in the formulation of Proposition 4.9 one can simply replace} $%
\partial _{h_{i}}^{\pm }f$ \emph{by} $\partial _{h_{i}}^{{}}f.\smallskip $

\textbf{Proof:} \ By the (IbP)-formula (4.57), for all $g\in
C_{0}^{1}(\Omega _{-p;k}^{R};h_{i})$%
\begin{equation*}
\int_{\Omega }g(\omega )\,\partial _{h_{i}}^{\pm }|\omega _{k}|_{C_{\beta
}^{{}}}d\mu (\omega )=-\int_{\Omega }|\omega _{k}|_{C_{\beta }^{{}}}\left[
\partial _{h_{i}}^{{}}g+gb_{i}\right] (\omega )\,d\mu (\omega ).
\end{equation*}%
Since the set $C_{0}^{1}(\Omega _{-p;k}^{R};h_{i})$ is dense in $L^{1}(\mu )$
(cf. Lemma 4.3 (iii)), this implies 
\begin{equation}
\mu \left( \omega \in \Omega \text{\thinspace }\left\vert \text{ }\partial
_{\varphi _{n}}^{+}|\omega _{k}|_{C_{\beta }^{{}}}=\partial _{\varphi
_{n}}^{-}|\omega _{k}|_{C_{\beta }^{{}}}\text{,\quad }\forall (k,n)\in 
\mathbb{Z}^{d+1}\right. \right) =1.  \tag{4.67}
\end{equation}%
We note that the identity in (4.65) $\partial _{h_{i}}^{+}f=\partial
_{h_{i}}^{-}f\ (\mu -$a.e.$)$ might also be derived from a general result in
[Ku82, Lemma 1.3] on the so-called stochastic $H-$G\^{a}teaux
differentiability of Lipshitz functions $f:B\rightarrow \mathbb{R}$ w.r.t.
any positive quasi-invariant measure $\mu $ on an abstract Wiener space $%
B\supset H\supset B^{\ast }$.

\noindent $\blacksquare \smallskip $

So, based on Propositions~4.8, instead of Euclidean Gibbs measures $\mu \in 
\mathcal{G}$ initially defined as random fields on the lattice $\mathbb{Z}%
^{d},$ we can just study probability measures on $\Omega $ satisfying the
(IbP)-formula (4.57) with the prescribed by (4.38) logarithmic derivatives $%
b_{i},\ i\in \mathbb{Z}^{d+1}.$ Let us stress that the $b_{i}$ only depend
on the given potentials $V$ and $W$ and hence are the same for all $\mu \in 
\mathcal{G}$ associated with the heuristic Hamiltonian (3.1). Solutions $\mu
\in \mathcal{M}^{b}$ to the (IbP)-formula (4.57) will also be called \emph{%
symmetrizing} measures. For further important connections to reversible
diffusion processes and Dirichlet operators in infinite dimensions we refer
e.g. to [AKR97a,b, AKRT01, BR01, BRW01] and to Subsect.\thinspace 5.3 below.

\subsection{Further discussion of the (IbP)-formula}

Here we present some modifications of the (IbP)-formula (4.57) to be useful
in applications. The domain $C_{0}^{1,\pm }(\Omega _{-p}^{R},h_{i}),$ on
which we have initially proved the (IbP)-characterization of Euclidean Gibbs
measures in Proposition 4.8, was chosen for reasons of technical convenience
only, and indeed it plays not more than an intermediate role. Among others,
the subsequent statement will show that as soon as we have the (IbP)-formula
(4.57) in a usual setting with $\partial _{h_{i}}f$ on some\emph{\
\textquotedblright minimal domain\textquotedblright } consisting of $h_{i}-$%
differentiable functions with bounded support (e.g., $C_{0}^{1}(\Omega
_{-p,k}^{R};h_{i})$) or of smooth cylinder functions (e.g., $\mathcal{F}%
C_{b}^{\infty }(\Omega ;\mathbb{Z}^{d}\times S_{\beta })$ or\emph{\ }$%
\mathcal{F}C_{b}^{\infty }(\Omega ;\mathcal{H})$), then we can always extend
it, by substituting $\partial _{h_{i}}^{\pm }f$ for $\partial _{h_{i}}f,$ to
a \emph{\textquotedblright maximal domain\textquotedblright\ }of definition
on functions which are differentiable along $h_{i}$ at least on the right or
left (e.g., $C_{b}^{1,\pm }(\Omega _{-p}^{R},h_{i})$). We remind that the
spaces of test functions used here have been already introduced respectively
in Subsections 4.2.1, 4.2.3 and by Definition 4.7. \smallskip \newline
\noindent \textbf{Proposition 4.11 (}Equivalent Domains for (IbP)-Formula) \ 
\emph{Let} $\mu (\Omega _{-p}^{R})=1$ \emph{with some }$p>d$,\emph{\ and let
us fix some direction }$h_{i}\in bas(\mathcal{H}),$ $i=(k,n)\in \mathbb{Z}%
^{d+1}.$\emph{\ Then the following assertions hold:\smallskip }

\textbf{(i)} \ \emph{Suppose that the (IbP)-formula (4.57)} $\emph{is}$ 
\emph{valid for} \emph{functions from }$C_{0}^{1,\pm }(\Omega
_{-p}^{R},h_{i}).$ \emph{Then} \emph{it extends by continuity to }$%
C_{dec}^{1,\pm }(\Omega _{-p}^{R},h_{i})$ \emph{and, in particular, to }$%
C_{0}^{1}(\Omega _{-p,k}^{R};h_{i})\subset C_{dec}^{1,\pm }(\Omega
_{-p}^{R},h_{i}).$ \emph{If, moreover, we know that} $b_{i}\in L^{1}(\mu ),$%
\emph{\ then (4.57) further extends from }$C_{0}^{1,\pm }(\Omega
_{-p}^{R},h_{i})$ \emph{to }$C_{b}^{1,\pm }(\Omega _{-p}^{R},h_{i}),$ \emph{%
and in particular to all cylinder functions from }$\mathcal{F}C_{b}^{\infty
}(\Omega ;\mathbb{Z}^{d}\times S_{\beta })$\emph{\ or }$\mathcal{F}%
C_{b}^{\infty }(\Omega ;\mathcal{H}).\smallskip $

\textbf{(ii) \ }\emph{The inverse statement to (i)\ is also true: Suppose
that the (IbP)-formula (4.57) is valid for functions from }$C_{0}^{1}(\Omega
_{-p,k}^{R};h_{i}),$ \emph{then it extends by continuity to all }$f\in
C_{0}^{1,\pm }(\Omega _{-p}^{R},h_{i}).$\emph{\ If }$b_{i}\in L^{1}(\mu ),$ 
\emph{then (4.57) extends from} $\mathcal{F}C_{b}^{\infty }(\Omega ;\mathbb{Z%
}^{d}\times S_{\beta })$ \emph{resp. }$\mathcal{F}C_{b}^{\infty }(\Omega ;%
\mathcal{H})$\emph{\ to }$C_{b}^{1,\pm }(\Omega _{-p}^{R},h_{i}).\smallskip $

\textbf{Proof:} \textbf{(i)} \ For any $f\in C_{b}^{1,\pm }(\Omega
_{-p}^{R},h_{i})$ ($\supset C_{dec}^{1,\pm }(\Omega _{-p}^{R},h_{i})$)$,$
let us take its approximation by $\{f^{(K)}\}_{K\in \mathbb{N}}\subset
C_{0}^{1,\pm }(\Omega _{-p}^{R},h_{i})$ of the form%
\begin{equation}
f^{(K)}(\omega ):=f(\omega )\chi _{K}^{{}}(||\omega ||_{-p,R;k}),  \tag{4.68}
\end{equation}%
where $\{\chi _{K}\}_{K\in \mathbb{N}}\subset C_{0}^{1}(\mathbb{R}%
_{+}\rightarrow \lbrack 0,1])$ is a cut-off sequence with the properties 
\begin{gather}
\chi _{K}^{{}}(s)=1\text{ for }s\in \lbrack 0,K],\ \chi _{K}^{{}}(s)=0\text{
for }s\in \lbrack K+1,\infty )  \notag \\
\text{and\quad }\chi _{K+1}^{{}}(s+1)=\chi _{K}^{{}}(s)\text{ for every }%
s\geq 0.  \tag{4.69}
\end{gather}%
The statement follows by Lebesgue's dominated convergence theorem applied,
as $K\rightarrow \infty ,$ to both sides of (4.57) with $f^{(K)}$ replacing $%
f.\smallskip $

\textbf{(ii)} \ By Lemma 4.3 (i) and Lebesgue's convergence theorem, (4.53)
extends either from $\mathcal{F}C_{b}^{\infty }(\Omega ;\mathbb{Z}^{d}\times
S_{\beta })$ or\emph{\ }$\mathcal{F}C_{b}^{\infty }(\Omega ;\mathcal{H})$ to
the sets $C_{b}^{1}(\Omega _{-p}^{R};h_{i})\supset C_{b}^{1}(\Omega
_{-p;k}^{R};h_{i})$ provided $b_{i}\in L^{1}(\mu ).$ On the other hand, by
Lemma 4.3 (i,ii) and again by Lebesgue's dominated convergence theorem,
(4.57) extends from $C_{0}^{1}(\Omega _{-p;k}^{R};h_{i})$ to $C_{0}^{1,\pm
}(\Omega _{-p}^{R},h_{i}).$ And finally, the passage from $C_{0}^{1,\pm
}(\Omega _{-p}^{R},h_{i})$ to $C_{b}^{1,\pm }(\Omega _{-p}^{R},h_{i})$ was
just performed in the proof of (i).

\noindent $\blacksquare \smallskip $

Now we present equivalent versions of the (IbP)-formula (4.57) for smooth
cylinder functions on $\Omega .$ Respectively, it would be important to
extend (4.57), when applied to the test functions $f\in \mathcal{F}%
C_{b}^{1}(\Omega ;\mathbb{Z}^{d}\times S_{\beta }),$ from the admissible
single spin directions $\varphi \in lin\{\varphi _{n}\}_{n\in \mathbb{Z}}$
to the Green function $\mathfrak{G}_{\tau }:=A_{\beta }^{-1}\delta _{\tau },$
$\tau \in S_{\beta },$ of the operator $A_{\beta }$ (for its definition and
properties see (3.29)--(3.31)). Heuristically, this is possible for the
following reason: even though $\varphi :=\mathfrak{G}_{\tau }\notin \mathcal{%
D}(A_{\beta }^{{}})$ and hence, in general, $d\mu (\omega +\theta h)\bot
d\mu (\omega )$ for $\mu \in \mathcal{M}_{t}^{b},$ $\theta \in \mathbb{R}$
and $h:=\varphi \otimes e_{k}\in \mathcal{H},$ we have the absolute
continuity for the corresponding space-time projections of these measures.
Technically, such extension will be performed through some standard \emph{%
approximation procedure} for $\mathfrak{G}_{\tau }$, which also will be
needed below in the proof of all our Main Theorems I--III. We also have to
assume here that\emph{\ }$\mu \in \mathcal{M}_{t}^{b}$ possesses the
following a priori integrability properties:

\begin{description}
\item[\textbf{(H}$_{\protect\mu }$\textbf{)}] \emph{For all }$i:=(k,n)\in 
\mathbb{Z}^{d+1}:$%
\begin{equation}
|\omega _{k}|_{C_{\beta }^{{}}},\text{ }|F_{k}^{V,W}(\omega )|_{L_{\beta
}^{1}}\in L^{1}(\mu )\quad \emph{and\ thus\ also\quad }b_{i}\in L^{1}(\mu ).
\tag{4.70}
\end{equation}
\end{description}

\noindent As will be shown later in Main Theorem II, Hypothesis $(\mathbf{H}%
_{\mu })$ \emph{always hold }under Assumptions ($\mathbf{W}$), ($\mathbf{J}$%
), ($\mathbf{V}$). Finally we note that local versions of the (IbP)-formula
analogous to (4.71) are of common use in the literature on quantum models
(see e.g. [GJ81]).\smallskip \newline
\textbf{Proposition 4.12 }(\textquotedblright Local\textquotedblright\
(IbP)-Formula)\textbf{\ \ }\emph{Let the measure }$\mu \in \mathcal{M}%
_{t}^{b}$ \emph{satisfy Hypothesis (}\textbf{H}$_{\mu }).$\emph{\ Then the
following (IbP)-formulas for cylinder functions from} $\mathcal{F}%
C_{b}^{1}(\Omega ;\mathbb{Z}^{d}\times S_{\beta })$ \emph{resp}. $\mathcal{F}%
C_{b}^{1}(\Omega ;\mathcal{H})$\emph{\ hold:} 
\begin{multline}
\text{\textbf{(i)\quad }}\sum_{l=1}^{L}\delta _{k-k_{l}}^{{}}\mathfrak{G}%
(\tau ,\tau _{l}^{{}})\int_{\Omega }\partial _{l}f_{L}\left( \omega
_{k}^{{}}(\tau _{1}),...,\omega _{k}^{{}}(\tau _{L}^{{}})\right) \,d\mu
(\omega )  \notag \\
=-\int_{\Omega }f_{L}\left( \omega _{k}^{{}}(\tau _{1}^{{}}),...,\omega
_{k}^{{}}(\tau _{L}^{{}})\right) \left( \omega _{k}+A_{\beta
}^{-1}F_{k}^{V,W}(\omega )\right) (\tau )d\mu (\omega );  \tag{4.71}
\end{multline}%
\begin{multline}
\text{\textbf{(ii)\quad }}\sum_{l=1}^{L}\delta _{k-k_{l}}(A_{\beta
}^{-1}\varphi _{n_{l}^{{}}}^{{}})(\tau )\int_{\Omega }\partial
_{l}f_{L}\left( <\omega ,h_{i_{1}}^{{}}>_{\mathcal{H}},...,<\omega
,h_{i_{L}}^{{}}>_{\mathcal{H}}^{{}}\right) d\mu (\omega )  \notag \\
=-\int_{\Omega }f_{L}\left( <\omega ,h_{i_{1}}>_{\mathcal{H}},...,<\omega
,h_{i_{L}}>_{\mathcal{H}}\right) \left( \omega _{k}+A_{\beta
}^{-1}F_{k}^{V,W}(\omega )\right) (\tau )\,d\mu (\omega )  \tag{4.72}
\end{multline}%
\emph{for any} $f_{L}\in C_{b}^{1}(\mathbb{R}^{L}),$ $h\in \lbrack L_{\beta
}^{2}]^{\mathbb{Z}_{0}^{d}},$ $n_{l}\in \mathbb{Z}$, $k,k_{l}\in \mathbb{Z}%
^{d},$ $\tau ,\tau _{l}\in S_{\beta },$ $1\leq l\leq L$ \emph{and} $L\in 
\mathbb{N}$.\smallskip 

\textbf{Proof: (i) \ }At first,\textbf{\ }we shall exploit the assumptions $%
|F_{k}^{V,W}(\omega )|_{L_{\beta }^{1}},$ $|\omega _{k}|_{L_{\beta }^{2}}\in
L^{1}(\mu )$ from Hypothesis $\mathbf{(H}_{\mu }\mathbf{).}$ According to
Proposition 4.11, (4.57) is then equivalent to the (IbP)-formula 
\begin{multline}
\int_{\Omega }\partial _{h_{i}}f(\omega )\,d\mu (\omega
)=\sum_{l=1}^{L}\delta _{k-k_{l}}^{{}}\varphi _{n}^{{}}(\tau
_{l}^{{}})\int_{\Omega }\partial _{l}f_{L}\left( \omega
_{k_{1}^{{}}}^{{}}(\tau _{1}^{{}}),...,\omega _{k_{L}^{{}}}^{{}}(\tau
_{L}^{{}})\right) \,d\mu (\omega )  \notag \\
=-\int_{\Omega }f_{L}\left( \omega _{k_{1}^{{}}}(\tau _{1}),...,\omega
_{k_{L}^{{}}}(\tau _{L})\right) \left[ (A_{\beta }\varphi _{n},\omega
_{k})_{H}^{{}}+(F_{k}^{V,W}(\omega ),\varphi _{n})_{H}\right] d\mu (\omega )
\tag{4.73}
\end{multline}%
valid for all cylinder functions $f\in \mathcal{F}C_{b}^{1}(\Omega ;\mathbb{Z%
}^{d}\times S_{\beta })$ of the form (4.30). Since $T_{\beta }:=lin\{\varphi
_{n}\}_{n\in \mathbb{Z}}$ is an essential domain in $L_{\beta }^{2}$ for the
self-adjoint operator $A_{\beta }$ with $\mathcal{D}(A_{\beta
}^{{}})=W_{\beta }^{2,2}$ (see Subsect.\thinspace 3.3.1 above), (4.73)
further extends by continuity to arbitrary $\varphi \in \mathcal{D}(A_{\beta
})$ as 
\begin{multline}
\sum_{l=1}^{L}\delta _{k-k_{l}}^{{}}\varphi (\tau _{l}^{{}})\int_{\Omega
}\partial _{l}f_{L}\left( \omega _{k_{1}^{{}}}^{{}}(\tau
_{1}^{{}}),...,\omega _{k_{L}^{{}}}^{{}}(\tau _{L}^{{}})\right) \,d\mu
(\omega )  \notag \\
=-\int_{\Omega }f_{L}\left( \omega _{k_{1}^{{}}}^{{}}(\tau
_{1}^{{}}),...,\omega _{k_{L}^{{}}}^{{}}(\tau _{L}^{{}})\right) \left[
(A_{\beta }^{{}}\varphi ,\omega _{k}^{{}})_{H}^{{}}+(F_{k}^{V,W}(\omega
),\varphi )_{H}^{{}}\right] d\mu (\omega ).  \tag{4.74}
\end{multline}

Now, fixing $\tau \in S_{\beta },$ we would like to replace $\varphi $ in
(4.74) by the Green function $\mathfrak{G}_{\tau }:=A_{\beta }^{-1}\delta
_{\tau }.$ With this purpose we construct the (so-called \emph{Yosida})
approximation of $\mathfrak{G}_{\tau }\in W_{\beta }^{2,1}\subset C_{\beta }$
by 
\begin{equation}
\varphi _{\tau }^{(K)}:=\left( 1+K^{-1}A_{\beta }\right) ^{-1}\mathfrak{G}%
_{\tau }\in \mathcal{D}(A_{\beta }),\quad K\in \mathbb{N}.  \tag{4.75}
\end{equation}%
Since $A_{\beta }^{{}}$ generates a Markovian $C_{0}$-semigroup on $C_{\beta
},$ its resolvent has the following properties: 
\begin{gather*}
\sup_{K\in \mathbb{N}}||\left( 1+K^{-1}A_{\beta }\right) ^{-1}||_{C_{\beta
}\rightarrow C_{\beta }}\leq 1, \\
\lim_{K\rightarrow \infty }|\left( 1+K^{-1}A_{\beta }^{{}}\right)
^{-1}\upsilon -\upsilon |_{C_{\beta }^{{}}}=0,\text{ }\forall \upsilon \in
C_{\beta }^{{}}.
\end{gather*}%
Therefore, $\lim_{K\rightarrow \infty }|\varphi _{\tau }^{(K)}-\mathfrak{G}%
_{\tau }|_{C_{\beta }^{{}}}=0$ and, moreover, 
\begin{gather}
|(A_{\beta }^{{}}\varphi _{\tau }^{(K)},\omega _{k})_{H}^{{}}|\leq |\omega
_{k}|_{C_{\beta }^{{}}},  \notag \\
\lim_{K\rightarrow \infty }(A_{\beta }^{{}}\varphi _{\tau }^{(K)},\omega
_{k})_{H}^{{}}=\lim_{K\rightarrow \infty }(\mathbb{\delta }_{\tau
}^{{}},\left( 1+K^{-1}A_{\beta }^{{}}\right) ^{-1}\omega
_{k})_{H}^{{}}=\omega _{k}(\tau ).  \tag{4.76}
\end{gather}%
On the other hand, $A_{\beta }^{-1}F_{k}^{V,W}:\Omega _{-p}^{R}\rightarrow
C_{\beta }^{{}}$ and hence 
\begin{gather}
\lim_{K\rightarrow \infty }(F_{k}^{V,W}(\omega ),\varphi _{\tau
}^{(K)})_{H}^{{}}=(A_{\beta }^{-1}F_{k}^{V,W}(\omega ))(\tau ),  \notag \\
|(F_{k}^{V,W}(\omega ),\varphi _{\tau }^{(K)})_{H}^{{}}|\leq
|F_{k}^{V,W}(\omega )|_{L_{\beta }^{1}}|\mathfrak{G}_{\tau }^{{}}|_{C_{\beta
}^{{}}}.  \tag{4.77}
\end{gather}%
Finally, taking into account (4.76), (4.77) and the next assumption $|\omega
_{k}|_{C_{\beta }^{{}}}\in L^{1}(\mu )$ from Hypothesis $\mathbf{(H}_{\mu }%
\mathbf{)},$ by Lebesgue's dominated convergence theorem (4.74) extends to
(4.71).\smallskip

\textbf{(ii) \ }The proof for the case of cylinder functions $f\in \mathcal{F%
}C_{b}^{1}(\Omega _{\beta }^{{}};\mathcal{H})$ of the form (4.28) is
analogous. Namely, under the assumptions $|F_{k}^{V,W}(\omega )|_{L_{\beta
}^{1}},$ $|\omega _{k}|_{L_{\beta }^{2}}\in L^{1}(\mu ),$ the (IbP)-formula
(4.57) is equivalent to%
\begin{multline}
\sum_{l=1}^{L}<(\mathbf{1\otimes }A_{\beta }^{-1}\mathbf{)}%
h,h_{i_{L}}^{{}}>_{\mathcal{H}}^{{}}\int_{\Omega }\partial _{l}f_{L}\left(
<\omega ,h_{i_{1}}^{{}}>_{\mathcal{H}}^{{}},...,<\omega ,h_{i_{L}}^{{}}>_{%
\mathcal{H}}^{{}}\right) d\mu (\omega )  \notag \\
\int_{\Omega }f_{L}\left( <\omega ,h_{i_{1}}^{{}}>_{\mathcal{H}%
}^{{}},...,<\omega ,h_{i_{L}}^{{}}>_{\mathcal{H}}^{{}}\right) \left\langle
\omega +(\mathbf{1\otimes }A_{\beta }^{-1}\mathbf{)}F^{V,W}(\omega
),h\right\rangle _{\mathcal{H}}^{{}}d\mu (\omega )  \tag{4.78}
\end{multline}%
and, if $|\omega _{k}|_{C_{\beta }^{{}}}\in L^{1}(\mu ),$ it further extends
by continuity to (4.72).

\noindent $\blacksquare $\smallskip

\noindent \textbf{Remark 4.13} \ For the particular QLS\ models (2.2) and
(2.21) with local interaction, the corresponding (IbP)-characterization (as
well as the flow characterization, cf. Remark 4.3) is available \emph{for all%
} (not necessarily tempered) $\mu \in \mathcal{G}.$ So, as was shown in
[AKPR03, Proposition 2], the (IbP)-formula (4.57) holds (without any a
priori information about global integrability of $b_{i})$ for all test
functions $f\in C_{b}^{1}(\Omega ;h_{i})$ satisfying the\ extra decay
condition (4.55).

\section{Applications of the (IbP)-formula}

\subsection{(IbP)-formula for probability kernels of the local specification}

As is immediately to see from the definitions (3.41)--(3.44) and for
classical lattice systems has been already mentioned in [Roy77], the
following flow description for the probability kernels\emph{\ }$\pi
_{\Lambda }^{{}}(d\omega |\xi )$ of the local specification $(\pi _{\Lambda
}^{{}})_{\Lambda \Subset \mathbb{Z}^{d}}$ is true:\smallskip

\noindent \textbf{Proposition 5.1.} \ \emph{Measures }$\pi _{\Lambda
}^{{}}(d\omega |\xi ),$ $\xi \in \Omega _{t},$ $\Lambda \Subset \mathbb{Z}%
^{d},$ \emph{(as well as the finite-volume Gibbs distributions }$\mu
_{\Lambda }^{{}}(d\omega _{\Lambda }|\xi _{\Lambda ^{c}})$\emph{)} \emph{are
quasi-invariant w.r.t. shifts} 
\begin{equation*}
\omega \longmapsto \omega +\theta h_{i},\emph{\ \ }\forall i=(k,n),\text{ }%
k\in \Lambda \Subset \mathbb{Z}^{d},\text{ }n\in \mathbb{Z}\text{ \emph{and} 
}\theta \in \mathbb{R},
\end{equation*}%
\emph{with the same Radon--Nikodym derivatives} \emph{as those for the
corresponding infinite-volume Gibbs measures }$\mu \in \mathcal{G}_{t}$%
.\smallskip\ 

More precisely, Proposition 5.1 means that for every such admissible shift $%
\theta h_{i}\in \mathcal{H}$%
\begin{equation}
\frac{d\pi _{\Lambda }(\omega +\theta h_{i}|\xi )}{d\pi _{\Lambda }(\omega
|\xi )}=a_{\theta h_{i}}^{{}}(\omega )>0,\text{ \ }\forall \omega \in \Omega
_{t}\text{ }(\pi _{\Lambda }^{{}}(d\omega |\xi )-a.e.)  \tag{5.1}
\end{equation}%
or, equivalently, for all $f\in L^{1}(\Omega ,$\thinspace $\pi _{\Lambda
}^{{}}(d\omega |\xi ))$ 
\begin{equation*}
\int_{\Omega }f(\omega )a_{\theta h_{i}}^{{}}(\omega )\pi _{\Lambda
}^{{}}(d\omega |\xi )=\int_{\Omega }f(\omega -\theta h_{i})\pi _{\beta
,\Lambda }^{{}}(d\omega |\xi ).
\end{equation*}%
A reasoning similar to that used in the proof of Propositions 4.8 and 4.11
(i) then gives the corresponding infinitesimal version of (5.1):\smallskip 

\noindent \textbf{Proposition 5.2} \ \emph{Measures }$\pi _{\Lambda
}^{{}}(d\omega |\xi ),$ $\xi \in \Omega _{t},$ $\Lambda \Subset \mathbb{Z}%
^{d},$ \emph{(as well as the finite-volume Gibbs distributions }$\mu
_{\Lambda }^{{}}(d\omega _{\Lambda }|\xi _{\Lambda ^{c}})$\emph{\ on }$%
\Omega _{\Lambda }$\emph{)} \emph{satisfy the} \emph{(IbP)-formula }%
\begin{gather}
\int_{\Omega }\partial _{h_{i}}^{\pm }f(\omega )\,\pi _{\Lambda
}^{{}}(d\omega |\xi )=-\int_{\Omega }f(\omega )b_{i}(\omega )\,\pi _{\Lambda
}(d\omega |\xi ),  \notag \\
\forall i=(k,n),\text{ }k\in \Lambda \Subset \mathbb{Z}^{d},\text{ }n\in 
\mathbb{Z},  \tag{5.2}
\end{gather}%
\emph{for all functions} $f\in \bigcup {}_{p>d}C_{0}^{1,\pm }(\Omega
_{-p}^{R};h_{i}).$ \smallskip

We note that for the proof of the above Propositions 5.1 and 5.2, as well as
of Proposition 5.3 below, it is essential\emph{\ }that (as was checked in
Lemma 4.1 and 4.4), 
\begin{equation*}
a_{\theta h_{i}},\text{ }b_{i}\in C_{b,loc}(\Omega _{-p;k}^{R}),\quad
\forall i=(k,n)\in \mathbb{Z}^{d+1},\text{ }\theta \in \mathbb{R}.
\end{equation*}%
Hence, if a sequence $\pi _{\Lambda ^{(K)}}^{{}}(d\omega |\xi ^{(K)}),$
where $\xi ^{(K)}\in \Omega _{-p}^{R},$ $\Lambda ^{(K)}\nearrow \mathbb{Z}%
^{d}$ as $K\rightarrow \infty ,$ weakly converges\emph{\ }on the Banach
space $\Omega _{-p;k}^{R}$ to some probability measure $\mu _{\ast },$ one
can also pass to the limit in both sides of (5.1) and (5.2). So, for $\mu
:=\mu _{\ast }$ we again have both the flow description (4.7) and the
(IbP)-formula (4.57), which hold in any direction $h_{i},$ $i=(k,n)\in 
\mathbb{Z}^{d+1},$ for all functions $f\in C_{0}^{1}(\Omega
_{-p;k}^{R};h_{i})\subset C_{0}^{1,\pm }(\Omega _{-p}^{R};h_{i})$. Combining
these properties of $\mu _{\ast }$ with Propositions 4.2, 4.8 and 4.11, we
have thus proved the following:\smallskip 

\noindent \textbf{Proposition 5.3} (Thermodynamic Limit Points are
Gibbs)\quad \emph{Fix any configuration space} $\Omega _{-p}^{R}$, $p>d$, 
\emph{and} \emph{consider a sequence of measures }$\pi _{\Lambda
^{(K)}}^{{}}(d\omega |\xi ^{(K)}),$ $K\in \mathbb{N},$\emph{\ where} $\xi
^{(K)}\in \Omega _{-p}^{R}$ $\emph{and}$ $\Lambda ^{(K)}\Subset \mathbb{Z}%
^{d},$ $\Lambda ^{(K)}\nearrow \mathbb{Z}^{d}$\emph{\ as }$K\rightarrow
\infty .$ \emph{Then each of its accumulation points} $\mu _{\ast }^{{}}$%
\emph{\ w.r.t. the topology of weak convergence of measures on the Polish
space }$\Omega _{-p}^{R},$ \emph{provided such exist,\ is Gibbs.\smallskip }

In this way, the alternative characterization of Euclidean Gibbs measures
enables us to study the existence problem for $\mu \in \mathcal{G}_{t}$ just
by showing some uniform estimates on measures $\pi _{\Lambda }^{{}}(d\omega
|\xi )$ which certainly imply their tightness.

\subsection{Proof of Main Theorems I--III under additional Hypotheses ($%
\mathbf{H}$) and ($\mathbf{H}_{loc}$)}

In this subsection we complete the proof of our Main Theorems I--III under
the \emph{additional} \emph{Hypotheses} ($\mathbf{H}$) and ($\mathbf{H}%
_{loc} $) concerning the uniform integrability of $|\omega _{k}|_{L_{\beta
}^{2}}$ and $|F_{k}(\omega )|_{L_{\beta }^{1}}$ w.r.t. tempered Gibbs
measures $\mu \in \mathcal{G}_{t}$ and their local specifications $\{\pi
_{\Lambda }\}_{\Lambda \Subset \mathbb{Z}^{d}}$. The \emph{verification of
these crucial} \emph{hypotheses} (which is one of the hardest parts of this
paper) we postpone to Subsect.\thinspace 7.3, after we have available the
necessary techniques based on the coercivity properties of the logarithmic
derivatives $b_{i}$ and developed in Sections 6 and 7 below. The respective
Theorems 7.5 and 7.6 will confirm that both ($\mathbf{H}$) and ($\mathbf{H}%
_{loc}$) \emph{always hold }for the system (3.1) as soon as the interaction
satisfies Assumptions ($\mathbf{W}$), ($\mathbf{J}$) and ($\mathbf{V}$) with
a proper relation between the parameters involved.\smallskip

So, within this subsection we assume:

\begin{description}
\item[$\mathbf{(H)}$] \emph{For every }$Q\geq 1$%
\begin{equation*}
I_{2Q}(F):=\sup_{\mu \in \mathcal{G}_{t}}\sup_{k\in \mathbb{Z}%
^{d}}\int_{\Omega }\left[ |\omega _{k}|_{L_{\beta }^{2}}+|F_{k}^{V,W}(\omega
)|_{L_{\beta }^{1}}\right] ^{2Q}d\mu (\omega )<\infty ;
\end{equation*}

\item[$\mathbf{(H}_{loc}\mathbf{)}$] \emph{For some fixed} $Q>1$ \emph{and} 
\emph{a boundary condition} $\xi \in \Omega _{t}$: 
\begin{equation*}
I_{2Q,\xi }(F):=\sup_{\Lambda \Subset \mathbb{Z}^{d}}\sup_{k\in \Lambda
}\int_{\Omega }\left[ |\omega _{k}|_{L_{\beta }^{2}}+|F_{k}^{V,W}(\omega
)|_{L_{\beta }^{1}}\right] ^{2Q}\pi _{\Lambda }(d\omega |\xi )<\infty .
\end{equation*}
\end{description}

\noindent Our nearest aim is to demonstrate that, taking into account the
regularity properties of the Green function $\mathfrak{G}(\tau ,\tau
^{\prime })$ of the operator $A_{\beta }^{{}}$, from the (IbP)-formulas
(4.57) resp. (5.2) and Hypotheses $\mathbf{(H)}$ resp. $\mathbf{(H}_{loc}%
\mathbf{)}$ one can already get uniform estimates on H\"{o}lder loop spaces $%
C_{\beta }^{\alpha }$ for tempered Gibbs measures $\mu \in \mathcal{G}_{t}$
as well as the existence of the latter. Recall that, according to our
previous agreement, $\mathcal{G}_{t}$ denotes either $\mathcal{G}_{(e)t}^{R}$
for the particular QLS\ models I--III or respectively $\mathcal{G}%
_{(s)t}^{R} $ for the general model (3.1).

\subsubsection{A priori estimates and support properties of Euclidean Gibbs
measures on H\"{o}lder loops}

The following result is the crucial step towards the proof of Main Theorem
II:\smallskip

\noindent \textbf{Lemma 5.4 }(Kolmogorov type moment estimates)\textbf{.} \ 
\emph{For given} $k\in \mathbb{Z}^{d},$ \emph{let a measure} $\mu \in 
\mathcal{M}_{t}$ \emph{satisfy the} \emph{(IbP)-formula (4.57)} \emph{in all
directions} $h_{i}$ \emph{with} $i=(k,n),$ $n\in \mathbb{Z}$, \emph{and,
moreover, obeys the integrability property}%
\begin{equation}
\int_{\Omega }\left[ |\omega _{k}|_{L_{\beta }^{2}}+|F_{k}^{V,W}(\omega
)|_{L_{\beta }^{1}}\right] ^{2Q}d\mu (\omega )\leq I_{2Q}  \tag{5.3}
\end{equation}%
\emph{with some fixed }$Q\in \mathbb{N}$ \emph{and} $I_{2Q}:=I_{2Q}(F_{k};%
\mu )\in (0,\infty )$\emph{.} \emph{Then} \emph{(5.3)} \emph{implies the
moment estimates}%
\begin{gather}
\int_{\Omega }\omega _{k}^{2Q}(\tau )\,d\mu (\omega ):=C_{2Q}<\infty , 
\tag{5.4} \\
\int_{\Omega }[\omega _{k}(\tau )-\omega _{k}(\tau ^{\prime })]^{2Q}\,d\mu
(\omega )\leq \Delta C_{2Q}\cdot \rho ^{Q}(\tau ,\tau ^{\prime })  \tag{5.5}
\end{gather}%
\emph{for all} $\tau ,\tau ^{\prime }\in S_{\beta }$ \emph{and with the
absolute} \emph{(i.e., independent on} $\mu $ \emph{and }$k$\emph{)
constants }$C_{2Q}$, $\Delta C_{2Q}\in (0,\infty ).$ \emph{Moreover}, $%
C_{2Q} $ \emph{and} $\Delta C_{2Q}$ \emph{themselves} \emph{are linear
functions} \emph{of} $I_{2Q}.$\smallskip

\textbf{Proof: (i)} \ Let us consider the following two families of test
functions indexed by $\tau ,\tau ^{\prime }\in S_{\beta }$ and $\varepsilon
>0$:%
\begin{gather}
f(\omega ):=f_{\tau ,\tau ^{\prime },\varepsilon }(\omega ):=\omega
_{k}^{2Q-1}(\tau )Z^{-1}(\omega _{k}),  \tag{5.6} \\
g(\omega ):=g_{\tau ,\tau ^{\prime },\varepsilon }(\omega ):=\left[ \Delta
\omega _{k}\right] ^{2Q-1}Z^{-1}(\omega _{k}),  \tag{5.7}
\end{gather}%
where we set for convenience 
\begin{equation*}
Z(\omega _{k}):=Z_{\varepsilon }(\omega _{k}):=1+\varepsilon |\omega
_{k}|_{C_{\beta }^{{}}}^{2Q},\text{ \ }\Delta \omega _{k}:=\Delta \omega
_{k}(\tau ,\tau ^{\prime }):=\omega _{k}(\tau )-\omega _{k}(\tau ^{\prime }).%
\text{\ }
\end{equation*}%
According to Definition 4.7, $f,g\in C_{b}^{1,\pm }(\Omega _{-p}^{R};h_{i})$
for any $p>d/2$ and $h_{i}:=e_{_{k}}\otimes \varphi _{n},$ $n\in \mathbb{Z}.$
Moreover, using (4.26) and (4.31), one can easily calculate the
corresponding derivatives along vectors $h:=e_{_{k}}\otimes \varphi $ with
arbitrary $\varphi \in C_{\beta }$%
\begin{equation*}
\partial _{h}^{\pm }f(\omega )=(2Q-1)\omega _{k}^{2Q-2}(\tau )\varphi (\tau
)Z^{-1}(\omega _{k})-2\varepsilon Q\omega _{k}^{2Q-1}(\tau )|\omega
_{k}|_{C_{\beta }^{{}}}^{2Q-1}\partial _{\varphi }^{\pm }|\omega
_{k}|_{C_{\beta }^{{}}}^{{}}Z^{-2}(\omega _{k}),
\end{equation*}%
\begin{multline*}
\partial _{h}^{\pm }g(\omega )=(2Q-1)\left[ \Delta \omega _{k}(\tau ,\tau
^{\prime })\right] ^{2Q-2}\left[ \varphi (\tau )-\varphi (\tau ^{\prime })%
\right] Z_{{}}^{-1}(\omega _{k}) \\
-2\varepsilon Q\left[ \Delta \omega _{k}\right] ^{2Q-1}|\omega
_{k}|_{C_{\beta }^{{}}}^{2Q-1}\partial _{\varphi }^{\pm }|\omega
_{k}|_{C_{\beta }}Z^{-2}(\omega _{k})
\end{multline*}%
and then majorate them by 
\begin{gather}
|\partial _{h}^{\pm }f(\omega )|\leq 4Q|\varphi |_{C_{\beta }^{{}}}\omega
_{k}^{2Q-2}(\tau )Z^{-1}(\omega _{k}),  \tag{5.8} \\
|\partial _{h}^{\pm }g(\omega )|\leq 8Q|\varphi |_{C_{\beta
}^{{}}}^{{}}[\Delta \omega _{k}]^{2Q-2}Z^{-1}(\omega _{k}).  \tag{5.9}
\end{gather}

Let $\mu \in \mathcal{G}_{t}=\mathcal{M}_{t}^{b}$ and $\mu (\Omega _{\beta
}^{-p,R_{W}})=1$ for some $p>d.$ Since by (5.3) certainly $b_{i}\in
L^{1}(\mu ),$ Proposition 4.12 enables us to apply to $f,g$ the
(IbP)-formula (4.57) in directions $h:=e_{_{k}}\otimes \varphi $ with $%
\varphi \in T_{\beta }:=lin\{\varphi _{n}\}_{n\in \mathbb{Z}}.$ Therefrom,
in view of (5.8) and (5.9), we get that 
\begin{multline}
\int_{\Omega }\omega _{k}^{2Q-1}(\tau )(A_{\beta }\varphi ,\omega
_{k})_{H}Z^{-1}d\mu   \notag \\
\leq |\varphi |_{C_{\beta }^{{}}}\int_{\Omega }\left[ 4Q\omega
_{k}^{2Q-2}(\tau )+|\omega _{k}(\tau )|^{2Q-1}|F_{k}^{V,W}(\omega
)|_{L_{\beta }^{1}}\right] Z^{-1}d\mu   \tag{5.10}
\end{multline}%
and 
\begin{multline}
\int_{\Omega }\left[ \Delta \omega _{k}\right] ^{2Q-1}(A_{\beta }\varphi
,\omega _{k})_{H}Z^{-1}d\mu   \notag \\
\leq |\varphi |_{C_{\beta }^{{}}}^{{}}\int_{\Omega }\left[ 8Q[\Delta \omega
_{k}]^{2Q-2}+|\Delta \omega _{k}|^{2Q-1}|F_{k}^{V,W}(\omega )|_{L_{\beta
}^{1}}\right] Z^{-1}d\mu .  \tag{5.11}
\end{multline}%
Moreover, due to the conditions $|\omega _{k}|_{L_{\beta }^{2}},$ $%
|F_{k}^{V,W}(\omega )|_{L_{\beta }^{1}}\in L^{1}(\mu )$ imposed in ($\mathbf{%
H}$), both these inequalities extend by continuity to all $\varphi \in 
\mathcal{D}(A_{\beta }^{{}})$. \smallskip 

\textbf{(ii) }\ Now we proceed analogously to the proof of Proposition 4.12.
Namely, using the Yosida approximation of the Green function $\mathfrak{G}%
_{\tau }=A_{\beta }^{-1}\delta _{\tau }$ by $(\varphi _{\tau }^{(K)})_{K\in 
\mathbb{N}}\subset \mathcal{D}(A_{\beta }^{{}})$ with the properties
(4.75)--(4.77), by Lebesgue's dominated convergence theorem we conclude from
(5.10) written for $\varphi ^{(K)}:=\varphi _{\tau }^{(K)}$ that 
\begin{equation}
\int_{\Omega }\omega _{k}^{2Q}(\tau )Z^{-1}d\mu \leq |\mathfrak{G}_{\tau
}|_{C_{\beta }^{{}}}\int_{\Omega }\left[ 4Q\omega _{k}^{2Q-2}(\tau )+|\omega
_{k}(\tau )|^{2Q-1}|F_{k}^{V,W}(\omega )|_{L_{\beta }^{1}}\right] Z^{-1}d\mu
.  \tag{5.12}
\end{equation}%
Recall in this respect that by (3.30) and (3.31) 
\begin{equation*}
\sup_{\tau \in S_{\beta }}|\mathfrak{G}_{\tau }|_{C_{\beta }^{{}}}\leq 
\mathfrak{g}:=\left[ 2a\sqrt{\mathfrak{m}}\left( 1-e^{-\frac{a}{\sqrt{%
\mathfrak{m}}}\beta }\right) \right] ^{-1}.
\end{equation*}%
Therefore by H\"{o}lder's inequality applied in the RHS of (5.12) 
\begin{equation*}
\left( \int_{\Omega }\omega _{k}^{2Q}(\tau )Z^{-1}d\mu \right) ^{1/Q}\leq 
\mathfrak{g}\left[ 4Q+\left( \int_{\Omega }\omega _{k}^{2Q}(\tau )Z^{-1}d\mu
\right) ^{1/2Q}\left( \int_{\Omega }|F_{k}^{V,W}|_{L_{\beta }^{1}}^{2Q}d\mu
\right) ^{1/2Q}\right] ,
\end{equation*}%
which obviously yields%
\begin{equation}
\int_{\Omega }\omega _{k}^{2Q}(\tau )Z^{-1}d\mu \leq C_{2Q}:=(4\mathfrak{g)}%
^{Q}\left[ (4Q)^{Q}+\mathfrak{g}^{Q}c_{2Q}(F)\right] .  \tag{5.13}
\end{equation}%
Letting $\varepsilon \searrow 0$ in (5.13), from Fatou's lemma we readily
obtain the required estimate (5.4).\smallskip 

\textbf{(iii)} \ In a similar way, taking the approximation of $\mathfrak{G}%
_{\tau }-\mathfrak{G}_{\tau ^{\prime }}^{{}}$ by $\varphi _{\tau ,\tau
^{\prime }}^{(K)}:=\varphi _{\tau }^{(K)}-\varphi _{\tau ^{\prime }}^{(K)}$,
we conclude from (5.11) that 
\begin{multline}
\int_{\Omega }\left[ \Delta \omega _{k}\right] ^{2Q}Z^{-1}d\mu \leq |%
\mathfrak{G}_{\tau }^{{}}-\mathfrak{G}_{\tau ^{\prime }}|_{C_{\beta }^{{}}} 
\notag \\
\times \int_{\Omega }\left[ 8Q[\Delta \omega _{k}]^{Q-2}+|\Delta \omega
_{k}(\tau ,\tau ^{\prime })|^{Q-1}|F_{k}^{V,W}(\omega )|_{L_{\beta }^{1}}%
\right] Z^{-1}d\mu .  \tag{5.14}
\end{multline}%
Then again by (3.31) and H\"{o}lder's inequality 
\begin{multline*}
\left( \int_{\Omega }\left[ \Delta \omega _{k}(\tau ,\tau ^{\prime })\right]
^{2Q}Z^{-1}d\mu \right) ^{1/Q} \\
\leq \frac{a\mathfrak{g}}{\sqrt{\mathfrak{m}}}\rho (\tau ,\tau ^{\prime })%
\left[ 8Q+\left( \int_{\Omega }\left[ \Delta \omega _{k}(\tau ,\tau ^{\prime
})\right] ^{2Q}Z^{-1}d\mu \right) ^{1/2Q}\left( \int_{\Omega
}|F_{k}^{V,W}|_{L_{\beta }^{1}}^{2Q}d\mu \right) ^{1/2Q}\right] ,
\end{multline*}%
and thus for all $\varepsilon >0$%
\begin{equation}
\int_{\Omega }\left[ \Delta \omega _{k}(\tau ,\tau ^{\prime })\right]
^{2Q}Z^{-1}d\mu \leq \Delta C_{2Q}\cdot \rho ^{Q}(\tau ,\tau ^{\prime }) 
\tag{5.15}
\end{equation}%
with the constant 
\begin{equation}
\Delta C_{2Q}:=\left( 4\mathfrak{g}\frac{a}{\sqrt{\mathfrak{m}}}\right) ^{Q}%
\left[ (8Q)^{Q}+\left( \mathfrak{g}\frac{a\beta }{\sqrt{\mathfrak{m}}}%
\right) ^{Q}I_{2Q}\right]  \tag{5.16}
\end{equation}%
Letting $\varepsilon \searrow 0$ in (5.15), by Fatou's lemma we obtain the
required estimate (5.5). Finally we note that, according to (5.23), (5.25),
both $C_{2Q}$ are $\Delta C_{2Q}$ are linear functions of $I_{2Q}$ and in no
other way depend on $k$ and $\mu .$

\noindent $\blacksquare \smallskip $

\noindent \textbf{Remark 5.5.} \textbf{(i) } Except for Hypothesis (\textbf{H%
}) about the uniform integrability of $|F_{k}^{V,W}(\omega )|_{L_{\beta
}^{1}}^{Q}$ w.r.t. $\mu \in \mathcal{G}_{t}$, no more assumptions on the
interaction potentials (like, e.g., $(\mathbf{V})$, $(\mathbf{J})$, $(%
\mathbf{W})$ in Definition 3.2) are at all needed for the proof of Lemma
5.4. In fact, the result is completely determined by the regularity
properties of the Green function $\mathfrak{G}_{\tau }$ of the elliptic
operator $A_{\beta }$ in the Hilbert space $H:=L_{\beta }^{{}}.\smallskip $

Having obtained for $\mu \in \mathcal{G}_{t}$ the estimates from Lemma 5.4,
we are ready to prove our Main Theorem II previously announced in
Subsect.\thinspace 2.3.\smallskip

\noindent \textbf{Corollary 5.6. \ }\emph{Let Hypothesis} $\mathbf{(H)}$ 
\emph{be fulfilled}$,$ \emph{then the set }$\mathcal{G}_{t}$\emph{\ of all
tempered Euclidean Gibbs measures has the following properties:\smallskip }

\textbf{(i) }(cf.\textbf{\ Theorem II) \ }\emph{The uniform a priori bound
(2.18), i.e., }%
\begin{equation*}
\sup_{\mu \in \mathcal{G}_{t}}\sup_{k\in \mathbb{Z}^{d}}\int_{\Omega
}|\omega _{k}|_{C_{\beta }^{\alpha }}^{Q}\,d\mu (\omega )\leq C_{Q,\alpha
}<\infty ,\text{ \ }\forall Q\geq 1,\text{ }\alpha \in \lbrack 0,\frac{1}{2}%
),
\end{equation*}%
\emph{holds for all }$\mu \in \mathcal{G}_{t},$ \emph{where }$\mathcal{G}%
_{t}:=\mathcal{G}_{(s)t}^{R}$ \emph{for} \emph{the general QLS model (3.1)
and respectively }$\mathcal{G}_{t}:=\mathcal{G}_{(e)t}^{R}$ \emph{(}$%
\supseteq \mathcal{G}_{(s)t}^{R}$\emph{) for the particular QLS models
I--III (with the pair interaction satisfying the stronger decay Assumption }$%
(\mathbf{J}_{\mathbf{0}})$(ii)\emph{).\smallskip }

\textbf{(ii) }(cf.\textbf{\ Corollary after Theorem II) \ }\emph{Every} $\mu
\in \mathcal{G}_{t}$ \emph{is supported by} 
\begin{equation}
\Omega _{\text{\texttt{supp}}}:=\bigcap\nolimits_{\text{ }0\leq \alpha <1/2,%
\text{ }p>d}\mathcal{C}_{-p}^{\alpha }\subset \bigcap\nolimits_{\text{ }%
0\leq \alpha <1/2}[C_{\beta }^{\alpha }]^{\mathbb{Z}^{d}}\subset \Omega 
\tag{5.17}
\end{equation}%
\emph{(where, recalling the definition (3.13),}%
\begin{equation*}
\mathcal{C}_{-p}^{\alpha }:=l^{2}(\gamma _{-p};C_{\beta }^{\alpha })\text{ \ 
\emph{with} \ }||\omega ||_{-p,\alpha }^{2}:=\sum_{k\in \mathbb{Z}%
^{d}}(1+|k|)^{-2p}|\omega _{k}|_{C_{\beta }^{\alpha }}^{2}<\infty \text{ 
\emph{)}}.
\end{equation*}%
\emph{Moreover,} \emph{the set} $\mathcal{G}_{t}$ \emph{is compact} \emph{%
(and hence, is a Choquet simplex) w.r.t. the topology of weak convergence of
measures on} \emph{any Banach space }$\mathcal{C}_{-p}^{\alpha }$ \emph{with}
$p>2d$ \emph{and} $0\leq \alpha <1/2.$ \emph{For the particular QLS models
I, II} \emph{(as such with the local interaction),} $\mathcal{G}_{t}$ \emph{%
is even compact} \emph{w.r.t. the topology of weak convergence of measures on%
} \emph{any of } \emph{product spaces }$[C_{\beta }^{\alpha }]^{\mathbb{Z}%
^{d}}$ \emph{with} $0\leq \alpha <1/2$\emph{\ (these latter are considered
as the Fr\'{e}chet spaces equipped by the system of norms} $|\omega
_{k}|_{C_{\beta }^{\alpha }\text{ }},$ $k\in \mathbb{Z}^{d}$\emph{)}.
\smallskip

\textbf{Proof: (i)} \ We employ a standard argument related to Kolmogorov's
continuity criterion. More precisely, using inequality (3.d) in [BY82]
(which in turn is a consequence of the well-known Garsia--Rodemich--Rumsey
lemma), one can deduce from (5.5) that 
\begin{equation}
\int_{\Omega }\sup_{\tau \neq \tau ^{\prime }\in S_{\beta }^{{}}}\left[ 
\frac{|\omega _{k}(\tau )-\omega _{k}(\tau ^{\prime })|}{\rho ^{\alpha
}(\tau ,\tau ^{\prime })}\right] ^{2Q}d\mu (\omega )\leq C_{2Q,\alpha
}^{\prime }<\infty   \tag{5.18}
\end{equation}%
for all $Q>1$ and $\alpha \in \lbrack 0,\frac{1}{2}-\frac{1}{2Q}).$ Here $%
C_{2Q,\alpha }^{\prime }$ is some universal constant, which is the same for
all measures on $\Omega $ satisfying (5.5) and can be calculated explicitly
(cf. [BY82]).\ Letting $Q\rightarrow \infty ,$ both (5.4) and (5.18) give us
the required bound (2.18), namely that%
\begin{equation*}
\sup_{\mu \in \mathcal{G}_{t}}\sup_{k\in \mathbb{Z}^{d}}\int_{\Omega
}|\omega _{k}|_{C_{\beta }^{\alpha }}^{Q}\,d\mu (\omega )\leq C_{Q,\alpha
}<\infty ,\text{ \ }\forall Q\geq 1,\text{ }\alpha \in \lbrack 0,\frac{1}{2}%
).
\end{equation*}

\textbf{(ii)} \ The support property (5.17) follows immediately from (2.18)
by the estimate%
\begin{equation}
\sup_{\mu \in \mathcal{G}_{t}}E_{\mu }||\omega ||_{-p,\alpha }\leq \left(
\sup_{\mu \in \mathcal{G}_{t}}\sup_{k\in \mathbb{Z}^{d}}E_{\mu }|\omega
_{k}|_{C_{\beta }^{\alpha }}^{2}\right) ^{1/2}\sum\nolimits_{k\in \mathbb{Z}%
^{d}}(1+|k|)^{-2p}<\infty .  \tag{5.19}
\end{equation}%
On the other hand, taking into account (3.4) and (3.9), it is easy to check
that the embedding $\mathcal{C}_{-p^{\prime }}^{\alpha ^{\prime }}$ $%
\underset{\longrightarrow }{\subset }\mathcal{C}_{-p}^{\alpha }$ is compact
as soon as $\alpha ^{\prime }>\alpha $ and $p>p^{\prime }+d.$ Thus by
Prokhorov's criterion we conclude from (5.19) that the family $\mathcal{G}%
_{t}$ is tight w.r.t. the topology of weak convergence in all $\mathcal{C}%
_{-p}^{\alpha }$ with $p>2d,\mathcal{\ }$and hence in all $[C_{\beta
}^{\alpha }]^{\mathbb{Z}^{d}}$ with $0\leq \alpha <1/2$ as well. And
finally, $\mathcal{G}_{t}$ is closed in the above topologies according to
Remark 3.12 (iii).

\noindent $\blacksquare $

\subsubsection{Existence of Euclidean Gibbs measures}

Respectively applying Lemma 5.4 to the probability kernels $\pi _{\Lambda
}^{{}}(d\omega |\xi )$ of the local specification $(\pi _{\Lambda
})_{\Lambda \Subset \mathbb{Z}^{d}}$, one gets for them the following
Kolmogorov type estimates:\smallskip \newline
\textbf{Lemma 5.7. \ }\emph{For given} $\Lambda \subseteq \mathbb{Z}^{d}$%
\emph{\ and boundary condition} $\xi \in \Omega _{t},$ \emph{let the measure 
}$\pi _{\Lambda }^{{}}(d\omega |\xi )$ \emph{obey the integrability property}%
\begin{equation}
\int_{\Omega }\left[ |\omega _{k}|_{L_{\beta }^{2}}+|F_{k}^{V,W}(\omega
)|_{L_{\beta }^{1}}\right] ^{2Q}\pi _{\Lambda }(d\omega |\xi )\leq I_{2Q,k},%
\text{ \ }k\in \Lambda ,  \tag{5.20}
\end{equation}%
\emph{with some} $Q\in \mathbb{N}$ \emph{and corresponding} $I_{2Q,k}\in
(0,\infty )$\emph{.} \emph{Then (5.20) implies the moment estimates for all} 
$\tau ,\tau ^{\prime }\in S_{\beta }$%
\begin{gather}
\sup_{\tau \in S_{\beta }^{{}}}\int_{\Omega }\omega _{k}^{2Q}(\tau )\pi
_{\Lambda }^{{}}(d\omega |\xi )\leq H_{2Q}(1+I_{2Q,k})=:C_{2Q},  \tag{5.21}
\\
\int_{\Omega }[\omega _{k}(\tau )-\omega _{k}(\tau ^{\prime })]^{2Q}\pi
_{\Lambda }(d\omega |\xi )\leq \Delta H_{2Q}(1+I_{2Q,k})\cdot \rho ^{Q}(\tau
,\tau ^{\prime })=:\Delta C_{2Q}  \tag{5.22}
\end{gather}%
\emph{with absolute (i.e., independent on} $k,$ $\Lambda $ \emph{and }$\xi $%
\emph{) constants }$H_{2Q}$, $\Delta H_{2Q}\in (0,\infty ).$\smallskip

>From Lemma 5.7 one standardly derives a priori bounds in the spin spaces $%
C_{\beta }^{\alpha }$:\smallskip 

\noindent \textbf{Corollary 5.8. \ }\emph{Under the assumptions of Lemma
5.7, for all }$\alpha \in \lbrack 0,\frac{1}{2}-\frac{1}{2Q})$\emph{\ holds\ 
}%
\begin{equation}
\int_{\Omega }|\omega _{k}|_{C_{\beta }^{\alpha }}^{2Q}\,\pi _{\Lambda
}(d\omega |\xi )\leq H_{2Q,\alpha }(1+I_{2Q,k})=:C_{2Q,\alpha }<\infty  
\tag{5.23}
\end{equation}%
\emph{with an absolute (i.e., independent on} $k,$ $\Lambda $ \emph{and }$%
\xi $\emph{)} \emph{constant }$H_{2Q,\alpha }\in (0,\infty )$\emph{%
.\smallskip }

\textbf{Proof: \ }By arguments similar to those used in the proof of
Corollary 5.6, the above estimates (5.21), (5.22) yield that for any $\alpha
\in \lbrack 0,\frac{1}{2}-\frac{1}{2Q})$:%
\begin{equation}
\int_{\Omega }\sup_{\tau \neq \tau ^{\prime }\in S_{\beta }^{{}}}\left[ 
\frac{|\omega _{k}(\tau )-\omega _{k}(\tau ^{\prime })|}{\rho ^{\alpha
}(\tau ,\tau ^{\prime })}\right] ^{2Q}\pi _{\Lambda }(d\omega |\xi )\leq
C_{2Q,\alpha }^{\prime }<\infty .  \tag{5.24}
\end{equation}%
An \emph{important} point here is that (as was calculated in [BY82]) the
constant $C_{2Q,\alpha }^{\prime }$ in (5.24) linearly depends on $\Delta
C_{2Q}$ in the RHS in (5.22). Hence, having regard to Lemma 5.7, the desired
estimate (5.23) holds with the constant $C_{2Q,\alpha }$ which lineary
depends on $I_{2Q}$ in the RHS in (5.20)$.$

\noindent $\blacksquare \smallskip $

In turn, Corollary 5.8 now readily implies the existence of $\mu \in 
\mathcal{G}_{t}$ as was announced in our Main Theorem I in
Subsect.\thinspace 2.3.\smallskip

\noindent \textbf{Corollary 5.9. \ }\emph{Let} \emph{Hypothesis} $\mathbf{(H}%
_{loc}\mathbf{)}$ \emph{be fulfilled} \emph{for some }$Q>1$ \emph{and} \emph{%
a given boundary condition} $\xi \in \Omega _{t}$ \emph{satisfying}%
\begin{equation}
\xi \in \lbrack C_{\beta }^{\alpha }]^{\mathbb{Z}^{d}}\text{ \ \emph{and} \ }%
\sup_{k\in \mathbb{Z}^{d}}|\xi _{k}|_{C_{\beta }^{\alpha }}<\infty . 
\tag{5.25}
\end{equation}%
\emph{for} \emph{some }$\alpha :=\alpha (\xi )>0.$ \emph{Then the following
statements hold:\smallskip }

\textbf{(i) }(cf. A Priori Estimate (2.19) in \textbf{Theorem III})\textbf{\
\ }\emph{For all }$0\leq \alpha ^{\prime }<\min \{\alpha ,\frac{1}{2}-\frac{1%
}{2Q}\}$%
\begin{equation}
\sup_{\Lambda \Subset \mathbb{Z}^{d}}\sup_{k\in \mathbb{Z}^{d}}\int_{\Omega
}|\omega _{k}|_{C_{\beta }^{\alpha ^{\prime }}}^{2Q}\,\pi _{\Lambda
}(d\omega |\xi )<\infty  \tag{5.26}
\end{equation}%
\emph{and hence the} \emph{family} $\left\{ \pi _{\Lambda }(d\omega |\xi
)\right\} _{\Lambda \Subset \mathbb{Z}^{d}}$\emph{\ is tight in all spaces} $%
\mathcal{C}_{-p}^{\alpha ^{\prime }}$ \emph{with} $p>2d.\smallskip $

\textbf{(ii) }(Existence of $\mu \in \mathcal{G}_{t};$ cf. \textbf{Theorem
I) \ }\emph{The set} $\mathcal{G}_{t}$ \emph{is not empty as such containing
each accumulation point as} $\Lambda \nearrow \mathbb{Z}^{d}$ \emph{for} $%
\left\{ \pi _{\Lambda }(d\omega |\xi )\right\} _{\Lambda \Subset \mathbb{Z}%
^{d}}.\smallskip $

\textbf{Proof: (i) }\ (5.26) follows immediately from (5.24) due to the
special choice (5.22) of the boundary condition $\xi $. On the other hand,
by (3.4), (3.9) and the definition (3.13) of the Banach spaces $\mathcal{C}%
_{-p}^{\alpha }:=l^{2}(\gamma _{-p};C_{\beta }^{\alpha })$ it is easy to
check that the embedding $\mathcal{C}_{-p^{\prime }}^{\alpha ^{\prime }}$ $%
\underset{\longrightarrow }{\subset }\mathcal{C}_{-p}^{\alpha }$ is compact
as soon as $\alpha ^{\prime }>\alpha $ and $p>p^{\prime }+d.$ Thus by
Prokhorov's criterion we conclude from (5.26) that the family of
distributions $\pi _{\Lambda ^{(K)}}(d\omega |\xi ),$ where $\Lambda
^{(K)}\nearrow \mathbb{Z}^{d}$ as $K\rightarrow \infty ,$ is tight in all $%
\mathcal{C}_{-p}^{\alpha }$ with $p>2d$ and $\alpha \in \lbrack 0,\frac{1}{2}%
-\frac{1}{2Q}).\mathcal{\ \smallskip }$

\textbf{(ii)} \ By (i) there exists a subsequence $\pi _{\Lambda
^{(K_{L})}}(d\omega |\xi ),$ $L\in \mathbb{N},$ which converges weakly to
some probability measure $\mu _{\ast }$ on $\mathcal{C}_{-p}^{\alpha }.$
Since $\mathcal{C}_{-p}^{\alpha }$ is continuously embedded into $\Omega
_{-p}^{R},$ this subsequence converges to $\mu _{\ast }$ also weakly on all
Polish spaces $\Omega _{-p}^{R}$ with $p>2d.$ This means by Proposition 5.3
that $\mu _{\ast }\in \mathcal{G}_{t}:=\mathcal{G}_{(s)t}^{R}$.

\noindent $\blacksquare \smallskip $

Next, we analyze in more detail the dependence on a boundary condition $\xi
\in \Omega _{t}$ in the a priori estimates for $\{\pi _{\Lambda }(d\omega
|\xi )\}_{\Lambda \Subset \mathbb{Z}^{d}}.\smallskip $

\noindent \textbf{Corollary 5.10. (i) \ }\emph{Let }$\xi \in \Omega
_{t}:=\Omega _{(s)t}^{R}$ \emph{and hence (in the notation of
Subsect.\thinspace 3.2.3) }%
\begin{equation}
||\xi ||_{\mathcal{L}_{-p}^{R}}^{2}:=\sum\nolimits_{k\in \mathbb{Z}%
^{d}}(1+|k|)^{-p}|\xi _{k}|_{L_{\beta }^{R}}^{2}<\infty \text{ \ }\emph{\ }%
\text{\emph{for some} \ }p=p(\xi )>d.  \tag{5.27}
\end{equation}%
\emph{For the general QLS\ model (3.1) then holds for} $Q\geq 1,$ $p^{\prime
}>pRQ+d$ \emph{and} $\alpha \in \lbrack 0,\frac{1}{2}-\frac{1}{Q}):$%
\begin{equation}
\sup_{\Lambda \Subset \mathbb{Z}^{d}}\sum\nolimits_{k\in \Lambda
}(1+|k|)^{-p^{\prime }}\int_{\Omega }|\omega _{k}|_{C_{\beta }^{\alpha
}}^{Q}\pi _{\Lambda }(d\omega |\xi )<\infty .  \tag{5.28}
\end{equation}%
\emph{\smallskip }

(\textbf{ii) }(cf. Estimate (2.20) in \textbf{Theorem III})\textbf{\ \ }%
\emph{Let }$\xi \in \Omega _{t}:=\Omega _{(e)t}^{R}$ \emph{and hence }%
\begin{equation}
||\xi ||_{\mathcal{L}_{-\delta }^{R}}^{2}:=\sum\nolimits_{k\in \mathbb{Z}%
^{d}}\exp (-\delta |k|)|\xi _{k}|_{L_{\beta }^{R}}^{2}<\infty \text{ \ }%
\emph{\ }\text{\emph{for all} \ }\delta >0.  \tag{5.29}
\end{equation}%
\emph{For the particular QLS\ models I, II (and also for the QLS\ model III
satisfying additionally Assumption }(\textbf{J}$_{\mathbf{0}})($ii$)$\emph{)}
\emph{then} \emph{holds for all} $Q\geq 1,$ $\delta >0$ \emph{and} $\alpha
\in \lbrack 0,\frac{1}{2}),$: 
\begin{equation}
\sup_{\Lambda \Subset \mathbb{Z}^{d}}\sum\nolimits_{k\in \Lambda }\exp
(-\delta |k|)\int_{\Omega }|\omega _{k}|_{C_{\beta }^{\alpha }}^{Q}\pi
_{\Lambda }(d\omega |\xi )<\infty .  \tag{5.30}
\end{equation}%
\emph{\smallskip }

\textbf{Proof: }The proof of both statements (i) and (ii) follows the same
line. Namely, having regard to Corollary 5.8, it suffices to check that
(5.27) resp. (5.29) implies for all $Q\in \mathbf{N}$ and $k\in \mathbb{Z}%
^{d}$%
\begin{equation}
\sup_{\Lambda \Subset \mathbb{Z}^{d}\backslash \{k\}}\int_{\Omega }\left[
|\omega _{k}|_{L_{\beta }^{2}}+|F_{k}^{V,W}(\omega )|_{L_{\beta }^{1}}\right]
^{Q}\pi _{\Lambda }(d\omega |\xi ):=I_{Q,k}  \tag{5.31}
\end{equation}%
with a nonnegative sequence $I_{Q}=(I_{Q,k})_{k\in \mathbb{Z}^{d}}$ such
that for all $p^{\prime }>pRQ+d$ 
\begin{equation}
\sum\nolimits_{k\in \mathbb{Z}^{d}}(1+|k|)^{-p^{\prime }}I_{Q,k}<\infty , 
\tag{5.32}
\end{equation}%
or respectively for all $\delta >0$%
\begin{equation}
\sum\nolimits_{k\in \Lambda }\exp (-\delta |k|)I_{Q,k}<\infty .  \tag{5.33}
\end{equation}%
The verification of these key conjections (together with Hypotheses (H) and
(H$_{loc}))$ will be performed in Subsect.\thinspace 7.3.

\noindent $\blacksquare \smallskip $

\noindent \textbf{Remark 5.11. (i) \ }For the\emph{\ translation invariant
systems}, by an obvious modification of the arguments used above one can
construct the so-called \emph{periodic }Euclidean Gibbs measures $\mu
_{per}\in \mathcal{G}_{t}.$ They are defined as accumulation points for the
family $\mu _{\mathfrak{T}(\Lambda )},$ $\Lambda \Subset \mathbb{Z}^{d},$ of
local Gibbs distributions with periodic boundary conditions and, hence, are
certainly $\mathbb{Z}_{0}^{d}-$translation invariant (cf. the related
discussion in Remark 3.10 (iv)). We emphasize that in doing so the crucial
estimate 
\begin{equation}
\sup_{\Lambda \Subset \mathbb{Z}^{d}}\sup_{k\in \mathbb{Z}^{d}}\int_{\Omega
_{\Lambda }}|\omega _{k}|_{C_{\beta }^{\alpha }}^{2Q}\text{\thinspace }\mu _{%
\mathfrak{T}(\Lambda )}(d\omega _{\Lambda })<\infty   \tag{5.34}
\end{equation}%
holds, provided we assume the following Hypothesis to be satisfied:

\begin{description}
\item[$\mathbf{(H}_{per}\mathbf{)}$] \emph{For some fixed} $Q>1$ 
\begin{equation*}
I_{2Q,per}(F):=\sup_{\Lambda \Subset \mathbb{Z}^{d}}\sup_{k\in \Lambda
}\int_{\Omega _{\Lambda }}\left[ |\omega _{k}|_{L_{\beta
}^{2}}+|F_{k}^{V,W}(\omega )|_{L_{\beta }^{1}}\right] ^{2Q}\mu _{\mathfrak{T}%
(\Lambda )}(d\omega _{\Lambda })<\infty .
\end{equation*}%
$\smallskip $
\end{description}

\textbf{(ii) }(\emph{Existence of Superstable Gibbs States})\textbf{\ \ }As%
\textbf{\ }already\textbf{\ }mentioned\textbf{\ }in Remark 2.4 (ii), $%
\mathcal{G}_{t}:=\mathcal{G}_{(s)t}^{R}$ contains a class $\mathcal{G}%
_{(ss)t}$ of the so-called Ruelle type \textquotedblleft \emph{superstable}%
\textquotedblright\ Gibbs measures, which for the considered quantum lattice
systems (in the particular case $R=2$ only) has been introduced in [PY94] by
the support condition 
\begin{gather}
\mathcal{G}_{(ss)t}:=\left\{ \mu \in \mathcal{G}\text{ }\left\vert \ \mu
(\Omega _{(ss)t})=1\right. \right\} ,  \notag \\
\Omega _{(ss)t}:=\left\{ \omega \in \Omega \left\vert \ \sup_{n\in \mathbb{N}%
}\left[ (1+2n)^{-d}\sum\limits_{|k|\leq n}|\omega _{k}|_{L_{\beta }^{2}}^{2}%
\right] <\infty \right. \right\} .  \tag{5.35}
\end{gather}%
But for any measure $\mu \in \mathcal{M}(\Omega ),$ which is \emph{%
translation invariant} and satisfies the a priori estimates (2.18), the
support condition (even much stronger than (5.35)) holds for all $Q\geq 1$
and $\alpha \in \lbrack 0,1/2),$ namely:%
\begin{equation}
\sup_{n\in \mathbb{N}}\left\{ (1+2n)^{-d}\sum\nolimits_{|k|\leq n}|\omega
_{k}|_{C_{\beta }^{\alpha }}^{Q}\right\} \leq C_{Q,\alpha }(\omega )<\infty
,\quad \forall \omega \in \Omega \text{ \ (}\mu -\text{a.e.).}  \tag{5.36}
\end{equation}%
The latter follows from the Birkhoff--Khinchin ergodic theorem (cf. e.g.
[DeuS89]) applied to the stationary process $\omega _{k},$ $k\in \mathbb{Z}%
^{d},$ on the probability space $(\Omega ,\mu ).$ Together with (i) this
means that we can refine the statement of Theorem I by claiming the
existence of $\mu _{per}\in \mathcal{G}_{(ss)t}$ with the additional support
property (5.36).

\subsection{Comparison with the stochastic dynamics method}

Here we would like to point out some advantages of our approach in
comparison to the stochastic dynamics employed in [AKRT01] to the quantum
lattice models (2.1).\smallskip 

In that paper we have restricted ourselves to the case of harmonic pair
interactions (described by a dynamical matrix $0\leq \mathbf{D}%
:=(a_{k,j})_{k,j\in \mathbb{Z}^{d}}\in \mathcal{L}(l^{2}(\mathbb{Z}^{d}))$,
see Remark 3.3 (iv)); in more generality the method could be applied to
many-particle interactions of \emph{at most quadratic growth}. Concerning
the one-particle potentials $V_{k}\in C_{b,loc}^{2}(\mathbb{R\rightarrow }%
\mathbb{R}),$ the following \emph{semi-monotonicity} 
\begin{equation}
(V_{k}^{\prime }(q_{1})-V_{k}^{\prime }(q_{2}))(q_{1}-q_{2})\geq
K_{5}^{-1}(q_{1}-q_{2})^{2}-L_{5}  \tag{5.37}
\end{equation}%
and \emph{at most polynomial growth} 
\begin{equation}
|V_{k}^{\prime }(q_{1})|\leq K_{6}(1+|q_{1}|)^{R}  \tag{5.38}
\end{equation}%
conditions with \emph{some fixed} $K_{5},$ $K_{6},$ $L_{5}>0$ and $R\geq 1$
are required to hold \emph{uniformly} for all $k\in \mathbb{Z}^{d}$ and $%
q_{1},$ $q_{2}\in \mathbb{R}.$

So, in [AKRT01] we have firstly constructed a Markov process $%
x_{t}=(x_{k,t})_{k\in \mathbb{Z}^{d}},$ $t\geq 0,$ which gives the unique
(generalized) solution to the (so-called \emph{Langevin}) stochastic
evolution equation with a drift term being the logarithmic gradient $%
b=(b_{k})_{k\in \mathbb{Z}^{d}}$ of the measures $\mu \in \mathcal{G}_{t}.$
More precisely, $x_{t},$ $t\geq 0,$ takes values in the Banach state space $%
X:=\mathcal{C}_{-p}=l^{2}(\mathbb{Z}^{d}\rightarrow C_{\beta };\gamma
_{-p})\subset \Omega $ (with large enough $p=p(R)>d)$ and satisfies the
following infinite system of stochastic partial differential equations
(SPDE's): 
\begin{equation}
\left\{ 
\begin{array}{c}
\frac{\partial }{\partial t}x_{k,t}=-\frac{1}{2}\left[ A_{\beta
}^{{}}x_{k,t}+\sum_{\text{ }j\in \mathbb{Z}^{d}}a_{k,j}x_{j,t}+V_{k}^{\prime
}(x_{k,t})\right] +\dot{w}_{k,t} \\ 
k\in \mathbb{Z}^{d}\quad (t>0,\text{ }\tau \in S_{\beta }^{{}}).%
\end{array}%
\right.  \tag{5.39}
\end{equation}%
Here $A_{\beta }^{{}}=-\mathfrak{m}\partial ^{2}/\partial \tau ^{2}+a^{2}%
\mathbf{1}$ is the self-adjoint operator in $L_{\beta }^{2}$ introduced in
Subsect.$\,$3.3.1 and $\dot{w}_{k,t}(\tau )$ is a Gaussian white noise on $%
\Omega \times \lbrack 0,\infty )$ (heuristically, $E\dot{w}%
_{k_{1},t_{1}}(\tau _{1})\times \dot{w}_{k_{2},t_{2}}(\tau _{2})=\delta
_{k_{1}-k_{2}}^{{}}\delta _{t_{1}-t_{2}}^{{}}\delta _{\tau _{1}-\tau
_{2}}^{{}}$). In the trivial case when $\mathbf{D}=0$ and $V=0,$ the
solution of (5.39), starting with initial data $g_{0}:=\zeta \in \mathcal{C}%
_{-p}^{{}},$ is explicitly given by the Ornstein--Uhlenbeck process $%
g_{t}=(g_{k,t})_{k\in \mathbb{Z}^{d}},$ $t\geq 0,$%
\begin{equation}
g_{k,t}:=e^{-tA_{\beta }^{{}}/2}\zeta
_{k}+\int\limits_{0}^{t}e^{(t-s)A_{\beta }^{{}}/2}dw_{k,s},\quad k\in 
\mathbb{Z}^{d},\text{ }t\geq 0.  \tag{5.40}
\end{equation}%
Taking into account the regularity properties (3.31) of the Green function $%
\mathfrak{G}_{\tau }=A_{\beta }^{-1}\delta _{\tau }$, one can deduce from
(5.40) that $g_{t},$ $t\geq 0,$ possesses a continuous modification in the
spaces of H\"{o}lder loops $\mathcal{C}_{-p}^{\alpha },$ $p>d,$ $0<\alpha
<1/2,$ and its polynomial moments are \emph{ultimately bounded}, i.e., $%
\forall Q\geq 1$%
\begin{equation}
\sup_{g_{0}\in \mathcal{C}_{-p}^{{}}}\sup_{k\in \mathbb{Z}%
^{d}}\limsup_{t\rightarrow \infty }E|g_{k,t}|_{C_{\beta }^{\alpha
}}^{Q}<\infty .  \tag{5.41}
\end{equation}%
\begin{equation}
\sup_{g_{0}\in \mathcal{C}_{-p}^{{}}}\sup_{k\in \mathbb{Z}%
^{d}}\limsup_{t\rightarrow \infty }E|g_{k,t}|_{C_{\beta }^{\alpha
}}^{Q}<\infty .  \tag{5.42}
\end{equation}%
Moreover, $g_{t},$ $t\geq 0,$ is ergodic with the unique invariant (and also
reversible) distribution $\gamma _{\text{\texttt{inv}}}(d\omega
):=\prod_{k\in \mathbb{Z}^{d}}\gamma _{\beta }(d\omega _{k}),$ i.e., the
corresponding laws $\mathcal{L}(g_{t})$ converge weakly, as $t\rightarrow
\infty ,$ to $\gamma _{\text{\texttt{inv}}}$. So it is reasonable to compare
the solution $x_{t},$ $t\geq 0,$ of the nonlinear problem (5.39) with the
Gaussian process $g_{t},$ $t\geq 0.$ If their initial values coincide, i.e., 
$x_{0}=g_{0}:=\zeta ,$ then for the deviation process $y_{t}:=x_{t}-g_{t},$%
\begin{equation*}
y_{k,t}=x_{k,t}-g_{k,t}=\int\limits_{0}^{t}e^{(t-s)A_{\beta }^{{}}/2}\left(
\sum_{j\in \mathbb{Z}^{d}}a_{k,j}x_{j,s}+V_{k}^{\prime }(x_{k,s})\right) ds,
\end{equation*}%
some helpful \emph{energy estimates} hold. So, under the above assumptions
(5.37) and (5.38), for all $Q\geq 1$:%
\begin{equation}
\sup_{\substack{ x_{0}\in \mathcal{C}_{-p}  \\ x_{0}=g_{0}}}\sup_{k\in 
\mathbb{Z}^{d}}\limsup_{t\rightarrow \infty }E||y_{k,t}||_{C_{\beta
}^{{}}}^{Q}<\infty ,\quad \sup_{\substack{ x_{0}\in \mathcal{C}_{-p}^{{}} 
\\ x_{0}=g_{0}}}\sup_{k\in \mathbb{Z}^{d}}\limsup_{t\rightarrow \infty }%
\frac{1}{t}\int\limits_{t}^{2t}E||y_{k,s}||_{W_{\beta }^{2,1}}^{Q}ds<\infty .
\tag{5.43}
\end{equation}%
Combined, (5.42) and (5.43) give us a crucial estimate for the process $%
x_{t},$ $t\geq 0,$ we shall use in the sequel: 
\begin{equation}
\sup_{x_{0}\in \mathcal{C}_{-p}^{{}}}\sup_{k\in \mathbb{Z}%
^{d}}\limsup_{t\rightarrow \infty }\frac{1}{t}\int%
\limits_{t}^{2t}E||x_{k,s}||_{C_{\beta }^{\alpha }}^{Q}ds:=I_{Q,\alpha
}^{\prime }<\infty .  \tag{5.44}
\end{equation}

As was further shown in [AKRT01], 
\begin{equation*}
(\mathbb{P}_{t}^{{}}f)(\omega ):=E\{f(x_{t}^{{}})\,|\,x_{0}^{{}}=\omega
\},\quad \omega \in \mathcal{C}_{\beta }^{-p},
\end{equation*}%
is a Feller transition semigroup $\mathbb{P}_{t},$ $t\geq 0,$ in the space $%
C_{b}(\mathcal{C}_{-p}^{{}})$ of all bounded continuous functions $f:%
\mathcal{C}_{-p}\rightarrow \mathbb{R}.$ Let $\mathcal{R}(\mathcal{C}_{-p})$
resp. $\mathcal{I}(\mathcal{C}_{-p}^{{}})$ denote the family of all \emph{%
reversible} resp. \emph{invariant} distributions for the Markov process $%
x_{t},$ $t\geq 0$. Then the following basic relation between Gibbs and
reversible distributions is true: 
\begin{equation}
\mu \in \mathcal{G}\text{\quad and\quad }\mu (\mathcal{C}_{-p}^{{}})=1%
\Longleftrightarrow \mu \in \mathcal{R}(\mathcal{C}_{-p})\subseteq \mathcal{I%
}(\mathcal{C}_{-p}^{{}}).  \tag{5.45}
\end{equation}%
(for the proof of (5.45) involving the Ito stochastic calculus and
(IbP)-formulas cf. [Fu91, KRZ96]). Moreover, in our situation one can
directly verify (cf. e.g. [Fu91]) that the finite volume Gibbs measures $\mu
_{\Lambda }(d\omega |\xi ),$ $\Lambda \Subset \mathbb{Z}^{d},$ $\xi \in 
\mathcal{C}_{-p}^{{}},$ are exactly the reversible distributions for the
corresponding cutoff dynamics $x_{t}^{\Lambda }=(x_{k,t}^{\Lambda })_{k\in
\Lambda }$ in $C_{\beta }^{\Lambda },$%
\begin{equation}
\left\{ 
\begin{array}{c}
\frac{\partial }{\partial t}x_{k,t}^{\Lambda }=-\frac{1}{2}\left[ A_{\beta
}^{{}}x_{k,t}^{\Lambda }+\sum_{\text{ }j\in \Lambda }a_{k,j}x_{j,t}^{\Lambda
}+V_{k}^{\prime }(x_{k,t}^{\Lambda })+\sum_{\text{ }j\in \Lambda
^{c}}a_{k,j}\xi _{j,t}^{\Lambda }\right] +\dot{w}_{k,t} \\ 
k\in \Lambda \quad (t>0,\text{ }\tau \in S_{\beta }^{{}}).%
\end{array}%
\right.   \tag{5.46}
\end{equation}%
An important point is that, for fixed boundary and initial conditions $\xi ,$
$\zeta \in \mathcal{C}_{-p},$ the analogues of (5.43) and (5.44) hold
uniformly in $\Lambda \Subset \mathbb{Z}^{d}$ for the corresponding
solutions $x_{t}^{\Lambda }$, $t\geq 0,$ starting from $x_{k,0}^{\Lambda
}=\zeta _{k},$ $k\in \Lambda ,$ e.g., 
\begin{equation}
\sup_{\Lambda \Subset \mathbb{Z}^{d},\text{ }k\in \Lambda }\sup_{t\geq 0}%
\frac{1}{t}\int\limits_{t}^{2t}E|x_{k,s}^{\Lambda }|_{C_{\beta }^{\alpha
}}^{2}\,ds:=I_{Q,\alpha }^{\prime }(\xi ,\zeta )<\infty .  \tag{5.47}
\end{equation}

Thus, in order to get the required information on $\mu \in \mathcal{G}%
_{t}\subseteq \mathcal{I}(\mathcal{C}_{-p}^{{}}$), one could apply standard
tools used in the theory of SDE's for the long-time analysis of diffusion
processes. So, by the ergodic theorem for invariant distributions, (5.47)
implies immediately that 
\begin{equation*}
\sup_{\substack{ \mu \in \mathcal{I}(\mathcal{C}_{-p})  \\ k\in \mathbb{Z}%
^{d} }}\int_{\Omega }|\omega _{k}|_{C_{\beta }^{\alpha }}^{Q}\,d\mu (\omega
)\leq \sup_{\substack{ x_{0}\in \mathcal{C}_{-p}  \\ k\in \mathbb{Z}^{d}}}%
\limsup_{t\rightarrow \infty }E|x_{k,t}|_{C_{\beta }^{\alpha
}}^{Q}=I_{Q,\alpha }^{\prime }<\infty ,
\end{equation*}%
which confirms the corresponding result in Theorem II in this special
situation. Moreover, the existence of invariant measures $\mu \in \mathcal{I}%
(\mathcal{C}_{-p}^{{}}$) is an\emph{\ }easy\emph{\ }consequence of (5.47)
and the usual Bogolyubov-Krylov argument. For the existence of $\mu \in 
\mathcal{G}=\mathcal{R}(\mathcal{C}_{-p}$) stated in Theorem I, by
Proposition 5.3 it suffices to prove the tightness in $\mathcal{C}_{-p}$of
the family $\{\pi _{\Lambda }(d\omega |0)\}_{\Lambda \Subset \mathbb{Z}^{d}}$
of local Gibbs distributions with fixed boundary condition $\xi =0.$ By
Prokhorov's criterion this is a consequence of the uniform boundedness of
their moments in $\mathcal{C}_{-p}^{\alpha },$ $0<\alpha <1/2,$%
\begin{equation}
\sup_{\Lambda \Subset \mathbb{Z}^{d},\text{ }k\in \Lambda }\int_{(C_{\beta
}^{{}})^{\Lambda }}|\omega _{k}|_{C_{\beta }^{\alpha }}^{2}\,d\mu _{\beta
,\Lambda }^{{}}(\omega _{\Lambda }|0):=I_{Q,\alpha }^{\prime }(\xi ,\zeta
)<\infty  \tag{5.48}
\end{equation}%
(see also Subsect.\thinspace 5.5.2). Since the finite volume dynamics (5.46)
are ergodic, (5.48) immediately follows from the estimates (5.47)
above.\smallskip

To summarize, a disadvantage of the stochastic dynamics method is (apart
from the fact that a lot of advanced technique it involved) that the
assumptions on the solvability of the related infinite dimensional evolution
equations are usually more restrictive than those under which the
(IbP)-formula (4.57) makes sense. It should be also emphasized that the
ergodicity problem for the interacting stochastic systems with unbounded
spins like (5.30) (except for the special cases of linear or strictly
dissipative ones, cf. [DPZ96, AKRT01]) is extremely difficult itself. Thus
in this paper, rather to discuss the processes and their invariant
(reversible) distributions, we start \emph{directly} from the definition of
Gibbs measures as solutions to the (IbP)-formulas.{\Huge \bigskip }

\begin{center}
{\Huge Part II: Abstract Setting\bigskip }
\end{center}

\section{Symmetrizing measures on Banach (e.g. loop) spaces}

\noindent In order to clarify the concept and stress the key ideas, we now
put the problem in a general framework of symmetrizing measures on Banach
lattices. In this section we discuss in detail the case of symmetrizing
measures $\mu \in \mathcal{M}^{b}(X)$ on a single Banach (in particular,
loop) space $X$. Assuming that the logarithmic derivative $b$ of these
measures has a linear component $A$ (being a positive selfadjoint operator
in some tangent space $H\supseteq X$) and the nonlinear one $F$ (possessing
certain coercivity properties w.r.t. $H$), we investigate the interplay
between the properties of the operators $A,$ $F$ and the integrability and
support properties of $\mu $. Apart from its origin in infinite dimensional
stochastic analysis, this setting may also be of independent interest
because of its possible applications to the study of time-reversible
distributions of stochastic evolution equations (in particular, SPDE's of
reaction-diffusion type).

\subsection{(IbP)-formula in a general setting}

Having in mind the concrete properties of the one-particle Euclidean
distributions from Subsect.\thinspace 3.3.1, one can formulate an abstract
setting of the problem as follows:

Let us fix a tangent space $H$ as a separable Hilbert space with inner
product $(\cdot ,\cdot )_{H}$ and norm $|\cdot |_{H}$. Let 
\begin{equation}
X\subset B\subset H\subset B^{\ast }  \tag{6.1}
\end{equation}%
be a rigging of $H$ by a locally convex vector (e.g., Banach) space $X$ and
by reflexive Banach spaces $B$ and its dual $B^{\ast }$ with (at least
outside zero) \emph{differentiable} norms $|\cdot |_{B}$ and $|\cdot
|_{B^{\ast }}$ respectively. The embeddings in (6.1) are dense and
continuous, which implies the relation 
\begin{equation}
\iota _{{}}^{-1}|\cdot |_{B^{\ast }}\leq |\cdot |_{H}\leq \iota
_{{}}^{{}}|\cdot |_{B}.  \tag{6.2}
\end{equation}%
with some finite constant $\iota >0.$ The duality between $B$ and $B^{\ast }$
is given by the inner product in $H$ and will also be denoted by $(\cdot
,\cdot )_{H}.$

Let $A$ be a \emph{positive self-adjoint }operator in $H$ with inverse\emph{%
\ }$A^{-1}$ \emph{of finite trace}. The operator $A$ has discrete spectrum $%
\{\lambda _{n}\}_{n\in \mathbb{N}}$ and a complete orthonormal system of
eigenvectors $\{\varphi _{n}\}_{n\in \mathbb{N}}:=bas(H),$%
\begin{equation}
A\varphi _{n}=\lambda _{n}\varphi _{n},\quad
Tr_{H}^{{}}A^{-1}=\sum\limits_{n\in \mathbb{N}}\lambda _{n}^{-1}<\infty . 
\tag{6.3}
\end{equation}%
Let $\mathbb{P}_{K}^{{}},$ $K\in \mathbb{N},$ be the finite-dimensional
projections generated by the first $K$ vectors of $bas(H)$, i.e., 
\begin{equation}
\mathbb{P}_{K}^{{}}x:=\sum\limits_{n=1}^{K}(x,\varphi _{n})_{H}^{{}}\varphi
_{n},\quad x\in H.  \tag{6.4}
\end{equation}%
The following is crucial: we assume that $\varphi _{n}\in X,\ \forall n\in 
\mathbb{N},$ and $\{\varphi _{n}\}_{n\in \mathbb{N}}$ is a \emph{Schauder
basis} in $B$, so that there exist finite constants $\kappa $ and $\varsigma 
$ such that 
\begin{equation}
\sup_{n\in \mathbb{N}}\{|\varphi _{n}|_{B}^{{}}\}\leq \kappa \text{,\quad }|%
\mathbb{P}_{K}^{{}}x|_{B}^{{}}\leq \varsigma |x|_{B}^{{}}\text{\quad for all 
}x\in B,\ K\in \mathbb{N}.  \tag{6.5}
\end{equation}

In accordance with the notation in Subsect.\thinspace 4.2.1, $%
C_{b}^{1}(X;\varphi )$ resp. $C_{b,loc}^{1}(X;\varphi )$ will be the spaces
of all functions $f:X\rightarrow \mathbb{R}$ which are continuous and
globally resp. locally bounded together with their partial derivative $%
\partial _{\varphi }f:X\rightarrow \mathbb{R}$ along the fixed direction $%
\varphi \in X.$ Then we define 
\begin{equation*}
C_{b}^{1}(X):=\bigcap_{\varphi \in X}C_{b}^{1}(X;\varphi )\text{\quad
resp.\quad }C_{b,loc}^{1}(X):=\bigcap_{\varphi \in X}C_{b,loc}^{1}(X;\varphi
),
\end{equation*}%
and their subspaces $\mathcal{F}C_{b}^{1}(X)$ resp. $\mathcal{F}%
C_{b,loc}^{1}(X)$ of cylinder functions w.r.t. $bas(H)$. To shorten the
notation, $\partial _{n}f:=\partial _{\varphi _{n}}f$ will stand for the
derivatives along the basis vectors $\varphi _{n}$, $n\in \mathbb{N}$.

Furthermore, let $F:X\rightarrow B^{\ast }$ be a nonlinear mapping which is 
\emph{continuous and locally bounded together with its partial derivatives}
in directions\emph{\ }$\varphi _{n}\in bas(H),$ i.e., $F\in \bigcap_{n\in 
\mathbb{N}}C_{b,loc}^{1}(X\rightarrow B^{\ast };\varphi ).$ We define a 
\emph{measurable vector} \emph{field} $X\ni x\rightarrow
b(x):=\{b_{n}(x)\}_{n\in \mathbb{N}}\in \mathbb{R}^{\mathbb{N}}$ by 
\begin{equation}
b_{n}(x):=-(A\varphi _{n},x)_{H}-(F(x),\varphi _{n})_{H},\quad n\in \mathbb{N%
},\ x\in X.  \tag{6.6}
\end{equation}%
>From the above assumptions 
\begin{equation}
b_{n}\in C_{b,loc}^{1}(X)\quad \text{and\quad }|b_{n}(x)|\leq \kappa
(\lambda _{n}|x|_{B^{\ast }}+|F(x)|_{B^{\ast }}).  \tag{6.7}
\end{equation}%
Let $\mathcal{M}^{b}(X)$ denote the family of all Borel probability measures 
$\mu $ on $X$ which satisfy the following \emph{(IbP)-formula} 
\begin{equation}
\int_{X}\partial _{n}f\,d\mu =-\int_{X}fb_{n}d\mu ,\qquad \forall n\in 
\mathbb{N},  \tag{6.8}
\end{equation}%
on the dense subset $C_{dec}^{1}(X)\subset C_{b}^{1}(X)$ of all functions $f$
having the following decay property: 
\begin{equation}
|f(x)|\leq C\left( 1+|x|_{B^{\ast }}+|F(x)|_{B^{\ast }}\right) ^{-1},\ x\in
X.  \tag{6.9}
\end{equation}%
with some finite $C=C(f)>0.$ Obviously, both integrals in (6.8) make sense
for such $f$. We will say that $\mu \in \mathcal{M}^{b}(X)$ are \emph{%
symmetrizing measures} in the sense that they have the given logarithmic
derivative $b.\smallskip $

The above set up in particular applies to the single loop spaces 
\begin{gather*}
X=C_{\beta }^{{}},\ H=L_{\beta }^{2}\quad \text{and\quad }B=L_{\beta }^{R},\
B^{\ast }=L_{\beta }^{R}, \\
\text{with\quad }2<R<\infty ,\ 1<R^{\prime }=R(R-1)^{-1}<R<\infty ,
\end{gather*}%
and the one-particle Euclidean measure $d\sigma _{\beta }(\upsilon )$
defined in Subsect.\thinspace 3.3.1. The (IbP)-formula (6.8) for them on the
domain $C_{dec}^{1}(X)$ with 
\begin{equation*}
b_{n}(x)=-(A_{\beta }\varphi _{n},x)_{H}-(V_{k}^{\prime }(x),\varphi
_{n})_{H},\quad x\in X,
\end{equation*}%
follows from the flow description (3.34) as a particular case of
Propositions 4.8 and 4.11 (i). 

\subsection{A priori integrability estimates for symmetrizing measures}

Supposing $\mathcal{M}^{b}(X)\neq \emptyset ,$ in this subsection we shall
derive sufficient conditions implying integrability of certain functionals
(logarithmic derivatives, polynomials, higher powers of norms etc.) w.r.t.
measures $\mu \in \mathcal{M}^{b}(X).$ It should be emphasized that to get
these a priori estimates we do not need to use the fact that $F$ is of
gradient type. Since the existence of Euclidean Gibbs measures readily
follows from the corresponding uniform integrability estimates for their
local specifications (see Subsect.\thinspace 5.2.2), here we do not discuss
at all the problem whether there exists any $\mu \in \mathcal{M}^{b}(X)$
(which is, of course, not the case for arbitrary $b$). Nevertheless, we note
that, developing in a proper way the ideas of [AKRT00, Theorem\thinspace
2.3] and [BR01, Theorem\thinspace 5.1], it is also possible to give an
abstract existence criterion for $\mu \in \mathcal{M}^{b}(X)$ based on their
approximation by finite dimensional measures $\mu _{n},$ $n\in \mathbb{N},$
on $\mathbb{R}^{n}$.

\subsubsection{The coercivity functional and its integrability properties}

In order to study the properties of $\mu \in $ $\mathcal{M}^{b}(X),$ \emph{%
conventionally} in this paper we introduce the following quantity related to
the vector field $F$ (i.e., the nonlinear part of the logarithmic derivative 
$b$ defined by (6.6)):\smallskip \newline
\textbf{Definition 6.1.} \ \emph{The functional } 
\begin{equation}
L_{H}^{F}:X\rightarrow \mathbb{R},\quad L_{H}^{F}(x):=(F(x),x)_{H}, 
\tag{6.10}
\end{equation}%
\emph{will be called coercivity functional for the vector field }$%
F:X\rightarrow B^{\ast }$\emph{\ w.r.t. the Hilbert space }$H$.\smallskip

The \emph{first }and already nontrivial step is to prove the integrability
of $L_{H}^{F}$.\medskip \newline
\textbf{Theorem 6.2 }(A Priori Integrability of the Coercivity Functional)%
\textbf{.} \ \emph{Let the following assumptions on the vector field }$%
F:X\rightarrow B^{\ast }$ \emph{hold} \emph{uniformly for all} $x\in X$ 
\emph{and} $n\in \mathbb{N}:$

\begin{description}
\item[\textbf{(}$\mathbf{F}$\textbf{)}] $\exists \mathcal{K}_{0}>0$ \emph{and%
} $\mathcal{L}_{0}\geq 0:\quad |\partial _{n}F(x)|_{B^{\ast }}\leq \mathcal{K%
}_{0}(|F(x)|_{B^{\ast }}+|x|_{B})+\mathcal{L}_{0};$
\end{description}

\noindent \emph{and respectively for its coercivity functional }$%
L_{H}^{F}:X\rightarrow \mathbb{R}$

\begin{description}
\item[\textbf{(}$\mathbf{L}_{1}$\textbf{)}] $\exists \mathcal{K}_{1}>0$ 
\emph{and} $\mathcal{L}_{1}\geq 0:$
\end{description}

$|(F(x),\varphi _{n})_{H}|+|(\partial _{n}F(x),\varphi
_{n})_{H}|+|(x,\varphi _{n})_{H}|\leq \mathcal{K}_{1}L_{H}^{F}(x)+\mathcal{L}%
_{1};$

\begin{description}
\item[\textbf{(}$\mathbf{L}_{2}$\textbf{)}] $\exists \mathcal{K}_{2}>0$ 
\emph{and} $\mathcal{L}_{2}\geq 0:\quad |(\partial _{n}F(x),x)_{H}|\leq 
\mathcal{K}_{2}L_{H}^{F}(x)+\mathcal{L}_{2}.$
\end{description}

\emph{If the above parameters satisfy the following relation for some }$%
Q\geq 1$%
\begin{equation}
\Xi _{Q-1}:=\mathcal{K}_{1}\left[ 1+(Q-1)\mathcal{K}_{2}\right]
Tr_{H}^{{}}A^{-1}<1  \tag{6.11}
\end{equation}%
\emph{(which can be achieved, keeping fixed the other parameters in (6.11),
by taking} $Tr_{H}A^{-1}$ \emph{or }$\mathcal{K}_{1}>0$\emph{\ }$\emph{small}
$\emph{\ enough),} \emph{then} 
\begin{equation}
\sup_{\mu \in \mathcal{M}^{b}(X)}\int_{X}|L_{H}^{F}(x)|^{Q}d\mu (x)\leq 
\mathcal{C}_{Q}^{\prime }<\infty .  \tag{6.12}
\end{equation}%
\smallskip

\textbf{Proof}: \ To start with, note that either ($\mathbf{L}_{1}$) or ($%
\mathbf{L}_{2}$) imply the \emph{global lower boundedness} \emph{of }$%
L_{H}^{F},$\emph{\ }i.e.,\emph{\ }that 
\begin{equation}
\inf_{x\in X}L_{H}^{F}(x)\geq -l:=-\min \{-\mathcal{L}_{1}\mathcal{K}%
_{1}^{-1},-\mathcal{L}_{2}\mathcal{K}_{2}^{-1}\}>-\infty .  \tag{6.13}
\end{equation}%
We would like to apply the (IbP)-formula (6.8) along the basis vectors $%
\varphi _{n}$ to the following family of test functions $f\in C_{dec}^{1}(X)$
indexed by $n\in \mathbb{N}$ and $0<\varepsilon \leq 1:$%
\begin{equation}
f(x):=f_{n,\varepsilon }(x):=[\tilde{L}(x)]^{Q-1}(F(x),\varphi
_{n})_{H}^{{}}Z^{-1}(x),\qquad x\in X.  \tag{6.14}
\end{equation}%
Here, for brevity, we have introduced the notation 
\begin{gather}
\varepsilon \leq \tilde{L}(x):=\tilde{L}_{\varepsilon
}(x):=L_{H}^{F}(x)+l+\varepsilon ,  \notag \\
1\leq Z(x):=Z_{\varepsilon }(x):=1+\varepsilon |F(x)|_{B^{\ast
}}^{2Q}+\varepsilon |x|_{B}^{2Q}.  \tag{6.15}
\end{gather}%
Then elementary calculations give for all $x\in X$%
\begin{multline}
\partial _{n}f_{n,\varepsilon }(x)=[\tilde{L}(x)]^{Q-1}(\partial
_{n}F(x),\varphi _{n})_{H}^{{}}Z^{-1}(x)  \notag \\
+(Q-1)[\tilde{L}(x)]^{Q-2}\left[ (\partial
_{n}F(x),x)_{H}^{{}}+(F(x),\varphi _{n})_{H}^{{}}\right] (F(x),\varphi
_{n})_{H}^{{}}Z^{-1}(x)  \notag \\
-[\tilde{L}(x)]^{Q-1}(F(x),\varphi _{n})_{H}^{{}}(Z^{-2}\partial _{n}Z)(x). 
\tag{6.16}
\end{multline}

To get the upper bound on the RHS in (6.16), let us first show that%
\begin{equation*}
\sup_{n\in N}\sup_{x\in X}\left\vert (Z^{-1}\partial _{n}Z)(x)\right\vert
<\infty .
\end{equation*}%
By the chain rule, for any $R>1$ 
\begin{equation}
\partial _{n}|F(x)|_{B^{\ast }}^{R}=\left\{ 
\begin{array}{cc}
R|F(x)|_{B^{\ast }}^{R-1}(\partial _{n}F(x),\partial |F(x)|_{B^{\ast }})_{H},
& F(x)\neq 0, \\ 
0, & F(x)=0,%
\end{array}%
\right.  \tag{6.17}
\end{equation}%
and 
\begin{equation}
\partial _{n}|x|_{B}^{R}=\left\{ 
\begin{array}{cc}
R|x|_{B^{\ast }}^{R-1}(\partial |x|_{B},\varphi _{n})_{H}, & x\neq 0, \\ 
0, & x=0,%
\end{array}%
\right.  \tag{6.18}
\end{equation}%
which can be readily estimated as 
\begin{gather}
|\partial _{n}|F(x)|_{B^{\ast }}^{2Q}|\leq 2Q|F(x)|_{B^{\ast
}}^{2Q-1}|\partial _{n}F(x)|_{B^{\ast }},  \notag \\
|\partial _{n}|x|_{B}^{2Q}|\leq 2Q|x|_{B}^{2Q-1}|\varphi _{n}|_{B},\quad
\forall x\in X.  \tag{6.19}
\end{gather}%
Herefrom, by (\textbf{F}), (6.5) and Young's inequality 
\begin{equation}
a^{R-1}b\leq \frac{R-1}{R}a^{R}+\frac{1}{R}b^{R},\text{ \ }\forall a,b\in 
\mathbb{R}_{+}\,,  \tag{6.20}
\end{equation}%
we have that uniformly for all $x\in X$%
\begin{multline}
\left\vert (Z^{-1}\partial _{n}Z)(x)\right\vert \leq 2Q\varepsilon \frac{%
|F(x)|_{B^{\ast }}^{2Q-1}|\partial _{n}F(x)|_{B^{\ast }}+\kappa
|x|_{B}^{2Q-1}}{1+\varepsilon |F(x)|_{B^{\ast }}^{2Q}+\varepsilon
|x|_{B}^{2Q}}  \notag \\
\leq \mathcal{Z}:=2Q[2\mathcal{K}_{0}+\varepsilon ^{\frac{1}{2Q}}\left( 
\mathcal{L}_{0}+\kappa \right) ].  \tag{6.21}
\end{multline}%
Then, substituting (6.21) into (6.16) and using both ($\mathbf{L}_{1}$) and (%
$\mathbf{L}_{2}$), we conclude that 
\begin{equation}
\partial _{n}f_{n,\varepsilon }(x)\leq \left\{ C_{Q,1}[\tilde{L}%
(x)]^{Q}+C_{Q,2}[\tilde{L}(x)]^{Q-1}|(F(x),\varphi _{n})_{H}|\right\}
Z^{-1}(x)  \tag{6.22}
\end{equation}%
with the constants 
\begin{equation*}
C_{Q,1}:=\mathcal{K}_{1}[1+(Q-1)\mathcal{K}_{2}],\quad
C_{Q,2}:=C_{Q,\varepsilon ,2}:=\mathcal{Z}+(Q-1)\mathcal{K}_{1}.
\end{equation*}

On the other hand, by the (IbP)-formula (6.8) 
\begin{equation}
\int_{X}[\tilde{L}(x)]^{Q-1}\left\{ \lambda _{n}(F(x),\varphi
_{n})_{H}(x,\varphi _{n})_{H}+(F(x),\varphi _{n})_{H}^{2}\right\}
Z^{-1}(x)d\mu (x)=\int_{X}\partial _{n}f_{n,\varepsilon }(x)d\mu (x), 
\tag{6.23}
\end{equation}%
and hence, combining (6.22) and (6.23), 
\begin{multline}
\int_{X}[\tilde{L}(x)]^{Q-1}(F(x),\varphi _{n})_{H}(x,\varphi
_{n})_{H}Z^{-1}(x)d\mu (x)  \notag \\
\leq \lambda _{n}^{-1}\int_{X}[\tilde{L}(x)]^{Q-1}\left\{ C_{Q,1}\tilde{L}%
(x)-\left[ (F(x),\varphi _{n})_{H}^{2}-C_{Q,2}|(F(x),\varphi _{n})_{H}|%
\right] \right\} Z^{-1}(x)d\mu   \notag \\
\leq \lambda _{n}^{-1}\int_{X}\left\{ C_{Q,1}[\tilde{L}(x)]^{Q}+\frac{1}{4}%
C_{Q,2}^{2}[\tilde{L}(x)]^{Q-1}\right\} Z^{-1}(x)\,d\mu (x).  \tag{6.24}
\end{multline}%
Now let us take the sum of the inequalities (6.24) over all $n\in \mathbb{N},
$ keeping in mind that $\sum_{n=1}^{\infty }\lambda
_{n}^{-1}=Tr_{H}^{{}}A^{-1}<\infty $ and that (due to (6.5)) uniformly for
all $K\in \mathbb{N}$%
\begin{equation*}
\left\vert \sum\nolimits_{n=1}^{K}(F(x),\varphi _{n})_{H}(x,\varphi
_{n})_{H}\right\vert =|(F(x),\mathbb{P}_{K}x)_{H}|\leq \varsigma
|(F(x)|_{B^{\ast }}^{{}}|x|_{B}^{{}}.
\end{equation*}%
Then Lebesgue's dominated convergence theorem yields that 
\begin{equation}
(1-C_{Q,1}Tr_{H}A^{-1})\int_{X}[\tilde{L}(x)]^{Q}Z^{-1}(x)\,d\mu (x)\leq (%
\frac{1}{4}C_{Q,2}^{2}+l+\varepsilon )\int_{X}\left[ \tilde{L}(x)\right]
^{Q-1}Z^{-1}(x)\,d\mu (x).  \tag{6.25}
\end{equation}%
But according to (6.11) we have that $\Xi _{Q-1}:=C_{Q,1}Tr_{H}^{{}}A^{-1}<1,
$ and thus H\"{o}lder's inequality together with (6.23) and (6.25) imply
that 
\begin{multline}
(1-\Xi _{Q-1})\left\{ \int_{X}[\tilde{L}(x)]^{Q}Z^{-1}(x)\,d\mu (x)\right\}
^{\frac{1}{Q}}  \notag \\
\leq \left\{ Q[2\mathcal{K}_{0}+\varepsilon ^{\frac{1}{2Q}}\left( \mathcal{L}%
_{0}+\kappa \right) ]+\frac{1}{2}(Q-1)\mathcal{K}_{1}\right\}
^{2}+l+\varepsilon .  \tag{6.26}
\end{multline}%
Finally, letting $\varepsilon \searrow 0$ in (6.26), from Fatou's lemma we
obtain that 
\begin{equation*}
(1-\Xi _{Q-1})\left\{ \int_{X}[L_{H}^{F}(x)+l]^{Q}d\mu (x)\right\} ^{\frac{1%
}{Q}}\leq \left[ 2Q\mathcal{K}_{0}+\frac{1}{2}(Q-1)\mathcal{L}_{1}\right]
^{2}+l,
\end{equation*}%
and thus%
\begin{equation}
\left\{ \int_{X}|L_{H}^{F}(x)|^{Q}d\mu (x)\right\} ^{\frac{1}{Q}}\leq
\left\{ \left[ 2Q\mathcal{K}_{0}+\frac{1}{2}(Q-1)\mathcal{K}_{1}\right]
^{2}+2l\right\} (1-\Xi _{Q-1})^{-1}<\infty ,  \tag{6.27}
\end{equation}%
which implies the required uniform integrability (6.12).

\noindent $\blacksquare $

\subsubsection{Integrability of partial logarithmic derivatives and
consequences}

Here we consider some direct applications of Theorem 6.2 proved above.
First, let us analyze in more detail its assumptions. Note that the partial
logarithmic derivatives (6.5) can be obviously decomposed with arbitrary $%
a\in \mathbb{R}$ as 
\begin{gather}
b_{n}(x):=-(\tilde{A}\varphi _{n},x)_{H}-(\tilde{F}(x),\varphi
_{n})_{H},\quad n\in \mathbb{N},\ x\in X,  \notag \\
\text{where\quad }\tilde{A}:=A+a^{2}\mathbf{1}\text{\quad and\quad }\tilde{F}%
(x):=F(x)-a^{2}x.  \tag{6.28}
\end{gather}%
Since $Tr_{H}^{{}}A^{-1}<\infty ,$ one gets $Tr_{H}^{{}}\tilde{A}%
^{-1}<\varepsilon $ for any given $\varepsilon >0$ by choosing a large
enough $a^{2}\geq a^{2}(\varepsilon )>0.$ On the other hand, as can easily
be seen from the definitions (6.2), (6.5) and (6.28), $\tilde{F}$ also
satisfies Assumption ($\mathbf{F}$) but with the constants $\mathcal{\tilde{K%
}}_{0}:=\mathcal{K}_{0}(1+a^{2}\iota ^{2})$ and $\mathcal{\tilde{L}}_{0}:=%
\mathcal{L}_{0}+a^{2}\kappa .$ Moreover, let us suppose (having in mind the
single loop case with Assumption ($\mathbf{V}_{\text{\textbf{iii}}}$) on the
one-particle potentials $V_{k}$, see Remark 6.1) that $L_{H}^{F}$ satisfies
additionally the following \emph{growth }condition at infinity:

\begin{description}
\item[($\mathbf{L}_{3}$)] $\exists R\geq 2,$ $\mathcal{K}_{3}>0$ and $%
\mathcal{L}_{3}\geq 0:\quad |x|_{B}^{R}\leq \mathcal{K}_{3}L_{H}^{F}(x)+%
\mathcal{L}_{3};$
\end{description}

\noindent which in turn, when $R>2,$ implies the standard \emph{coercivity}
property:

\begin{description}
\item[($\mathbf{L}_{4}$)] $\forall \mathcal{K}_{4}>0$ $\exists \mathcal{L}%
_{4}\geq 0:\quad |x|_{H}^{2}\leq \mathcal{K}_{4}L_{H}^{F}(x)+\mathcal{L}%
_{4}. $
\end{description}

\noindent From ($\mathbf{L}_{4}$) obviously for any $\delta _{1},\delta
_{2}>0$ and $0<\mathcal{K}_{4}\leq \delta _{1}(a^{2}+\delta _{2})^{-1}$%
\begin{multline}
L_{H}^{\tilde{F}}(x):=(\tilde{F}(x),x)_{H}^{{}}=(1-\delta
_{1})L_{H}^{F}(x)+(\delta _{1}L_{H}^{F}(x)-a^{2}|x|_{H}^{2})  \notag \\
\geq (1-\delta _{1})L_{H}^{F}(x)+\delta _{2}|x|_{H}^{2}-\delta _{1}\mathcal{L%
}_{4}\mathcal{K}_{4}^{-1},  \tag{6.29}
\end{multline}%
so that $L_{H}^{\tilde{F}}$ also satisfies Assumptions ($\mathbf{L}_{1}$)
and ($\mathbf{L}_{2}$) with corresponding constants 
\begin{equation*}
\mathcal{\tilde{K}}_{1}:=\mathcal{K}_{1}(1-\delta _{1})^{-1}\text{ \ and \ }%
\mathcal{\tilde{K}}_{2}:=\mathcal{K}_{2}(1-\delta _{1})^{-1}.
\end{equation*}%
Thus, for any given $Q\geq 1$, by a proper decomposition (4.28) with $%
a^{2}\geq a^{2}(Q,\mathcal{K}_{1},\mathcal{K}_{2})>0$ we can always achieve
that 
\begin{equation*}
\tilde{\Xi}_{Q-1}:=\mathcal{\tilde{K}}_{1}\left[ 1+(Q-1)\mathcal{\tilde{K}}%
_{2}\right] Tr_{H}\tilde{A}^{-1}<1.
\end{equation*}%
All together, this yields the following modification of Theorem
6.2:\smallskip \newline
\textbf{Theorem 6.2%
%TCIMACRO{\U{b4}}%
%BeginExpansion
\'{}%
%EndExpansion
} (Refinement of Theorem 6.2)\textbf{.} \ \emph{Along with }($\mathbf{F}$)%
\emph{,} ($\mathbf{L}_{1}$)\emph{\ and }($\mathbf{L}_{2}$)\emph{, let also
Assumption} ($\mathbf{L}_{4}$)\emph{\ hold. Then the coercivity functional} $%
L_{H}^{F}$\emph{\ possesses the integrability property} \emph{(4.12) for all}
$Q\geq 1.\smallskip $

The next statement is readily apparent from Theorem 6.2 and its
proof.\smallskip

\noindent \textbf{Corollary 6.3} (Integrability of Polynomials and Partial
Logarithmic Derivatives)\textbf{. \ }\emph{Let the assumptions of Theorem
6.2 hold} \emph{and additionally suppose that} 
\begin{equation*}
\Xi _{Q-1}<1\quad \emph{for\ some}\text{\emph{\ fixed} }Q\geq 1.
\end{equation*}%
\emph{Then} 
\begin{gather}
\sup_{\mu \in \mathcal{M}^{b}(X)}\sup_{n\in \mathbb{N}}\int_{X}|(x,\varphi
_{n})_{H}^{{}}|^{Q}d\mu (x)\}=:\mathcal{C}_{Q}<\infty ,  \notag \\
\sup_{\mu \in \mathcal{M}^{b}(X)}\int_{X}|b_{n}(x)|^{Q}d\mu (x)=:\mathcal{C}%
_{n,Q}<\infty ,\quad \forall n\in \mathbb{N},  \tag{6.30}
\end{gather}%
\emph{and hence the (IbP)-formula (6.8) can be extended to all functions }$%
f\in C_{b,loc}^{1}(X)$ \emph{satisfying the extra growth condition:} 
\begin{equation}
\forall n\in \mathbb{N\ \exists }C=C(f,n)>0:\ |f|+|\partial _{n}f|\leq
C(1+L_{H}^{F})^{Q-1}.  \tag{6.31}
\end{equation}

\textbf{Proof:} \ The integrability of $b_{n}^{Q}(x)$ and $(x,\varphi
_{n})_{H}^{Q}$ uniformly w.r.t. all $\mu \in \mathcal{M}^{b}(X)$ is evident
from ($\mathbf{L}_{1}$) and (6.12). In order to extend the (IbP)-formula
(6.8) to functions satisfying the growth condition (6.31), we construct a
proper approximation of any such $f$ by $f_{\varepsilon }\in C_{dec}^{1}(X),$
$0<\varepsilon \leq 1.$ Namely, define $f_{\varepsilon
}(x)=f(x)Z_{\varepsilon }^{-1}(x),$ where $Z_{\varepsilon }^{{}}$ is given
by (6.15). Using the upper bound (6.20) for $|Z_{\varepsilon }^{-1}\partial
_{n}Z_{\varepsilon }^{{}}|$, it is easy to check that 
\begin{equation*}
f_{\varepsilon }^{{}}\rightarrow f,\ \partial _{n}f_{\varepsilon
}^{{}}\rightarrow \partial _{n}f,\ \varepsilon \searrow 0,\text{ \ and \ }%
|f_{\varepsilon }^{{}}|+|\partial _{n}f_{\varepsilon }^{{}}|\leq
C_{f,n}(1+L_{H}^{F})^{Q-1}.
\end{equation*}%
Then the validity of the (6.8) for this $f$ can be shown in a standard way
by Lebesgue's dominated convergence theorem applied to both sides of the
identity $\int_{X}\partial _{n}f_{\varepsilon }^{{}}d\mu
=-\int_{X}f_{\varepsilon }^{{}}b_{n}d\mu $ when $\varepsilon \searrow 0.$

\noindent $\blacksquare \medskip $

\noindent \textbf{Remark 6.4. }It would be worth to compare our method with
a usual Lyapunov function approach to the study of invariant distributions
of stochastic differential equations (for its rigorous implementation in
infinite dimensional spaces we refer to [BR01]). Namely, keeping fixed the
assumptions of Subsect.\thinspace 6.1, let us consider the SDE 
\begin{equation}
dx_{t}=-\frac{1}{2}[Ax_{t}+F(x_{t})]dt+dw_{t},\quad t>0,  \tag{6.32}
\end{equation}%
where $w_{t}$, $t\geq 0,$ is a standard Wiener process with the identity
correlation operator in $H.$ Let a diffusion process $x_{t}$, $t\geq 0,$
taking values in a Polish space $X\subset H,$ be a solution to (6.32) and
let $(\mathbb{H},\mathcal{D}(\mathbb{H}))$ be a generator of its associated
Feller transition semigroup $\mathbb{P}_{t},$ $t\geq 0,$ in $C_{b}(X).$ Then
all invariant measures $\mu $ for $x_{t}$, $t\geq 0,$ are also
infinitesimally invariant, that is $\int_{X}\mathbb{H}u\,d\mu =0,$ $\forall
u\in \mathcal{D}(\mathbb{H})$ (cf. Subsect.\thinspace 5.3). For this reason
in [BR01] one deals directly with measures $\mu $ on $(X,\mathcal{B}(X))$
(including all those satisfying the (IbP)-formulas (6.8)) defined as
solutions of the equation%
\begin{equation}
\int \mathbb{H}u\,d\mu =-\sum_{n\in \mathbb{N}}\int (\partial
_{n}^{{}}+b_{n})\partial _{n}^{{}}u\,d\mu =0  \tag{6.33}
\end{equation}%
on certain classes of test functions $u$ from $\mathcal{D}(\mathbb{H})$. If
there exists a Lyapunov function $\phi :X\rightarrow \mathbb{R}$ such that $%
\mathbb{H}\phi \leq -g+c,$ then one readily gets the uniform estimate $\int
gd\mu \leq c$ for all distributions $\mu $ solving (6.33). The latter would
mean, in the proof of Theorem 6.2, a special choice of $(f_{n}^{{}})_{n\in 
\mathbb{N}}$ being of a gradient type, i.e., $f_{n}^{{}}=\partial
_{n}^{{}}\phi ,$ $\forall n\in \mathbb{N}$. But this is not the case for the
test functions we have used in the proof and thus our method is not directly
applicable to invariant measures. To this end, we mention that for general
nonlinear drifts $F:X\rightarrow X^{\ast }$ (except the case of $F$ having
at most linear growth, cf. [BR01, Sect.\thinspace 7]) by now there are no
analytic proofs (i.e., without using the SDE (6.34)) of the existence and a
priori estimates for the invariant measures $\mu $ as solutions of (6.33).

\subsubsection{Finiteness of polynomial (exponential) moments}

In this subsection we derive sufficient conditions for the integrability of $%
|x|_{B}^{Q},$ $\forall Q\geq 1,$ (and, moreover, $\exp \lambda |x|_{B},$ $%
\forall \lambda \in \mathbb{R}$) w.r.t. $\mu \in \mathcal{M}%
^{b}(X).\smallskip $

\textbf{Theorem 6.5 }(A Priori Polynomial Moment Estimate)\textbf{.}\emph{\
\ Let Assumptions }($\mathbf{F}$)\emph{, }($\mathbf{L}_{1}$)\emph{\ and }($%
\mathbf{L}_{2}$) \emph{of the Theorem 6.2 hold and let us suppose that} 
\begin{equation*}
\Xi _{0}:=\mathcal{K}_{1}Tr_{H}^{{}}A^{-1}<1.
\end{equation*}%
\emph{If, in addition, the asymptotic growth assumption on the coercivity
functional}

\begin{description}
\item[($\mathbf{L}_{B}$)] $\lim_{x\in X,\text{ }|x|_{B}\rightarrow \infty
}L_{H}^{F}(x)=+\infty ,$
\end{description}

\noindent \emph{holds, then the following moment estimate is satisfied for
all }$Q\geq 1$%
\begin{equation}
\sup_{\mu \in \mathcal{M}^{b}(X)}\int_{X}|x|_{B}^{Q}d\mu (x)\leq \mathcal{C}%
_{Q}^{\prime \prime }<\infty .  \tag{6.34}
\end{equation}%
\emph{Moreover, in this case the (IbP)-formula (6.8) can be extended to all }%
$f\in C_{b,loc}^{1}(X)$ \emph{satisfying the following polynomial growth
condition for every fixed }$n\in \mathbb{N}$ \emph{with corresponding }$%
Q:=Q(f,n)\geq 1$ \emph{and} $C:=C(f,n)>0$\emph{:}

\begin{equation}
|f(x)|+|\partial _{n}f(x)|\leq C(1+|x|_{B})^{Q-1},\ x\in X.  \tag{6.35}
\end{equation}

\textbf{Proof:} \ Firstly we note that $\Xi _{0}<1$ implies by Theorem 6.2 
\begin{equation}
\sup_{\mu \in \mathcal{M}^{b}(X)}\int_{X}|L_{H}^{F}(x)|d\mu (x)\leq \mathcal{%
C}_{1}^{\prime }<\infty .  \tag{6.36}
\end{equation}%
Below we will also use the notation and the relevant estimates
(6.17)--(6.27) from the proof of Theorem 6.2.

So, fix $Q>2$ and $n\in \mathbb{N},$ and consider the following family of
test functions $g\in C_{dec}^{1}(X),$ $0<\varepsilon \leq 1,$%
\begin{equation}
g(x):=g_{n,\varepsilon }(x):=|x|_{B}^{Q-1}(F(x),\varphi
_{n})_{H}^{{}}\,Z^{-1}(x),\qquad x\in X,  \tag{6.37}
\end{equation}%
with the cut-off term $Z(x):=Z_{\varepsilon }(x)$ given by (6.15). Their
partial derivatives $\partial _{n}g$ in the direction $\varphi _{n}$ can be
estimated for all $x\in X$ as follows: 
\begin{multline*}
\partial _{n}g_{n,\varepsilon }(x)\leq |x|_{B}^{Q-2}\left\{ |x|_{B}(\partial
_{n}F(x),\varphi _{n})_{H}+(Q-1)|(F(x),\varphi _{n})_{H}|\cdot \,|\varphi
_{n}|_{B}^{{}}\right\} Z^{-1} \\
+\mathcal{Z}|x|_{B}^{Q-1}|(F(x),\varphi _{n})_{H}|\,Z^{-1}.
\end{multline*}%
Thus, by ($\mathbf{L}_{1}$) and the (IbP)-formula (6.8) 
\begin{multline*}
\int_{X}|x|_{B}^{Q-1}(F(x),\varphi _{n})_{H}(x,\varphi
_{n})_{H}Z^{-1}(x)\,d\mu (x) \\
\leq \lambda _{n}^{-1}\int_{X}\left\{ |x|_{B}^{Q-2}\left[ |x|_{B}+(Q-1)%
\kappa \right] \,\left[ \mathcal{K}_{1}L_{H}^{F}(x)+\mathcal{L}_{1}\right]
\right.  \\
\left. -|x|_{B}^{Q-1}\left[ (F(x),\varphi _{n})_{H}^{2}-\mathcal{Z}%
|(F(x),\varphi _{n})_{H}|\right] \right\} Z^{-1}(x)\,d\mu (x)
\end{multline*}%
and herefrom, summing over all $n\in \mathbb{N}$, we obtain 
\begin{multline}
(1-\mathcal{K}_{1}Tr_{H}^{{}}A^{-1})%
\int_{X}|x|_{B}^{Q-1}L_{H}^{F}(x)Z^{-1}(x)\,d\mu (x)  \notag \\
\leq Tr_{H}^{{}}A^{-1}\left( \mathcal{L}_{1}+\frac{1}{4}\mathcal{Z}%
^{2}\right) \int_{X}|x|_{B}^{Q-1}Z^{-1}(x)\,d\mu (x)  \notag \\
+Tr_{H}^{{}}A^{-1}(Q-1)\kappa \int_{X}|x|_{B}^{Q-2}\left[ \mathcal{K}%
_{1}L_{H}^{F}(x)+\mathcal{L}_{1}\right] Z^{-1}(x)\,d\mu (x).  \tag{6.38}
\end{multline}%
Then an elementary application of Young's inequality in the last line of
(6.38) gives that for any fixed but small enough $\delta >0$ (more
precisely, such that $\delta <(Tr_{H}A^{-1})^{-1}-\mathcal{K}_{1}$) 
\begin{multline}
\int_{X}|x|_{B}^{Q-1}L_{H}^{F}(x)Z^{-1}(x)\,d\mu (x)  \notag \\
\leq C_{Q,1}\int_{X}|x|_{B}^{Q-1}Z^{-1}(x)\,d\mu
(x)+C_{Q,2}\int_{X}(|L_{H}^{F}(x)|+1)\,d\mu (x)  \tag{6.39}
\end{multline}%
with positive constants 
\begin{equation*}
C_{Q,1}:=\frac{\left( \mathcal{L}_{1}+\frac{1}{4}\mathcal{Z}^{2}+\delta
\right) Tr_{H}^{{}}A^{-1}}{(1-(\mathcal{K}_{1}+\delta )Tr_{H}^{{}}A^{-1})}%
\quad \text{and\quad }C_{Q,2}:=C_{Q,2}(\mathcal{K}_{1},\mathcal{L}_{1}).
\end{equation*}%
Taking into account (6.36) and choosing any $\rho >0$ from Assumption ($%
\mathbf{L}_{B}$) large enough so that 
\begin{equation*}
L_{H}^{F}(x)\geq C_{Q,1}+1\quad \text{when\quad }|x|_{B}\geq \rho ,
\end{equation*}%
and hence obviously 
\begin{equation}
|x|_{B}^{Q-1}L_{H}^{F}(x)\geq (C_{Q,1}+1)|x|_{B}^{Q-1}-(C_{Q,1}+1-l)\rho
^{Q-1},\quad \forall x\in X,  \tag{6.40}
\end{equation}%
we get that 
\begin{equation}
\sup_{0<\varepsilon \leq 1}\int_{X}|x|_{B}^{Q-1}Z_{\varepsilon
}^{-1}(x)\,d\mu (x)\leq (C_{Q,1}+1-l)\rho ^{Q-1}+C_{Q,2}(\mathcal{C}%
_{1}^{\prime }+1)<\infty .  \tag{6.41}
\end{equation}%
Finally, letting $\varepsilon \searrow 0$ in (6.41), from Fatou's lemma we
obtain the required estimate (6.34).

\noindent $\blacksquare \medskip $\newline
\textbf{Corollary 6.3%
%TCIMACRO{\U{b4}}%
%BeginExpansion
\'{}%
%EndExpansion
} (Refinement of Corollary 6.3)\textbf{.}\quad \emph{Suppose that the
partial logarithmic derivatives} $b_{n}$ \emph{(or, even} \emph{more, the
coercivity functional} $L_{H}^{F}$\emph{)}$,$ \emph{have at most polynomial
growth, i.e.,}

\begin{description}
\item[($\mathbf{P}_{B}$)] $\exists \mathcal{K}_{n}^{\prime },\mathcal{L}%
_{n}^{\prime },R_{n}^{\prime }>0$\quad \emph{(resp}. $\mathcal{K}^{\prime },%
\mathcal{L}^{\prime },R^{\prime }>0$\emph{)}$:$%
\begin{equation*}
|b_{n}(x)|\leq \mathcal{K}_{n}^{\prime }|x|_{B}^{R_{n}^{\prime }}+\mathcal{L}%
_{n}^{\prime }\text{\quad \emph{(resp.} }|L_{H}^{F}(x)|\leq \mathcal{L}%
^{\prime }|x|_{B}^{R^{\prime }}+\mathcal{L}^{\prime }\text{\emph{)}},\quad
\forall x\in X.
\end{equation*}
\end{description}

\noindent \emph{Then under the assumptions of Theorem 6.5 for all} $Q\geq 1$%
\begin{equation*}
\sup_{\mu \in \mathcal{M}^{b}(X)}\int_{X}|b_{n}(x)|^{Q}d\mu (x)<\infty \text{%
\quad \emph{(resp.} }\sup_{\mu \in \mathcal{M}^{b}(X)}%
\int_{X}|L_{H}^{F}(x)|^{Q}d\mu (x)<\infty \text{\emph{)}}.
\end{equation*}%
\smallskip

Suppose that the coercivity functional $L_{H}^{F}$ satisfies also Assumption
($\mathbf{L}_{3}$) with some fixed $R>2$ and $\mathcal{K}_{3},\mathcal{L}%
_{3}>0,$ and thus surely ($\mathbf{L}_{B}$). Let us use again the
decomposition (6.28) of the logarithmic derivative $b$ with $a^{2}>0$ large
enough. Then, as is readily seen from (6.29), ($\mathbf{L}_{1}$) and ($%
\mathbf{L}_{B}$) respectively imply 
\begin{equation*}
\tilde{\Xi}_{0}:=\mathcal{\tilde{K}}_{1}Tr_{H}^{{}}\tilde{A}^{-1}<1\text{%
\quad and\quad }\lim_{|x|_{B}^{{}}\rightarrow \infty }L_{H}^{\tilde{F}%
}(x)=\infty .
\end{equation*}%
Thus, we get the following modification of Theorem 6.5:\smallskip \newline
\textbf{Theorem 6.5%
%TCIMACRO{\U{b4}}%
%BeginExpansion
\'{}%
%EndExpansion
} (Refinement of Theorem 6.5)\textbf{.} \ \ \emph{Along with }($\mathbf{F}$)%
\emph{, }($\mathbf{L}_{1}$)\emph{\ and }($\mathbf{L}_{2}$)\emph{, let
Assumption }($\mathbf{L}_{3}$)\emph{\ hold with some} \emph{fixed }$R>2.$ 
\emph{Then the moment estimate} \emph{(6.36) is satisfied for all} $Q\geq 1.$%
\medskip \newline
\textbf{Remark 6.6.} \textbf{(i)} \ In order to get that 
\begin{equation}
\sup_{\mu \in \mathcal{M}^{b}(X)}\int_{X}|x|_{H}^{Q}d\mu (x)<\infty ,\quad
\forall Q\geq 1,  \tag{6.42}
\end{equation}%
it suffices to just replace Assumption ($\mathbf{L}_{B}$) in the formulation
of Theorem 6.5 by a weaker one, namely:

\begin{description}
\item[($\mathbf{L}_{H}$)] $\lim_{x\in X,\text{ }|x|_{H}^{{}}\rightarrow
\infty }L_{H}^{F}(x)=\infty .$
\end{description}

\textbf{(ii) \ }In fact, under the assumptions of Theorem 6.5, an a priori
exponential moment estimate holds 
\begin{equation}
\sup_{\mu \in \mathcal{M}^{b}(X)}\int_{X}\exp \lambda |x|_{B}^{{}}d\mu
(x)\leq \mathcal{C}_{\lambda }^{\prime \prime }<\infty ,\quad \forall
\lambda \in \mathbb{R}.  \tag{6.43}
\end{equation}%
The proof of estimate (6.43) is similar to that of (6.34), but with the
family of test functions $g\in C_{b}^{1}(X\backslash \{0\}),$\ indexed by $%
n\in \mathbb{N}$ and $0<\varepsilon \leq 1,$%
\begin{equation}
g(x):=g_{n,\varepsilon }(x):=\frac{(F(x),\varphi _{n})_{H}^{{}}\,\exp
\lambda |x|_{B}}{1+\varepsilon |F(x)|_{B^{\ast }}^{2}+\varepsilon \exp
2\lambda |x|_{B}},\qquad x\in X.  \tag{6.44}
\end{equation}

\subsubsection{Support properties of symmetrizing measures}

Below we study a relation between the support properties of measures $\mu
\in \mathcal{M}^{b}(X)$ and spectral properties of the operator $A$ in $H.$
For this purpose we introduce the scale of Hilbert spaces defined in terms
of the powers of $A:$ 
\begin{equation}
H^{\alpha }:=\left\{ x\in H\,\left\vert |x|_{H^{\alpha }}:=\left[
\sum\limits_{n\in \mathbb{N}}\lambda _{n}^{\alpha }(x,\varphi _{n})_{H}^{2}%
\right] ^{1/2}<\infty \right. \right\} ,\ \alpha \geq 0;\quad H^{0}:=H. 
\tag{6.45}
\end{equation}%
\medskip \newline
\textbf{Theorem 6.7 }(A Priori Moment Estimates for Sobolev Norms)\textbf{.}
\ \emph{Let the assumptions of Theorem 6.2} \emph{hold and suppose that} 
\begin{equation}
\Xi _{2Q-1}<1\quad \emph{for\ some\quad }Q\geq 1.  \tag{6.46}
\end{equation}%
\emph{If, in addition, the following assumption on the linear part }$A$ 
\emph{of the logarithmic derivative holds:}

\begin{description}
\item[($\mathbf{T}_{\protect\alpha }^{{}}$)] $Tr_{H}^{{}}A^{\alpha
-1}<\infty $\quad \emph{for some}\quad $\alpha \geq 0,$
\end{description}

\noindent \emph{then} $\mu (X\cap H_{{}}^{\alpha })=1$ \emph{and, moreover,} 
\begin{equation}
\sup_{\mu \in \mathcal{M}^{b}(X)}\int_{X}|x|_{H_{{}}^{\alpha }}^{Q+1}d\mu
\leq \mathcal{C}_{Q,\alpha }<\infty .  \tag{6.47}
\end{equation}

\textbf{Proof:} \ Let for instance $Q>2$ and let us apply the (IbP)-formula
(6.8) along basis vectors $\varphi _{n},$ $n\in \mathbb{N},$ to the
following cylinder test functions $g\in \mathcal{F}C^{1}(X)$ indexed by $%
K,n\in \mathbb{N}$, $K\geq n,$%
\begin{equation}
g(x):=g_{K,n}(x):=(x,\varphi _{n})_{H}|P_{K}x|_{H_{{}}^{\alpha
}}^{Q-1},\quad x\in X,  \tag{6.48}
\end{equation}%
(which is allowed in this case by Corollary 6.3; if $Q\leq 2$ one obviously
takes $\partial _{n}^{\pm }g_{n}$ instead of $\partial _{n}g_{n}$). Firstly,
by virtue of (6.18), 
\begin{equation*}
\exists \partial _{n}g_{K,n}(x)=\left\{ 
\begin{array}{cc}
|P_{K}x|_{H_{{}}^{\alpha }}^{Q-1}+(Q-1)|P_{K}x|_{H_{{}}^{\alpha
}}^{Q-3}(P_{K}x,\varphi _{n})_{H}(P_{K}x,\varphi _{n})_{H_{{}}^{\alpha }}, & 
P_{K}x\neq 0,\smallskip \smallskip \\ 
0, & P_{K}x=0,%
\end{array}%
\right.
\end{equation*}%
and thus 
\begin{equation}
|\partial _{n}g_{K,n}(x)|\leq Q|P_{K}x|_{H_{{}}^{\alpha }}^{Q-1},\quad
\forall x\in X.  \tag{6.49}
\end{equation}%
Hence the (IbP)-formula yields for $1\leq n\leq K$%
\begin{equation*}
\int_{X}\left[ \lambda _{n}(x,\varphi _{n})_{H}^{2}+(F(x),\varphi
_{n})_{H}(x,\varphi _{n})_{H}\right] \,|P_{K}x|_{H_{{}}^{\alpha }}^{Q-1}d\mu
(x)\leq Q\int_{X}|P_{K}x|_{H_{{}}^{\alpha }}^{Q-1}d\mu (x),
\end{equation*}%
and herefrom by Young's inequality 
\begin{multline}
\lambda _{n}\int_{X}(x,\varphi _{n})^{2}\,|P_{K}x|_{H_{{}}^{\alpha
}}^{Q-1}d\mu (x)  \notag \\
\leq \int_{X}\left[ |P_{K}x|_{H_{{}}^{\alpha
}}^{Q}+Q|P_{K}x|_{H_{{}}^{\alpha }}^{Q-1}\right] d\mu (x)+\sup_{n\in \mathbb{%
N}}\int_{X}|(F(x),\varphi _{n})(x,\varphi _{n})|^{Q}d\mu (x).  \tag{6.50}
\end{multline}%
Summing over $1\leq n\leq K$ and using ($\mathbf{T}_{\alpha }$), (6.12) and
(6.46), we conclude from (6.50) that for all $K\in \mathbb{N}$ and $0<\delta
\leq 1$%
\begin{multline}
\int_{X}|P_{K}x|_{H^{\alpha }}^{Q+1}d\mu (x)  \notag \\
\leq Tr_{H}^{{}}A^{\alpha -1}\left\{ \delta \int_{X}|P_{K}x|_{H_{{}}^{\alpha
}}^{Q+1}d\mu (x)+\int_{X}(\mathcal{K}_{1}L_{H}^{F}(x)+\mathcal{L}%
_{1})^{2Q}d\mu (x)+C(Q,\delta )\right\}  \notag \\
\leq Tr_{H}^{{}}A^{\alpha -1}\left\{ \delta \int_{X}|P_{K}x|_{H_{{}}^{\alpha
}}^{Q+1}d\mu (x)+C_{Q}(\mathcal{K}_{1},\mathcal{L}_{1},\mathcal{C}%
_{2Q}^{\prime },\delta )\right\} .  \tag{6.51}
\end{multline}%
Finally, choosing any $0<\delta <(Tr_{H}^{{}}A^{\alpha -1})^{-1}$ and then
letting $K\rightarrow \infty $ in (6.51), from Fatou's lemma we obtain the
required estimate (6.47).

\noindent $\blacksquare \medskip $

Using again the decomposition (6.28) of the logarithmic derivative $b$, we
get the corresponding modification of Theorem 6.7:\smallskip \newline
\textbf{Theorem 6.7%
%TCIMACRO{\U{b4}}%
%BeginExpansion
\'{}%
%EndExpansion
} (Refinement of Theorem 6.7)\textbf{.} \ \emph{Along with }($\mathbf{F}$)%
\emph{, }($\mathbf{L}_{1}$)\emph{\ and }($\mathbf{L}_{2}$)\emph{, let
Assumptions }($\mathbf{L}_{4}$) \emph{and }($\mathbf{T}_{\alpha }$) \emph{%
hold. Then for all }$Q\geq 1$ \emph{the uniform integrability estimate} 
\emph{(5.49) for the Sobolev norms }$|\cdot |_{H^{\alpha }}$\emph{\
holds.\smallskip }\newline

\section{Symmetrizing measures on Banach (e.g. loop) lattices}

Having in mind applications to the Euclidean Gibbs measures on the
\textquotedblright loop lattice\textquotedblright\ $\Omega :=[C_{\beta }]^{%
\mathbb{Z}^{d}}$, we further enrich the abstract setting of the previous
section by adding a lattice structure. So, our aim here will be to develop
an abstract framework for symmetrizing measures $\mu \in \mathcal{M}^{b}(%
\mathcal{X})$ on Banach lattices $\mathcal{X}:=X^{\mathbb{Z}^{d}}.$ Such
measures will be defined in Subsect.\thinspace 7.1 as the solutions to the
(IbP)-formula (7.16). The required a priori estimates on $\mu \in \mathcal{M}%
^{b}(\mathcal{X})$ will be formulated in Subsect.\thinspace 7.2 as Theorems
7.1 and 7.3. In Subsect.\thinspace 7.3 we come back to the Euclidean Gibbs
measures and lastly, on the basis of the abstract results obtained, verify
the validity of Hypotheses (H) and (H$_{loc}$) for them. Recall that these
hypotheses were crucial for the proof in Sect.\thinspace 5 of our Main
Theorems I--III describing the properties of $\mu \in \mathcal{G}_{t}.$ Let
us also emphasize that here we do not touch at all such problems as
existence and (what is an especially difficult and completely open problem)
uniqueness for $\mu \in \mathcal{M}^{b}(\mathcal{X})$, since in applications
to the Euclidean Gibbs measures we prefer more standard methods for their
investigation (cf. Subsect.\thinspace 2.6.1).

\subsection{(IbP)-formula on Banach lattices}

\subsubsection{Support spaces for symmetrizing measures}

In order to include the case of Gibbs measures on loop lattices (for their
(IbP)-description see Subsect.\thinspace 4.4), one can modify the abstract
setting of Subsect.\thinspace 6.1 as follows:\smallskip

As before, let 
\begin{equation}
X\subset B\subset H\subset B^{\ast }  \tag{7.1}
\end{equation}%
be a rigging of the Hilbert space $(H,$ $(\cdot ,\cdot )_{H})$ by a locally
convex space $X$ and by reflexive Banach spaces $(B,\ |\cdot |_{B})$ and its
dual $(B^{\ast },\ |\cdot |_{B^{\ast }})$ with the properties (6.1)--(6.4).
We only recall that in $H$ we fix an orthonormal basis $bas(H)=\{\varphi
_{n}\}_{n\in \mathbb{N}},$ indexed by any countable set, $\mathbb{N}$ say,
and consisting of the eigenvectors $\varphi _{n}$ of the self-adjoint
operator $A>0,$ i.e., $A\varphi _{n}=\lambda _{n}\varphi _{n},$ with $%
Tr_{H}A^{-1}=\sum_{n\in \mathbb{N}}\lambda _{n}^{-1}<\infty $.\smallskip

On the other hand, let us be given one more rigging 
\begin{equation}
E_{+}\subset E_{0}\subset E_{-}  \tag{7.2}
\end{equation}%
of the Hilbert space $(E_{0},$ $(\cdot ,\cdot )_{0})$ by Hilbert spaces $%
(E_{+},$ $(\cdot ,\cdot )_{+})$ and its dual $(E_{-},$ $(\cdot ,\cdot
)_{-}). $ Again, as in (7.1), all the spaces in (7.2) are separable and all
the embeddings are dense and continuous. Without loss of generality, let 
\begin{equation*}
|\cdot |_{-}\leq |\cdot |_{0}\leq |\cdot |_{+}\text{ .}
\end{equation*}%
What is important, the operators $O^{-}:E_{0}\rightarrow E_{-}$ and $%
O^{+}:E_{+}\rightarrow E_{0}$ are supposed to have finite Hilbert--Schmidt
norm 
\begin{equation}
\mathbf{||}O^{\pm }\mathbf{||}_{HS}:=||O^{-}:E_{0}\rightarrow
E_{-}||_{HS}=||O^{+}:=E_{+}\rightarrow E_{0}||_{HS}<\infty .  \tag{7.3}
\end{equation}%
The duality between $E_{+}$ and $E_{-}$ is given by the scalar product $%
(\cdot ,\cdot )_{0}$ in $E_{0}$ and by Riesz' representation theorem there
is a canonical isometry $\mathbf{I}:E_{-}\rightarrow E_{+}$ defined by 
\begin{equation}
(e,e^{\prime })_{0}=(e^{\prime },e)_{0}=(\mathbf{I}e,e^{\prime })_{+}=(e,%
\mathbf{I}^{-1}e^{\prime })_{-},\quad e\in E_{-},\ e^{\prime }\in E_{+}. 
\tag{7.4}
\end{equation}%
Let $\{e_{k}\}_{k\in \mathbb{Z}}$, $\{e_{k}^{+}\}_{k\in \mathbb{Z}}$ and $%
\{e_{k}^{-}\}_{k\in \mathbb{Z}}$ be orthonormal bases in $E_{0},$ $E_{+}$
and $E_{-}$ respectively, indexed by any countable set, with abelian group
structure, say $\mathbb{Z}$. Due to (7.3) and (7.4), they can be always
taken so that 
\begin{gather}
e_{k}\in E_{+}\ \text{and }e_{k}^{+}=\gamma _{k}^{-1/2}e_{k},\
e_{k}^{-}=\gamma _{k}^{1/2}e_{k}\text{ ,}  \notag \\
\gamma _{k}\geq 1,\text{ \ }\sum\nolimits_{k\in \mathbb{Z}}\gamma _{k}^{-1}=|%
\mathbf{|}O^{\pm }\mathbf{||}_{HS}^{2}\leq 1,  \tag{7.5}
\end{gather}%
with some weight sequence $\gamma =\{\gamma _{k}\}_{k\in \mathbb{Z}}.$ The
latter means that one has the natural isomorphism $e\leftrightarrow
\{(e,e_{k})_{0}\}_{k\in \mathbb{N}}$ between the spaces $E_{0},$ $E_{+}$, $%
E_{-}$ and the coordinate spaces $l^{2}(\mathbb{Z}),\ l^{2}(\mathbb{Z}%
;\gamma ),$ $l^{2}(\mathbb{Z};\gamma ^{-1})$ (cf. Subsect.\thinspace 3.2.2)
respectively. To make things technically easier, we impose one more specific
assumption on the rigging (7.2) 
\begin{equation}
\gamma _{k-j}=\gamma _{j-k}\leq \gamma _{k}\,\gamma _{j}\quad \forall k,j\in 
\mathbb{Z},  \tag{7.6}
\end{equation}%
which in particular implies that the shift operator 
\begin{equation}
T_{j}\{e_{k}\}_{k\in \mathbb{Z}}:=\{e_{k+j}\}_{k\in \mathbb{Z}},\quad
||T_{j}||_{l^{2}}=1,\quad j\in \mathbb{Z},  \tag{7.7}
\end{equation}%
is also bounded in all $l^{q}(\mathbb{Z};\gamma ^{-q}),$ $q\geq 1,$ and $%
||T_{j}||_{l^{q}(\gamma ^{-q})}\leq \gamma _{j}.$ This obviously yields the
inclusion 
\begin{equation}
l^{2}(\mathbb{Z};\gamma ^{-1})\subset l^{1}(\mathbb{Z};\gamma
^{-1})\subseteq l^{q}(\mathbb{Z};\{\gamma _{k-j}^{-q}\}_{k\in \mathbb{N}%
}),\quad q\geq 1,\ j\in \mathbb{Z},  \tag{7.8}
\end{equation}%
with the following relations between the corresponding norms for arbitrary
vectors $e\in l^{2}(\mathbb{Z};\gamma ^{-1})$:%
\begin{equation}
\gamma _{j}^{-1}|e|_{l^{q}(\{\gamma _{k-j}^{-q}\}_{k\in \mathbb{N}})}\leq
|e|_{l^{q}(\gamma ^{-q})}\leq |e|_{l^{1}(\gamma ^{-1})}\leq
|e|_{l^{2}(\gamma ^{-1})}|\mathbf{|}O^{\pm }\mathbf{||}_{HS}\leq
|e|_{l^{2}(\gamma ^{-1})}.  \tag{7.9}
\end{equation}

We stress once more that the riggings (7.1) resp. (7.2) include the cases of
single loop spaces resp. spaces of scalar sequences over $\mathbb{Z}^{d}$ in
Subsect.\thinspace 3.2.1 resp. Subsect.\thinspace 3.2.2 as special
cases.\smallskip

Thereafter, we define (having in mind, in particular, the \emph{spaces of
loop sequences} from Subsect.\thinspace 3.2.3):

\begin{itemize}
\item the tangent Hilbert space $(\mathcal{H},$ $<\cdot ,\cdot >_{\mathcal{H}%
}),$%
\begin{multline}
\mathcal{H}:=l^{2}(\mathbb{Z}\rightarrow H)  \notag \\
:=\left\{ x=(x_{k})_{k\in \mathbb{Z}}\in H^{\mathbb{Z}}\left\vert ||x||_{%
\mathcal{H}}:=<x,x>_{\mathcal{H}}^{1/2}:=\left[ \sum\limits_{k\in \mathbb{Z}%
}|x_{k}|_{H}^{2}\right] ^{1/2}<\infty \right. \right\}  \notag \\
\cong E_{0}\otimes H:=\left\{ x:=\sum\limits_{k\in \mathbb{Z}}e_{k}\otimes
x_{k}\left\vert \,||x||_{E_{0}\otimes H}^{{}}=:||x||_{0}<\infty \right.
\right\} ;  \tag{7.10}
\end{multline}

\item the reflexive Banach space with smooth norm $(\mathcal{B},$ $||\cdot
||_{\mathcal{B}}),$%
\begin{multline}
\mathcal{B}:=l^{2}(\mathbb{Z}\rightarrow B;\gamma ^{-1})  \notag \\
:=\left\{ x=(x_{k})_{k\in \mathbb{Z}}\in B^{\mathbb{Z}}\left\vert ||x||_{%
\mathcal{B}}:=\left[ \sum\limits_{k\in \mathbb{Z}}\gamma
_{k}^{-1}|x_{k}|_{B}^{2}\right] ^{1/2}\right. <\infty \right\}  \notag \\
\cong E_{-}\otimes B:=\left\{ x:=\sum\limits_{k\in \mathbb{Z}}e_{k}\otimes
x_{k}\left\vert \,||x||_{E_{-}\otimes B}=:||x||_{\mathcal{B}}<\infty \right.
\right\} ;  \tag{7.11}
\end{multline}

\item the locally convex Polish space $\mathcal{X}:=\mathcal{B}\cap X^{%
\mathbb{Z}},$ as a support space for measures to be considered, with the
metric 
\begin{equation}
\rho _{\mathcal{X}}(x,x^{\prime }):=\left[ \sum\limits_{k\in \mathbb{Z}%
}\gamma _{k}^{-1}\left( |x_{k}-x_{k}^{\prime }|_{B}^{2}+\frac{%
|x_{k}-x_{k}^{\prime }|_{X}^{2}}{1+|x_{k}-x_{k}^{\prime }|_{X}^{2}}\right) %
\right] ^{1/2}.  \tag{7.12}
\end{equation}
\end{itemize}

\subsubsection{Smooth functions and measures}

Again, as in Subsects. 4.2.1 and 6.1, we define the spaces $C_{b}^{1}(%
\mathcal{X};h)$ resp. $C_{b,loc}^{1}(\mathcal{X};h)$ of all functions $f:%
\mathcal{X}\rightarrow \mathbb{R}$ which are continuous and globally resp.
locally bounded together with their partial derivative $\partial _{h}f:%
\mathcal{X}\rightarrow \mathbb{R}$ along a given direction $h\in \mathcal{X}%
. $ Analogously, 
\begin{equation*}
C_{b}^{1}(\mathcal{X}):=\bigcap_{h\in \mathcal{X}}C_{b}^{1}(\mathcal{X}%
;h)\quad \text{resp.\quad }C_{b,loc}^{1}(\mathcal{X}):=\bigcap_{h\in 
\mathcal{X}}C_{b,loc}^{1}(\mathcal{X};h).
\end{equation*}%
We fix an orthonormal basis $bas(\mathcal{H}):=\{h_{i}\}_{i\in \mathbb{Z}%
^{d+1}}\mathcal{\text{ }}$in $\mathcal{H},$%
\begin{equation}
h_{i}:=\{\delta _{k,j}\varphi _{n}\}_{j\in \mathbb{N}}\cong e_{k}\otimes
\varphi _{n},\quad i=(k,n)\in \mathbb{Z\times N}.  \tag{7.13}
\end{equation}%
The corresponding subspaces of cylinder functions will be denoted by $%
\mathcal{F}C_{b}^{1}(\mathcal{X})$ resp. $\mathcal{F}C_{b,loc}^{1}(\mathcal{X%
})$. For shortness, $\partial _{i}f:=\partial _{(k,n)}f$ will denote the
derivative along the basis vector $h_{i}=h_{(k,n)}$.\smallskip

Given any $k\in \mathbb{Z},$ let 
\begin{equation*}
F_{k}^{0}:X\rightarrow B^{\ast }\quad \text{and\quad }G_{k}:\mathcal{X}%
_{-}\rightarrow B^{\ast }
\end{equation*}%
be some nonlinear mappings which are \emph{continuous and} \emph{locally
bounded} \emph{together with their partial derivatives} in all directions $%
\varphi _{n}\in bas(H)$ resp. $e_{k}\otimes \varphi _{n}\in bas(\mathcal{H}%
), $ $n\in \mathbb{N}.$ Having regard to (4.38), (4.39) and (6.6), we define
a measurable vector field $\mathcal{X}\ni x\rightarrow
b(x):=\{b_{i}(x)\}_{i\in \mathbb{Z\times N}}\in \mathbb{R}^{\mathbb{Z\times N%
}}$ (the so-called \emph{\ logarithmic gradient}) by 
\begin{gather}
b_{i}(x)=b_{(k,n)}(x):=-(A\varphi _{n},x_{k})_{H}-(F_{k}(x),\varphi
_{n})_{H},  \notag \\
F_{k}(x):=F_{k}^{0}(x_{k})+G_{k}(x),\text{ \ \ }i=(k,n)\in \mathbb{Z\times N}%
.  \tag{7.14}
\end{gather}%
Here, along with the operator $A$ as the \emph{linear} component of $b,$ it
is convenient to separate its\emph{\ nonlinear diagonal} components $%
F_{k}^{0}(x_{k})$ and resp. \emph{nondiagonal} ones $G_{k}(x).$ From
assumptions (6.3)--(6.5) on the eigenvectors $\varphi _{n},$ $n\in \mathbb{N}%
,$ of the operator $A$ it follows that 
\begin{equation}
b_{(k,n)}\in C_{b,loc}^{1}(\mathcal{X})\quad \text{and\quad }%
|b_{(k,n)}(x)|\leq \kappa (\lambda _{n}|x_{k}|_{B^{\ast
}}+|F_{k}(x)|_{B^{\ast }}).  \tag{7.15}
\end{equation}%
Let $\mathcal{M}^{b}(\mathcal{X})$ denote the family of all Borel
probability measures $\mu $ on $\mathcal{X},$ which satisfy for any $%
i=(k,n)\in \mathbb{Z\times N}$ the \emph{(IbP)-formula} 
\begin{equation}
\int_{\mathcal{X}}\partial _{(k,n)}f(x)\,d\mu (x)=-\int_{\mathcal{X}%
}f(x)b_{(k,n)}(x)d\mu (x)  \tag{7.16}
\end{equation}%
on the corresponding dense subset $C_{dec}^{1}(\mathcal{X};h_{i})\subset
C_{b}^{1}(\mathcal{X};h_{i})$ of all functions $f$ possessing the additional%
\emph{\ }decay property 
\begin{equation}
\sup_{x\in \mathcal{X}}\left\{ |f(x)|\left( 1+|x_{k}|_{B^{\ast
}}+|F_{k}(x)|_{B^{\ast }}\right) \right\} <\infty .  \tag{7.17}
\end{equation}

\subsection{A priori integrability properties}

We start with a priori integrability properties of $\mu \in \mathcal{M}^{b}(%
\mathcal{X})$ supposing that such for sure exist, i.e., $\mathcal{M}^{b}(%
\mathcal{X})\neq \emptyset .$

\subsubsection{Main theorem}

Theorem 7.1 presented below is an extension both of Theorems 6.2 and 6.5 to
the case of measures on Banach (e.g. loop) lattices. Having in mind concrete
applications to quntum lattice systems like (3.1), from the very beginning
we impose here the polynomial boundedness of the nonlinear nondiagonal terms 
$G_{k}(x)$ and strong enough coercivity properties of the nonlinear diagonal
terms $F_{k}^{0}(x_{k})$ in the presentation (7.14) for the logarithmic
derivatives $b_{(k,n)}$ (see Assumptions ($\mathbb{G}_{1}$),\ ($\mathbb{G}%
_{2}$) and resp. ($\mathbb{L}_{3}$) below).\smallskip

According to Definition 6.1, we introduce the \emph{coercivity functionals}
corresponding to the vector fields\emph{\ }$F_{k}^{0}:X\rightarrow B^{\ast
}, $ $k\in \mathbb{Z},$ w.r.t. the\emph{\ }tangent Hilbert space\emph{\ }$%
\mathcal{H}$ by 
\begin{equation}
L_{k}:=L_{H}^{F_{k}^{0}}:X\rightarrow \mathbb{R},\quad
L_{H}^{k}(x_{k}):=(F_{k}^{0}(x_{k}),x_{k})_{H},\quad x_{k}\in X.  \tag{7.18}
\end{equation}%
\smallskip \newline
\textbf{Theorem 7.1} (A Priori Moment Estimates and Integrability of
Coercivity Functionals) \ \emph{Fix some }$R\geq 1$\emph{\ and a weight
system }$\gamma =\{\gamma _{k}\}_{k\in \mathbb{Z}}$ \emph{with the
properties (7.5)--(7.9).} \emph{Let} $\mathcal{J}=\{\mathcal{J}%
_{k,j}\}_{k,j\in \mathbb{N}}$ \emph{be an infinite symmetric matrix} \emph{%
with nonnegative entries} $\mathcal{J}_{k,j}=\mathcal{J}_{j,k}\geq 0$ \emph{%
such that}

\begin{description}
\item[($\mathbb{J}$)] $|||\mathcal{J}|||:=\sup_{k\in \mathbb{Z}}\sum_{j\in 
\mathbb{Z}}\mathcal{J}_{k,j}\gamma _{k-j}^{R}<\infty .$
\end{description}

\emph{Furthermore, let the following assumptions on the vector fields }$%
F_{k}^{0}:X\rightarrow B^{\ast }\ $\emph{and}$\ G_{k}:\mathcal{X}%
_{-}\rightarrow B^{\ast }$ \emph{hold} \emph{uniformly for all} $x\in 
\mathcal{X}$ \emph{and} $i=(k,n)\in \mathbb{Z\times N}:$

\begin{description}
\item[\textbf{(}$\mathbb{F}_{0}$\textbf{)}] $\exists \mathcal{K}_{0}>0$ 
\emph{and} $\mathcal{L}_{0}\geq 0:\quad |\partial
_{n}F_{k}^{0}(x_{k})|_{B^{\ast }}\leq \mathcal{K}%
_{0}(|F_{k}^{0}(x_{k})|_{B^{\ast }}+|x_{k}|_{B}^{R})+\mathcal{L}_{0};$

\item[($\mathbb{G}_{1}$)] $\exists \mathcal{M}_{1}>0$\emph{\ and }$\mathcal{N%
}_{1}\geq 0:$

$|(G_{k}(x),\varphi _{n})_{H}|+|(\partial _{i}G_{k}(x),\varphi
_{n})_{H}|\leq \mathcal{M}_{1}\sum_{j\in \mathbb{Z}}\mathcal{J}%
_{k,j}|x_{j}|_{B}^{R}+\mathcal{N}_{1};$

\item[($\mathbb{G}_{2}$)] $\exists \mathcal{M}_{2}>0$ \emph{and} $\mathcal{N}%
_{2}\geq 0:$

$|G_{k}(x)|_{B^{\ast }}+|(G_{k}(x),x_{k})_{H}|+|(\partial
_{i}G_{k}(x),x_{k})_{H}|\leq \mathcal{M}_{2}\sum_{j\in \mathbb{Z}}\mathcal{J}%
_{k,j}|x_{j}|_{B}^{R}+\mathcal{N}_{2};$
\end{description}

\noindent \emph{and respectively for the coercivity functionals }$%
L_{k}:X\rightarrow \mathbb{R}$

\begin{description}
\item[\textbf{(}$\mathbb{L}_{1}$\textbf{)}] $\exists \mathcal{K}_{1}>0$ 
\emph{and} $\mathcal{L}_{1}\geq 0:$

$|(F_{k}^{0}(x_{k}),\varphi _{n})_{H}|+|(\partial
_{n}F_{k}^{0}(x_{k}),\varphi _{n})_{H}|\leq \mathcal{K}_{1}L_{k}(x_{k})+%
\mathcal{L}_{1};$

\item[\textbf{(}$\mathbb{L}_{2}$\textbf{)}] $\exists \mathcal{K}_{2}>0$ 
\emph{and} $\mathcal{L}_{2}\geq 0:\quad |(\partial
_{n}F_{k}^{0}(x_{k}),x_{k})_{H}|\leq \mathcal{K}_{2}L_{k}(x_{k})+\mathcal{L}%
_{2};$

\item[($\mathbb{L}_{3}$)] $\exists \mathcal{K}_{3}>0$ \emph{and} $\mathcal{L}%
_{3}\geq 0:\quad |x_{k}|_{B}^{R}\leq \mathcal{K}_{3}L_{k}(x_{k})+\mathcal{L}%
_{3}.$
\end{description}

\emph{If, in addition, the following two relations between the parameters
are satisfied:} 
\begin{gather}
\Xi _{0}:=\mathcal{K}_{1}Tr_{H}^{{}}A^{-1}<1,  \tag{7.19} \\
\Theta _{0}:=|||\mathcal{J}|||\cdot \mathcal{K}_{3}\frac{\mathcal{M}_{2}+%
\mathcal{M}_{1}Tr_{H}^{{}}A^{-1}}{1-\mathcal{K}_{1}Tr_{H}^{{}}A^{-1}}<1, 
\tag{7.20}
\end{gather}%
\emph{then} 
\begin{equation}
\sup_{\mu \in \mathcal{M}^{b}(\mathcal{X})}\sup_{k\in \mathbb{Z}}\int_{%
\mathcal{X}}|L_{k}(x_{k})|d\mu (x)\leq \mathcal{C}_{1}^{\prime }<\infty 
\tag{7.21}
\end{equation}%
\emph{and} \emph{for all }$Q\geq 1$%
\begin{gather}
\sup_{\mu \in \mathcal{M}^{b}(\mathcal{X}_{-})}\sup_{k\in \mathbb{Z}}\int_{%
\mathcal{X}_{-}}|x_{k}|_{B}^{Q}d\mu (x)\leq \mathcal{C}_{Q}^{\prime \prime
}<\infty ,  \tag{7.22} \\
\sup_{\mu \in \mathcal{M}^{b}(\mathcal{X}_{-})}\sup_{k\in \mathbb{Z}}\int_{%
\mathcal{X}_{-}}\left( \sum_{j\in \mathbb{Z}}\mathcal{J}%
_{k,j}|x_{j}|_{B}^{R}\right) ^{Q}d\mu (x)\leq \mathcal{C}_{Q}^{\prime \prime
\prime }<\infty .  \tag{7.23}
\end{gather}%
\emph{If, in addition, for some }$Q\geq 1$%
\begin{equation}
\Xi _{Q-1}:=Tr_{H}^{{}}A^{-1}\mathcal{K}_{1}\left[ 1+(Q-1)\mathcal{K}_{2}%
\right] <1,  \tag{7.24}
\end{equation}%
\emph{then} 
\begin{equation}
\sup_{\mu \in \mathcal{M}^{b}(\mathcal{X}_{-})}\sup_{k\in \mathbb{Z}}\int_{%
\mathcal{X}_{-}}|L_{k}(x_{k})|^{Q}d\mu (x)\leq \mathcal{C}_{Q}^{\prime
}<\infty .  \tag{7.25}
\end{equation}%
\smallskip

\textbf{Proof:} \ To get the estimates (7.21)--(7.23) which are uniform in $%
k\in \mathbb{Z}$, in view of (7.6)--(7.9) we can endow the space $\mathcal{X}%
:=l^{2}(\mathbb{Z}\rightarrow B;\gamma ^{-1})\cap X^{\mathbb{Z}}\subset
l^{R}(\mathbb{Z}\rightarrow B;\gamma ^{-R})$ with the family of (mutually
equivalent) norms $||\cdot ||_{R,k_{0}}\,,$ $k_{0}\in \mathbb{Z},$%
\begin{multline}
||x||_{R,k_{0}}:=||T_{k_{0}}x||_{l^{R}(\gamma ^{-R})}:=\left[
\sum\nolimits_{k\in \mathbb{Z}}\gamma _{k-k_{0}}^{-R}|x_{k}|_{B}^{R}\right]
^{1/R}  \notag \\
\leq \gamma _{j}||x||_{l^{R}(\gamma ^{-R})}\leq \gamma _{j}|\mathbf{|}O^{\pm
}\mathbf{||}_{HS}||x||_{\mathcal{B}}\leq \gamma _{j}||x||_{\mathcal{B}}. 
\tag{7.26}
\end{multline}%
An important observation is that the matrix $\mathcal{J}$ is uniformly
bounded in all $l^{1}(\mathbb{Z};\{\gamma _{k-k_{0}}^{-R}\}_{k\in \mathbb{Z}%
})$, that is 
\begin{equation}
|||\mathcal{J}|||_{k_{0}}:=\sup_{k\in \mathbb{Z}}\left\{ \sum_{j\in \mathbb{Z%
}}\mathcal{J}_{k,j}\gamma _{k-k_{0}}^{-R}\gamma _{j-k_{0}}^{R}\right\} \leq
\sup_{k\in \mathbb{Z}}\sum_{j\in \mathbb{Z}}\mathcal{J}_{k,j}\gamma
_{k-j}^{R}=|||\mathcal{J}|||<\infty .  \tag{7.27}
\end{equation}%
To prove the theorem, we perform induction on $Q$ and proceed in several
steps.\smallskip

\textbf{Step 1:} \ As already mentioned in (7.13), Assumptions ($\mathbb{L}%
_{1}$), ($\mathbb{L}_{2}$) imply the \emph{uniform lower boundedness} of $%
L_{k}$, i.e., that 
\begin{equation}
\inf_{k\in \mathbb{Z}}\inf_{x_{k}\in X}L_{k}(x_{k})\geq -l:=-\min \{\mathcal{%
L}_{1}\mathcal{K}_{1}^{-1},\mathcal{L}_{2}\mathcal{K}_{2}^{-1}\}>-\infty . 
\tag{7.28}
\end{equation}%
For the sake of convenience we introduce the following functionals on $%
\mathcal{X}:$%
\begin{gather}
1\leq \tilde{L}_{k}(x):=L_{k}(x_{k})+(G_{k}(x),x_{k})_{H}+\mathcal{M}%
_{2}\sum_{j\in \mathbb{Z}}\mathcal{J}_{k,j}|x_{j}|_{B}^{R}+\mathcal{N}%
_{2}+l+1,  \notag \\
1\leq Z_{Q,k}(x):=Z_{Q,k,\varepsilon ,\varepsilon ^{\prime
}}(x):=1+\varepsilon ^{\prime }(|F_{k}^{0}(x_{k})|_{B^{\ast
}}^{2Q}+|x_{k}|_{B}^{2QR})+\varepsilon ||x||_{-}^{2QR},\smallskip   \notag \\
Z_{Q}(x):=1+\varepsilon ||x||_{-}^{2QR},  \tag{7.29}
\end{gather}%
where $Q\geq 1,$ $k\in \mathbb{Z}$ and $0<\varepsilon ,\varepsilon ^{\prime
}\leq 1.$ Fix any basis vector $h_{i},\ i=(k,n)$ and integrate by parts
along direction $h_{i}$ the following two families of test functions on $%
\mathcal{X}$ simultaneously: 
\begin{gather}
f_{i}(x):=f_{Q,i,\varepsilon ,\varepsilon ^{\prime }}(x):=[\tilde{L}%
_{k}(x)]^{Q-1}(F_{k}(x),\varphi _{n})_{H}^{{}}\,Z_{Q,k,\varepsilon
,\varepsilon ^{\prime }}^{-1}(x),\smallskip   \tag{7.30} \\
g_{i}(x):=g_{Q^{\prime },i,k_{0},\varepsilon ,\varepsilon ^{\prime
}}(x):=||x||_{R,k_{0}}^{R(Q^{\prime }-1)}(F_{k}(x),\varphi
_{n})_{H}^{{}}\,Z_{Q^{\prime },k,\varepsilon ,\varepsilon ^{\prime
}}^{-1}(x),  \tag{7.31}
\end{gather}%
where $Q\geq 1,\ Q^{\prime }>2,$ $k_{0}\in \mathbb{Z}$ and $0<\varepsilon
,\varepsilon ^{\prime }\leq 1.$ Note that by construction $f_{i},\ g_{i}\in
C_{b}^{1}(\mathcal{X})$ satisfy the growth condition (7.17). Then elementary
calculations (analogous to (6.16)--(6.22)) give us that for all $x\in 
\mathcal{X}$%
\begin{multline}
\partial _{i}f_{i}(x)\leq \lbrack \tilde{L}_{k}]^{Q-1}(\partial
_{i}F_{k},\varphi _{n})_{H}Z_{Q,k}^{-1}+[\tilde{L}_{k}]^{Q-1}|(F_{k},\varphi
_{n})_{H}|\cdot |Z_{Q,k}^{-2}\partial _{i}Z_{Q,k}|  \notag \\
+(Q-1)[\tilde{L}_{k}]^{Q-2}|(F_{k},\varphi _{n})_{H}|Z_{Q,k}^{-1}\left[
(\partial _{i}F_{k},x_{k})_{H}^{{}}+(F_{k},\varphi _{n})_{H}^{{}}+R\mathcal{M%
}_{2}\sum\nolimits_{j\in \mathbb{Z}}\mathcal{J}_{k,j}|x_{j}|_{B}^{R-1}|%
\varphi _{n}|_{B}\right]   \tag{7.32}
\end{multline}%
and respectively 
\begin{multline}
\partial _{i}g_{i}(x)\leq ||x||_{R,k_{0}}^{R(Q^{\prime }-1)}(\partial
_{i}F_{k},\varphi _{n})_{H}Z_{Q^{\prime
},k}^{-1}+||x||_{R,k_{0}}^{R(Q^{\prime }-1)}|(F_{k},\varphi _{n})_{H}|\cdot
|Z_{Q^{\prime },k}^{-2}\partial _{i}Z_{Q^{\prime },k}|  \notag \\
+R(Q^{\prime }-1)||x||_{R,k_{0}}^{R(Q^{\prime
}-1)-1}||h_{i}||_{R,k_{0}}|(F_{k},\varphi _{n})|_{H}Z_{Q^{\prime },k}^{-1} 
\tag{7.33}
\end{multline}%
with the uniform bound%
\begin{gather}
|Z_{Q,k}^{-1}\partial _{i}Z_{Q,k}|\leq 2Q\frac{\varepsilon ^{\prime
}(|F_{k}^{0}|_{B^{\ast }}^{2Q-1}|\partial _{n}F_{k}^{0}|_{B^{\ast
}}+R|x_{k}|_{B}^{2QR-1}|\varphi _{n}|_{B})+\varepsilon |\varphi
_{n}|_{B}R||x||_{\mathcal{H}}^{2QR-1}}{1+\varepsilon ^{\prime
}(|F_{k}^{0}|_{B^{\ast }}^{2Q}+|x_{k}|_{B}^{2QR})+\varepsilon ||x||_{%
\mathcal{H}}^{2QR}}  \notag \\
\leq 2Q[2\mathcal{K}_{0}+(\varepsilon ^{\prime })^{\frac{1}{2Q}}\left( 
\mathcal{L}_{0}+\kappa \right) +\varepsilon ^{\frac{1}{2QR}}R\kappa ]=:%
\mathcal{Z}_{Q,\varepsilon ,\varepsilon ^{\prime }}^{\prime }=:\mathcal{Z}%
_{Q}^{\prime },\quad \mathcal{Z}_{Q,\varepsilon ,\varepsilon ^{\prime
}}^{\prime }\left\vert _{\varepsilon ^{\prime }=0}\right. =:\mathcal{Z}_{Q}.
\tag{7.34}
\end{gather}%
Thus, 
\begin{multline}
\partial _{i}f_{i}\leq \lbrack \tilde{L}_{k}]^{Q-2}\left[ \tilde{L}_{k}+(Q-1)%
\hat{L}_{k,2}\right] \Gamma _{k,1}Z_{Q,k}^{-1}  \notag \\
+[\tilde{L}_{k}]^{Q-2}\left\vert (F_{k},\varphi _{n})_{H}^{{}}\right\vert
\cdot \left[ \mathcal{Z}_{Q}^{\prime }\tilde{L}_{k}+(Q-1)\Gamma _{k,1}\right]
Z_{Q,k}^{-1}  \tag{7.35}
\end{multline}%
and respectively 
\begin{multline}
\partial _{i}g_{i}\leq ||x||_{R,k_{0}}^{R(Q^{\prime }-1)-1}\left[
||x||_{R,k_{0}}+R(Q^{\prime }-1)\gamma _{k-k_{0}}^{-2}\kappa \right] \Gamma
_{k,1}Z_{Q^{\prime },k}^{-1}  \notag \\
+\mathcal{Z}_{Q^{\prime }}^{\prime }||x||_{R,k_{0}}^{R(Q^{\prime }-1)}\cdot
\left\vert (F_{k},\varphi _{n})_{H}^{{}}\right\vert \,Z_{Q^{\prime },k}^{-1},
\tag{7.36}
\end{multline}%
where we denote 
\begin{align}
\Gamma _{k,1}(x)& :=\mathcal{K}_{1}L_{k}(x_{k})+\mathcal{M}_{1}\sum_{j\in 
\mathbb{Z}}\mathcal{J}_{k,j}|x_{j}|_{B}^{R}+\mathcal{L}_{1}+\mathcal{N}_{1},
\notag \\
\Gamma _{k,2}(x)& :=\mathcal{K}_{2}L_{k}(x_{k})+\mathcal{M}_{2}(1+\kappa
R)\sum_{j\in \mathbb{Z}}\mathcal{J}_{k,j}|x_{j}|_{B}^{R}+\mathcal{L}_{2}+%
\mathcal{N}_{2}.  \tag{7.37}
\end{align}%
On the other hand, by the (IbP)-formula (7.16) 
\begin{multline}
\int_{\mathcal{X}}[\tilde{L}_{k}(x)]^{Q-1}(F_{k}(x),\varphi
_{n})_{H}^{{}}(x_{k},\varphi _{n})_{H}^{{}}Z_{Q,k}^{-1}(x)d\mu (x)  \notag \\
=\lambda _{n}^{-1}\int_{\mathcal{X}}\left\{ \partial _{i}f(x)-[\tilde{L}%
_{k}(x)]^{Q-1}(F_{k}(x),\varphi _{n})_{H}^{2}Z_{Q,k}^{-1}(x)\right\} d\mu (x)
\tag{7.38}
\end{multline}%
and respectively 
\begin{multline}
\int_{\mathcal{X}}||x||_{R,k_{0}}^{R(Q^{\prime }-1)}(F_{k}(x),\varphi
_{n})_{H}^{{}}(x_{k},\varphi _{n})_{H}^{{}}Z_{Q^{\prime },k}^{-1}(x)d\mu (x)
\notag \\
=\lambda _{n}^{-1}\int_{\mathcal{X}}\left\{ \partial
_{i}g(x)-||x||_{R,k_{0}}^{R(Q^{\prime }-1)}(F_{k}(x),\varphi
_{n})_{H}^{2}Z_{Q^{\prime },k}^{-1}(x)\right\} d\mu (x).  \tag{7.39}
\end{multline}%
Taking the sum of the inequalities (7.38) and (7.39) each over $n\in \mathbb{%
N}$ and using the same arguments as in the proof of Theorem 6.2, one arrives
at the estimates%
\begin{multline}
\int_{\mathcal{X}}[\tilde{L}_{k}(x)]^{Q}Z_{Q,k}^{-1}(x)d\mu (x)  \notag \\
\leq Tr_{H}^{{}}A^{-1}\sup_{n\in \mathbb{N}}\int_{\mathcal{X}}\left\{
\partial _{i}f(x)-[\tilde{L}_{k}(x)]^{Q-1}(F_{k}(x),\varphi
_{n})_{H}^{2}Z_{Q,k}^{-1}(x)\right\} d\mu (x)  \notag \\
+\int_{\mathcal{X}}[\tilde{L}_{k}(x)]^{Q-1}\left[ \mathcal{M}%
_{2}\sum\nolimits_{j\in \mathbb{Z}}\mathcal{J}_{k,j}|x_{j}|_{B}^{R}+\mathcal{%
N}_{2}+l+1)\right] Z_{Q,k}^{-1}(x)d\mu (x)  \tag{7.40}
\end{multline}%
and respectively 
\begin{multline}
\int_{\mathcal{X}}||x||_{R,k_{0}}^{R(Q^{\prime }-1)}L_{k}(x_{k})Z_{Q^{\prime
},k}^{-1}(x)d\mu (x)  \notag \\
\leq Tr_{H}^{{}}A^{-1}\sup_{n\in \mathbb{N}}\int_{\mathcal{X}}\left\{
\partial _{i}g(x)-||x||_{R,k_{0}}^{R(Q^{\prime }-1)}(F_{k}(x),\varphi
_{n})_{H}^{2}Z_{Q^{\prime },k}^{-1}(x)\right\} d\mu (x)  \notag \\
+\int_{\mathcal{X}}||x||_{R,k_{0}}^{R(Q^{\prime }-1)}\left[ \mathcal{M}%
_{2}\sum\nolimits_{j\in \mathbb{N}}\mathcal{J}_{k,j}|x_{j}|_{B}^{R}+\mathcal{%
N}_{2}+l+1\right] Z_{Q^{\prime },k}^{-1}(x)d\mu (x).  \tag{7.41}
\end{multline}%
\smallskip 

\textbf{Step 2: }$\mathbf{Q=1}$ \ We have from (7.35) and (7.40) that more
precisely 
\begin{multline}
\int_{\mathcal{X}}\left[ L_{k}(x_{k})+(G_{k}(x),x_{k})_{H}\right]
Z_{1,k}^{-1}(x)d\mu (x)\leq Tr_{H}^{{}}A^{-1}\int_{\mathcal{X}}\Gamma
_{k,1}(x)Z_{1,k}^{-1}(x)d\mu (x)  \notag \\
+Tr_{H}^{{}}A^{-1}\sup_{n\in \mathbb{N}}\int_{\mathcal{X}}\left\{ \mathcal{Z}%
_{1}^{\prime }\left\vert (F_{k}(x),\varphi _{n})_{H}^{{}}\right\vert
-(F_{k}(x),\varphi _{n})_{H}^{2}\right\} Z_{1,k}^{-1}(x)d\mu (x)  \tag{7.42}
\end{multline}%
and furthermore due to ($\mathbb{G}_{2}$), ($\mathbb{L}_{1}$) and ($\mathbb{L%
}_{3}$) 
\begin{multline}
\int_{\mathcal{X}}|x_{k}|_{B}^{R}Z_{1,k}^{-1}(x)d\mu (x)-\mathcal{L}_{3}\leq
\int_{\mathcal{X}}\mathcal{K}_{3}L_{k}(x_{k})Z_{1,k}^{-1}(x)d\mu (x)  \notag
\\
\leq \mathcal{K}_{3}\frac{\mathcal{M}_{2}+\mathcal{M}_{1}Tr_{H}^{{}}A^{-1}}{%
1-\Xi _{0}}\int_{\mathcal{X}}\left[ \sum\nolimits_{j\in \mathbb{Z}}\mathcal{J%
}_{k,j}|x_{j}|_{B}^{R}\right] Z_{1,k}^{-1}(x)d\mu (x)  \notag \\
+\mathcal{K}_{3}\frac{\mathcal{N}_{2}+\left( \mathcal{L}_{1}+\mathcal{N}_{1}+%
\frac{1}{4}[\mathcal{Z}_{1}^{\prime }]^{2}\right) Tr_{H}^{{}}A^{-1}}{1-\Xi
_{0}}.  \tag{7.43}
\end{multline}%
Letting $\varepsilon ^{\prime }\searrow 0$ in (7.43) and then taking the sum
over $k\in \mathbb{Z}$ with the weights $\gamma _{k-k_{0}}^{-R}$, by (7.27)
and Lebesgue's dominated convergence theorem we get that 
\begin{multline}
\int_{\mathcal{X}}||x||_{R,k_{0}}^{R}(x)Z_{1}^{-1}d\mu (x)  \notag \\
\leq \left( 1-\Theta _{0}\right) ^{-1}\left( 1-\Xi _{0}\right) ^{-1}||O^{\pm
}||_{HS}^{2}\mathcal{K}_{3}\left\{ \mathcal{N}_{2}+\left( \mathcal{L}_{1}+%
\mathcal{N}_{1}+\frac{1}{4}[\mathcal{Z}_{1}]^{2}\right) Tr_{H}^{{}}A^{-1}+%
\mathcal{K}_{3}^{-1}\mathcal{L}_{3}\right\} .  \tag{7.44}
\end{multline}%
Note that in doing so we used that by (7.5) 
\begin{equation*}
\sum_{k\in \mathbb{Z}}\gamma _{k-k_{0}}^{-R}\leq \sum_{k\in \mathbb{Z}%
}\gamma _{k}^{-1}=||O^{\pm }||_{HS}^{2}
\end{equation*}%
and respectively by (7.20)%
\begin{equation}
0\leq \Xi _{0}:=\mathcal{K}_{1}Tr_{H}^{{}}A^{-1}<1,\text{ \ }0\leq \Theta
_{0}:=|||\mathcal{J}|||\cdot \mathcal{K}_{3}\frac{\mathcal{M}_{2}+\mathcal{M}%
_{1}Tr_{H}^{{}}A^{-1}}{1-\Xi _{0}}<1.  \tag{7.45}
\end{equation}%
Herefrom, letting $\varepsilon \searrow 0$, by Fatou's lemma we obtain that 
\begin{equation}
\sup_{k_{0}\in \mathbb{Z}}\int_{\mathcal{X}}||x||_{R,k_{0}}^{R}(x)d\mu
(x)\leq \mathcal{C}_{R}<\infty  \tag{7.46}
\end{equation}%
with a constant $\mathcal{C}_{R}:=\mathcal{C}_{R}(\mathcal{K}_{0},\mathcal{K}%
_{1},...,\mathcal{N}_{2})$ which equals the RHS in (7.44) with $\mathcal{Z}%
_{1}:=\mathcal{Z}_{1,\varepsilon =0}=4\mathcal{K}_{0}$. Hence due to (7.26)
and (7.27) 
\begin{gather}
\sup_{\mu \in \mathcal{M}^{b}(\mathcal{X})}\sup_{k\in \mathbb{Z}}\int_{%
\mathcal{X}}|x_{k}|_{B}^{R}d\mu (x)\leq \mathcal{C}_{R}^{\prime \prime
}:=\gamma _{0}^{R}\mathcal{C}_{R}<\infty ,  \tag{7.47} \\
\sup_{\mu \in \mathcal{M}^{b}(\mathcal{X})}\sup_{k\in \mathbb{Z}}\int_{%
\mathcal{X}}\sum\nolimits_{j\in \mathbb{Z}}\mathcal{J}_{k,j}|x_{j}|_{B}^{R}d%
\mu (x)\leq \mathcal{C}_{1}^{\prime \prime \prime }:=|||\mathcal{J}|||\cdot
\gamma _{0}^{R}\mathcal{C}_{R}<\infty .  \tag{7.48}
\end{gather}%
>From (7.43) and (7.48) by Fatou's lemma we conclude in turn that also 
\begin{equation}
\sup_{\mu \in \mathcal{M}^{b}(\mathcal{X})}\sup_{k\in \mathbb{Z}}\int_{%
\mathcal{X}}|L_{k}(x_{k})|d\mu (x)\leq \mathcal{C}_{1}^{\prime }<\infty . 
\tag{7.49}
\end{equation}%
\smallskip

\textbf{Step} \textbf{3:} $\mathbf{Q}^{\prime }\mathbf{>2}$ \ Now we return
to the general case of $Q^{\prime }>2.$ From (7.36), (7.41) and Young's
inequality it readily follows that for every $k\in \mathbb{N}$ and $0<\delta
,\varepsilon ,\varepsilon ^{\prime }\leq 1$%
\begin{multline}
\int_{\mathcal{X}}|x_{k}|_{B}^{R}||x||_{R,k_{0}}^{R(Q^{\prime
}-1)}Z_{Q^{\prime },k}^{-1}(x)d\mu (x)-\mathcal{L}_{3}\leq \mathcal{K}%
_{3}\int_{\mathcal{X}}L_{k}(x_{k})||x||_{R,k_{0}}^{R(Q^{\prime
}-1)}Z_{Q^{\prime },k}^{-1}(x)d\mu (x)  \notag \\
\leq \mathcal{K}_{3}\frac{\mathcal{M}_{2}+\mathcal{M}_{1}Tr_{H}^{{}}A^{-1}+%
\delta }{1-\Xi _{0}-\delta }\int_{\mathcal{X}}||x||_{R,k_{0}}^{R(Q^{\prime
}-1)}\sum\nolimits_{j\in \mathbb{Z}}\mathcal{J}_{k,j}|x_{j}|_{B}^{R}Z_{Q^{%
\prime },k}^{-1}(x)d\mu (x)  \notag \\
+C_{Q^{\prime }}^{\prime }\int_{\mathcal{X}}\left[ ||x||_{R,k_{0}}^{R(Q^{%
\prime }-1)}+L_{k}(x_{k})+\sum\nolimits_{j\in \mathbb{Z}}\mathcal{J}%
_{k,j}|x_{j}|_{B}^{R}+1\right] Z_{Q^{\prime },k}^{-1}(x)d\mu (x).  \tag{7.50}
\end{multline}%
with $C_{Q^{\prime }}^{\prime }:=C_{Q^{\prime }}(\mathcal{K}_{1},...,%
\mathcal{N}_{2};\delta ,\mathcal{Z}_{Q^{\prime }}^{\prime })\in (0,\infty )$
continuously depending, among the other parameters, on $Q^{\prime }\geq 2$.
Letting $\varepsilon ^{\prime }\searrow 0$ (and hence $Z_{Q^{\prime
},k}\searrow Z_{Q^{\prime }}$) in (7.50) and then summing over $k\in \mathbb{%
Z}$ with the weights $\gamma _{k-k_{0}}^{-R},$ by (7.27) and Lebesgue's
convergence theorem we get that for all $Q^{\prime }\geq 2$ 
\begin{multline}
\left[ 1-|||\mathcal{J}|||\cdot \mathcal{K}_{3}\frac{\mathcal{M}_{2}+%
\mathcal{M}_{1}Tr_{H}^{{}}A^{-1}+\delta }{1-\Xi _{0}-\delta }\right] \int_{%
\mathcal{X}}||x||_{R,k_{0}}^{RQ^{\prime }}Z_{Q^{\prime }}^{-1}(x)d\mu (x) 
\notag \\
\leq C_{Q^{\prime }}^{^{\prime \prime }}\int_{\mathcal{X}}\left[
||x||_{R,k_{0}}^{R(Q^{\prime }-1)}+1\right] Z_{Q^{\prime }}^{-1}(x)d\mu (x) 
\tag{7.51}
\end{multline}%
with some new constant $C_{Q^{\prime }}^{\prime \prime }:=C_{Q^{\prime }}(%
\mathcal{K}_{1},...,\mathcal{N}_{2};\delta ,\mathcal{Z}_{Q^{\prime
}}^{\prime })\in (0,\infty )$ (even though $||x||_{R,k_{0}}$ is not
differentiable at $x=0$ and we cannot directly apply the (IbP)-formula
(7.16) when $Q^{\prime }=2$). Note that in doing so we took into account the
estimates (7.47)--(7.49) proved above. Suppose that we already know that 
\begin{equation*}
\sup_{k_{0}\in \mathbb{Z}}\int_{\mathcal{X}}||x||_{R,k_{0}}^{R(Q^{\prime
}-1)}(x)d\mu (x)\leq \mathcal{C}_{R(Q^{\prime }-1)}<\infty
\end{equation*}%
(as it was the case for $Q^{\prime }=2$ in (7.46)). Fix $\delta >0$ small
enough so that 
\begin{equation}
\Theta _{0}<\Theta _{\delta }:=|||\mathcal{J}|||\cdot \mathcal{K}_{3}\frac{%
\mathcal{M}_{2}+\mathcal{M}_{1}Tr_{H}^{{}}A^{-1}+\delta }{1-\Xi _{0}-\delta }%
<1.  \tag{7.52}
\end{equation}%
Then by Fatou's lemma, letting $\varepsilon \searrow 0$ (and thus $%
Z_{Q^{\prime }}\searrow 1$) in (7.51), we obtain that 
\begin{equation}
\sup_{k_{0}\in \mathbb{Z}}\int_{\mathcal{X}}||x||_{R,k_{0}}^{RQ^{\prime
}}d\mu (x)\leq \mathcal{C}_{R,Q^{\prime }}<\infty  \tag{7.53}
\end{equation}%
with a proper constant $\mathcal{C}_{R,Q^{\prime }}:=\mathcal{C}%
_{R,Q^{\prime }}(\mathcal{K}_{0},\mathcal{K}_{1},...,\mathcal{N}_{2}).$ So,
by induction, the estimate (7.53) is valid for all $Q^{\prime }\geq 2.$
Since 
\begin{equation*}
\sum_{j\in \mathbb{Z}}\mathcal{J}_{k,j}|x_{j}|_{B}^{R}\leq |||\mathcal{J}%
|||\cdot ||x||_{R,k}^{R},
\end{equation*}%
(7.53) immediately implies the desired estimates (7.22) and (7.23).\smallskip

\textbf{Step 4: }$\mathbf{Q>1}$ \ A similar reasoning as in the proof of
Theorem 6.2 shows by (7.35), (7.40) and Young's inequality that for every $%
k\in \mathbb{Z}$ and $0<\delta ,\varepsilon ,\varepsilon ^{\prime }\leq 1$%
\begin{multline}
\left\{ 1-\mathcal{K}_{1}(1+(Q-1)\mathcal{K}_{2})]Tr_{H}^{{}}A^{-1}-\delta
\right\} \int_{\mathcal{X}}[\tilde{L}_{k}(x_{k})]^{Q}Z_{Q,k}^{-1}(x)d\mu (x)
\notag \\
\leq C_{Q}\int_{\mathcal{X}}[\tilde{L}_{k}(x_{k})]^{Q-1}\left(
\sum\nolimits_{j\in \mathbb{Z}}\mathcal{J}_{k,j}|x_{j}|_{B}^{R}+1\right)
Z_{Q,k}^{-1}(x)d\mu (x)  \tag{7.54}
\end{multline}%
with some constant $C_{Q}:=C_{Q}(\mathcal{K}_{1},...,\mathcal{N}_{2};\delta ,%
\mathcal{Z}_{Q,\varepsilon ,\varepsilon ^{\prime }}^{\prime })\in (0,\infty
).$ If $1<Q\leq 2$ and thus $0<[\tilde{L}_{k}(x_{k})]^{Q-2}\leq 1,$ we
continue the estimate (7.54) in a trivial way to 
\begin{multline}
\left( 1-\Xi _{Q-1}-\delta \right) \int_{\mathcal{X}}[\tilde{L}%
_{k}(x_{k})]^{Q}Z_{Q,k}^{-1}(x)d\mu (x)  \notag \\
\leq C_{Q}\left\{ \int_{\mathcal{X}}\left( \sum\nolimits_{j\in \mathbb{Z}}%
\mathcal{J}_{k,j}|x_{j}|_{B}^{R}\right) d\mu (x)+1\right\} ,  \tag{7.55}
\end{multline}%
or, if otherwise $Q>2,$ then respectively by H\"{o}lder's inequality to 
\begin{multline}
\left( 1-\Xi _{Q-1}-\delta \right) \left( \int_{\mathcal{X}}[\tilde{L}%
_{k}]^{Q}Z_{Q,k}^{-1}d\mu \right) ^{\frac{1}{Q}}  \notag \\
\leq C_{Q}\left\{ \left( \int_{\mathcal{X}}\left( \sum\nolimits_{j\in 
\mathbb{Z}}\mathcal{J}_{k,j}|x_{j}|_{B}^{R}\right) ^{Q}d\mu (x)\right) ^{%
\frac{1}{Q}}+1\right\} .  \tag{7.56}
\end{multline}%
Fix $\delta >0$ small enough so that $\Xi _{Q-1}<\Xi _{Q-1}+\delta <1.$Then
by Fatou's lemma, letting $\varepsilon ,\varepsilon ^{\prime }\searrow 0$
and hence $Z_{Q,k}\searrow 1$, we conclude from (7.23), (7.55) and (7.56)
that 
\begin{equation*}
\sup_{\mu \in \mathcal{M}^{b}(\mathcal{X})}\sup_{k\in \mathbb{Z}}\int_{%
\mathcal{X}}[\tilde{L}_{k}(x_{k})]^{Q}(x)d\mu (x)\leq \tilde{C}_{Q}<\infty .
\end{equation*}%
Thus 
\begin{equation*}
\sup_{\mu \in \mathcal{M}^{b}(\mathcal{X})}\sup_{k\in \mathbb{Z}}\int_{%
\mathcal{X}}|L_{k}(x_{k})|^{Q}d\mu (x)\leq \mathcal{C}_{Q}^{\prime }<\infty ,
\end{equation*}%
which completes the proof of Theorem 7.1. Finally we note that of course all
the constants in the above estimates (i.e., $\mathcal{C}_{Q}^{\prime },$ $%
\mathcal{C}_{Q}^{\prime \prime },$ $\mathcal{C}_{Q}^{\prime \prime \prime })$
can be calculated explicitly when needed.

\noindent $\blacksquare \smallskip $

Moreover, let us suppose that the functionals $L_{k}$ satisfy additionally
the standard \emph{coercivity }property:\smallskip

\begin{description}
\item[($\mathbb{L}_{4}$)] $\forall \mathcal{K}_{4}>0$ $\exists \mathcal{L}%
_{4}\geq 0:\quad |x_{k}|_{H}^{2}\leq \mathcal{K}_{4}L_{k}(x_{k})+\mathcal{L}%
_{4}$

\emph{uniformly for all} $x_{k}\in X$ \emph{and} $k\in \mathbb{Z}.\smallskip 
$
\end{description}

As in Subsect.\thinspace 6.2, we write the decomposition of $b$ with
arbitrary $a\in \mathbb{R}$%
\begin{gather}
b_{i}(x)=b_{(k,n)}(x):=-(\tilde{A}\varphi _{n},x_{k})_{H}-(\tilde{F}%
_{k}^{0}(x_{k}),\varphi _{n})_{H}-(G_{k}(x),\varphi _{n})_{H},  \notag \\
\text{where }\tilde{A}:=A+a^{2}\mathbf{1,}\text{\quad }\tilde{F}%
_{k}^{0}(x_{k}):=F_{k}^{0}(x_{k})-a^{2}x_{k},  \notag \\
i=(k,n)\in \mathbb{Z\times N},\ x\in \mathcal{X}.  \tag{7.57}
\end{gather}%
Then choosing $a^{2}\geq a^{2}(Q,\varepsilon )>0$ large enough, we can
always achieve that $Tr_{H}\tilde{A}^{-1}<\varepsilon $ and $\Xi
_{Q-1}<\varepsilon $ for any given $\varepsilon >0.$ Analogously, by (7.29)
for any $\delta >0$ and $0<\mathcal{K}_{4}<\delta a^{-2}$%
\begin{equation*}
\tilde{L}_{k}(x_{k}):=(\tilde{F}_{k}^{0}(x_{k}),x_{k})_{H}\geq (1-\delta
)L_{k}^{{}}(x_{k})-\delta \mathcal{L}_{4}\mathcal{K}_{4}^{-1},
\end{equation*}%
and thus all $\tilde{L}_{k}^{{}}$satisfy ($\mathbb{L}_{3}$) with the same
constant $\mathcal{\tilde{K}}_{3}=\mathcal{K}_{3}(1-\delta )^{-1}.$
Therefore, instead of $\Theta _{0}<1$ in the formulation of Theorem 7.1, it
suffices to assume that $\mathcal{K}_{3}\mathcal{M}_{2}\cdot |||\mathcal{J}%
|||<1.$ Altogether this gives the following modification of Theorem
7.1:\smallskip \newline
\textbf{Theorem 7.1}$^{\prime }$ (Refinement of Theorem 7.1)\textbf{.} \ 
\emph{Let Assumptions }($\mathbb{J}$)\emph{, }($\mathbb{F}_{0}$)\emph{, }($%
\mathbb{G}_{1,2}$)\emph{\ and }($\mathbb{L}_{1-4})$ \emph{hold. Then the
moment estimates (7.22), (7.23) and (7.25) are satisfied for all }$Q\geq 1$ 
\emph{provided} 
\begin{equation}
\Theta _{0}^{\prime }:=\mathcal{K}_{3}\mathcal{M}_{2}\cdot |||\mathcal{J}%
|||<1.  \tag{7.58}
\end{equation}

\subsubsection{Integrability of logarithmic derivatives and Sobolev norms}

The next statement gives a corresponding generalization of Corollary 6.3. It
follows immediately from the proof of our main Theorem 7.1.\smallskip 
\newline
\textbf{Corollary 7.2} (Integrability of Logarithmic Derivatives)\textbf{\
(i) \ }\emph{Suppose that under the assumptions of Theorem 7.1, additionally,%
} 
\begin{equation}
\Theta _{0}<1\quad \text{\emph{and}}\quad \Xi _{Q-1}<1\quad \text{\emph{for
some}}\emph{\quad }Q\geq 1.  \tag{7.59}
\end{equation}%
\emph{Then} 
\begin{equation}
\sup_{\mu \in \mathcal{M}^{b}(\mathcal{X})}\sup_{k\in \mathbb{Z}}\int_{%
\mathcal{X}}|b_{(k,n)}|^{Q}d\mu \leq \mathcal{C}_{Q,n}<\infty ,\quad n\in 
\mathbb{N},  \tag{7.60}
\end{equation}%
\emph{and hence the (IbP)-formula (7.16) can be extended to all }$f\in
C_{b,loc}^{1}(\mathcal{X})$ \emph{satisfying the polynomial growth
condition: }$\forall n\in \mathbb{N}$\ $\ \mathbb{\exists }C:\mathbb{=}%
C(f,k,n)>0$ \emph{and} $Q:\mathbb{=}Q(f,k,n)\geq 1$ \emph{such that} $%
\forall x\in \mathcal{X}$%
\begin{equation}
|f(x)|+|\partial _{(k,n)}f(x)|\leq C(1+|x|_{l^{R}(\mathbb{Z}\rightarrow
B;\gamma ^{-1})})^{Q}.  \tag{7.61}
\end{equation}

\textbf{(ii)} \ \emph{If, moreover, the partial logarithmic derivatives} $%
b_{(k,n)}$ \emph{(or, even, the coercivity functionals} $L_{k}$\emph{)} 
\emph{have at most polynomial growth at the infinity, i.e., }$\exists 
\mathcal{K}_{n}^{\prime },\mathcal{L}_{n}^{\prime },R_{n}^{\prime }>0$ \emph{%
(resp.} $\mathcal{K}^{\prime },\mathcal{L}^{\prime },R^{\prime }>0$\emph{)
such that }$\forall x\in \mathcal{X}$%
\begin{gather}
|b_{(k,n)}(x)|\leq \mathcal{K}_{n}^{\prime }|x_{k}|_{B}^{R_{n}^{\prime }}+%
\mathcal{L}_{n}^{\prime }+\mathcal{M}_{1}\sum_{j\in \mathbb{Z}}\mathcal{J}%
_{k,j}|x_{j}|_{B}^{R}+\mathcal{N}_{1}  \notag \\
\text{\emph{(resp.}\quad }|L_{k}(x)|\leq \mathcal{K}^{\prime
}|x_{k}|_{B}^{R^{\prime }}+\mathcal{L}^{\prime }+\mathcal{M}_{2}\sum_{j\in 
\mathbb{Z}}\mathcal{J}_{k,j}|x_{j}|_{B}^{R}+\mathcal{N}_{2}\emph{)}, 
\tag{7.62}
\end{gather}%
\emph{then (7.57) (or, even stronger,} \emph{the estimate} 
\begin{equation*}
\sup_{\mu \in \mathcal{M}^{b}(\mathcal{X})}\sup_{k\in \mathbb{Z}}\int_{%
\mathcal{X}}|L_{k}(x_{k})|^{Q}d\mu (x)<\infty \text{\emph{)}}
\end{equation*}%
\emph{holds for all} $Q\geq 1.\smallskip $\newline

The next two statements give corresponding generalization of Theorems 6.7
and 6.7$^{\prime }$.\smallskip

\noindent \textbf{Theorem 7.3 }(A Priori Moment Estimates for Sobolev Norms)
\ \emph{Suppose that under the assumptions of Theorem 7.1, additionally,} 
\begin{equation*}
\Theta _{0}<1\quad \text{\emph{and}}\quad \Xi _{2Q-1}<1\text{ \emph{for some
given }}Q\geq 1.
\end{equation*}%
\emph{Furthermore, assume that the linear part of the logarithmic derivative 
}$b$ \emph{satisfies:}

\begin{description}
\item[($\mathbb{T}_{\protect\alpha }^{{}}$)] $Tr_{H}^{{}}A^{\alpha
-1}<\infty $\quad \emph{for some}\quad $\alpha \geq 0.$
\end{description}

\noindent \emph{Then} 
\begin{equation}
\mu (\mathcal{X}\cap l^{2}(\mathbb{Z}\rightarrow H^{\alpha };\gamma ^{-1}))=1
\tag{7.63}
\end{equation}%
\emph{and, moreover,} 
\begin{equation}
\sup_{\mu \in \mathcal{M}^{b}(\mathcal{X})}\sup_{k\in \mathbb{Z}}\int_{%
\mathcal{X}_{-}}|x_{k}|_{H_{{}}^{\alpha }}^{Q+1}d\mu \leq C_{Q+1,\alpha
}<\infty .  \tag{7.64}
\end{equation}

\textbf{Proof:} This proof is completely similar to the proof of Theorem
6.5. Due to Corollary 7.2 one can apply the (IbP)-formula (7.16) in the
given direction $h_{i}\in bas(\mathcal{H}),$ $i=(k,n)\in \emph{Z}^{d+1},$ to
the following cylinder functions 
\begin{gather*}
g_{K,n}(x):=g_{K,n,k}(x_{k}):=(x_{k},\varphi
_{n})_{H}|P_{K}x_{k}|_{H_{{}}^{\alpha }}^{Q-1}, \\
|g_{K,n}(x)|\leq |P_{K}x_{k}|_{H_{{}}^{\alpha }}^{Q},\text{\quad }|\partial
_{(k,n)}g_{k,n}(x)|\leq Q|P_{K}x_{k}|_{H_{{}}^{\alpha }}^{Q-1},\quad x\in 
\mathcal{X}_{-},
\end{gather*}%
where $Q>2,$ $k\in \mathbb{Z}$ and $K\geq n.$ Thus, as already calculated in
(6.50) and (6.51), we get 
\begin{multline*}
\lambda _{n}\int_{\mathcal{X}}(x_{k},\varphi
_{n})^{2}|P_{K}x_{k}|_{H_{{}}^{\alpha }}^{Q-1}d\mu (x) \\
\leq \int_{\mathcal{X}}\left[ |P_{K}x|_{H_{{}}^{\alpha
}}^{Q}+Q|P_{K}x|_{H_{{}}^{\alpha }}^{Q-1}\right] d\mu (x)+\sup_{n\in \mathbb{%
N}}\int_{\mathcal{X}}|(F_{k}(x_{k})+G_{k}(x),\varphi _{n})_{H}(x_{k},\varphi
_{n})_{H}|^{Q}d\mu (x).
\end{multline*}%
Herefrom, summing over $1\leq n\leq K$ and using ($\mathbb{L}_{1}^{\prime }$%
), ($\mathbb{G}_{1}$), (7.23) and (7.25), we conclude that uniformly for all 
$K\in \mathbb{N}$%
\begin{multline}
(1-\delta Tr_{H}^{{}}A^{\alpha -1})\int_{\mathcal{X}}|P_{K}x_{k}|_{H_{{}}^{%
\alpha }}^{Q+1}d\mu (x)  \notag \\
\leq Tr_{H}^{{}}A^{\alpha -1}\left\{ \int_{\mathcal{X}}(\mathcal{K}%
_{1}L_{H}^{F}(x)+\mathcal{M}_{1}\sum_{j\in \mathbb{Z}}\mathcal{J}%
_{k,j}|x_{j}|_{B}^{R}+\mathcal{L}_{1}+\mathcal{N}_{1})^{2Q}d\mu
(x)+C(Q,\delta )\right\} <\infty  \tag{7.65}
\end{multline}%
as soon as $0<\delta <(Tr_{H}^{{}}A^{\alpha -1})^{-1}\leq 1.$ Finally,
letting $K\rightarrow \infty $ in (7.65), from Fatou's lemma we obtain the
desired estimates (7.63) and (7.64).

\noindent $\blacksquare \smallskip $\newline
\textbf{Theorem 7.3}$^{\prime }$\textbf{\ }(Refinement of Theorem 7.3)%
\textbf{.} \ \emph{Let Assumptions }($\mathbb{J}$)\emph{, }($\mathbb{F}_{0}$)%
\emph{, }($\mathbb{G}_{1,2}$)\emph{, }($\mathbb{L}_{1-4}$)\emph{\ }$\emph{and%
}$ $(\mathbb{T}_{\alpha }$) \emph{hold. Then the moment estimates (7.61) for
the Sobolev norms }$|\cdot |_{H^{\alpha }}$\emph{\ are satisfied for all }$%
Q\geq 1$ \emph{provided in (7.58) }$\Theta _{0}^{\prime }<1$.

\subsection{Applications to Euclidean Gibbs states: proof of Hypotheses (%
\textbf{H}) and (\textbf{H}$_{\text{\texttt{loc}}}$)}

Here we come back to the Euclidean Gibbs measures and, on the basis of the
abstract Theorems 7.1$^{\prime }$ and 7.3$^{\prime }$, verify for them
Hypothesis (\textbf{H}) and (\textbf{H}$_{loc)},$ which were left as the
crucial steps for our proof of Main Theorems I--III in Subsect.\thinspace
5.2.\smallskip

\noindent \textbf{Theorem 7.5\ }(cf. \textbf{Hypothesis} ($\mathbf{H}$): A
Priori Estimates for Tempered Gibbs Measures on Lebesgue and Sobolev Loops)%
\emph{\ }\textbf{(i)}\emph{\ \ Let Assumptions }$(\mathbf{W}),\mathbf{(J)}$ 
\emph{and} $(\mathbf{V})$ \emph{on the interaction potentials of the general
QLS (3.1) hold with some fixed, but small enough (e.g., satisfying (7.74)
below)}%
\begin{equation}
0<K_{3}<K_{3}^{0}(R,||\mathbf{J}||_{0}).  \tag{7.66}
\end{equation}%
\emph{Then} \emph{every} $\mu \in \mathcal{G}_{t}:=\mathcal{G}_{(s)t}^{R}$ 
\emph{such that} 
\begin{equation*}
\mu (\Omega _{-p}^{R})=1\quad \text{\emph{for some}\ }p=p(\mu )>d
\end{equation*}%
$\emph{is\ actually}$\emph{\ sup}$\emph{ported\ by}$ $\emph{the\ set}$%
\begin{equation*}
\bigcap\limits_{0\leq \alpha <1/2,\text{ }p>d}l^{2}(\mathbb{Z}%
^{d}\rightarrow W_{\beta }^{2,\alpha };\gamma _{-p})
\end{equation*}%
\emph{(}$\emph{where}$ $\gamma _{p}=\{\gamma _{p,k}\}_{k\in \mathbb{Z}^{d}}$ 
\emph{is the weight sequence with }$\gamma _{p,k}:=(1+|k|)^{p}$ \emph{and} $%
W_{\beta }^{2,\alpha }$ \emph{is the space of Sobolev loops with the norm }$%
||\omega _{k}||_{W_{\beta }^{2,\alpha }}^{2}:=(A_{\beta }^{\alpha }\omega
_{k},\omega _{k})_{L_{\beta }^{2}}$\emph{;} \emph{cf. Subsect.\thinspace 3.2)%
}$.$ \emph{Moreover, for all} $Q\geq 1$ \emph{and }$\alpha \in \lbrack 0,1/2)
$%
\begin{equation}
\begin{array}{cc}
(\text{\textbf{i}})\quad  & \sup\limits_{\mu \in \mathcal{G}%
_{t}}\sup\limits_{k\in \mathbb{Z}^{d}}{\displaystyle \int\limits_{\Omega }}|\omega
_{k}|_{W_{\beta }^{2,\alpha }}^{Q}\,d\mu (\omega )<\infty , \\ 
(\text{\textbf{ii}})\quad  & \sup\limits_{\mu \in \mathcal{G}%
_{t}}\sup\limits_{k\in \mathbb{Z}^{d}}{\displaystyle \int\limits_{\Omega }}|V_{k}^{\prime
}(\omega _{k})\cdot \omega _{k}|^{Q}\,d\mu (\omega )<\infty , \\ 
(\text{\textbf{iii}})\quad  & \sup\limits_{\mu \in \mathcal{G}%
_{t}}\sup\limits_{k\in \mathbb{Z}^{d}}{\displaystyle \int\limits_{\Omega
}}|F_{k}^{V,W}(\omega )|_{L_{\beta }^{1}}^{Q}\,d\mu (\omega )<\infty .%
\end{array}
\tag{7.67}
\end{equation}%
\smallskip 

\textbf{(ii) \ }\emph{For the particular QLS models I, II from Section 2
(and also for the model III satisfying additionally Assumption} (\textbf{J}$%
_{\mathbf{0}}$) (ii)\emph{) a priori estimates (7.67) always hold even for
all} $\mu \in \mathcal{G}_{t}:=\mathcal{G}_{(e)t}^{R}\supseteq \mathcal{G}%
_{(s)t}^{R}.\smallskip $

\textbf{Proof:} \textbf{(i) \ }Let $p>d$ and let us take any $\mu \in 
\mathcal{G}_{t}=\mathcal{M}_{t}^{b}$ supported by $\mathcal{X}:=\Omega
_{-p}^{R}.$ For this $\mu $ we check the validity of the assumptions of
Theorems 7.1$^{\prime }$ and 7.4$^{\prime }$.\ Recall (cf.
Subsect.\thinspace 3.2) that in this concrete set up we have the following
single loop spaces over $S_{\beta }$%
\begin{equation}
X:=C_{\beta }^{{}},\text{ \ }H:=L_{\beta }^{2},\text{ \ }B:=L_{\beta }^{R},%
\text{ \ }B^{\ast }:=L_{\beta }^{R^{\prime }}  \tag{7.68}
\end{equation}%
and the spaces of scalar sequences over $\mathbb{Z}^{d}$%
\begin{equation}
E_{0}:=l^{2}(\mathbb{Z}^{d}),\text{ \ }E_{\pm }:=l^{2}(\mathbb{Z}^{d};\gamma
_{\pm p}).  \tag{7.69}
\end{equation}%
Actually, instead of the initial norm $||\cdot ||_{\mathcal{L}%
_{-p}^{R}}:=||\cdot ||_{-p,R}$ given by (3.12), we shall at once endow $%
\mathcal{B}:=\mathcal{L}_{-p}^{R}:=l^{2}(\mathbb{Z}^{d}\rightarrow L_{\beta
}^{R};\gamma _{-p})$ with a continuous system of \emph{mutually equivalent
norms} given by 
\begin{equation}
||\omega ||_{-p,R,\sigma }:=\left[ \sum\nolimits_{k\in \mathbb{Z}%
^{d}}(1+\sigma |k|)^{-p}|\omega _{k}|_{L_{\beta }^{R}}^{2}\right]
^{1/2},\quad 0<\sigma \leq 1.  \tag{7.70}
\end{equation}%
Note that each of the corresponding weight sequences $\gamma _{p,\sigma
}:=\{(1+\sigma |k|)^{p}\}_{k\in \mathbb{Z}^{d}}$ satisfies the required
properties (7.5)--(7.9). Besides, as can be easily verified, for any fastly
decreasing (i.e., satisfying Assumption $(\mathbf{J})$ with $M=2$) matrix $%
\mathcal{J}=\{\mathcal{J}_{k,j}\}_{k,j\in \mathbb{Z}^{d}}$ we have:%
\begin{gather}
|||\mathcal{J}|||_{\sigma }:=\sup_{k\in \mathbb{Z}^{d}}\sum\nolimits_{j\in 
\mathbb{Z}^{d}}\mathcal{J}_{k,j}\gamma _{p,\sigma ,k-j}^{2R}<\infty \text{,}
\notag \\
\lim_{\sigma \rightarrow +0}|||\mathcal{J}|||_{\sigma }=|||\mathcal{J}%
|||_{0}:=\sup_{k\in \mathbb{Z}^{d}}\sum\nolimits_{j\in \mathbb{Z}^{d}}%
\mathcal{J}_{k,j}.  \tag{7.71}
\end{gather}%
Then, as described in Subsect.\thinspace 4.3, the partial logarithmic
derivatives $b_{i}=b_{(k,n)}$ of $\mu $ has the form (7.14) with the
self-adjoint linear operator (cf. Subsect.\thinspace 2.3.1) 
\begin{equation*}
A=A^{\ast }:=A_{\beta }>0\quad \text{such that \ \ }Tr_{H}^{{}}A^{\alpha
-1}<\infty ,\text{ }\forall \alpha <1/2,
\end{equation*}%
and with the smooth nonlinear components (cf. Lemma 4.6) 
\begin{gather*}
F_{k}^{0}:=F_{k}^{V}\in \bigcap_{n\in \mathbb{N}}C_{b,loc}^{1}(X\rightarrow
B^{\ast };\varphi _{n}), \\
G_{k}=F_{k}^{W}\in \bigcap_{i\in \mathbb{Z}\times \mathbb{N}}C_{b,loc}^{1}(%
\mathcal{X}\rightarrow B^{\ast };h_{i}).
\end{gather*}%
By Proposition 4.12 (i), the (IbP)-formula (7.16) holds for all $f\in
C_{b}^{1}(\mathcal{X})$ satisfying the extra decay condition (7.17).

Furthermore, note that the required Assumptions ($\mathbb{F}_{0}$) on the
vector fields $F_{k}^{0}$ and resp. ($\mathbb{L}_{1-4}$) on their coercivity
functionals $L_{k}$ have been already obtained in (4.50) and (4.51). In our
case 
\begin{equation*}
\mathcal{K}_{3}:=K_{3},\text{ }\mathcal{K}_{4}:=K_{4}
\end{equation*}%
and the constant $\mathcal{K}_{4}>0$ could be taken arbitrary small since so
is $K_{4}$ in the initial Assumption ($\mathbf{V}_{\text{\textbf{iv}}}$) on
the self-interaction potentials $V_{k}.$ Analogously, the estimate (4.52)
implies ($\mathbb{G}_{1}$) and ($\mathbb{G}_{2}$) with 
\begin{equation*}
\mathcal{M}_{2}:=3\cdot 2^{R}
\end{equation*}%
and the fastly decreasing (cf. ($\mathbf{J}$) and Lemma 3.4 (i)) matrix 
\begin{equation*}
\mathcal{J=}\{\mathcal{J}_{k,j}\}_{k,j\in \mathbb{Z}^{d}},\quad \mathcal{J}%
_{k,j}:=\left\{ 
\begin{array}{cc}
\tilde{J}_{k,j}, & k\neq j \\ 
||J||_{0}, & k=j%
\end{array}%
\right. ,\quad |||\mathcal{J}|||_{0}\leq 2||J||_{0}.
\end{equation*}%
So, we can apply Theorems 7.1$^{\prime }$ and 7.3$^{\prime }$, provided for
some fixed $0<\sigma \leq 1$%
\begin{equation}
\Theta _{\sigma }^{\prime }:=3K_{3}2^{R}|||\mathcal{J}|||_{\sigma }<1. 
\tag{7.72}
\end{equation}%
But, by the continuity property (7.71), the relation (7.72) always holds for 
$\delta \in \lbrack 0,\delta _{0})$ as soon as 
\begin{equation}
\Theta _{0}^{\prime }:=6K_{3}2^{R}||\mathbf{J}||_{0}<1.  \tag{7.73}
\end{equation}%
In turn, (7.73) as well as (3.55) (the latter is sufficient for the
well-definedness of $\mu \in \mathcal{G}_{t},$ cf. Lemma 3.6) can be
achieved by choosing small enough $K_{3}>0$ such that at least 
\begin{equation}
K_{3}^{-1}>3R2^{R}||\mathbf{J}||_{0}^{{}}.  \tag{7.74}
\end{equation}%
Taking into account the upper bounds (4.50) and (4.52) on $%
|F_{k}^{V,W}(\omega )|_{L_{\beta }^{1}}$, all this together gives us the
following estimates for all $Q\geq 1$ and $0\leq \alpha <1/2$%
\begin{gather}
\sup_{k\in \mathbb{Z}^{d}}\int_{\Omega }|\omega _{k}|_{W_{\beta }^{2,\alpha
}}^{Q}\,d\mu (\omega )\leq \mathcal{C}_{Q,\alpha }^{{}}(p),  \notag \\
\sup_{k\in \mathbb{Z}^{d}}\int_{\Omega }|V_{k}^{\prime }(\omega _{k})\omega
_{k}|^{Q}\,d\mu (\omega )\leq \mathcal{C}_{Q}^{\prime }(p),  \notag \\
\sup_{k\in \mathbb{Z}^{d}}\int_{\Omega }|F_{k}^{V,W}(\omega )|_{L_{\beta
}^{1}}^{Q}\,d\mu (\omega )\leq \mathcal{C}_{Q}^{\prime \prime }(p), 
\tag{7.75}
\end{gather}%
which are \emph{uniform for all} $\mu \in \mathcal{M}^{b}(\Omega _{-p}^{R}).$
Moreover, the first estimate in (7.75) obviously implies that any $\mu \in 
\mathcal{G}_{t}=\mathcal{M}_{t}^{b}$ is in fact supported on every $\Omega
_{-p}^{R}$ as long as $p>d.$ Hence (7.75) implies (7.67), and further by the
embedding theorem (3.4) that 
\begin{equation}
\sup_{\mu \in \mathcal{G}_{t}}\sup_{k\in \mathbb{Z}^{d}}\int_{\Omega
}|\omega _{k}|_{L_{\beta }^{R^{\prime }}}^{Q}\,d\mu (\omega )<\infty ,\quad
\forall Q,R^{\prime }\geq 1.  \tag{7.76}
\end{equation}%
\smallskip 

\textbf{(ii)} \ The proof is analogous to that of (i), but with the
following obvious modification. Since any $\mu \in \mathcal{G}_{t}:=\mathcal{%
G}_{(e)t}^{R}$ is supported by $\Omega _{-\delta }^{R},$ $\forall \delta >0,$
in the concrete set up of Theorems 7.1$^{\prime }$ and 7.3$^{\prime }$ one
should put $\mathcal{B}:=\mathcal{L}_{-\delta }^{R}:=l^{2}(\mathbb{Z}%
^{d}\rightarrow L_{\beta }^{R};\gamma _{-\delta }),$ $\mathcal{X}:=\Omega
_{-\delta }^{R}:=\Omega \cap \mathcal{L}_{-\delta }^{R}$ and respectively $%
\gamma _{\delta }:=\{\exp \delta |k|\}_{k\in \mathbb{Z}^{d}}$ (cf.
definitions (3.18)--(3.20)). Then it is easy to check up that all conditions
of these theorems are satisfied with finite $|||\mathcal{J}|||:=\sup_{k\in 
\mathbb{Z}}\sum_{j\in \mathbb{Z}}\mathcal{J}_{k,j}\gamma _{\delta
,k-j}^{R}<\infty $ and with arbitrary small $\mathcal{K}_{3}$, $\mathcal{K}%
_{4}>0.$

\noindent $\blacksquare \smallskip $\newline
\textbf{Theorem 7.6\ }(cf. \textbf{Hypothesis}\emph{\ }($\mathbf{H}_{loc}$):
Uniform Estimates for Local Gibbs Specifications on Lebesgue and Sobolev
Loops)\textbf{.}\emph{\ \ Fix any boundary condition} 
\begin{equation}
\xi \in \Omega _{t}\text{\emph{\quad }}\emph{with\quad }\sup_{k\in \mathbb{Z}%
^{d}}|\xi _{k}|_{L_{\beta }^{R}}^{{}}<\infty .  \tag{7.77}
\end{equation}%
\emph{Then under the assumptions of Theorem 7.5 the following uniform
estimates on the local specifications }$\pi _{\Lambda }^{{}}(d\omega |\xi )$ 
\emph{hold for all} $Q\geq 1$ \emph{and} $0\leq \alpha <1/2:$%
\begin{equation}
\begin{array}{cc}
(\text{i})\quad  & \sup\limits_{\Lambda \Subset \mathbb{Z}%
^{d}}\sup\limits_{k\in \Lambda }{\displaystyle \int_{\Omega }}|\omega _{k}|_{W_{\beta
}^{2,\alpha }}^{Q}\,\pi _{\Lambda }^{{}}(d\omega |\xi )\leq \mathcal{C}%
_{Q,\alpha }(\xi )<\infty , \\ 
(\text{ii})\quad  & \sup\limits_{\Lambda \Subset \mathbb{Z}%
^{d}}\sup\limits_{k\in \Lambda }{\displaystyle \int_{\Omega }}|V_{k}^{\prime }(\omega
_{k})\omega _{k}|_{L_{\beta }^{1}}^{Q}\pi _{\Lambda }^{{}}(d\omega |\xi
)\leq \mathcal{C}_{Q}^{\prime }(\xi )<\infty , \\ 
(\text{iii})\quad  & \sup\limits_{\Lambda \Subset \mathbb{Z}%
^{d}}\sup\limits_{k\in \Lambda }{\displaystyle \int_{\Omega }}|F_{k}^{V,W}(\omega
)|_{L_{\beta }^{1}}^{Q}\,\pi _{\Lambda }^{{}}(d\omega |\xi )\leq \mathcal{C}%
_{Q}^{\prime \prime }(\xi )<\infty .%
\end{array}
\tag{7.78}
\end{equation}%
\smallskip 

\textbf{Proof: \ }Actually, as follows from Theorems 6.2 and 6.5 applied to
the finite volume Gibbs distributions $\nu _{\Lambda }^{{}}(d\omega
_{\Lambda }^{{}}|\xi _{\Lambda ^{c}})$, for any $\Lambda \Subset \mathbb{Z}%
^{d}$ and $\xi \in \Omega _{t}$ the corresponding integrals in (7.67) are
finite. We keep the concrete setting (7.68)--(7.70) already used in the
proof of Theorem 7.5 for the loop lattice $\mathcal{X}=\Omega _{-p}^{R}$
with some $p>d.$ In order to get the required bounds uniformly for all $%
\Lambda \Subset \mathbb{Z}^{d},$ let us go step by step through the scheme
of the proof of the Theorems 7.1 and 7.3. Namely, setting $\varepsilon
=\varepsilon ^{\prime }=0,$ let us apply the (IbP)-formula (4.57) along
directions $h_{i},$ $i=(k,n),$ when $k\in \Lambda ,$ to the test functions $%
f_{i},g_{i}$ given by (7.30), (7.31). Note that, in doing so, all the
estimates (7.32)--(7.42) are still valid for any $k\in \Lambda ,$ $\xi \in 
\mathcal{L}_{-p}^{R}$ and $\mu =\pi _{\Lambda }(d\omega |\xi ).$ Since $%
\omega _{\Lambda ^{c}}=\xi _{\Lambda ^{c}}$ $(\pi _{\Lambda }^{{}}(d\omega
|\xi )-$a.e.$)$, taking in the both sides of (7.43) and (7.50) the weighted
sums with $(1+\sigma |k-k_{0}|)^{-pR}$over $k\in \Lambda $ and then adding
to them in the obvious way the constants $(1+\sigma |k-k_{0}|)^{-pR}|\xi
_{k}|_{L_{\beta }^{R}}^{R}$ when $k\notin \Lambda ,$ we conclude that%
\begin{equation}
\begin{tabular}{ll}
$(\text{i})\quad $ & $\int_{\Omega }|||\omega |||_{k,\sigma }^{RQ}\pi
_{\Lambda }(d\omega |\xi )\leq \mathcal{C}_{Q,\sigma }\left( 1+|||\xi
|||_{k,\sigma }^{RQ}\right) ,$ \\ 
$(\text{ii})\quad $ & $\int_{\Omega }\left( \sum\nolimits_{j\in \mathbb{Z}%
^{d}}\mathcal{J}_{k,j}|\omega _{j}|_{L_{\beta }^{R}}^{R}\right) ^{Q}\pi
_{\Lambda }(d\omega |\xi )\leq \mathcal{C}_{Q,\sigma }^{\prime }\left(
1+|||\xi |||_{k,\sigma }^{RQ}\right) ,$%
\end{tabular}
\tag{7.79}
\end{equation}%
where for the shorthand we denote%
\begin{equation}
|||\xi |||_{k,\sigma }:=\left[ \sum\nolimits_{j\in \mathbb{Z}^{d}}(1+\sigma
|j-k|)^{-pR}|\xi _{j}|_{L_{\beta }^{R}}^{R}\right] ^{1/R}\leq C_{\sigma
}(1+|k|)^{p}||\xi ||_{\mathcal{L}_{-p}^{R}}<\infty .  \tag{7.80}
\end{equation}%
Note that we have again use the same trick as in the proof of Theorem 7.5 by
fixing a small enough value of $\sigma >0$ such that $\Xi _{0},\Theta
_{\delta }^{\prime }<1.$ As soon as (7.79) is proved, from here on we can
just repeat all the subsequent arguments from the proof of Theorems 7.1, 7.3
resp. Theorems 7.1$^{\prime }$, 7.3$^{\prime }$. In the final analysis,
taking into account the upper bounds (4.50) and (4.52) on $%
|F_{k}^{V,W}(\omega )|_{L_{\beta }^{1}},$ we get that%
\begin{equation}
\begin{tabular}{ll}
$(\text{i})\quad $ & $\int_{\Omega }|\omega _{k}|_{W_{\beta }^{2,\alpha
}}^{Q}\pi _{\Lambda }^{{}}(d\omega |\xi )\leq \mathcal{C}_{Q,\alpha ,\sigma
}\left( 1+|||\xi |||_{k,\sigma }^{2RQ}\right) <\infty ,$ \\ 
$(\text{ii})\quad $ & $\int_{\Omega }|V_{k}^{\prime }(\omega _{k_{0}})\omega
_{k_{0}}|_{L_{\beta }^{1}}^{Q}\pi _{\Lambda }(d\omega |\xi )\leq \mathcal{C}%
_{Q,\sigma }^{\prime \prime }\left( 1+|||\xi |||_{k,\sigma }^{RQ}\right)
<\infty .$ \\ 
$(\text{iii})\quad $ & $\int_{\Omega }|F_{k}^{V,W}(\omega )|_{L_{\beta
}^{1}}^{Q}\pi _{\Lambda }^{{}}(d\omega |\xi )\leq \mathcal{C}_{Q,\sigma
}^{\prime \prime \prime }\left( 1+|||\xi |||_{k,\sigma }^{RQ}\right) <\infty
.$%
\end{tabular}
\tag{7.81}
\end{equation}%
Moreover, it is important that all the constants $\mathcal{C}_{Q,\sigma
},...,\mathcal{C}_{R,Q,\sigma }^{\prime \prime \prime }$ in the RHS in
(7.83) and (7.85) are universal, i.e., do not depend on $k_{0}$, $\Lambda $
and $\xi $. Assuming now that the components of the loop sequence $\xi \in 
\mathcal{L}_{-p}^{R}$ are uniformly bounded, i.e.,%
\begin{equation*}
\sup_{k\in \mathbb{Z}^{d}}|\xi _{k}|_{L_{\beta }^{R}}^{{}}<\infty \text{ \
and hence \ }\sup_{k\in \mathbb{Z}^{d}}|||\xi |||_{k,\sigma }^{Q}<\infty ,
\end{equation*}%
we get the desired estimates (7.78), which are also uniform in $k\in \Lambda 
$ and $\Lambda \Subset \mathbb{Z}^{d}$.

\noindent $\blacksquare \smallskip $

\noindent \textbf{Remark 7.7 (i) }\ Now we are in a position to complete the
proof of Corollary 5.10 and hence of our Main Theorem III. Having use of
(3.9), for the proof of statement (i) in Corollary 5.10 it suffices to show
that $||\xi ||_{\mathcal{L}_{-p}^{R}}<\infty $ implies%
\begin{equation}
\sup_{\Lambda \Subset \mathbb{Z}^{d}}\sup_{k\in \Lambda }\left\{
(1+|k|)^{-pRQ}\int_{\Omega }\left[ |\omega _{k}|_{L_{\beta
}^{2}}+|F_{k}^{V,W}(\omega )|_{L_{\beta }^{1}}\right] ^{Q}\pi _{\Lambda
}(d\omega |\xi )\right\} <\infty .  \tag{7.82}
\end{equation}%
But (7.82) readily follows from (7.79) (i), (7.80) and (7.81) (iii). On the
other hand, modifying the concrete set up for the particular QLS models
I--III as was described in the proof of Theorem 7.5 (ii) above, in much the
same way one can show that $||\xi ||_{\mathcal{L}_{-\delta }^{R}}<\infty $
implies 
\begin{equation}
\sup_{\Lambda \Subset \mathbb{Z}^{d}}\sup_{k\in \Lambda }\left\{ e^{-\delta
RQ|k|}\int_{\Omega }\left[ |\omega _{k}|_{L_{\beta
}^{2}}+|F_{k}^{V,W}(\omega )|_{L_{\beta }^{1}}\right] ^{Q}\pi _{\Lambda
}(d\omega |\xi )\right\} <\infty .  \tag{7.83}
\end{equation}%
In turn, as was pointed out in the proof of statement (ii) in Corollary
5.10, (7.87) and (3.9) yield for all $Q\geq 1,$ $\delta ^{\prime }>\delta RQ$
and $\alpha \in \lbrack 0,\frac{1}{2}-\frac{1}{Q})$ 
\begin{equation}
\sup_{\Lambda \Subset \mathbb{Z}^{d}}\left\{ \sum\nolimits_{k\in \Lambda
}e^{-\delta ^{\prime }|k|}\int_{\Omega }|\omega _{k}|_{C_{\beta }^{\alpha
}}^{Q}\pi _{\Lambda }(d\omega |\xi )\right\} <\infty .  \tag{7.84}
\end{equation}%
Taking in (7.84) arbitrary small $\delta >0$ resp. large $Q>1$, we get the
desired estimate (2.20) in Main Theorem III.\smallskip 

\textbf{(ii) }\ Similarly\textbf{\ }one also can verify another Hypothesis 
\textbf{(H}$_{per}$\textbf{)} from Subsect. 5.2.2, which (according to
Remark 5.11 (i), (ii)) provides us an alternative way to show the existence
of $\mu \in \mathcal{G}_{t}.$

\section{Appendix: \ Euclidean approach and reconstruction theorem}

According to the standard algebraic approach (cf. [BrRo81]), equilibrium
states in quantum statistical mechanics are defined as those normal states
over a $C^{\ast }$-algebra $\mathcal{A}$ of quasi-local observables which
satisfy the so-called KMS (Kubo--Martin--Schwinger) condition w.r.t. to the
one-parameter group $\alpha _{t},\ t\in \mathbb{R},$ of time evolution
automorphisms on $\mathcal{A}$. This approach is especially applicable for
spin models with finite dimensional physical Hilbert spaces corresponding to
every single particle. But, unfortunately, the quantum lattice systems
considered in this paper do not fit principally into the framework of such
algebraic approach, since there is no infinite volume dynamics for such
models. Thus we use the following scheme for constructing their equilibrium
states:\smallskip

We begin with finite volumes $\Lambda \Subset \mathbb{Z}^{d}$ and local
Hamiltonians 
\begin{multline}
\mathbb{H}_{\Lambda }:=-\frac{1}{2\mathfrak{m}}\sum_{k\in \Lambda }\frac{%
d^{2}}{dx_{k}^{2}}+\frac{a^{2}}{2}\sum_{k\in \Lambda }x_{k}^{2}  \notag \\
+\sum\limits_{k\in \Lambda
}V_{k}(x_{k})+\sum_{M=2}^{N}\sum\limits_{\{k_{1},...,k_{M}\}\subset \Lambda
}W_{\{k_{1},...,k_{M}\}}(x_{k_{1}},...,x_{k_{M}})  \tag{8.1}
\end{multline}%
which are well defined as self-adjoint operators in $\mathcal{H}_{\Lambda
}:=L_{2}(\mathbb{R}^{\Lambda },\times _{k\in \Lambda }dx_{k})$. Then on the
corresponding algebras of bounded linear operators $\mathcal{A}_{\Lambda }=%
\mathcal{L}(\mathcal{H}_{\Lambda })$ we have the local automorphism groups\ $%
\alpha _{t,\Lambda }(A)=e^{itH_{\Lambda }}Ae^{-itH_{\Lambda }}$ (the so
called Heisenberg dynamics), and the local Gibbs states at fixed inverse
temperature $\beta >0$%
\begin{equation}
G_{\beta ,\Lambda }(A):=Tr(Ae^{-\beta H_{\Lambda }})/Tr(e^{-\beta H_{\Lambda
}}),\quad A\in \mathcal{A}_{\Lambda }.  \tag{8.2}
\end{equation}%
Starting from these finite volume objects, one would like to construct the
infinite volume automorphism group $\alpha _{t}(A)=\lim_{\Lambda \nearrow 
\mathbb{Z}^{d}}\alpha _{t,\Lambda }(A)$ and the limit Gibbs states 
\begin{equation}
G_{\beta ,\Lambda }(A)=\lim_{\Lambda \nearrow \mathbb{Z}^{d}}G_{\beta
,\Lambda }(A),\quad A\in \mathcal{A}_{loc}:=\bigcup_{\Lambda \Subset \mathbb{%
Z}^{d}}\mathcal{A}_{\Lambda }.  \tag{8.3}
\end{equation}%
However, it should be strongly emphasized that for the systems under
consideration it is impossible to control the thermodynamic limit (8.3)
using the operator technique alone.\smallskip

In order to overcome this principally difficulty, in [AH-K75] an approach to
the study of Gibbs states for quantum lattice systems, using Euclidean (%
\emph{rigorously defined}) path space integrals, was initiated. Conceptually
this approach is analogous to the well-known Euclidean strategy in quantum
field theory (see e.g., [Si74, Fr\"{o}77, GJ81]). More precisely, the
Euclidean method transforms the problem of giving a proper meaning to a
quantum Gibbs state $G_{\beta }^{{}}$ of the lattice system (8.1) at inverse
temperature $\beta >0$ into the problem of constructing some (classical)
Gibbs measure $\mu $ on the \emph{temperature loop lattice} $\Omega
=[C_{\beta }]^{^{\mathbb{Z}^{d}}}$ (already defined by (2.5)). Due to this
fact, various probabilistic techniques become available for the description
of equilibrium properties of quantum infinite-particle systems. Here we only
briefly outline these deep relations between quantum statistical mechanics
and stochastic processes according to the initial paper [AH-K75] (see also
[AKKR02] for the extended and up-to-date presentation):\smallskip

For any finite set of multiplication operators $\{A_{0},...,A_{n}\}\subset
L^{\infty }(\mathbb{R}^{\Lambda }),$ the spectral properties of the local
Hamiltonians $\mathbb{H}_{\Lambda }$ enable us to define the so called \emph{%
Matsubara }(or \emph{temperature} \emph{Euclidean} \emph{Green}) \emph{%
functions} 
\begin{gather}
\Gamma _{A_{0},...,A_{n}}^{\beta ,\Lambda }(\tau _{0},...,\tau
_{n}):=Tr(e^{-(\beta -(\tau _{n}-\tau _{0}))H_{\Lambda }}A_{n}...e^{-(\tau
_{1}-\tau _{0})H_{\Lambda }}A_{0})/Tr(e^{-\beta H_{\Lambda }}),  \notag \\
0\leq \tau _{0}\leq ...\leq \tau _{n}\leq \beta ,  \tag{8.4}
\end{gather}%
which have analytic continuations to the complex domain$\,0<\func{Re}%
z_{0}<...<\func{Re}z_{n}<\beta $ with the uniform bound $|\Gamma
_{A_{0},...,A_{n}}^{\beta ,\Lambda }(z_{0},...,z_{n})|\leq
||A_{0}||_{L^{\infty }}...||A_{n}||_{L^{\infty }}.$ Since on the boundary 
\begin{equation}
\Gamma _{A_{0},...,A_{n}}^{\beta ,\Lambda }(-i\tau _{0},...,-i\tau
_{n})=G_{\beta ,\Lambda }(\alpha _{\tau _{0},\Lambda }(A_{0})...\alpha
_{\tau _{n},\Lambda }(A_{n})),\quad \tau _{0},...,\tau _{n}\in \mathbb{R}, 
\tag{8.5}
\end{equation}%
and since (by the H\o egh-Krohn theorem [AH-K75]) the algebra spanned by the
operators $\alpha _{\tau ,\Lambda }(A)$ is a \emph{state} \emph{detemining
set}, the Matsubara functions (8.4) fully determine the Gibbs state $%
G_{\beta ,\Lambda }.$ \smallskip 

A \emph{crucial} \emph{observation} of the Euclidean method is that Green
functions (8.4) may be obtained as the moments of certain probability
measures. More precisely, let $\gamma _{\beta }^{{}}$ be a Gaussian measure
on $C_{\beta }$ with correlation operator $A_{\beta }^{-1}$ (cf. (3.32)).
Then by the Feynman--Kac formula 
\begin{equation}
\Gamma _{A_{0},...,A_{n}}^{\beta ,\Lambda }(\tau _{0},...,\tau
_{n})=\int_{C_{\beta }^{\Lambda }}A_{0}(\omega _{\Lambda }(\tau
_{0}))...A_{n}(\omega _{\Lambda }(\tau _{n}))d\mu _{\Lambda }(\omega
_{\Lambda }),  \tag{8.6}
\end{equation}%
where the probability measure $\mu _{\Lambda }$ is defined by 
\begin{equation}
d\mu _{\Lambda }(\omega _{\Lambda }):=\frac{1}{Z_{\Lambda }}\exp \left\{
-I_{\Lambda }(\omega _{\Lambda })\right\} \times _{k\in \Lambda }d\gamma
_{\beta }(\omega _{k}),  \tag{8.7}
\end{equation}%
with $\omega _{\Lambda }:=(\omega _{k})_{k\in \Lambda }\in C_{\beta
}^{\Lambda }$ and 
\begin{equation*}
I_{\Lambda }(\omega _{\Lambda }):=\sum\limits_{k\in \Lambda }\int_{S_{\beta
}}V_{k}(\omega _{k_{1}}(\tau ))d\tau
+\sum\limits_{M=2}^{N}\sum\limits_{\{k_{1},...,k_{M}\}\subset \Lambda
}\int_{S_{\beta }}W_{\{k_{1},...,k_{M}\}}(\omega _{k_{1}}(\tau ),...,\omega
_{k_{m}}(\tau ))d\tau .
\end{equation*}%
Suppose that a sequence $\{\mu _{\Lambda ^{(K)}}\}_{K\in \mathbb{N}}$
converges (in the local weak sense) to some measure $\mu _{\infty }$ on $%
\Omega $. As follows from (6.6), this implies the existence of the \emph{%
limit} \emph{temperature Green functions} 
\begin{gather}
\Gamma _{A_{0},...,A_{n}}^{\beta }(\tau _{0},...,\tau _{n})=\lim_{\Lambda
^{(K)}\nearrow \mathbb{Z}^{d}}\Gamma _{A_{0},...,A_{n}}^{\beta ,\Lambda
^{(K)}}(\tau _{0},...,\tau _{n}),  \notag \\
\forall A_{0},...,A_{n}\subset L_{\infty }(\mathbb{R}^{\Lambda }),\quad
\Lambda ^{(K)}\Subset \mathbb{Z}^{d},  \tag{6.8}
\end{gather}%
which satisfy the desired properties, such as \emph{analyticity }and the so
called \emph{reflection positivity on semicircle}. As was analyzed on the
axiomatic level in [AH-K75, KL81, BF02], there exists a correspondence
(analogously to the Osterwalder-Schrader \emph{reconstruction theorem} in
Euclidean field theory [Sim74, GJ81]), between equilibrium states,
temperature Green functions and stochastic processes. Namely, from $\Gamma
_{A_{0},...,A_{n}}^{\beta }(\tau _{0},...,\tau _{n})$ it is possible to
construct (up to unitary equivalence) a Hilbert space $\mathcal{H}^{\Gamma }$
with a self-adjoint Hamiltonian $H,$ a representation $\pi $ of the algebra
of observables $\mathcal{A}$ on $\mathcal{H}^{\Gamma }$ with a cyclic vector 
$\psi _{\beta },$ and the von Neumann algebra $\mathcal{B}\subset L(\mathcal{%
H}^{\Gamma })$ generated by $\pi (\mathcal{A}),$ such that the state $%
G_{\beta }(B):=(\psi _{\beta },B\psi _{\beta })_{\mathcal{H}^{\Gamma }}$ is
a KMS state on $\mathcal{B}$ w.r.t. the dynamics $\alpha
_{t}(B):=e^{itH}Be^{-itH},$ $B\in \mathcal{B}.$ On the other hand, due to
their properties the functions $\Gamma _{A_{0},...,A_{n}}^{\beta }(\tau
_{0},...,\tau _{n})$ uniquely determine a measure $\mu $ on $\Omega $ such
that 
\begin{equation}
E_{\mu }(A_{0}(\omega (\tau _{0}))...A_{n}(\omega (\tau _{n})))=\Gamma
_{A_{0},...,A_{n}}^{\beta }(\tau _{0},...,\tau _{n}).  \tag{7.9}
\end{equation}

In fact, measures $\mu _{\Lambda }$ in (7.7) correpond to the Gibbs
distributions in the finite volumes $\Lambda $ with the empty boundary
conditions of a lattice system on $\mathbb{Z}^{d}$ with the single spin
space $C_{\beta }$. Due to this observation, it would appear reasonable to
extend a class of limiting states to all Gibbs measures $\mu $ on $\Omega $
with the given Euclidean action functional $(I_{\Lambda })_{\Lambda \Subset 
\mathbb{Z}^{d}}$. In general, the set of all Gibbs measures $\mathcal{G}%
=\{\mu \}$ is large than the set of all accumulation points $\mu _{\infty }$
of $\{\mu _{\Lambda }\}$ when $\Lambda \nearrow \mathbb{Z}^{d}$.
Nevertheless, it is important to note that from any such Gibbs measure $\mu $
we are able to reconstruct a state $G_{\beta }$ with the Euclidean Green
functions (8.9). For the above reasons the measures $\mu \in \mathcal{G}%
_{\beta }$ will be called \emph{Euclidean Gibbs states (in the temperature
loop space representation) }for the quantum lattice system (3.1). Thus the
Euclidean approach provides us not only a \emph{constuctive way} to verify
the existence of limiting Gibbs states in the traditional scheme of
[BrRo81], but also enables us to \emph{extend substantially a class of states%
} for the considered quantum lattice systems.\smallskip

Nevertheless, despite a lot of papers dealing with the reconstruction
theorem on an abstract level or in some concrete models (see e.g. [AH-K75, Fr%
\"{o}77, GJO94a,b, BF02, ect.]), its consistent and mathematically rigorous
implementation in the case of temperature Gibbs states for quantum lattice
systems seems not have been performed in the existing literature. We
consider this as a very important open problem of high current interest,
which we shall work on in the future.\bigskip \bigskip \bigskip

\textsc{Acknowledgments}{\small :\quad Financial support by the DFG through
the Schwerpunktprogramm "Interacting Stochastic Systems of High Complexity"
(Research Projects AL 214/17-2 and RO 1195/5) and by the Lise-Meitner
Habilitation Programm is gratefully acknowledged.}

\bigskip

\bigskip

{\small S. A. :\quad Institut f\"{u}r Angewandte Mathematik, Universit\"{a}t
Bonn, D-53155 Bonn, Germany; BiBoS Research Centre, Bielefeld, Germany; and
CERFIM, Locarno, Switzerland.}

{\small Yu. K. :\quad Fakult\"{a}t f\"{u}r Mathematik and BiBoS Research
Centre, Bielefeld Universit\"{a}t, D-33615 Bielefeld, Germany; and Institute
of Mathematics, NASU, Kiev, Ukraine.}

{\small M. R. :\quad Fakult\"{a}t f\"{u}r Mathematik and BiBoS Research
Centre, Bielefeld Universit\"{a}t, D-33615 Bielefeld, Germany;}

{\small T. P. \ :\quad BiBoS Research Centre, Bielefeld Universit\"{a}t,
D-33615 Bielefeld, Germany;}

\bigskip

{\small E-mail:\quad albeverio@uni-bonn.de;
kondrat@mathematik.uni-bielefeld.de; }

{\small roeckner@mathematik.uni-bielefeld.de; pasurek@physik.uni-bielefeld.de%
}

\end{document}